\documentclass{scrartcl}

\usepackage{amsmath}
\usepackage[abbrev]{amsrefs}
\usepackage{amssymb}
\usepackage{amsthm}

\usepackage{1151arxiv}
\newtheorem{thm}{Theorem}[section]

\newtheorem{lem}[thm]{Lemma}
\newtheorem{prop}[thm]{Proposition}
\theoremstyle{definition}
\newtheorem{defn}[thm]{Definition}
\theoremstyle{remark}

\numberwithin{equation}{section}

\title{On a Simultaneous Approach to the Even and Odd Truncated Matricial Hamburger Moment Problems}
\author{Bernd Fritzsche \and Bernd Kirstein \and Conrad M\"{a}dler}
\dedication{Dedicated to Lev Aronovich Sakhnovich on the occasion of his 80th~birthday}
\begin{document}

\maketitle

\begin{abstract}
	The main goal of this paper is to achieve a simultaneous treatment of the even and odd truncated matricial Hamburger moment problems in the most general case. In the odd case, these results are completely new for the matrix case, whereas the scalar version was recently treated by V.~A.~Derkach, S.~Hassi and H.~S.~V.~de~Snoo~\cite{2011arXiv1101.0162D}. The even case was studied earlier by G.-N.~Chen and Y.-J.~Hu~\cite{MR1624548}. Our approach is based on Schur analysis methods. More precisely, we use two interrelated versions of Schur-type algorithms, namely an algebraic one and a function-theoretic one. The algebraic version was worked out in a former paper of the authors. It is an algorithm which is applied to finite or infinite sequences of complex matrices. The construction and investigation of the function-theoretic version of our Schur-type algorithm is a central theme of this paper. This algorithm will be applied to relevant subclasses of holomorphic matrix-valued functions of the Herglotz-Nevanlinna class. Using recent results on the holomorphicity of the Moore-Penrose inverse of matrix-valued Herglotz-Nevanlinna functions, we obtain a complete description of the solution set of the moment problem under consideration in the most general situation.
\end{abstract}
\tableofcontents
\section{Introduction}
The investigation of matrix versions of the Hamburger moment problem was initiated in the second half of the 1940's by M.~G.~Krein (see~\cites{MR0034964,MR0026759}). Mainly, based on new approaches (V.~P.~Potapov's method of fundamental matrix inequalities, reproducing kernel Hilbert spaces, Schur analysis) a renaissance of this topic could be observed in the last three decades (see e.\,g.\ the papers Kovalishina~\cite{MR703593}, Katsnelson~\cites{MR645308MR701996MR734686MR738449,MR1473259}, Dym~\cite{MR1018213}, Bolotnikov~\cite{MR1395706}, Chen/Hu~\cite{MR1624548} and the monographs Sakhnovich~\cite{MR1631843} and Bakonyi/Woerdeman~\cite{MR2807419}).

This paper continues the authors' investigations on matricial versions of the truncated Hamburger moment problem (see~\cites{MR2570113,MR2805417,103}). It was inspired by the recent work by V.~A.~Derkach, S.~Hassi and H.~S.~V.~de~Snoo on indefinite versions of truncated moment problems in the class of scalar generalized Nevanlinna functions (see~\cite{2011arXiv1101.0162D}). As a byproduct they obtained a result (see~\cite{2011arXiv1101.0162D}*{\ccoro{5.2}}) which seems to be new even for the odd truncated scalar Hamburger moment problem. More precisely, if $n$ is a non-negative integer and if $\seq{\su{j}}{j}{0}{2n+1}$ is a sequence of real numbers such that the Hankel matrix $\Hu{n}\defg\matauuo{\su{j+k}}{j,k}{0}{n}$ is positive Hermitian, then the set of \tStit{s} of all solutions of the Hamburger moment problem corresponding to the sequence $\seq{\su{j}}{j}{0}{2n+1}$ could be parametrized with the aid of a linear fractional transformation  the generating matrix-valued function of which is a \taaa{2}{2}{matrix} polynomial built from the sequence $\seq{\su{j}}{j}{0}{2n+1}$. The role of the set of parameter functions was taken by a particular subclass of the class of all functions which are holomorphic in the open upper half plane of the complex plane and have a non-negative imaginary part at each point. The study of such classes of holomorphic functions was initiated by I.~S.~Kats in~\cite{MR0080745} and then continued by S.~Hassi, H.~S.~V.~de~Snoo and A.~D.~I.~Willemsma (see~\cites{MR1376588,MR1451805}). Using several results originating in the just mentioned papers, V.~A.~Derkach, S.~Hassi and H.~S.~V.~de~Snoo~\cite{2011arXiv1101.0162D} were able to derive new Hamburger-Nevanlinna type theorems which played an important role in their approach to the moment problems for generalized Nevanlinna functions. Another cornerstone in the approach of V.~A.~Derkach, S.~Hassi and H.~S.~V.~de~Snoo~\cite{2011arXiv1101.0162D} is the application of the Schur-Chebyshev recursion algorithm, studied in the nondegenerate situation by M.~Derevyagin~\cite{MR1999775} (see also~\cite{MR2069282}). The main goal of this paper is to generalize the assertion of~\cite{2011arXiv1101.0162D}*{\ccoro{5.2}} to the matrix case without assuming some conditions of nondegeneracy. Roughly speaking, our strategy to realize this plan is modelled in basic steps along the lines of the approach in~\cite{2011arXiv1101.0162D}. We will use a lot of results from the recent paper~\cite{112} on various classes of holomorphic matrix-valued functions, which are generalizations of the scalar theory obtained in~\cites{MR0080745,MR1376588,MR1451805}. A further basic element in our approach is to use new Hamburger-Nevanlinna type theorems for matrix-valued holomorphic functions which will be derived in \rsect{S1510}. An other central feature of our approach is the use of Schur analysis methods. We will apply two interrelated versions of matricial Schur type algorithms, namely an algebraic one and a function-theoretic one. The algebraic version was worked out in~\cite{103}. It is an algorithm which is applied to a finite or infinite sequences of complex \tpqa{matrices}. An essential feature of this algorithm is that it preserves several properties of block Hankel matrices built from the sequences of complex \tqqa{matrices} under consideration (see \rsect{S1505}). The construction and investigation of the function-theoretic version of our Schur type algorithm is a central theme of this paper. This algorithm will be applied to relevant classes of holomorphic matrix-valued functions in the open upper half plane. In the scalar case, this algorithm coincides with the classical algorithm which was constructed by R.~Nevanlinna~\cite{Nev22} in adaptation of the classical Schur algorithm for bounded holomorphic functions, which is due to I.~Schur~\cite{Sch1718}. The idea how to build the algorithm in the matrix case was inspired by some constructions in the paper~\cite{MR1624548} by Chen and Hu. An essential point of our approach is an intensive use of the interplay between the function-theoretic and algebraic versions of our matricial Schur type algorithms. Both algorithms are formulated in terms of Moore-Penrose inverses of matrices. What concerns the function-theoretic version, it can be said that its effectiveness is mostly caused by recent results from~\cites{FKLR,112} on the holomorphicity of the Moore-Penrose inverse of special classes of holomorphic matrix-valued functions.

In order to describe more concretely the central topics studied in this paper, we introduce some notation. Throughout this paper, let $p$ and $q$ be positive integers. Let $\N$, $\NO$, $\Z$, $\R$, and  $\C$ be the set of all positive integers, the set of all non-negative integers, the set of all integers, the set of all real numbers, and the set of all complex numbers, respectively. For every choice of $\rho, \kappa \in \R \cup\{-\infty, +\infty\}$, let $\mn{\rho}{\kappa}\defg \setaa{k\in\Z}{\rho\leq k \leq \kappa}$\index{Z_,@$\mn{\rho}{\kappa}$}. We will write $\Cpq$, $\CHq$, $\Cggq$, and $\Cgq$ for the set of all complex \tpqa{matrices}, the set of all Hermitian complex \tqqa{matrices}, the set of all non-negative Hermitian complex \tqqa{matrices}, and the set of all positive Hermitian complex \tqqa{matrices}, respectively. 

We will use $\BAR$ to denote the $\sigma$\nobreakdash-algebra of all Borel subsets of $\R$. For all $\Omega \in \BAR\setminus\{\emptyset\}$, let $\BA{\Omega} \defg \BAR\cap \Omega$\index{B@$\BA{\Omega}$}. Furthermore, we will write $\MggqR$ to designate the set of all non-negative Hermitian \tqqa{measures} defined on $\BAR$, i.\,e., the set of $\sigma$\nobreakdash-additive mappings $\mu\colon\BAR \to \Cggq$. We will use the integration theory with respect to non-negative Hermitian \tqqa{measures} which was worked out independently by I.~S.~Kats~\cite{MR0080280} and M.~Rosenberg~\cite{MR0163346}. For all $j\in \NO$, we will use $\MggquR{j}$\index{M_>=,^()@$\MggquR{j}$} to denote the set of all $\sigma \in \MggqR$ such that the integral
\begin{equation}	\label{suo}
	\suo{j}{\sigma}
	\defg \int_\R t^j \sigma(\dif t)
\end{equation}
\index{_^[]@$\suo{j}{\sigma}$}exists. Furthermore, we set $\MgguqR{\infty}\defg\bigcap_{j=0}^\infty\MggquR{j}$\index{M_>=,^()@$\MggquR{\infty}$}.

\bremal{R1523}
	If $k,l\in\NO$ and $k<l$, then it can be verified, as in the scalar case, that the inclusion $\MgguqR{l}\subseteq\MgguqR{k}$ holds true.
\erema

The central problem in this paper is the truncated version of the following power moment problem of Hamburger type:
\begin{prob}[\mproblem{\R}{\kappa}{=}]
	Let $\kappa\in \NO\cup\set{+\infty}$ and let $\seq{\su{j}}{j}{0}{\kappa}$ be a sequence of complex \tqqa{matrices}. Describe the set $\MggqRag{\seq{\su{j}}{j}{0}{\kappa}}$\index{M_>=^[R;()_^,=]@$\MggqRag{\seq{\su{j}}{j}{0}{\kappa}}$} of all $\sigma \in \MggquR{\kappa}$ for which  $\su{j}^{[\sigma]}=\su{j}$ is fulfilled for all $j \in \mn{0}{\kappa}$.
\end{prob}\index{M[R;,=]@\mproblem{\R}{\kappa}{=}}
There is a further matricial version of the truncated Hamburger moment problem:
\begin{prob}[\mproblem{\R}{2n}{\leq}]
	Let $n\in\NO$ and let $\seq{\su{j}}{j}{0}{2n}$ be a sequence of complex \tqqa{matrices}. Describe the set $\MggqRakg{\seq{\su{j}}{j}{0}{2n}}$\index{M_>=^[R;()_^2,<=]@$\MggqRakg{\seq{\su{j}}{j}{0}{2n}}$} of all $\sigma\in\MggquR{2n}$ for which  $\suo{j}{\sigma}=\su{j}$ is satisfied for each $j\in\mn{0}{2n-1}$ whereas the matrix $\su{2n}-\suo{2n}{\sigma}$ is non-negative Hermitian.
\end{prob}
\index{M[R;^2,<=]@\mproblem{\R}{2n}{\leq}}The first investigation of \rprob{\mproblem{\R}{2n}{\leq}} goes back to I.~V.~Kovalishina~\cite{MR703593} who used V.~P.~Potapov's method of fundamental matrix inequalities. In the nondegenerate case, she obtained in~\cite{MR703593}*{\ctheo{$\mathcal{H}$}} a complete description of the Stieltjes transforms of the solution set of \rprob{\mproblem{\R}{2n}{\leq}} in terms of a linear fractional transformation. An extension of V.~P.~Potapov's method to degenerate situations was worked out by V.~K.~Dubovoj (see~\cites{Dub83,MR751390}) for the case of the matricial Schur problem. Using a modification of V.~K.~Dubovoj's method, \rprob{\mproblem{\R}{2n}{\leq}} could be handled by V.~A.~Bolotnikov (see~\cite{MR1395706}*{\ctheo{4.6}}) in the degenerate case. A common method of solving simultaneously the nondegenerate and degenerate versions of \rprob{\mproblem{\R}{2n}{\leq}} was presented by Chen and Hu in~\cite{MR1624548}. Their method is based on the use of a matricial Schur type algorithm involving matrix-valued continued fractions.

What concerns the investigation of interrelations between the two moment problems under consideration, we refer the reader to the papers~\cites{MR1624548,MR2570113}. A detailed treatment of the history of these two moment problems is contained in the introduction to the paper~\cite{MR2570113}.

In order to state a necessary and sufficient condition for the solvability of each of the above formulated moment problems, we have to recall the notion of two types of sequences of matrices. If $n \in \NO$ and if $\seq{\su{j}}{j}{0}{2n}$ is a sequence of complex \tqqa{matrices}, then $\seq{\su{j}}{j}{0}{2n}$ is called \emph{Hankel non-negative definite} if the block Hankel matrix
\[
	\Hu{n}
	\defg \matauuo{\su{j+k}}{j,k}{0}{n}
\]
\index{H_@$\Hu{n}$}is non-negative Hermitian. Let $\seq{\su{j}}{j}{0}{\infty}$ be a sequence of complex \tqqa{matrices}. Then $\seq{\su{j}}{j}{0}{\infty}$ is called \emph{\tHnnd{}} if $\seq{\su{j}}{j}{0}{2n}$ is \tHnnd{} for all $n\in\NO$. For all $\kappa\in\NO\cup\set{+\infty}$, we will write  $\Hggqu{2\kappa}$\index{H_,2^>=@$\Hggqu{2\kappa}$} for the set of all Hankel non-negative definite sequences $\seq{\su{j}}{j}{0}{2\kappa}$ of complex \tqqa{matrices}. Furthermore, for all $n \in \NO$, let $\Hggequ{2n}$\index{H_,2^>=,e@$\Hggequ{2n}$} be the set of all sequences $\seq{\su{j}}{j}{0}{2n}$ of complex \tqqa{matrices} for which there exist complex \tqqa{matrices} $\su{2n+1}$ and $\su{2n+2}$ such that $\seq{\su{j}}{j}{0}{2(n+1)} \in\Hggqu{2(n+1)}$, whereas $\Hggequ{2n+1}$\index{H_,2+1^>=,e@$\Hggequ{2n+1}$} stands for the set of all sequences $\seq{\su{j}}{j}{0}{2n+1}$ of complex \tqqa{matrices} for which there exist some $\su{2n+2}\in\Cqq$ such that $\seq{\su{j}}{j}{0}{2(n+1)} \in\Hggqu{2(n+1)}$. For each $m\in\NO$, the elements of the set $\Hggequ{m}$ are called \emph{\tHnnde{} sequences}. For technical reason, we set $\Hggeqinf\defg\Hggqinf$\index{H_,8^>=,e@$\Hggeqinf$}. Now we can characterize the situations that the mentioned problems have a solution:

\begin{thm}[\cite{MR1624548}*{\ctheo{3.2}}]\label{T1615}
	Let $n\in\NO$ and let $\seq{\su{j}}{j}{0}{2n}$ be a sequence of complex \tqqa{matrices}. Then $\MggqRakg{\seq{\su{j}}{j}{0}{2n}} \neq \emptyset$ if and only if $\seq{\su{j}}{j}{0}{2n} \in\Hggqu{2n}$.
\end{thm}

For an extension of \rtheo{T1615}, we refer the reader to~\cite{MR2570113}*{\ctheo{4.16}}. This extension says that $\MggqRakg{\seq{\su{j}}{j}{0}{2n}}\neq\emptyset$ if and only if this set contains a molecular measure (i.\,e., a measure which is concentrated on a finite set of real numbers). Now we characterize the solvability of Problem~\mproblem{\R}{m}{=}.

\begin{thm}[\cite{MR2805417}*{\ctheo{6.6}}]\label{1-417}
	Let $\kappa\in\NO\cup\set{+\infty}$ and let $\seq{\su{j}}{j}{0}{\kappa}$ be a sequence of complex \tqqa{matrices}. Then $\MggqRag{\seq{\su{j}}{j}{0}{\kappa}} \neq \emptyset$ if and only if $\seq{\su{j}}{j}{0}{\kappa} \in\Hggequ{\kappa}$.
\end{thm}

Note that, in the case of an even integer $\kappa$, \rtheo{1-417} was proved in Chen/Hu~\cite{MR1624548}*{Theorem~3.1}. For the case that $\kappa=2n$ with some non-negative integer $n$, Chen and Hu also stated a parametrization of the solution set of Problem \mproblem{\R}{\kappa}{=} in the language of the Stieltjes transforms (see~\cite{MR1624548}*{\ctheo{4.1}}). The goal of our paper here is to  give a parametrization of the solution set of Problem \mproblem{\R}{\kappa}{=} in the case that $\kappa$ is an odd integer (see \rtheoss{T1452}{T1640}). As already mentioned above in the scalar case $q=1$, such a description of the solution set of Problem \mproblem{\R}{m}{=} with odd integer $m$ was given by Derkach, Hassi, and de Snoo~\cite{2011arXiv1101.0162D}*{\ccoro{5.2}}. During the work at the odd case we observed that a slight modification of our approach also works for the even case. For this reason, have worked out a simultaneous approach to the odd and even versions of the problem. In this way, we will also alternately prove the corresponding results in the case of an arbitrary even integer $m$ (see \rtheoss{T1425}{T1620}).

\section{The Class $\RqP$}\label{S1425}
For all $A\in\Cpq$, let $\Kerna{A}$\index{N()@$\Kerna{A}$} be the null space of $A$ and $\Bilda{A}$\index{R()@$\Bilda{A}$} be the column space of $A$. If $A\in\Cqq$, then let $\re A\defg\frac{1}{2}(A+A^\ad)$\index{R@$\re A$} and $\im A\defg\frac{1}{2\I}(A-A^\ad)$\index{I@$\im A$} be the real part and the imaginary part of $A$, respectively. Let
\begin{align}\label{RIugg}
	\Rqgg&\defg\setaa{A\in\Cqq}{\re A\in\Cggq}&
	&\text{and}&
	\Iqgg&\defg\setaa{A\in\Cqq}{\im A\in\Cggq}.
\end{align}\index{R_>=@$\Rqgg$}\index{I_>=@$\Iqgg$}If $A$ and $B$ are Hermitian complex \tqqa{matrices}, then we will write $A\geq B$ or $B\leq A$ to indicate that the matrix $A-B$ is non-negative Hermitian, and $A>B$ means that $A-B$ is positive Hermitian.
As usual, for all $A\in\Cpq$, let $A^\MP$ be the Moore-Penrose pseudoinverse of $A$. If $A\in \Cpq$, then $\norma{A}$ is the operator norm of $A$. The notation $\Iq$ (or short $\EM$) stands for the unit matrix which belongs to $\Cqq$ and $\Opq$ (or short $\NM$) designates the null matrix which belongs to $\Cpq$. If $\nu$ is a non-negative Hermitian \tqqa{measure} on a measurable space $(\Omega,\mathfrak{A})$, then we will use $\LaaaC{\Omega}{\mathfrak{A}}{\nu}$\index{L^1(,,;C)@$\LaaaC{\Omega}{\mathfrak{A}}{\nu}$} to denote the space of all Borel-measurable functions $f\colon\Omega\to\C$ for which the integral $\int_\Omega f\dif\nu$ exists. Let $\ohe\defg\setaa{z\in\C}{\im z\in (0,+\infty)}$\index{Pi_+@$\ohe$} be the open upper half plane of $\C$. Of central importance to this paper is the class $\RqP$\index{R_()@$\RqP$} of all matrix-valued functions $F\colon\ohe\to\Cqq$ which are holomorphic in $\ohe$ and which satisfy $F(z)\in\Iqgg$ for all $z\in\ohe$. The elements of $\RqP$ are called \emph{$\xx{q}{q}$-matrix-valued Herglotz-Nevanlinna functions}. For a comprehensive survey on the class $\RqP$, we refer the reader to~\cite{MR1784638} and~\cite{MR2222521}*{\cSect{8}}. Herglotz-Nevanlinna functions admit a well-known integral representation. To state this, we observe that, for all $\nu\in\MggqR$ and each $z\in\C\setminus\R$, the function $\hu{z}\colon\R\to\C$\index{h_@$\hu{z}$} defined by $\hua{z}{t}\defg\frac{1+tz}{t-z}$ belongs to $\LaaaC{\R}{\BAR}{\nu}$.

\btheol{T1554}
	\benui
		\item Let $F\in\RqP$. Then there are unique matrices $\alpha\in\CHq$ and $\beta\in\Cggq$ and a unique non-negative Hermitian measure $\nu\in\MggqR$ such that
		\begin{align}\label{B1-1}
			F(z)
			&=\alpha+z\beta+\int_\R\frac{1+tz}{t-z}\nu(\dif t)&
			\text{for each }z&\in\ohe.
		\end{align}
		\il{T1554.b} If $\alpha\in\CHq$, $\beta\in\Cggq$ and $\nu\in\MggqR$, then $F\colon\ohe\to\Cqq$ defined by \eqref{B1-1} belongs to the class $\RqP$.
	\eenui
\etheo

For all $F\in\RqP$, the unique triple $(\alpha,\beta,\nu) \in \CHq\times\Cggq\times\MggqR$ for which the representation \eqref{B1-1} holds true is called the \emph{Nevanlinna parametrization of $F$} and we will also write $(\alphaF,\betaF,\nuF)$\index{alpha_@$\alphaF$}\index{beta_@$\betaF$}\index{nu_@$\nuF$} for $(\alpha,\beta,\nu)$. In particular, $\nuF$ is said to be the \emph{Nevanlinna measure of $F$}

\bexaml{E1338}
	Let $F\colon\ohe\to\Cqq$ be defined by $F(z)\defg\Oqq$. In view of \rtheo{T1554}, then $F\in\NFq$ and $(\alphaF,\betaF,\nuF)=(\Oqq,\Oqq,\zmq)$, where $\zmq\colon\BAR\to\Cqq$ is given by $\zmqa{B}\defg\Oqq$.
\eexam

\bremal{R1355}
	Let $F\in\NFq$ with \tNp{} $(\alphaF,\betaF,\nuF)$. If the matrix $\beta\in\Cqq$ satisfies $\betaF+\beta\in\Cggq$, then \rtheo{T1554} shows that the function $G\colon\ohe\to\Cqq$ defined by $G(z)\defg F(z)+z\beta$ belongs to $\NFq$ and satisfies $(\alphau{G},\betau{G},\nuu{G})=(\alphaF,\betaF+\beta,\nuF)$.
\erema

\bremal{R0824}
	Let $n\in\N$ and $\seq{p_k}{k}{1}{n}$ be a sequence of positive integers. For each $k\in\mn{1}{n} $ let $F_k\in\NFu{p_k}$ and $A_k\in\Coo{p_k}{q}$. Then it can be easily verified that $F\defg\sum_{k=1}^nA_k^\ad F_kA_k$ belongs to $\NFq$ and satisfies
	\[
		(\alphaF,\betaF,\nuF)
		=\lrk\sum_{k=1}^nA_k^\ad\alphau{F_k}A_k,\sum_{k=1}^nA_k^\ad\betau{F_k}A_k,\sum_{k=1}^nA_k^\ad\nuu{F_k}A_k\rrk.
	\]
\erema

\begin{prop}[see~\cite{MR1784638}*{\ctheo{5.4~(iv)}}]\label{P1513}
	Let $F\in\NFq$ with \tNp{} $(\alphaF,\betaF,\nuF)$. Then $\alphaF=\re[F(\I)]$ and $\betaF=\lim_{y\to+\infty}[\frac{1}{\I y}F(\I y)]$.
\end{prop}

We are particularly interested in the structure of the values of functions belonging to $\NFq$.

\blemml{L1444}
	Let $F\in\NFq$. For each $z\in\ohe$, then $\Bilda{[F(z)]^\ad}=\Bilda{F(z)}$, $\Kerna{[F(z)]^\ad}=\Kerna{F(z)}$, and $[F(z)][F(z)]^\MP=[F(z)]^\MP[F(z)]$.
\elemm
\bproof
	This follows from \rlemm{L1412} and \rprop{P1645}.
\eproof

\begin{prop}[\cite{112}*{\cprop{3.7}}]\label{P1356}
	Let $F\in\NFq$ with \tNp{} $(\alphaF,\betaF,\nuF)$. For all $z\in\ohe$, then
	\begin{gather*}
		\KernA{F(z)}
		=\Kerna{\alphaF}\cap\Kerna{\betaF}\cap\KernA{\nuFa{\R}},\\
		\BildA{F(z)}
		=\Bilda{\alphaF}+\Bilda{\betaF}+\BildA{\nuFa{\R}},
	\end{gather*}
	and
	\begin{align*}
		\KernA{\im\lek F(z)\rek}
		&=\Kerna{\betaF}\cap\KernA{\nuFa{\R}},&
		\BildA{\im\lek F(z)\rek}
		&=\Bilda{\betaF}+\BildA{\nuFa{\R}}.
	\end{align*}
\end{prop}
\rprop{P1356} contains essential information on the class $\NFq$. It indicates that, for an arbitrary function $F$ belonging to $\NFq$, the null space $\Kerna{F(z)}$ and the column space $\Bilda{F(z)}$ are independent of the concrete point $z\in\ohe$ and, furthermore, in which way these linear subspaces of $\Cq$ are determined by the \tNp{} $(\alphaF,\betaF,\nuF)$ of $F$.

In the sequel, we will sometimes meet situations where interrelations of the null space (resp.\ column space) of a function $F\in\NFq$ to the null space (resp.\ column space) of a given matrix $A\in\Cpq$ are of interest. More precisely, we will frequently apply the following observation.

\blemml{L1030}
	Let $A\in\Cpq$ and let $F\in\NFq$ with \tNp{} $(\alphaF,\betaF,\nuF)$. Then:
	\benui
		\il{L1030.a} The following statements are equivalent:
		\baeqii{0}
			\il{L1030.i} $\Kerna{A}\subseteq\Kerna{\alphaF}\cap\Kerna{\betaF}\cap\Kerna{\nuFa{\R}}$.
			\il{L1030.ii} $\Kerna{A}\subseteq\Kerna{F(z)}$ for all $z\in\ohe$.
			\il{L1030.iii} There exists a number $z\in\ohe$ with $\Kerna{A}\subseteq\Kerna{F(z)}$.
			\il{L1030.iv} $FA^\MP A=F$.
			\il{L1030.v} $\Bilda{\alphaF}+\Bilda{\betaF}+\Bilda{\nuFa{\R}}\subseteq\Bilda{A^\ad}$.
			\il{L1030.vi} $\Bilda{F(z)}\subseteq\Bilda{A^\ad}$ for all $z\in\ohe$.
			\il{L1030.vii} There exists a number $z\in\ohe$ with $\Bilda{F(z)}\subseteq\Bilda{A^\ad}$.
			\il{L1030.viii} $A^\MP AF=F$.
		\eaeqii
		\il{L1030.b} Let~\ref{L1030.i} be satisfied. If $p=q$ and if $A\in\CEPq$, then $AA^\MP F=F$ and $\Bilda{F(z)}\subseteq\Bilda{A}$ for all $z\in\ohe$.
	\eenui
\elemm
\bproof
	\eqref{L1030.a} \impl{L1030.ii}{L1030.iii} and \impl{L1030.vi}{L1030.vii}: These implications hold true obviously.
	
	\impl{L1030.i}{L1030.ii} and \impl{L1030.iii}{L1030.i}: Use \rprop{P1356}.
	
	\baeq{L1030.ii}{L1030.iv}
		This equivalence follows from \rpart{L1622.a} of \rrema{L1622}.
	\eaeq
	
	\baeq{L1030.ii}{L1030.vi}
		Because of \rlemm{L1444} we have $\Bilda{F(z)}=[\Kerna{F(z)}]^\bot$ for all $z\in\ohe$. Hence,~\ref{L1030.ii} and~\ref{L1030.vi} are equivalent.
	\eaeq
	
	\impl{L1030.v}{L1030.vi} and \impl{L1030.vii}{L1030.v}: Use \rprop{P1356}.
	
	\baeq{L1030.vi}{L1030.viii}
		Use \rprop{P1619} and \rpart{R1553.c} of \rrema{R1553}.
	\eaeq
	
	\eqref{L1030.b} In view of $A\in\CEPq$, we have $\Bilda{A^\ad}=\Bilda{A}$ and, taking \rprop{P1645} into account, furthermore, $AA^\MP=A^\MP A$. Thus,~\eqref{L1030.b} follows from~\eqref{L1030.a}.
\eproof

%
%
%

A generic application of \rlemm{L1030} will be concerned with situations where the matrix $A$ even belongs to $\Cggq$.

In our subsequent considerations we will very often use the Moore-Penrose inverse of functions belonging to the class $\NFq$. In this connection, the following result turns out to be of central importance.

\begin{prop}[\cite{112}*{\cprop{3.8}}]\label{P1511}
	Let $F\in\NFq$ with \tNp{} $(\alphaF,\betaF,\nuF)$. Then $-F^\MP\in\NFq$ and $\alphau{-F^\MP}=-[F(\I)]^\MP\alphaF([F(\I)]^\MP)^\ad$.
\end{prop}


\section{On Some Subclasses of $\NFq$}\label{S1516}
An essential feature of our subsequent considerations is the use of a whole variety of different subclasses of $\NFq$. In this section, we summarize basic facts about these subclasses under the special orientation of this paper. The first part of these subclasses concerns objects which are already well-studied, whereas the larger remaining part deals with subclasses of $\NFq$, which were introduced and studied very recently by the authors in~\cite{112}. The latter collection of subclasses are characterized by growth properties on the positive imaginary axis. It should be mentioned that the scalar versions of the function classes were introduced and studied in the papers~\cites{MR0080745,MR1376588,MR1451805}. An important subclass of the class $\RqP$ is the set $\RqI$ of all $F\in\RqP$  for which the function $g\colon\R\to\R$\index{g@$g$} defined by $g(t)\defg t^2+1$ belongs to $\LaaaR{\R}{\BAR}{\nuF}$, where $\nuF$ is taken from the \tNpa{F}. In view of \rrema{R1523}, a member $F$ of the class $\RqP$ belongs to $\RqI$ if and only if  $\nuF\in\MggquR{2}$. For all $F\in\RqI$, then the mapping $\sigmaF\colon\BAR\to\Cqq$\index{sigma_@$\sigmaF$} given by
\bgl{sigmaua}
	\sigmaFa{B}
	\defg\int_B(t^2+1)\nuFa{\dif t}
\eg
is a well-defined non-negative Hermitian measure belonging to $\MggqR$. The measure $\sigmaF$ is called the \emph{\tsma{F}}. In this paper, we encounter mostly situations in which, for a given function $F\in\NFprimeq$, the \tsm{} $\sigmaF$ plays a more important role than the Nevanlinna measure $\nuF$. Observe that the functions of the class $\RqI$ also admit a special integral representation (see~\cite{112}*{Theorem~4.3}).

Now we will see that each measure $\sigma\in\MggqR$ generates by a special integration procedure, in a natural way, a function belonging to $\NFprimeq$.

\bpropl{P0836}
	Let $\sigma\in\MggqR$.
	\benui
		\il{P0836.a} Let $\nu_{[\sigma]}\colon\BAR\to\Cqq$ be defined by
		\[
			B
			\mapsto\int_B\frac{1}{t^2+1}\sigma(\dif t).
		\]
		Then $\nu_{[\sigma]}\in\MgguqR{2}$ and
		\begin{align*}
			\suo{j}{\nu_{[\sigma]}}&=\int_\R\frac{t^j}{t^2+1}\sigma(\dif t)&\text{for each }j\in&\set{1,2}.
		\end{align*}
		\il{P0836.b} Let $\Stitu{\sigma}\colon\ohe\to\Cqq$ be defined by
		\bgl{Stitua}
			\Stitua{\sigma}{z}
			\defg\int_\R\frac{1}{t-z}\sigma(\dif t).
		\eg
		Then $\Stitu{\sigma}$ is a matrix-valued function belonging to the class $\NFprimeq$ with \tNp{} $(\alphau{\Stitu{\sigma}},\betau{\Stitu{\sigma}},\nuu{\Stitu{\sigma}})=(\suo{1}{\nu_{[\sigma]}},\Oqq,\nu_{[\sigma]})$ and \tsm{} $\sigmau{\Stitu{\sigma}}=\sigma$.
	\eenui
\eprop
\bproof
	\eqref{P0836.a} This follows immediately from~\cite{112}*{\cprop{B.5}}.
	
	\eqref{P0836.b} Taking into account that, for each $z\in\ohe$ and each $t\in\R$, the identity
	\[
		\frac{1+tz}{t-z}
		=\lrk\frac{1}{t-z}-\frac{t}{1+t^2}\rrk(1+t^2)
	\]
	holds true, the assertion of~\eqref{P0836.b} follows by direct computations, using \rrema{R1523} and~\cite{112}*{\cprop{B.5}}.
\eproof

Let $\sigma\in\MggqR$. Then the function $\Stitu{\sigma}$ defined in \eqref{Stitua} is called the \emph{\tStita{\sigma}}. \rPart{P0836.b} of \rprop{P0836} shows that each $\sigma\in\MggqR$ occurs as \tsm{} of an appropriately chosen function belonging to $\NFprimeq$.
%
%
Now we want to characterize the set of all \tStit{s} of measures belonging to $\MggqR$ by the matricial generalization of a classical result due to R.~Nevanlinna~\cite{Nev22}.

\btheol{T0933}
	Let
	\bgl{RFtO}
		\RFtOq
		\defg\SetaA{F\in\RqP}{\sup_{y\in[1,+\infty)}y\normA{F(\I y)}<+\infty}.
	\eg
	\index{R^~_0,()@$\RFtOq$}Then $\RFtOq=\setaa{\Stitu{\sigma}}{\sigma\in\MggqR}$ and the mapping $\sigma\mapsto\Stitu{\sigma}$ is a bijective correspondence between $\MggqR$ and $\RFtOq$.
\etheo

For each $F\in\RFtOq$ the unique measure $\sigma\in\MggqR$ satisfying $\Stitu{\sigma}=F$ is called the \emph{\tStima{F}}. \rtheo{T0933} indicates that the \tStit{} $\Stitu{\sigma}$ of a measure $\sigma\in\MggqR$ is characterized by a particular mild growth on the positive imaginary axis.

In view of \rtheo{T0933}, Problem \mproblem{\R}{\kappa}{=} can be given a first reformulation as an equivalent problem in the class $\RtOq$ as follows:

\begin{prob}[\iproblem{\R}{\kappa}{=}]
 Let $\kappa\in\NO\cup\set{+\infty}$ and let $\seq{\su{j}}{j}{0}{\kappa}$ be a sequence of complex \tqqa{matrices}. Describe the set of all matrix-valued functions $S \in \RtOq$, the Stieltjes measure of which belongs to  $\MggqRag{\seq{\su{j}}{j}{0}{\kappa}}$.
\end{prob}\index{R[R;,=]@\iproblem{\R}{\kappa}{=}}

In \rsect{S1510}, we will state a reformulation of the original power moment problem \mproblem{\R}{\kappa}{=} as an equivalent problem of finding a prescribed asymptotic expansion in a sector of the open upper half plane $\ohe$. Furthermore, we will see that a detailed analysis of the behaviour of the concrete functions of $F\in\NFq$ under study on the positive imaginary axis is extremely useful. For this reason, we turn now our attention to some subclasses of $\NFq$ which are described in terms of their growth on the positive imaginary axis. First we consider the set
\bgl{Ruo[-2]}
	\NFuq{-2}
	\defg\SetaA{F\in\NFq}{\lim_{y\to+\infty}\lek\frac{1}{y}\normA{F(\I y)}\rek=0}.
\eg
\index{R_^[-2]()@$\RqK{-2}$}


In the following, we will use the symbol $\Leb$ to denote the Lebesgue measure defined on $\BAR$.

\bexaml{E1047*}
	Let $A\in\Iqgg$ and let $F\colon\ohe\to\Cqq$ be given by $F(z)\defg A$. Then it is immediately checked that $F\in\NFuq{-2}$. Let $\mu\colon\BAR\to[0,+\infty)$ be defined by
	\[
		\mu(B)
		\defg\frac{1}{\pi}\int_B\frac{1}{1+u^2}\Leb(\dif u).
	\]
	Using the residue theorem, it can be verified then by direct computation that
	\[
		(\alphaF,\betaF,\nuF)
		=\lrk\re A,\Oqq,(\im A)\mu\rrk.
	\]
\eexam

\brema
	From \eqref{RFtO} and \eqref{Ruo[-2]} it is obvious that $\RFtOq\subseteq\NFuq{-2}$.
\erema

\bremal{R1538}
	In view of \eqref{Ruo[-2]} and \rprop{P1513}, we have
	\[
		\NFuq{-2}
		=\SetAa{F\in\NFq}{\betaF=\Oqq}.
	\]
\erema

\bremal{R1535}
	Let $n\in\N$ and let $\seq{p_k}{k}{1}{n}$ be a sequence from $\N$. For all $k\in\mn{1}{n} $, let $F_k\in\NFuu{-2}{p_k}$ and let $A_k\in\Coo{p_k}{q}$. In view of \rrema{R1538} and~\cite{112}*{\crema{3.4}}, then $\sum_{k=1}^nA_k^\ad F_kA_k\in\NFuq{-2}$.
\erema

Now we state modifications of \rprop{P1356} and \rlemm{L1030} for $\NFuq{-2}$.

\bpropl{L1608}
	Let $F\in\NFuq{-2}$ and let $(\alphaF,\betaF,\nuF)$ the \tNpa{F}. For all $z\in\ohe$, then
	\begin{align*}
		\KernA{F(z)}&
		=\Kerna{\alphaF}\cap\KernA{\nuFa{\R}},&
		\BildA{F(z)}
		&=\Bilda{\alphaF}+\BildA{\nuFa{\R}}\\
	\intertext{and}
		\KernA{\im\lek F(z)\rek}
		&=\KernA{\nuFa{\R}},&
		\BildA{\im\lek F(z)\rek}
		&=\BildA{\nuFa{\R}}.
	\end{align*}
\eprop
\bproof
	Combine \rprop{P1356} and \rrema{R1538}.
\eproof

\blemml{L1028}
	Let $A\in\Cpq$ and let $F\in\NFuq{-2}$. Then the statements
	\baeqi{8}
		\il{L1028.ix} $\Kerna{A}\subseteq\Kerna{\alphaF}\cap\Kerna{\nuFa{\R}}$.
	\eaeqi
	and
	\baeqi{9}
		\item $\Bilda{\alphaF}+\Bilda{\nuFa{\R}}\subseteq\Bilda{A^\ad}$.
	\eaeqi
	are equivalent. Furthermore,~\ref{L1028.ix} is equivalent to each of the statements~\ref{L1030.i}--\ref{L1030.viii} stated in \rlemm{L1030}.
\elemm
\bproof
	Combine \rrema{R1538} and \rlemm{L1030}.
\eproof

In the following, the notation $\tilde\Leb$\index{lambda^~@$\tilde\Leb$} stands for the Lebesgue measure defined on $\BA{[1,+\infty)}$. We now recall some subclasses of $\NFq$, which were introduced and studied in~\cite{112}. A further subclass of $\RqP$, which we will need in the following, is the set
\bgl{Ruo[-1]}
	\RqK{-1}
	\defg\SetaA{F \in \RqP}{\int_{[1,+\infty)} \frac{1}{y}\normA{\im F(\I y)}\tilde{\Leb}(\dif y)
	<+\infty}.
\eg
\index{R_^[-1]()@$\RqK{-1}$}

\bexam
	Let $A\in\Iqgg$ and let $F\colon\ohe\to\Cqq$ be defined by $F(z)\defg A$. Then from \rexam{E1047*} and \eqref{Ruo[-1]} it is obvious that $F\in\NFuq{-1}$ if and only if $\im A=\Oqq$.
\eexam

\bremal{R1425}
	From~\cite{112}*{\clemm{5.1}} we get more information about the \tNp{} of the functions belonging to $\NFuq{-1}$, namely for all $F \in \RqK{-1}$, we have $\betaF=\Oqq$, $\nuF\in\MggquR{1}$, and the function $h\colon\R\to\R$\index{h@$h$} defined by $h(t)\defg\frac{t^2+1}{\abs{t}+1}$ belongs to $\LaaaR{\R}{\BAR}{\nuF}$. This implies that, for all $F \in \RqK{-1}$, the mapping  $\muF\colon\BAR\to\Cqq$\index{mu_@$\muF$} given by
	\bgl{muua}
		\muFa{B}
		\defg\int_B\frac{t^2+1}{\abs{t}+1}\nuFa{\dif t}
	\eg
	is a well-defined non-negative Hermitian measure belonging to $\MggqR$ and, in view of \eqref{suo} and~\cite{112}*{Remark~B.4}, that the matrix
	\bgl{gammau}
		\gammaF
		\defg\alphaF-\suo{1}{\nuF}
	\eg
	\index{gamma@$\gammaF$}satisfies $(\gammaF)^\ad=\gammaF$.
\erema

\brema
	From \rremass{R1538}{R1425} we see that $\NFuq{-1}\subseteq\NFuq{-2}$.
\erema

The next result indicates that functions which belong to $\RqK{-1}$ admit a particular characterization in terms of a constant Hermitian matrix and an integral representation:

\begin{thm}[\cite{112}*{\ctheo{5.6}}]\label{T1347}
	\benui
		\item Each matrix-valued function $F$ belonging to $\RqK{-1}$ admits, for all $z\in\ohe$, the representation
		\[
			F(z)
			=\gammaF+\int_\R\frac{\abs{t}+1}{t-z}\muFa{\dif t},
		\]
		where $\gammaF$ and $\muF$ are given via \eqref{gammau} and \eqref{muua}, respectively.
		\il{T1347.b} Let $\gamma\in\CHq$ and let $\mu\in\MggqR$. Then $F\colon\ohe\to\Cqq$ defined by
		\[
			F(z)
			\defg\gamma+\int_\R\frac{\abs{t}+1}{t-z}\mu(\dif t)
		\]
		belongs to $\RqK{-1}$ and $(\gammaF,\muF)=(\gamma,\mu)$ holds true.
	\eenui
\end{thm}

\bexaml{E1350}
	Let $F\colon\ohe\to\Cqq$ be defined by $F(z)\defg\Oqq$. In view of \rpart{T1347.b} of \rtheo{T1347}, then $F\in\NFuq{-1}$ and $(\gammaF,\muF)=(\Oqq,\zmq)$, where $\zmq\colon\BAR\to\Cqq$ is given by $\zmqa{B}\defg\Oqq$.
\eexam

The next result describes the asymptotics of the functions belonging to $\NFuq{-1}$ on the positive imaginary axis.

\begin{prop}[\cite{112}*{\cprop{5.8}}]\label{P1401}
	Let $F\in\NFuq{-1}$. Then
	\begin{align*}
		\lim_{y\to+\infty}\re F(\I y)&=\gammaF,&
		\lim_{y\to+\infty}\im F(\I y)&=\Oqq
	\end{align*}
	and
	\bgl{N15-3}
		\lim_{y\to+\infty} F(\I y)
		=\gammaF.
	\eg
\end{prop}

In our scale of subclasses of $\NFq$, we next consider the class
\bgl{Ruo[0]}
	\RqK{0}
	\defg\SetaA{F\in\RqP}{\sup_{y\in[1,+\infty)}y\normA{\im F(\I y)}<+\infty}.
\eg
\index{R_^[0]()@$\RqK{0}$}In view of~\cite{112}*{\cprop{6.4}}, we have
\begin{equation}\label{FSG}
	\RqK{0}
	=\RqK{-1}\cap\RqI,
\end{equation}
where $\RqK{-1}$ is given via \eqref{Ruo[-1]}. Now we recall a special characterization for the functions of the class $\RqK{0}$. 

\begin{thm}[\cite{112}*{\ctheo{6.3}}]\label{T1446}
	\benui
		\il{T1446.a} Let $F\in\RqK{0}$. Then $F$ belongs to the class $\NFprimeq$ and, if $\gammaF$ and $\sigmaF$ are given via \eqref{gammau} and \eqref{sigmaua}, respectively, then
		 \[ F(z)=\gammaF+\int_\R\frac{1}{t-z}\sigmaFa{\dif t}\]
		 for all $z\in\ohe$.
		\item For all $\gamma\in\CHq$ and each $\sigma\in\MggqR$, the function $F\colon\ohe\to\Cqq$ given by
		\[
			F(z)
			\defg\gamma+\int_\R\frac{1}{t-z}\sigma(\dif t)
		\]
		belongs to $\RqK{0}$ and satisfies $(\gammaF,\sigmaF)=(\gamma,\sigma)$.
	\eenui
\end{thm}

For each $\kappa\in\set{-2,-1,0}$, the class $\RqK{\kappa}$ is already defined. In view of \eqref{FSG}, we have $\NFuq{0}\subseteq\NFprimeq$. Thus, the functions $F$ belonging to $\NFuq{0}$ have a well-defined \tsm{} $\sigmaF$, which is given via \eqref{sigmaua}. So, for all $\kappa \in  \N \cup \{+\infty\}$, the class
\bgl{Ruo}
	\RqK{\kappa}
	\defg\SetaA{{F\in\RqK{0}}}{\sigmaF\in\MggqkappaR}
\eg
\index{R_^[]()@$\RqK{\kappa}$}is well-defined. We recall that then the following result holds true:

\begin{prop}\label{P0601}
	Let $\kappa\in\mn{-2}{+\infty}\cup\set{+\infty}$. Then
	\[
		\RqK{\kappa}
		=\SetaA{F\in\RqP}{\betaF=\Oqq\text{ and }\nuF\in\MggquR{\kappa+2}},
	\]
	where $\betaF$ and $\nuF$ are taken from the \tNpa{F}.
\end{prop}
\bproof
	Use \rrema{R1538} in the case $\kappa=-2$ and~\cite{112}*{\cprop{7.3}} in the case $\kappa\geq-1$.
\eproof

\bremal{R1215}
	Observe that \rprop{P0601} shows that the proper inclusions
	\[
		\RqK{+\infty}
		\subsetneq\RqK{l}
		\subsetneq\RqK{k}
		\subsetneq\RqP
	\]
	are fulfilled for all $l\in\NO$ and all $k\in\mn{-1}{l-1}$.
\erema


For all $\kappa\in\mn{-1}{+\infty}\cup\set{+\infty}$, we now consider the class
\bgl{Ruu}
	\RFkappaq
	\defg\SetAa{F\in \RqK{\kappa}}{\gammaF=\Oqq}.
\eg\index{R_,()@$\RFkappaq$}

\brema
	From \eqref{Ruu} and \rrema{R1215} one sees that
	\[
		\RFuq{+\infty}
		\subsetneq\RFuq{l}
		\subsetneq\RFuq{k}
		\subsetneq\NFq
	\]
	holds true for all $l\in\NO$ and all $k\in\mn{-1}{l-1}$.
\erema

\bexaml{E1351}
	Let $F\colon\ohe\to\Cqq$ be defined by $F(z)\defg\Oqq$. In view of \rexam{E1350}, then $F\in\RFuq{-1}$.
\eexam

The following results complement \rprop{P1356}. More precisely, we state now some modifications of \rprop{P1356} for various subclasses of $\NFq$.

\begin{lem}\label{L1406}
	Let $\kappa\in\mn{-1}{+\infty}\cup\set{+\infty}$ and let $F\in\NFkappaq$. For each $z\in\ohe$, then
	\begin{align*}
		\BildA{F(z)}&=\Bilda{\gammaF}+\BildA{\muFa{\R}},&
		\KernA{F(z)}&=\Kerna{\gammaF}\cap\KernA{\muFa{\R}}\\
	\intertext{and}
		\BildA{\im\lek F(z)\rek}&=\BildA{\muFa{\R}},&
		\KernA{\im\lek F(z)\rek}&=\KernA{\muFa{\R}}.
	\end{align*}
\end{lem}
\bproof
	Use~\cite{112}*{\clemm{5.4}} and \rrema{R1215}.
\eproof

\begin{lem}[\cite{112}*{\clemm{8.1}}]\label{R1410}
	Let $\kappa\in\mn{-1}{+\infty}\cup\set{+\infty}$ and let $F\in\RFkappaq$. For each $z\in\ohe$, then
	\begin{align*}
		\BildA{F(z)}
		&=\BildA{\muFa{\R}}
		=\BildA{\im\lek F(z)\rek}\\
	\intertext{and}
		\KernA{F(z)}
		&=\KernA{\muFa{\R}}
		=\KernA{\im\lek F(z)\rek}.
	\end{align*}
\end{lem}


\begin{lem}[\cite{112}*{\clemm{8.2}}]\label{R1414}
	Let $\kappa\in\NO\cup\set{+\infty}$ and let $F\in\RFkappaq$. For each $z\in\ohe$, then
	\begin{align*}
		\BildA{F(z)}
		&=\BildA{\sigmaFa{\R}}
		=\BildA{\im\lek F(z)\rek}\\
	\intertext{and}
		\KernA{F(z)}
		&=\KernA{\sigmaFa{\R}}
		=\KernA{\im\lek F(z)\rek}.
	\end{align*}
\end{lem}


\bremal{R1437}
	Let  $\kappa\in\mn{-1}{+\infty}\cup\set{+\infty}$ and $F\in\RFkappaq$. Then from \eqref{Ruu} and \rprop{P1401} we see that
	\bgl{R1437.B1}
		\lim_{y\to+\infty}F(\I y)
		=\Oqq.
	\eg
\erema

\bremal{R1433}
	In view of \rprop{P0601}, \eqref{gammau} and \eqref{Ruu}, for all $\kappa\in\mn{-1}{+\infty}\cup\set{+\infty}$, we have 
	\begin{multline*}
		\RFkappaq\\
		=\SetaA{F\in \RqP}{\betaF=\Oqq\text{ and }\nuF\in\MggquR{\kappa+2}\text{ and }\alphaF=\suo{1}{\nuF}}.
	\end{multline*}
	In particular, for all $\kappa \in \NO\cup\set{+\infty}$ and each $m\in\mn{-1}{\kappa-1}$, the class $\RFkappaq $ is a proper subset of $\RFuu{m}{q}$.
\erema

\blemml{R1041}
	Let $A\in\Cpq$, let $\kappa\in\mn{-1}{+\infty}\cup\set{+\infty}$, and let $F\in\RFuq{\kappa}$. Then the statements
	\baeqi{10}
		\il{R1041.xi} $\Kerna{A}\subseteq\Kerna{\muFa{\R}}$.
	\eaeqi
	and
	\baeqi{11}
		\item $\Bilda{\muFa{\R}}\subseteq\Bilda{A^\ad}$.
	\eaeqi
	are equivalent. Furthermore,~\ref{R1041.xi} is equivalent to each of the statements~\ref{L1030.i}--\ref{L1030.viii} in \rlemm{L1030}.
\elemm
\bproof
	Combine \rrema{R1433} and \rlemmss{R1410}{L1030}.
\eproof

\brema
	From~\cite{112}*{Proposition~6.4} we know that 
	\[
		\ROq
		=\RFuq{-1} \cap \RqI.
	\]
\erema

\bremal{P1514}
	From \eqref{Ruu}, \eqref{Ruo[-1]} and \rprop{P1401} we see that
	\begin{multline*}
		\RFuq{-1}=\\
		\SetaA{F\in\NFq}{\int_{[1,+\infty)}\frac{1}{y}\normA{\im\lek F(\I y)\rek}\tilde{\Leb}(\dif y)
	<+\infty\text{ and }\lim_{y\to+\infty}\normA{F(\I y)}=0}.
	\end{multline*}
\erema

Now we get that the classes given in \eqref{Ruu} and \eqref{RFtO} coincide.

%

\bpropl{P1417}
	$\RFOq=\RFtOq$.
\eprop
\bproof
	Combine \rtheoss{T0933}{T1446}.
\eproof

\bcorol{C1403}
	Let $\kappa\in\NO\cup\set{+\infty}$. Then $\RFuq{\kappa}\subseteq\RFtOq$.
\ecoro
\bproof
	Combine \eqref{Ruu}, \rrema{R1215}, and \rprop{P1417}.
\eproof

\bremal{R1455}
	Let $n\in\N$, $\kappa\in\mn{-1}{+\infty}\cup\set{+\infty}$, and $\seq{p_k}{k}{1}{n}$ be a sequence from $\N$. For $k\in\mn{1}{n}$ let $F_k\in\RFuu{\kappa}{p_k}$ and $A_k\in\Coo{p_k}{q}$. Then~\cite{112}*{\crema{8.3}} shows that $\sum_{k=1}^nA_k^\ad F_kA_k\in\RFuq{\kappa}$.
\erema

\section{On the Classes $\PFevenqa{A}$ and $\PFoddqa{A}$}\label{S1431}
In this section, we consider particular subclasses of the classes $\NFuq{-2}$ and $\RFuq{-1}$, which were introduced in \eqref{Ruo[-2]} and \eqref{Ruu}, respectively. We have seen in \rprop{P1356} that, for an arbitrary function $F\in\NFq$, the null space of $F(z)$ is independent from the concrete choice of $z\in\ohe$. For the cases $F\in\NFuq{-2}$ or $F\in\RFuq{-1}$, a complete description of this constant null space was given in \rprop{L1608} and \rlemm{R1410}, respectively. Against to this background, we single out now special subclasses of $\NFuq{-2}$ and $\RFuq{-1}$ which are characterized by the interrelation of the constant null space to the null space of a prescribed matrix $A\in\Cpq$. More precisely, for all $A\in\Cpq$, let
\begin{align}\label{Pevenua}
	\PFevenqa{A}
	&\defg\SetAa{F\in\NFuq{-2}}{\Kerna{A}\subseteq\Kerna{\alphaF}\cap\KernA{\nuFa{\R}}}\\
\intertext{and let}\label{Poddua}
	\PFoddqa{A}
	&\defg\SetaA{F\in\RFuq{-1}}{\Kerna{A}\subseteq\KernA{\muFa{\R}}},
\end{align}
\index{Peven@$\PFevenqa{A}$}\index{Podd@$\PFoddqa{A}$}where $\muF$ is given via \eqref{muua}. The choice of the terminology is caused by the role which the sets introduced in \eqref{Pevenua} and \eqref{Poddua} will later play in the framework of the even and odd version of the truncated matricial Hamburger moment problem, respectively. The role of the matrix $A$ will be taken then by matrices which are generated from the sequence of data matrices of the problem of consideration via a Schur type algorithm.

\bremal{R0827}
	If $A\in\Cpq$ satisfies $\Kerna{A}=\set{\Ouu{q}{1}}$, then $\PFevenqa{A}=\NFuq{-2}$ and $\PFoddqa{A}=\RFuq{-1}$. In particular, this situation arises in the case that $p=q$ and $\det A\neq0$ are fulfilled.
\erema

\bexaml{E1354}
	Let $A\in\Cpq$ and let $F\colon\ohe\to\Cqq$ be defined by $F(z)\defg\Oqq$. In view of $\Oqq\in\Iqgg$, \rexamss{E1047*}{E1338}, and \eqref{Pevenua}, then $F\in\PFevenqa{A}$, and, in view of \rexamss{E1351}{E1350} and \eqref{Poddua}, furthermore $F\in\PFoddqa{A}$.
\eexam


\blemml{L1407}
	Let $A\in\Cpq$ and $F\colon\ohe\to\Cqq$. Then:
	\benui
		\il{L1407.a} The following statements are equivalent:
		\baeqii{0}
			\item $F\in\PFevenqa{A}$.
			\item $F\in\NFuq{-2}$ and $\Kerna{A}\subseteq\Kerna{F(z)}$ for all $z\in\ohe$.
			\item $F\in\NFq$, $\lim_{y\to+\infty}[\frac{1}{y}\norma{F(\I y)}]=0$, and $\Kerna{A}\subseteq\Kerna{F(z)}$ for all $z\in\ohe$.
			\item $F\in\NFq$, $\betaF=\Oqq$, and $\Kerna{A}\subseteq\Kerna{\alphaF}\cap\KernA{\nuFa{\R}}$.
		\eaeqii
		\il{L1407.b} The following statements are equivalent:
		\baeqii{4}
			\item $F\in\PFoddqa{A}$.
			\item $F\in\RFuq{-1}$ and $\Kerna{A}\subseteq\Kerna{F(z)}$ for all $z\in\ohe$.
			\item $F\in\NFq$, $\int_{[1,+\infty)}\frac{1}{y}\norma{\im[F(\I y)]}\tilde{\Leb}(\dif y)<+\infty$, $\lim_{y\to+\infty}\norma{F(\I y)}=0$, and $\Kerna{A}\subseteq\Kerna{F(z)}$ for all $z\in\ohe$.
			\item $F\in\NFq$, $\betaF=\Oqq$, $\nuF\in\MggquR{1}$, $\alphaF=\suo{1}{\nuF}$, and $\Kerna{A}\subseteq\Kerna{\nuFa{\R}}$.
		\eaeqii
	\eenui
\elemm
\bproof
	\eqref{L1407.a} This follows from \eqref{Pevenua}, \rlemm{L1028}, \eqref{Ruo[-2]}, and \rrema{R1538}.
	
	\eqref{L1407.b} Use \eqref{Poddua}, \rlemm{R1041}, \rremass{P1514}{R1433}, \eqref{Ruu}, and~\cite{112}*{\crema{5.3}}.
\eproof

\brema
	Let $z_0\in\ohe$. From \rlemm{L1407}, \rprop{L1608} and \rlemm{R1410} one can easily see then that:
	\benui
		\item If $F\in\NFuq{-2}$, then $F\in\PFevenqa{F(z_0)}$.
		\item If $F\in\RFuq{-1}$, then $F\in\PFoddqa{F(z_0)}$.
	\eenui
\erema
%

\bremal{R1543}
	Let $A\in\Cpq$ and let $B\in\Coo{r}{q}$ with $\Kerna{A}\subseteq\Kerna{B}$. In view of \eqref{Pevenua} and \eqref{Poddua}, then $\PFevenqa{B}\subseteq\PFevenqa{A}$ and $\PFoddqa{B}\subseteq\PFoddqa{A}$.
\erema

\bpropl{L1551}
	Let $A\in\Cpq$.
	\benui
		\il{L1551.a} $\PFevenqa{A}=\setaa{A^\MP AFA^\MP A}{F\in\NFuq{-2}}$.
		\il{L1551.b} $\PFoddqa{A}=\setaa{A^\MP AFA^\MP A}{F\in\RFuq{-1}}$.
	\eenui
\eprop
\bproof
	\eqref{L1551.a} Let $F\in\PFevenqa{A}$. In view of \eqref{Pevenua}, we have then
	\begin{align}\label{L1551.1}
		F&\in\NFuq{-2}&
		&\text{and}&
		\Kerna{A}&\subseteq\Kerna{\alphaF}\cap\KernA{\nuFa{\R}}.
	\end{align}
	Taking \eqref{L1551.1} into account, \rlemm{L1028} yields for $z\in\ohe$ now $F(z)A^\MP A=F(z)$ and $A^\MP AF(z)=F(z)$. Consequently, $A^\MP AFA^\MP A=F$. Combining this with \eqref{L1551.1}, we infer
	\bgl{L1551.5}
		\PFevenqa{A}
		\subseteq\SetaA{A^\MP AFA^\MP A}{F\in\NFuq{-2}}.
	\eg
	Conversely, let us consider an arbitrary
	\bgl{L1551.6}
		F
		\in\NFuq{-2}.
	\eg
	In view of \eqref{L1551.6} and $(A^\MP A)^\ad=A^\MP A$, from \rrema{R1535}, we see that $G\defg A^\MP AFA^\MP A$ fulfills
	\bgl{L1551.8}
		G
		\in\NFuq{-2}
	\eg
	and $GA^\MP A=(A^\MP AFA^\MP A)A^\MP A=G$. Thus, \rpart{L1622.a} of \rrema{L1622} gives $\Kerna{A}\subseteq\Kerna{G(z)}$ for each $z\in\ohe$. Combining this with \eqref{L1551.8} and applying \rpart{L1407.a} of \rlemm{L1407}, we conclude $G\in\PFevenqa{A}$. Thus, $\PFevenqa{A}\supseteq\setaa{A^\MP AFA^\MP A}{F\in\NFuq{-2}}$. This inclusion and \eqref{L1551.5} show that \rpart{L1551.a} holds.
	
	\eqref{L1551.b} Let $F\in\PFoddqa{A}$. In view of \eqref{Poddua}, we have then
	\bgl{L1551.11}
		F
		\in\RFuq{-1}
	\eg
	and $\Kerna{A}\subseteq\KernA{\muFa{\R}}$. Consequently, from \rlemm{R1041} we see that $A^\MP AFA^\MP A=F$ holds true. Thus, \eqref{L1551.11} yields
	\bgl{L1551.15}
		\PFoddqa{A}
		\subseteq\SetaA{A^\MP AFA^\MP A}{F\in\RFuq{-1}}.
	\eg
	Conversely, we now consider an arbitrary $F\in\RFuq{-1}$. From \rrema{R1455} and $(A^\MP A)^\ad=A^\MP A$ we get then that $G\defg A^\MP AFA^\MP A$ belongs to $\RFuq{-1}$ and fulfills $\Kerna{A}\subseteq\Kerna{G(z)}$ for each $z\in\ohe$. Thus, the application of \rpart{L1407.b} of \rlemm{L1407} yields $G\in\PFoddqa{A}$. Consequently,
	\[
		\PFoddqa{A}
		\supseteq\SetaA{A^\MP AFA^\MP A}{F\in\RFuq{-1}}.
	\]
	Combining this with \eqref{L1551.15} completes the proof of \rpart{L1551.b}.
\eproof


\bcoro
	Let $A\in\Cpq$ and let $F\in\PFevenqa{A}\cup\PFoddqa{A}$. Then $FA^\MP A=F$ and $A^\MP AF=F$.
\ecoro


\blemm
	Let $A\in\Cpq$. Then $\PFoddqa{A}\subseteq\PFevenqa{A}$.
\elemm
\bproof
	Let $F\in\PFoddqa{A}$. In view of \eqref{Poddua}, we have then $F\in\RFuq{-1}$. From \rrema{R1433} we get then $F\in\NFq$ and $\betaF=\Oqq$. Hence, \rrema{R1538} yields $F\in\NFuq{-2}$. Furthermore, we obtain
	\[
		\Kerna{A}
		\subseteq\KernA{\muFa{\R}}
		=\KernA{F(\I)}
		=\Kerna{\alphaF}\cap\KernA{\nuFa{\R}},
	\]
	where the inclusion is due to $F\in\PFoddqa{A}$ and \eqref{Poddua}, the 1st~equation is due to $F\in\RFuq{-1}$ and \rlemm{R1410}, and the 2nd~equation is due to $F\in\NFuq{-2}$ and \rprop{L1608}. In view of \eqref{Pevenua}, we get then $F\in\PFevenqa{A}$.
\eproof

%
%

The following result contains essential information on the structure of the sets $\PFevenqa{A}$ and $\PFoddqa{A}$, where $A\in\Cpq$.

\bpropl{L1249*}
	Let $A\in\Cpq$. Then:
	\benui
		\il{L1249*.a} If $A=\Opq$, then $\PFevenqa{A}=\set{F}$ and $\PFoddqa{A}=\set{F}$, where $F\colon\ohe\to\Cqq$ is defined by $F(z)\defg\Oqq$.
		\il{L1249*.b} Suppose that $r\defg\rank A$ fulfills $r\geq1$. Let $u_1,u_2,\dotsc,u_r$ be an orthonormal basis of $\Bilda{A^\ad}$ and let $U\defg\brow u_1,u_2,\dotsc,u_r\erow$. Then
		\begin{align*}
			\PFevenqa{A}&=\SetaA{UfU^\ad}{f\in\NFuu{-2}{r}}\\
		\intertext{and}
			\PFoddqa{A}&=\SetaA{UfU^\ad}{f\in\RFuu{-1}{r}}.
		\end{align*}
		\il{L1249*.c} If $f,g\in\NFuu{-2}{r}\cup\RFuu{-1}{r}$ are such that $UfU^\ad=UgU^\ad$, then $f=g$.
	\eenui
\eprop
\bproof
	\eqref{L1249*.a} This follows from \rprop{L1551} and \rexam{E1354}.
	
	\eqref{L1249*.b} Let $G\in\PFevenqa{A}$ (resp.\ $G\in\PFoddqa{A}$). In view of \rpart{L1551.a} (resp.\ \rpart{L1551.b}) of \rprop{L1551}, there exists an $F\in\NFuq{-2}$ (resp.\ $F\in\RFuq{-1}$) such that $G=A^\MP AFA^\MP A$. Let $f\defg U^\ad FU$. Because of \rrema{R1535} (resp.\ \rrema{R1455}), then $f\in\NFuu{-2}{r}$ (resp.\ $f\in\RFuu{-1}{r}$). In view of \rrema{R1526}, we have $UU^\ad=A^\MP A$. Thus, $G=UU^\ad FUU^\ad=UfU^\ad$. Hence,
	\begin{align*}
		\PFevenqa{A}&\subseteq\SetaA{UfU^\ad}{f\in\NFuu{-2}{r}}\\
	\intertext{(resp.}
		\PFoddqa{A}&\subseteq\SetaA{UfU^\ad}{f\in\RFuu{-1}{r}}\text{).}
	\end{align*}
	Conversely, let $f\in\NFuu{-2}{r}$ (resp.\ $f\in\RFuu{-1}{r}$). In view of \rrema{R1535} (resp.\ \rrema{R1455}), then
	\ba
		UfU^\ad&\in\NFuq{-2}&(\text{resp.\ }UfU^\ad&\in\RFuq{-1}).
	\ea
	Now we consider an arbitrary $x\in\Kerna{A}$. In view of the construction of $U$ and the relation $[\Kerna{A}]^\bot=\Bilda{A^\ad}$, we get $U^\ad x=\Ouu{r}{1}$. Thus, $x\in\Kerna{U^\ad}$. Consequently, for each $z\in\ohe$ we get $\Kerna{A}\subseteq\Kerna{(UfU^\ad)(z)}$. The application of \rpart{L1407.a} (resp.\ \rpart{L1407.b}) of \rlemm{L1407} yields now
	\begin{align*}
		\SetaA{UfU^\ad}{f\in\NFuu{-2}{r}}&\subseteq\PFevenqa{A}\\
	\intertext{(resp.}
		\SetaA{UfU^\ad}{f\in\RFuu{-1}{r}}&\subseteq\PFoddqa{A}\text{).}
	\end{align*}
	This completes the proof of~\eqref{L1249*.b}.
	
	\eqref{L1249*.c} In view of \rrema{R1526}, we have $U^\ad U=\Iu{r}$. Thus $UfU^\ad=UgU^\ad$ implies $f=g$.
\eproof

\section{The Classes  $\Rqkj$ and $ \Rkqj$}\label{S1438}
In this section, we consider particular subclasses of the class $\RqK{\kappa}$, which was introduced in \eqref{Ruo[0]} for $\kappa=0$ and in \eqref{Ruo} for $\kappa\in\N\cup\set{+\infty}$. Because of \eqref{FSG} and \eqref{Ruo}, we have the inclusion
\begin{align}\label{G1004}
	\NFuq{\kappa}&\subseteq\NFprimeq&\text{for each }\kappa&\in\NO\cup\set{+\infty}.
\end{align}
In view of \eqref{G1004}, for each function $F$ belonging to one of the classes $\NFuq{\kappa}$ with some $\kappa\in\NO\cup\set{+\infty}$, the spectral measure $\sigmaF$ given by \eqref{sigmaua} is well-defined. Now taking \rrema{R1523} into account, we turn our attention to subclasses of functions $F\in\NFuq{\kappa}$ with prescribed parameter $\gammaF$ and prescribed first $\kappa+1$ power moments of the \tsm{} $\sigmaF$.

Taking \eqref{sigmaua} and \eqref{gammau} into account, for all $\kappa \in \NO\cup\set{+\infty}$ and each sequence $\seq{\su{j}}{j}{-1}{\kappa}$ of complex \tqqa{matrices}, now we consider the class
\begin{multline} \label{Ruoa}
	\RqoS{\kappa}\\
	\defg\SetaA{F \in \RqK{\kappa}}{\gammaF=-\su{-1} \text{ and } \sigmaF\in \MggqRAg{\seq{\su{j}}{j}{0}{\kappa}}}
\end{multline}
\index{R_^[][;]@$\Rqos{\kappa}$}and, for all $\kappa \in \NO\cup\set{+\infty}$ and each sequence $\seq{\su{j}}{j}{0}{\kappa}$ from $\Cqq$, furthermore 
\begin{equation} \label{Ruua}
	\RuqS{\kappa}
	\defg \SetaA{F \in \RKq}{\sigmaF\in \MggqRAg{\seq{\su{j}}{j}{0}{\kappa}}},
\end{equation}
\index{R_,[;]@$\Ruqs{\kappa}$}where $\RKq$ is defined in \eqref{Ruu}.

\bremal{R1404}
	Let $\kappa\in\NO\cup\set{+\infty}$ and let $\seq{\su{j}}{j}{0}{\kappa}$ be a sequence from $\Cqq$. Let $t_{-1}\defg\Oqq$ and let $t_j\defg\su{j}$ for all $j\in\mn{0}{\kappa}$. Then taking \eqref{Ruoa}, \eqref{Ruua}, and \eqref{Ruu} into account we see that
	\[
		\RFuqA{\kappa}{\seq{\su{j}}{j}{0}{\kappa}}
		=\NFuqA{\kappa}{\seq{t_j}{j}{-1}{\kappa}}.
	\]
\erema

\brema
	Let $\kappa\in\NO\cup\set{+\infty}$ and let $\seq{\su{j}}{j}{-1}{\kappa}$ be a sequence from $\Cqq$. If $\iota\in\NO\cup\set{+\infty}$ with $\iota\leq\kappa$, then \eqref{Ruoa} and \rrema{R1215} show that
	\[
		\NFuqA{\kappa}{\seq{\su{j}}{j}{-1}{\kappa}}
		=\bigcap_{m=0}^\kappa\NFuqA{m}{\seq{\su{j}}{j}{-1}{m}}.
	\]
\erema

\bremal{R0938*}
	Let $\kappa\in\NO\cup\set{+\infty}$ and let $\seq{\su{j}}{j}{0}{\kappa}$ be a sequence from $\Cqq$. If $\iota\in\NO\cup\set{+\infty}$ with $\iota\leq\kappa$ then \eqref{Ruua}, \eqref{Ruu}, and \rrema{R1215} show that
	\[
		\RFuqA{\kappa}{\seq{\su{j}}{j}{0}{\kappa}}
		\subseteq\RFuqA{\iota}{\seq{\su{j}}{j}{0}{\iota}}.
	\]
	In particular,
	\[
		\RFuqA{\kappa}{\seq{\su{j}}{j}{0}{\kappa}}
		=\bigcap_{m=0}^\kappa\RFuqA{m}{\seq{\su{j}}{j}{0}{m}}.
	\]
\erema



Now we characterize those sequences for which the sets defined in \eqref{Ruoa} and \eqref{Ruua} are non-empty.

\bpropl{P1407}
	Let $\kappa\in\NO\cup\set{+\infty}$ and let $\seq{\su{j}}{j}{-1}{\kappa}$ be a sequence from $\Cqq$. Then:
	\benui
		\il{P1407.a} $\NFuqa{\kappa}{\seq{\su{j}}{j}{-1}{\kappa}}\neq\emptyset$ if and only if $\seq{\su{j}}{j}{0}{\kappa}\in\Hggequ{\kappa}$ and $\su{-1}\in\CHq$.
		\il{P1407.b} $\RFuqa{\kappa}{\seq{\su{j}}{j}{0}{\kappa}}\neq\emptyset$ if and only if $\seq{\su{j}}{j}{0}{\kappa}\in\Hggequ{\kappa}$.
	\eenui
\eprop
\bproof
	\eqref{P1407.a} Combine \eqref{Ruoa}, \eqref{Ruo}, and \rtheoss{T1446}{T1554}.
	
	\eqref{P1407.b} Combine \eqref{Ruua}, \eqref{Ruu}, and \rtheoss{T1446}{T1554}.
\eproof


Now we state a useful characterization of the set of functions given in \eqref{Ruua}.

\bpropl{R0947*}
	Let $\kappa\in\NO\cup\set{+\infty}$ and let $\seq{\su{j}}{j}{0}{\kappa}$ be a sequence of complex \tqqa{matrices}. Then
	\[
		\RuqS{\kappa}
		=\SetaA{F\in\RFtOq}{\sigmaF\in\MggqRAg{\seq{\su{j}}{j}{0}{\kappa}}}.
	\]
\eprop
\bproof
	In view of \eqref{Ruua} and \rcoro{C1403}, we have
	\[
		\RuqS{\kappa}
		\subseteq\SetaA{F\in\RFtOq}{\sigmaF\in\MggqRAg{\seq{\su{j}}{j}{0}{\kappa}}}.
	\]

	Conversely, now let $F\in\RFtOq$ be such that
	\bgl{Nr.2-3-1}
		\sigmaF
		\in\MggqRAg{\seq{\su{j}}{j}{0}{\kappa}}. 
	\eg
	From \eqref{Nr.2-3-1} we get
	\bgl{R0947*.2}
		\sigmaF
		\in\MgguqR{\kappa},
	\eg
	whereas $F\in\RFtOq$ and \rprop{P1417} imply $F\in\RFOq$. Hence, \eqref{Ruu} yields $F\in\NFuq{0}$ and $\gammaF=\Oqq$. From $F\in\NFuq{0}$, \eqref{R0947*.2}, and \eqref{Ruo} we see that $F\in\NFuq{\kappa}$. Combining this with $\gammaF=\Oqq$, we infer from \eqref{Ruu} that $F\in\RFuq{\kappa}$. Because of \eqref{Nr.2-3-1} and \eqref{Ruua}, then it follows $F\in\RuqS{\kappa}$. Thus, the proof is complete.
%
\eproof


\bcorol{R1103}
	Let $\kappa\in\NO\cup\set{+\infty}$ and let $\seq{\su{j}}{j}{0}{\kappa}$ be a sequence of complex \tqqa{matrices}. In view of \eqref{Stitua} then
	\[
		\RuqS{\kappa}
		=\SetaA{\Stitu{\sigma}}{\sigma\in\MggqRAg{\seq{\su{j}}{j}{0}{\kappa}}}.
	\]
\ecoro
\bproof
	Combine \rprop{R0947*}, \rtheo{T0933}, and \rpart{P0836.b} of \rprop{P0836}.
\eproof

\rcoro{R1103} shows that $\RFuqa{\kappa}{\seq{\su{j}}{j}{0}{\kappa}}$ coincides with the solution set of Problem~\iproblem{\R}{\kappa}{=}, which is via \tStit{} equivalent to the original Problem~\mproblem{\R}{\kappa}{=}. Thus, the investigation of the set $\RFuqa{\kappa}{\seq{\su{j}}{j}{0}{\kappa}}$ is a central theme of our further considerations. The next result contains essential information on the functions belonging to this set.

\bpropl{L0908}
	Let $\kappa\in\NO\cup\set{+\infty}$, let $\seq{\su{j}}{j}{0}{\kappa}\in\Hggequ{\kappa}$, and let $F\in\RFuqa{\kappa}{\seq{\su{j}}{j}{0}{\kappa}}$. Then:
	\benui
		\il{L0908.a} For each $z\in\ohe$, the equations $\Kerna{F(z)}=\Kerna{\su{0}}$ and $\Bilda{F(z)}=\Bilda{\su{0}}$ hold.
		\il{L0908.b} For each $z\in\ohe$, the equations $[F(z)][F(z)]^\MP=\su{0}\su{0}^\MP$ and $[F(z)]^\MP[F(z)]=\su{0}^\MP\su{0}$ hold.
		\il{L0908.c} The function $F$ belongs to the class $\NFq$ and its \tNp{} $(\alphaF,\betaF,\nuF)$ satisfies
		\begin{align}\label{N18-1}
			\Kerna{\su{0}}
			&=\Kerna{\alphaF}\cap\Kerna{\betaF}\cap\KernA{\nuFa{\R}}\\
		\intertext{and}\label{N18-2}
			\Bilda{\su{0}}
			&=\Bilda{\alphaF}+\Bilda{\betaF}+\BildA{\nuFa{\R}}.
		\end{align}
	\eenui
\eprop
\bproof
	\eqref{L0908.a} In view of the choice of $F$, we get from \eqref{Ruua} that
	\bgl{L0908.1}
		F
		\in\RFuq{\kappa}
	\eg
	and $\sigmaF\in\MggqRag{\seq{\su{j}}{j}{0}{\kappa}}$. Thus, we have $\suo{0}{\sigmaF}=\su{0}$. Otherwise, from \eqref{suo} we have $\suo{0}{\sigmaF}=\sigmaFa{\R}$. Hence, $\sigmaFa{\R}=\su{0}$. Combining this with \eqref{L0908.1}, we obtain from \rlemm{R1414} all assertions of~\eqref{L0908.a}.
	
	\eqref{L0908.b} The assertions of~\eqref{L0908.b} follow from~\eqref{L0908.a} by application of \rrema{L1622}.
	
	\eqref{L0908.c} From \eqref{L0908.1} and \rrema{R1433} we get $F\in\NFq$. Now the assertions of~\eqref{L0908.c} follow by combination of \rprop{P1356} with~\eqref{L0908.a}.
\eproof


The next result establishes a connection to the class $\PFoddqa{\su{0}}$ introduced in \rsect{S1431}.

\blemml{L1553}
	Let $\kappa\in\NO\cup\set{+\infty}$ and let $\seq{\su{j}}{j}{0}{\kappa}$ be a sequence of complex \tqqa{matrices}. Then $\RFuqa{\kappa}{\seq{\su{j}}{j}{0}{\kappa}}\subseteq\PFoddqa{\su{0}}$.
\elemm
\bproof
	Let $F\in\RFuqa{\kappa}{\seq{\su{j}}{j}{0}{\kappa}}$. From \eqref{Ruua} we get then $F\in\RFuq{\kappa}$. Thus, from \rrema{R1215} and \eqref{Ruu} we infer $F\in\RFuq{-1}$. Let $z\in\ohe$. In view of $F\in\RFuq{-1}$ \rlemm{R1410} yields $\Kerna{F(z)}=\Kerna{\muFa{\R}}$, whereas \rpart{L0908.c} of \rprop{L0908} gives $\Kerna{F(z)}=\Kerna{\su{0}}$. Hence, $\Kerna{\su{0}}=\Kerna{\muFa{\R}}$. In view of \eqref{L1551.11} and \eqref{Poddua}, this implies $F\in\PFoddqa{\su{0}}$.
\eproof

Now we want to discuss the asymptotic behaviour of functions belonging to $\NFq$. For this reason, we need a particular construction, which will be introduced now.

\bremal{R1013}
	Let $\kappa\in\mn{-1}{+\infty}\cup\set{+\infty}$, let  $\seq{\su{j}}{j}{-1}{\kappa}$ be a sequence of complex \tpqa{matrices}, let $\mG$ be a non-empty subset of $\C$, and let $F\colon \mG\to \Cpq$ be a matrix-valued function. For all $k\in\mn{-1}{\kappa}$, let then $\trauo{F}{k}{s}\colon \mG\to \Cpq$\index{_^<>@$\trauo{F}{k}{s}$} be defined by 
	\begin{equation} \label{Fuoa}
		\trauoa{F}{k}{s}{z}
		\defg z^{k+1}F(z)+\sum_{j=0}^{k+1}z^{k+1-j}\su{j-1}.	
	\end{equation}
	For every choice of integers $k$ and $l$ with $-1 \leq k<l\leq \kappa$ and each $z \in \mG$, then it is immediately checked that
	\bgl{R1013.1}
		\trauoa{F}{l}{s}{z}
		=z^{l-k}\trauoa{F}{k}{s}{z}+\sum_{j=1}^{l-k}z^{l-k-j}\su{k+j}
	\eg
	and in the case $z\neq 0$ furthermore	 
	\begin{equation} \label{N4-12}
		\trauoa{F}{k}{s}{z}=z^{k-l}\lek\trauoa{F}{l}{s}{z}-\sum_{j=0}^{l-k-1}z^j\su{l-j}\rek.
	\end{equation}
\erema

In the following, we will often use the construction of \rrema{R1013} for the case that $\mathcal{G}=\ohe$ and that the function $F$ belongs to particular subclasses of $\RqP$. We start with the case that $F$ belongs to the class introduced in \eqref{Ruoa} and investigate the associated sequence $\seq{\trauo{F}{k}{s}}{k}{-1}{\kappa}$. First we show that these functions admit integral representations with respect to $\sigmaF$.

\bpropl{L1605}
	Let $\kappa\in\NO\cup\set{+\infty}$, let $\seq{\su{j}}{j}{-1}{\kappa}$ be a sequence of complex \tqqa{matrices}, and let $F\in \Rqos{\kappa}$. For each  $k\in\mn{-1}{\kappa}$ and each $z\in\ohe$, then
	\begin{equation} \label{N4-3}
		\trauoa{F}{k}{s}{z}
		=\int_\R\frac{t^{k+1}}{t-z}\sigmaF(\dif t),
	\end{equation}
	where $\sigmaF$ is given via \eqref{sigmaua}.
\eprop
\bproof
	Since $F$ belongs to $ \Rqos{\kappa}$, we have 
	\begin{align} \label{N4-5}
		F&\in \RqK{\kappa},&
		\gammaF&=-\su{-1}&
		&\text{and}&
		\sigmaF&\in\MggqRAg{\seq{\su{j}}{j}{0}{\kappa}}.
	\end{align}
	The last relation in \eqref{N4-5} implies
	\begin{align}\label{L1605.1}
		\int_\R t^j\sigmaF(\dif t)&=\su{j}&\text{for each }j&\in\mn{0}{\kappa}.
	\end{align}
	From the definition \eqref{Ruo} of the class $\RqK{\kappa}$, we see that $F$ belongs in particular to $\RqK{0}$. Thus, for all $z \in \ohe$, \rpart{T1446.a} of \rtheo{T1446} and the second relation in \eqref{N4-5} yield
	\begin{equation} \label{N4-1}
			F(z)=-\su{-1} + \int_\R\frac{1}{t-z}\sigmaF(\dif t) 
	\end{equation}
	and, in view of \eqref{Fuoa} and \eqref{N4-1}, then, in particular, that \eqref{N4-3} holds true for $k=-1$. Because of \eqref{Fuoa} and \eqref{N4-1}, for all $z \in \ohe$, we obtain
	\bsp
		\trauoa{F}{0}{s}{z}
		&=z\lek F(z)+\su{-1}\rek+s_{0}
		= \int_\R\frac{z}{t-z}\sigmaF(\dif t) +\int_\R 1\sigmaF(\dif t) \\
		&= \int_\R \lrk\frac{z}{t-z} +1\rrk\sigmaF(\dif t)
		=\int_\R\frac{t^{0+1}}{t-z}\sigmaF(\dif t)
	\esp
	and, consequently, \eqref{N4-3} for $k=0$. It remains  to consider the case $k \in \mn{1}{\kappa}$.  For all $z \in \ohe$ and each $t \in \R$	, it is readily checked that 
	\begin{equation} \label{N5-1}
			\frac{z^{k+1}}{t-z} + \sum_{j=1}^{k+1}z^{k+1-j} t^{j-1} = \frac{t^{k+1}}{t-z}
	\end{equation}
	is true. Using \eqref{Fuoa}, \eqref{N4-1}, \eqref{L1605.1}, and \eqref{N5-1}, for all $z \in \ohe$, we get then
	\bsp
		\trauoa{F}{k}{s}{z}
		&=z^{k+1}\lek-\su{-1}+\int_\R\frac{1}{t-z}\sigmaF(\dif t)\rek+\sum_{j=0}^{k+1}z^{k+1-j} \su{j-1}\\
		&=z^{k+1} \int_\R\frac{1}{t-z}\sigmaF(\dif t)
		+  \sum_{j=1}^{k+1}z^{k+1-j} \int_\R t^{j-1}\sigmaF(\dif t)\\
		&= \int_\R \lrk\frac{z^{k+1}}{t-z} + \sum_{j=1}^{k+1}z^{k+1-j} t^{j-1}\rrk\sigmaF(\dif t)
		=  \int_\R\frac{t^{k+1}}{t-z}\sigmaF(\dif t).
	\esp
	Thus, \eqref{N4-3} is also proved in the case  $k \in \mn{1}{\kappa}$.
\eproof

\bremal{R1050}
	Let $\kappa\in\NO\cup\set{+\infty}$, let $\seq{\su{j}}{j}{-1}{\kappa}$ be a sequence of complex \tqqa{matrices},  let $k\in\mn{0}{\kappa}$, and let $F\in\Rqos{\kappa}$. Then \eqref{R1013.1} and \rprop{L1605} show that
	\[
		\trauoa{F}{l}{s}{z}
		=\sum_{j=0}^{l-k-1}z^j\su{l-j}+z^{l-k}\int_\R\frac{t^{k+1}}{t-z}\sigmaF(\dif t)
	\]
	holds true for all $l\in\mn{k+1}{\kappa}$ and each $z\in\ohe$. 	
\erema

Now we consider the functions $\trauo{F}{k}{s}$ occurring in \rprop{L1605} especially for odd numbers $k$.

\bpropl{R1640}
	Let  $\kappa\in\NO\cup\set{+\infty}$ and let  $\seq{\su{j}}{j}{-1}{\kappa}$ be a sequence  of complex \tqqa{matrices}. Furthermore let $F\in\Rqos{\kappa}$ and let $\sigmaF$ be given via \eqref{sigmaua}. Then:
	\benui
		\il{R1640.a} For all $n\in\NO$ with $2n-1\leq\kappa$, the function $\trauo{F}{2n-1}{s}$ belongs to $\RFuq{-1}$ and the measures $\nuu{\trauo{F}{2n-1}{s}}$ and $\mu_{\trauo{F}{2n-1}{s}}$ admit, for all $B\in\BAR$, the representations
		\begin{equation} \label{N6-1-1}
			\nuu{\trauo{F}{2n-1}{s}}(B)
			=\int_B\frac{t^{2n}}{t^2+1}\sigmaFa{\dif t}
		\end{equation}
		and
		\begin{equation} \label{N6-1}
			\mu_{\trauo{F}{2n-1}{s}}(B)
			=\int_B\frac{t^{2n}}{\abs{t}+1}\sigmaFa{\dif t}.
		\end{equation}
		\il{R1640.b} For all $n\in\NO$ with $2n\leq\kappa$, it holds $\trauo{F}{2n-1}{s}\in\Ruua{\kappa-2n}{q}{\seq{\su{2n+j}}{j}{0}{\kappa-2n}}$ and,  for all $B\in\BAR$, furthermore
		\begin{equation} \label{N6-3}
			\sigmau{\trauo{F}{2n-1}{s}}(B)=\int_B t^{2n}\sigmaFa{\dif t}.
		\end{equation}
	\eenui
\eprop
\bproof
	Since $F$ belongs to $\Rqos{\kappa}$, we have \eqref{N4-5}. In particular, $\sigmaF \in \MggquR{\kappa}$. \rprop{L1605} shows that, for all $z \in \ohe$,  furthermore
	\begin{equation} \label{N6}
		\trauoa{F}{2n-1}{s}{z}
		=\int_\R\frac{t^{2n}}{t-z}\sigmaF(\dif t).		
	\end{equation}
	
	\eqref{R1640.a} Let $n\in\NO$ with $2n-1\leq\kappa$. If $n \geq 1$, then $0\leq \frac{t^{2n}}{t^2+1} \leq t^{2n-2}$ and $\abs{t \frac{t^{2n}}{t^2+1}} \leq \abs{t^{2n-1}}$ for all $t \in \R$, whereas if $n =0$, then $0\leq \frac{t^{2n}}{t^2+1} \leq 1$ and $\abs{t \frac{t^{2n}}{t^2+1}} \leq 1$ for all $t \in \R$. Thus, from $\sigmaF \in \MggquR{\kappa}$ and~\cite{112}*{\cprop{B.5}} we get that the following statement holds true:
	\bAeqi{0}
		\il{R1640.I} The mapping $\nu \colon \BAR \to \Cqq$ given by
		\begin{equation} \label{N7-1}
			\nu(B) \defg \int_B\frac{t^{2n}}{t^2+1}\sigmaFa{\dif t}
		\end{equation}
		is a well-defined non-negative Hermitian measure which belongs to $\MggquR{1}$ and which fulfills 
		\[
			\suo{1}{\nu}
			=\int_\R t \frac{t^{2n}}{t^2+1}\sigmaF(\dif t).		
		\]
	\eAeqi
	For all $z \in \ohe$, from \eqref{N6} and~\ref{R1640.I} we then get 
	\bgl{R1640.1}
		\trauoa{F}{2n-1}{s}{z}
		=\suo{1}{\nu}+ \int_\R\frac{t^{2n}}{t-z}\sigmaF(\dif t) -\int_\R t \frac{t^{2n}}{t^2+1}\sigmaF(\dif t).		
	\eg
	Since~\ref{R1640.I} and~\cite{112}*{\cprop{B.5}} provide us 
	\bsp
		\int_\R\frac{t^{2n}}{t-z}\sigmaF(\dif t) -\int_\R t \frac{t^{2n}}{t^2+1}\sigmaF(\dif t)
		&=\int_\R\lrk\frac{t^{2n}}{t-z}- t \frac{t^{2n}}{t^2+1}\rrk\sigmaF(\dif t)\\
		&= \int_\R \frac{1+tz}{t-z}\frac{t^{2n}}{t^2+1}\sigmaF(\dif t)\\
		&=\int_\R \frac{1+tz}{t-z} \nu(\dif t)
	\esp
	for all $z \in \ohe$, from \eqref{R1640.1} we conclude that 
	\[
	\trauoa{F}{2n-1}{s}{z}
		=\suo{1}{\nu}+ z\cdot\Oqq + \int_\R \frac{1+tz}{t-z} \nu(\dif t)
	\]
	for all $z \in \ohe$. Since \eqref{suo} and~\cite{112}*{\cprop{B.4}} show that $(\suo{1}{\nu})^\ad=\suo{1}{\nu}$, \rtheo{T1554} then yields that $\trauo{F}{2n-1}{s}\in \RqP$ with \tNp{}
	\bgl{R1640.5}
		(\alphau{\trauo{F}{2n-1}{s}}, \betau{\trauo{F}{2n-1}{s}}, \nuu{\trauo{F}{2n-1}{s}})
		= (\suo{1}{\nu}, \Oqq, \nu).
	\eg
	From \eqref{R1640.5} we get
	\bgl{R1640.2}
		\alphau{\trauo{F}{2n-1}{s}}
		= \suo{1}{\nuu{\trauo{F}{2n-1}{s}}}
	\eg
	and, taking \eqref{N7-1} into account, also that formula \eqref{N6-1-1} is true for all $B\in\BAR$ and that
	\bgl{R1640.6}
		\nuu{\trauo{F}{2n-1}{s}}
		\in\MggquR{1}.
	\eg
	In view of \eqref{R1640.5}, we have $\betau{\trauo{F}{2n-1}{s}}=\Oqq$. Combining this with \eqref{R1640.2} and \eqref{R1640.6}, we see from \rrema{R1433} that $\trauo{F}{2n-1}{s}\in\RFuq{-1}$ and, in view of \eqref{Ruu}, especially $\trauo{F}{2n-1}{s}\in\NFuq{-1}$.
	Taking \eqref{muua} and \eqref{N6-1-1} into account and using~\cite{112}*{\cprop{B.5}}, it follows
	\[
	\mu_{\trauo{F}{2n-1}{s}}(B)
		= \int_B \frac{t^{2}+1}{\abs{t}+1} \nuu{\trauo{F}{2n-1}{s}}(\dif t)
		= \int_B \frac{t^{2n}}{\abs{t}+1} \sigmaF(\dif t)
	\]
	for all $B\in \BAR$. Thus, \eqref{N6-1} is verified and the proof of \rpart{R1640.a} is complete.

	\eqref{R1640.b} In view of $F\in\Rqos{\kappa}$, we get from \eqref{Ruoa} that
	\bgl{R1640.3}
		\sigmaF
		\in\MgguqR{\kappa}
	\eg
	and
	\begin{align}\label{R1640.4}
		\suo{j}{\sigmaF}&=\su{j}&\text{for each }j&\in\mn{0}{\kappa}.
	\end{align}
	Now we assume that $n\in \NO$ is such that $2n\leq\kappa$. Then \eqref{R1640.3} implies $\sigmaF \in \MggquR{2n}$. Therefore, from~\cite{112}*{\cprop{B.5}} we see that the mapping $\sigma\colon \BAR \to \Cqq$ given by
	\begin{equation} \label{N8-2}
		\sigma(B)\defg \int_B t^{2n}\sigmaF(\dif t)
	\end{equation}
	is a well-defined non-negative Hermitian measure belonging to $\MggqR$ and that
	\begin{equation} \label{N8-1}
		\int_B \frac{1}{t-z} \sigma(\dif t) = \int_B \frac{t^{2n}}{t-z} \sigmaF(\dif t)
	\end{equation}
	for all $z \in \ohe$. Combining \eqref{N6} and \eqref{N8-1}, we obtain
	\[
		\trauoa{F}{2n-1}{s}{z}
		=\Oqq + \int_\R \frac{1}{t-z} \sigma(\dif t)
	\]
	for all $z \in \ohe$, which, in view of \rtheo{T1446}, implies 
	\begin{align} \label{N8}
		\trauo{F}{2n-1}{s}&\in \RqK{0}&
		&\text{and}&
		(\gammau{\trauo{F}{2n-1}{s}} , \sigmau{\trauo{F}{2n-1}{s}})
		=( \Oqq, \sigma).
	\end{align}
	Consequently, from \eqref{N8-2} we get \eqref{N6-3} for all $B \in \BAR$. For all $m \in \mn{0}{\kappa-2n}$, we have $m+2n \leq \kappa$, hence \eqref{R1640.3} and \rrema{R1523} imply $\sigmaF \in \MggquR{m+2n}$, and, because of \eqref{N8-2} and~\cite{112}*{\cprop{B.5}}, furthermore, we infer $\sigma \in \MggquR{m}$ and 
	\[
		\int_\R t^{m+2n} \sigmaF(\dif t)
		= \int_\R t^m t^{2n} \sigmaF(\dif t)
		=\int_\R t^m\sigma(\dif t).
	\]
	Thus, $\sigma \in \MggquR{\kappa-2n}$. In view of \eqref{N8}, this means $\sigmau{\trauo{F}{2n-1}{s}} \in \MggquR{\kappa-2n}$. Thus, since \eqref{N8} implies $\gammau{\trauo{F}{2n-1}{s}}=\Oqq$, from \eqref{Ruu}, we get  $\trauo{F}{2n-1}{s} \in \Rkq {\kappa -2n}$. For all $m \in \mn{0}{\kappa-2n}$, we have $m+2n \leq \kappa$, so that \eqref{N8}, \eqref{N8-2},~\cite{112}*{\cprop{B.5}}, \eqref{suo}, and \eqref{R1640.4} imply
	\[
		\int_\R t^{m} \sigmau{\trauo{F}{2n-1}{s}}(\dif t)
		= \int_\R t^m \sigma(\dif t)
		=\int_\R t^m t^{2n} \sigmaF(\dif t)
		=\suo{2n+m}{\sigmaF}
		=\su{2n+m}.
	\]
	Consequently, $\sigmau{\trauo{F}{2n-1}{s}} \in \MggqRag{\seq{\su{2n+j}}{j}{0}{\kappa-2n}}$. Because of \eqref{Ruua}, then we see that $\trauo{F}{2n-1}{s}$ belongs to $\RkqSjj{\kappa-2n}{2n+j}{0}$.
\eproof 

It should be mentioned that in the scalar case the membership of $\trauo{F}{2n-1}{s}$ to $\RFuq{-1}$, which is contained in \rpart{R1640.a} of \rprop{R1640}, was already obtained in~\cite{MR1451805}*{\ctheo{3.2}}.

The following results complement the theme of \rprop{R1640}. They will play an important role in the proof of \rtheo{P0611-2}.

\bpropl{R1340}
	Let $n\in\NO$ and let $\seq{\su{j}}{j}{-1}{2n+1}$ be a sequence of complex \tqqa{matrices}. Furthermore, let $F\in \RqKJ{2n}{-1}$ and let $\sigmaF$ be given via \eqref{sigmaua}. Then:
	\benui
		\il{R1340.a} For all $z\in\ohe$, the matrix-valued function ${\trauo{F}{2n+1}{s}}$ can be represented via
		\begin{equation} \label{N9-1}
			\trauoa{F}{2n+1}{s}{z}
			=\su{2n+1}+z\int_\R\frac{t^{2n+1}}{t-z}\sigmaF(\dif t).
		\end{equation} 
		\il{R1340.b} Suppose $\su{2n+1}^\ad=\su{2n+1}$. Then $\trauo{F}{2n+1}{s} \in \RqP$, and the \tNp{} $(\alphau{\trauo{F}{2n+1}{s}},\betau{\trauo{F}{2n+1}{s}},\nuu{\trauo{F}{2n+1}{s}})$ of $\trauo{F}{2n+1}{s}$ is given by
		\[
			\alphau{\trauo{F}{2n+1}{s}}
			=\su{2n+1}-\int_\R\frac{t^{2n+1}}{t^2+1}\sigmaFa{\dif t},
		\]
		$\betau{\trauo{F}{2n+1}{s}}=\Oqq$, and
		\begin{align*}
			\nuua{\trauo{F}{2n+1}{s}}{B}
			&=\int_B\frac{t^{2n+2}}{t^2+1}\sigmaFa{\dif t}&\text{for each }B&\in\BAR.
		\end{align*}
		If $\trauo{F}{2n+1}{s}$ belongs to $\RqK{-1}$, then $F\in \RqK{2n+1}$ and
		\[
			\suo{2n+1}{\sigmaF}
			=\su{2n+1}-\gammau{\trauo{F}{2n+1}{s}},
		\]
		where $\gammau{\trauo{F}{2n+1}{s}}$ is given by \eqref{gammau}.
	\eenui
\eprop
\bproof
	\eqref{R1340.a} Formula \eqref{N9-1} immediately follows from \rrema{R1050}.

	\eqref{R1340.b} Since $F$ belongs to $\RqKJ{2n}{-1}$, in view of \eqref{Ruoa}, we have 
	\begin{align} \label{N9-2}
			F&\in \RqK{2n},&
			\gammaF&= -\su{-1},&
			&\text{and}&
			\sigmaF&\in \MggqRAg{\seq{\su{j}}{j}{0}{2n}}.
	\end{align} 
	Because of $0 \leq \frac{t^{2n+2}}{t^2+1}\leq t^{2n}$ for all $t \in \R$ and the third relation in \eqref{N9-2}, from~\cite{112}*{\cprop{B.5}} we then get that $\nu\colon\BAR\to\Cqq$ given by 
	\begin{equation} \label{N9-3}
		\nu(B)
		\defg\int_B\frac{t^{2n+2}}{t^2+1}\sigmaFa{\dif t}
	\end{equation} 
is a well-defined non-negative Hermitian measure belonging to $\MggqR$ for which the identity
	\begin{align} \label{N9}
		 \int_\R \frac{1+tz}{t-z}\nu(\dif t)
		 &=\int_\R\frac{1+tz}{t-z}\frac{t^{2n+2}}{t^2+1}\sigmaFa{\dif t},&z&\in\ohe,
	\end{align} 
	holds true. If $n \geq 1$, then $\abs{\frac{t^{2n+1}}{t^2+1}} \leq \abs{t^{2n-1}}$ is fulfilled for all $t \in \R$. If $n=0$, then $\abs{\frac{t^{2n+1}}{t^2+1}} \leq 1$  for all $t \in \R$. Thus, we see from~\cite{112}*{\clemm{B.1}} and \eqref{N9-2} that the integral $\int_\R \frac{t^{2n+1}}{t^2+1}\sigmaFa{\dif t}$ exists. For every choice of $z $ in $\ohe$, from \eqref{N9-1} and \eqref{N9} we conclude
	\[
		\begin{split}
			\trauoa{F}{2n+1}{s}{z}-\su{2n+1}
			&= z  \int_\R\frac{t^{2n+1}}{t-z} \sigmaFa{\dif t} +\int_\R \frac{t^{2n+1}}{t^2+1}\sigmaFa{\dif t} - \int_\R \frac{t^{2n+1}}{t^2+1}\sigmaFa{\dif t} \\
			&= \int_\R \lrk\frac{zt^{2n+1}}{t-z} +\frac{t^{2n+1}}{t^2+1} \rrk\sigmaFa{\dif t} - \int_\R \frac{t^{2n+1}}{t^2+1}\sigmaFa{\dif t} \\
			&= \int_\R\frac{1+tz}{t-z}\frac{t^{2n+2}}{t^2+1}\sigmaFa{\dif t} -\int_\R \frac{t^{2n+1}}{t^2+1}\sigmaFa{\dif t} \\
			&= \int_\R\frac{1+tz}{t-z} \nu(\dif t) -\int_\R \frac{t^{2n+1}}{t^2+1}\sigmaFa{\dif t}
		\end{split}
	\]
and, consequently,
	\begin{equation} \label{N10-1}
			\trauoa{F}{2n+1}{s}{z}
			=\su{2n+1} -\int_\R \frac{t^{2n+1}}{t^2+1}\sigmaFa{\dif t}+z\cdot\Oqq +\int_\R\frac{1+tz}{t-z} \nu(\dif t).
	\end{equation} 
Thanks to $\su{2n+1}^\ad=\su{2n+1}$ and~\cite{112}*{\crema{B.4}}, the matrix $\su{2n+1}-\int_\R\frac{t^{2n+1}}{t^2+1}\sigmaFa{\dif t}$ is Hermitian. Thus, in view of \eqref{N10-1}, applying \rtheo{T1554} yields $\trauo{F}{2n+1}{s} \in \RqP$ and  
	\[
		(\alphau{\trauo{F}{2n+1}{s}},\betau{\trauo{F}{2n+1}{s}},\nuu{\trauo{F}{2n+1}{s}})
		=\lrk\su{2n+1}-\int_\R\frac{t^{2n+1}}{t^2+1}\sigmaFa{\dif t}, \Oqq, \nu\rrk.
	\]
	Now suppose that
	\bgl{R1340.1}
		\trauo{F}{2n+1}{s}
		\in\RqK{-1}.
	\eg
	In view of \eqref{N9-2} and the definition \eqref{Ruo} of the class $\RqK{2n}$, we also have $\trauo{F}{2n+1}{s} \in \RqK{0}$. Furthermore, \rrema{R1425} yields $\nuu{\trauo{F}{2n+1}{s}} \in \MgguqR{1}$. Thus, in view of $\nuu{\trauo{F}{2n+1}{s}}=\nu$, the integral  $\int_\R t \nu(\dif t)$ exists. Because of \eqref{N9-3} and~\cite{112}*{\cprop{B.5}}, then we see that the integral $\int_\R \frac{t^{2n+3}}{t^2+1}\sigmaFa{\dif t}$ exists. Hence, since $\abs{t^{2n+1}} \leq \abs{\frac{t^{2n+3}}{t^2+1}}$ holds for all $t \in \R$, from~\cite{112}*{\clemm{B.1}} we get that  $\sigmaF \in \MgguqR{2n+1}$. In view of \eqref{Ruo} this shows that $F$ belongs to $ \RqK{2n+1}$. Because of \eqref{R1340.1}, we infer from \rprop{P1401} that 
 	\begin{equation} \label{N10-2}
 			\lim_{k\to+\infty}\re \trauoa{F}{2n+1}{s}{\I k}
 			=\gammau{\trauo{F}{2n+1}{s}}.
	\end{equation} 
Using \eqref{N9-1}, $\su{2n+1}^\ad=\su{2n+1}$, and~\cite{112}*{\crema{B.4}}, for all $k \in \N$, we conclude
 	\begin{equation} \label{N10-3}
 			\re \trauoa{F}{2n+1}{s}{\I k}
 			= \su{2n+1}+ \int_\R \re \lrk\frac{\I kt^{2n+1}}{t-\I k}\rrk\sigmaFa{\dif t}
 			= \su{2n+1}- \int_\R \frac{k^2t^{2n+1}}{t^2+k^2}\sigmaFa{\dif t}.
	\end{equation} 
Thus, from \eqref{N10-2} and \eqref{N10-3} we get
  	\begin{equation} \label{N11-1}
 			 \su{2n+1}-\gammau{\trauo{F}{2n+1}{s}}
 			 = 	\lim_{k\to+\infty} \int_\R \frac{k^2t^{2n+1}}{t^2+k^2}\sigmaFa{\dif t}.
	\end{equation} 
For all $t \in \R$, we have 
  	\begin{equation} \label{N11-2}
 			 \lim_{k\to+\infty} \frac{k^2t^{2n+1}}{t^2+k^2}
 			 =t^{2n+1}.
	\end{equation} 
Since $\abs{\frac{k^2t^{2n+1}}{t^2+k^2}}\leq \abs {t^{2n+1}}$ holds for all $k \in \N$ and all $t \in \R$, from  $\sigmaF \in \MgguqR{2n+1}$, \eqref{N11-2},~\cite{112}*{\clemm{B.1}} and Lebesgue's dominated convergence theorem then 
  	\begin{equation} \label{N11-3}
 			 \lim_{k\to+\infty} \int_\R \frac{k^2t^{2n+1}}{t^2+k^2}\sigmaFa{\dif t}
 			 =\int_\R t^{2n+1}\sigmaFa{\dif t}
 			 =\suo{2n+1}{\sigmaF}
	\end{equation} 
follows. Combining \eqref{N11-1} and \eqref{N11-3} yields $\suo{2n+1}{\sigmaF} =\su{2n+1}-\gammau{\trauo{F}{2n+1}{s}}$.
\eproof

\section{On Hamburger-Nevanlinna Type Results for $\NFuqa{\kappa}{\seq{\su{j}}{j}{-1}{\kappa}}$}\label{S1510}
Let $n\in\NO$ and let $\seq{\su{j}}{j}{0}{2n}$ be a sequence from $\Cqq$. Then the moment problem \mproblem{\R}{2n}{=} can be reformulated as a problem of a prescribed asymptotic expansion for functions in $\RFtOq$. This is a consequence of a matricial version of a classical result due to Hamburger and Nevanlinna. This matricial version can be found in~\cite{MR703593}*{p.~47} and~\cite{MR1624548}*{\clemm{2.1}}, where it was stated without proof. It can be proved along the lines of the proof of the scalar result which was given in~\cite{MR0184042}*{Ch.~3, Sect.~2}. Before formulating the result, we introduce some notation. For all $r\in (0,+\infty)$ and each $\delta \in(0,\frac{\pi}{2}]$, let  
\bgl{Sigmauu}
	\ohesecuu{r}{\delta}
	\defg\SetaA{z \in \C}{ \abs{z}\geq r \text{ and } \delta \leq \arg r \leq \pi -\delta}.
\eg\index{Sigma@$\ohesecuu{r}{\delta}$}

Taking \eqref{Ruua}, \rcoro{C1403} and \rtheo{T0933} into account, now we can reformulate the matricial version of the Hamburger-Nevanlinna theorem:

\btheol{T1515}
	Let $n\in\NO$ and let $\seq{\su{j}}{j}{0}{2n}$ be a sequence of complex \tqqa{matrices}.
	\benui
		\il{T1515.a} Let $F\in\RFuqa{2n}{\seq{\su{j}}{j}{0}{2n}}$. Then
		\[
			\lim_{r\to+\infty} \sup_{z \in \ohesecuu{r}{\delta}} \normA{z^{2n+1}\lek F(z)+\sum_{j=0}^{2n}\frac{1}{z^{j+1}}\su{j}\rek}
			=0.
		\]
		\il{T1515.b} Let $\seq{\su{j}}{j}{0}{2n}$ be a sequence from $\CHq$, and let $F\in\NFq$ be such that
		\[
			\lim_{y\to+\infty}\normA{(\I y)^{2n+1}\lek F(\I y)+\sum_{j=0}^{2n}\frac{1}{(\I y)^{j+1}}\su{j}\rek}
			=0.
		\]
		Then $F\in\RFuqa{2n}{\seq{\su{j}}{j}{0}{2n}}$.
	\eenui
\etheo

\rPart{T1515.b} of \rtheo{T1515} will be often applied in the following. It contains a sufficient condition which implies that a function $F\in\NFq$ is the \tStit{} of a solution $\sigma$ of Problem~\mproblem{\R}{2n}{=}.

The main goal of this section is to find appropriate generalizations of \rtheo{T1515} for the class $\NFuqa{\kappa}{\seq{\su{j}}{j}{-1}{\kappa}}$ with arbitrary $\kappa\in\NO\cup\set{+\infty}$. In the scalar case, this theme was treated by~\cite{MR1451805}. As we will see soon, that similar as in~\cite{MR1451805}, the essential tool in our strategy is the use of the construction introduced in \rrema{R1013}. In this connection, it should be mentioned that, in view of \rrema{R1404}, in the case of an affirmative answer to the generalization of \rpart{T1515.b} of \rtheo{T1515}, we would also obtain a sufficient condition for a function $F\in\NFq$ to be the \tStit{} of a solution $\sigma$ of Problem~\mproblem{\R}{2n+1}{=}. The following result, which in the scalar case goes back to~\cite{MR1451805}*{\ctheo{3.2}}, meets our above formulated goal concerning \rpart{T1515.a} of \rtheo{T1515}.


\btheol{P1518}
	Let $\kappa\in\NO\cup\set{+\infty}$, let  $\seq{\su{j}}{j}{-1}{\kappa}$ be a sequence  of complex \tqqa{matrices}, and let $F\in \RqkJ{-1}$. For all $k\in\mn{-1}{\kappa}$ and each $\delta \in (0,\frac{\pi}{2}]$, then
	\begin{equation} \label{N12-1}
		\lim_{r\to+\infty} \sup_{z \in \ohesecuu{r}{\delta}} \normA{ \trauoa{F}{k}{s}{z}}
		=0
	\end{equation}
	holds true, where $\ohesecuu{r}{\delta}$ is given by \eqref{Sigmauu}.
\etheo
\bproof
	The strategy of our proof is inspired by the proof which was given in~\cite{MR0184042}*{Ch.~3, Sect.~2} for the scalar case $q=1$ of \rpart{T1515.a} of \rtheo{T1515}. Since the function $F$ belongs to $\RqkJ{-1}$, we have  
  \begin{align} \label{N12-2}
		F&\in \RqK{\kappa},&
		\gammaF&=-\su{-1},&
		&\text{and}&
		\sigmaF&\in \MggqRAg{\seq{\su{j}}{j}{0}{\kappa}}.
	\end{align} 
	From \rprop{L1605} we know  that \eqref{N4-3} holds true for all $z \in \ohe$. Now we let $k\in\mn{-1}{\kappa}$, let $\delta \in (0,\frac{\pi}{2}]$, and let $\epsilon \in (0, +\infty)$. We consider an arbitrary $u \in \Cq$. In view of~\cite{112}*{\clemm{B.3}}, then $\rho_u \defg u^\ad\sigmaF u$ is a finite measure on $(\R, \BAR)$, which belongs to $\MgguqR{\kappa}$. Using additionally \eqref{N4-3}, for all $z \in \ohe$, we then obtain
  \begin{equation} \label{N12-53}
		\Abs{u^\ad\trauoa{F}{k}{s}{z}u}
		=\Abs{\int_\R \frac{t^{k+1}}{t-z}\rho_u(\dif t)}
		\leq \int_\R \frac{\abs{t^{k+1}}}{\abs{t-z}}\rho_u (\dif t).
	\end{equation} 
	First we now consider the case $k\in\mn{0}{\kappa}$. Then the mapping $\mu\colon\BAR \to \C$ given by 
  \begin{equation} \label{N12-9}
		\mu(B)
		=\int_B \abs{t^k}\rho_u(\dif t)
	\end{equation} 
	is a well-defined finite measure, i.\,e., $\mu$ belongs to $\Mggoaa{1}{\R}$. Obviously,
  \begin{equation} \label{N12-3}
 		\lim_{n\to+\infty} \mu\lrk\R \setminus [-n,n]\rrk
 		=0.
	\end{equation} 
	Thanks to \eqref{N12-3}, there is an $N \in \N$ such that
  \begin{equation} \label{N12-4}
 	 \mu\lrk\R \setminus [-N,N]\rrk
 	 < \frac{\epsilon}{2} \sin \delta.
	\end{equation} 
	Clearly, if we set
	\[
		R
		\defg \frac{2 N^{k+1} \rho_u\lrk[-N,N]\rrk}{\epsilon \sin \delta}  +1,
	\]
	then $R$ belongs to $[1,+\infty)$ and we have 
  \begin{equation} \label{N12-22}
 		\frac{N^{k+1} \rho_u\lrk[-N,N]\rrk}{ R \sin \delta}
 		< \frac{\epsilon}{2}.
	\end{equation} 
	We consider an arbitrary $r \in [R, +\infty)$ and an arbitrary $z \in \ohesecuu{r}{\delta}$. Then 
  \begin{align} \label{N12-5}
 		\abs{z}&\in [r, + \infty)&
 		&\text{and}&
 		\arg z&\in [\delta, \pi-\delta].
	\end{align} 
	For all $t \in \R$, we get furthermore
  \begin{equation} \label{N12-31}
		\abs{t-z}
		\geq \Abs{\im (t-z)}
		=\abs{\im z}
		= \abs{z}\Abs{\sin(\arg z)}
		\geq \abs{z}\sin \delta.
	\end{equation} 
	For all $t \in [-N,N]$, we have	
  \begin{equation} \label{N12-26}
		0
		\leq \abs{t^{k+1}}
		\leq N^{k+1}
	\end{equation} 
	and, because of \eqref{N12-31}, we conclude
  \begin{equation} \label{N12-32}
		\abs{t-z}
		\geq \abs{z}\sin \delta
		\geq r \sin \delta
		\geq R \sin \delta.
	\end{equation} 
	If $t \in \R \setminus [-n,n]$, then
  \bsp
		\abs{t-z}
		\geq \abs{\e^{-\I  \arg z}}\Abs{t - \abs{z}\e^{\I  \arg z}} 
		&= \Abs{t \e^{-\I  \arg z}- \abs{z}}
		\geq \Abs{ \im \lrk t \e^{\I  \arg z }- \abs{z}\rrk} \\
		&= \abs{t}\Abs{\sin (\arg z)}
		\geq \abs{t} \sin \delta
		>0
	\esp
	and, consequently,
  \begin{equation} \label{N12-36}
		\frac{\abs{t^k}}{\sin \delta}
		= \frac{\abs{t^{k+1}}}{\abs{t}\sin \delta} 
		\geq \frac{\abs{t^{k+1}}}{\abs{t -z}}.
	\end{equation} 
	Using \eqref{N12-53}, \eqref{N12-26}, \eqref{N12-32}, \eqref{N12-36}, \eqref{N12-9}, \eqref{N12-22}, and \eqref{N12-4}, we get then
  \begin{equation} \label{N12-39}
  	\begin{split}
			\Abs{u^\ad\trauoa{F}{k}{s}{z}u}
			&\leq \int_{[-N,N]} \frac{\abs{t^{k+1}}}{\abs{t-z}}\rho_u(\dif t)+\int_{\R \setminus[-N,N]} \frac{\abs{t^{k+1}}}{\abs{t-z}}\rho_u(\dif t)   \\
			&\leq \int_{[-N,N]} \frac{N^{k+1}}{R \sin \delta }\rho_u(\dif t)  
			+\int_{\R \setminus[-N,N]} \frac{\abs{t^{k}}}{\sin \delta}\rho_u(\dif t)  \\
			&= \frac{N^{k+1}}{R \sin \delta }\rho_u\lrk[-N,N]\rrk  
			+ \frac{1}{\sin \delta}\mu\lrk\R \setminus[-N,N]\rrk
			< \epsilon.
  	\end{split}
	\end{equation} 
	Now we consider the case $k=-1$. Then 
  \begin{equation} \label{N12-131}
		R
		\defg \frac{\rho_u(\R)}{\epsilon \sin \delta } +1
	\end{equation} 
	belongs to $[1, +\infty)$. Let $r \in [R, +\infty)$ and  $z \in \ohesecuu{r}{\delta}$. Hence the relations in \eqref{N12-5} are true and  the inequalities in \eqref{N12-31} follow again. Consequently, \eqref{N12-32} is fulfilled. Because of  \eqref{N12-53}, $k=-1$, \eqref{N12-32}, and \eqref{N12-131}, we get then 
	\bspl{P1518.1}
		\Abs{u^\ad\trauoa{F}{k}{s}{z}u}
		\leq \int_{\R} \frac{1}{\abs{t-z}}\rho_u(\dif t)  
		\leq \int_{\R} \frac{1}{R \sin \delta }\rho_u(\dif t)  
		= \frac{\rho_u(\R)}{R \sin \delta }  
		< \epsilon.
	\esp
	Thus, in view of \eqref{N12-39} and \eqref{P1518.1}, we proved that, for all $k\in\mn{-1}{\kappa}$, for all $\delta \in (0,\frac{\pi}{2}]$, for all $\epsilon \in (0, +\infty)$, and each $u \in \Cq$, there is an $R \in [0, +\infty)$ such that $\abs{u^\ad\trauoa{F}{k}{s}{z}u} < \epsilon$ for all $r \in [R, + \infty)$ and each $z \in \ohesecuu{r}{\delta}$. In other words, for all $k\in\mn{-1}{\kappa}$, each $\delta \in (0,\frac{\pi}{2}]$, and each $u \in \Cq$, we have
	\[
 		\lim_{r\to+\infty} \sup_{z \in \ohesecuu{r}{\delta}}\Abs{u^\ad\trauoa{F}{k}{s}{z}u}
 		=0.
 	\]
	A standard argument of linear algebra (see, e.\,g.~\cite{MR1152328}*{\crema{1.1.1}}) yields then \eqref{N12-1}.
\eproof

\bcorol{C1014*}
	Let $\kappa\in\NO\cup\set{+\infty}$, let $\seq{\su{j}}{j}{-1}{\kappa}$ be a sequence of complex \tqqa{matrices}, and let $F\in\RqkJ{-1}$. For all $k\in\mn{-1}{\kappa}$, then
	\begin{equation} \label{N15-1}
		\lim_{y\to+\infty} \trauoa{F}{k}{s}{\I y}
		=\Oqq
	\end{equation} 
	and
	\bgl{N15-2*}
		-\su{k}
		=
		\begin{cases}
			\lim_{y\to+\infty}F(\I y)\incase{k=-1}\\
			\lim_{y\to+\infty}(\I y)^{k+1}[F(\I y)+\sum_{j=0}^k(\I y)^{-j}\su{j-1}]\incase{k\geq0}
		\end{cases}.
	\eg
\ecoro
\bproof
	Because of $F \in \RqkJ{-1}$, we have \eqref{N12-2}. In particular, $F$ belongs to $\NFuq{\kappa}$. Thus, \eqref{Ruo} gives $F\in\NFuq{0}$. Taking \eqref{FSG} into account, then $F\in\NFuq{-1}$ follows. Consequently, \rprop{P1401} yields \eqref{N15-3}. Let $k\in\mn{-1}{\kappa}$. The limit \eqref{N15-1} is an immediate consequence of \rtheo{P1518}. If $k=-1$, then \eqref{N12-2} and \eqref{N15-3} imply
	\[
		-\su{-1}
		=\gammaF
		=\lim_{y\to+\infty}F(\I y).
	\]
	Thus, \eqref{N15-2*} holds true for $k=-1$. If $k\in\mn{0}{\kappa}$, then \eqref{N15-1} and \eqref{Fuoa} imply
	\[
		-\su{k}
		=\lim_{y\to+\infty}\lek\trauoa{F}{k}{s}{\I y} -\su{k}\rek
		=\lim_{y\to+\infty}(\I y)^{k+1}\lek F(\I y)+\sum_{j=0}^k(\I y)^{-j}\su{j-1}\rek.\qedhere
	\]
\eproof

Now we are going to show that \rpart{T1515.a} of \rtheo{T1515} is an immediate consequence of \rtheo{P1518}.

\begin{proofof}{\rpart{T1515.a} of \rtheo{T1515}}
	Let $t_{-1}\defg\Oqq$ and let $t_j\defg\su{j}$ for all $j\in\mn{0}{2n}$. Then \rrema{R1404} yields $F\in\NFuqa{2n}{\seq{t_j}{j}{-1}{2n}}$. Thus, \rtheo{P1518} implies
	\bgl{T1515.1}
		\lim_{r\to+\infty}\sup_{z \in \ohesecuu{r}{\delta}}\normA{ \trauoa{F}{2n}{t}{z}}
		=0.
	\eg
	Because of $t_{-1}=\Oqq$ we see from \eqref{Fuoa} that
	\bgl{T1515.2}
		\trauoa{F}{2n}{t}{z}
		=z^{2n+1}\lek F(z)+\sum_{j=0}^{2n}\frac{1}{z^{j+1}}\su{j}\rek.
	\eg
	Now the combination of \eqref{T1515.1} and \eqref{T1515.2} completes the proof of \rpart{T1515.a} of \rtheo{T1515}.
\end{proofof}

Now we state a corresponding generalization of the second part of \rtheo{T1515}. It should be mentioned that in the scalar case the result goes back to~\cite{MR1451805}*{\ctheo{3.3}}.


\btheol{P0611-1}
	Let  $n \in \NO$,  let  $\seq{\su{j}}{j}{-1}{2n}$ be a sequence of Hermitian complex \tqqa{matrices},  and let $F\in \RqP$ be such that 
	\begin{equation} \label{N16-1}
		\lim_{y\to+\infty} \trauoa{F}{2n}{s}{\I y}
		=\Oqq.
	\end{equation} 
	Then $F$ belongs to $\RqKJ{2n}{-1}$.
\etheo
\bproof
	Let $G\colon \ohe \to \Cqq$ defined by $G(z)\defg \su{-1}$, and let $S\defg F+G$. Since $\su{-1}^\ad=\su{-1}$ holds, we see from  \rtheo{T1554} and~\cite{112}*{\crema{3.4}} that $G$ and $S$ belong to $\RqP$. Using \eqref{N16-1} and \eqref{Fuoa}, we obtain then 
  \bsp
 		\Oqq
 		&=\lim_{y\to+\infty}(\I y)^{2n+1}\lek F(\I y)+\sum_{j=0}^{2n+1}(\I y)^{-j}\su{j-1}\rek \\
 		&=\lim_{y\to+\infty}(\I y)^{2n+1}\lek S(\I y)+\sum_{k=0}^{2n}(\I y)^{-k-1}\su{k}\rek.
 	\esp
	Hence, \rpart{T1515.b} of \rtheo{T1515} shows that
	\bgl{P0611-1.1}
		S
		\in\RFuqA{2n}{\seq{\su{j}}{j}{0}{2n}}.
	\eg
	Thus, taking \eqref{P0611-1.1}, \eqref{Ruo}, and \eqref{Ruu} into account we see that $S$ belongs to $\RqK{0}$ and that $\gammau{S}=\Oqq$ holds. Since \rtheo{T1446} shows that $-G \in \RqK{0}$, that $\gammau{-G}=-\su{-1}$, and that $\sigmau{-G}$ is the zero measure in $\MggqR$, we see from \eqref{FSG} that $S$ and $-G$ both belong to $\RqK{-1} \cap \RqI$. Thus,~\cite{112}*{\crema{4.4}} yields $F \in \RqI$ and $\sigmaF =\sigmau{S}$. In particular, from \eqref{P0611-1.1} we infer then $\sigmaF \in \MggqRag{\seq{\su{j}}{j}{0}{2n}}$  and  $F \in \RqK{2n}$. Furthermore,~\cite{112}*{\crema{5.7}} provides us $\gammaF = \gammau{S} +\gammau{-G}=-\su{-1}$. Consequently, $F$ belongs to $\RqKJ{2n}{-1}$.
\eproof

Our next aim can be described as follows. Let $k\in\NO$ and let $\seq{\su{j}}{j}{-1}{k}$ be a sequence of Hermitian complex \tqqa{matrices}. Then we are looking for appropriate descriptions of the set $\NFuqa{k}{\seq{\su{j}}{j}{-1}{k}}$. First we consider the case of an even number $k$.

\bprop
	Let $n\in\NO$ and let  $\seq{\su{j}}{j}{-1}{2n}$ be a sequence of Hermitian complex \tqqa{matrices}. Then 
	\[
		\RqoS{2n}
		=\SetaA{F\in \RqP }{\lim_{y\to+\infty}\trauoa{F}{2n}{s}{\I y}=\Oqq}.
	\]
\eprop
\bproof
	Combine \eqref{N15-1} and  \rtheo{P0611-1}.
\eproof

Now we treat the case of a sequence $\seq{\su{j}}{j}{-1}{k}$ from $\CHq$ with odd number $k$.

The following result contains a useful sufficient condition which guarantees that a function $F\in\NFq$ belongs to the set $\NFuqa{2n+1}{\seq{\su{j}}{j}{-1}{2n+1}}$. In the scalar case, the result goes back to~\cite{MR1451805}*{\ctheo{3.3}}.

\btheol{P0611-2}
	Let  $n \in \NO$,  let  $\seq{\su{j}}{j}{-1}{2n+1}$ be a sequence of Hermitian complex \tqqa{matrices},  and let $F\in \RqP$ be such that 
  \begin{equation} \label{N17-1}
 		\lim_{y\to+\infty} \trauoa{F}{2n+1}{s}{\I y}
 		=\Oqq.
	\end{equation} 
	Then $\trauo{F}{2n+1}{s}$ belongs to $\RqP$. If $\trauo{F}{2n+1}{s}$ even belongs to $\RqK{-1}$, then the function $F$ belongs to $\RqKJ{2n+1}{-1}$.
\etheo
\bproof
	For all $y \in (0,+\infty)$, we get from formula \eqref{N4-12} in \rrema{R1013} that 
	\[
		\trauoa{F}{2n}{s}{\I y}
		=(\I y)^{-1}\lek\trauoa{F}{2n+1}{s}{\I y} -\su{2n+1}\rek.
	\]
	In view of \eqref{N17-1}, this implies \eqref{N16-1}. Thus, \rtheo{P0611-1} yields   $F \in \RqKJ{2n}{-1}$. \rPart{R1340.b} of \rprop{R1340} shows then that $F$ belongs to $\RqP$. Now we additionally suppose that $	\trauo{F}{2n+1}{s} \in \RqK{-1}$. Then \rpart{R1340.b} of \rprop{R1340} provides us $F \in \RqK{2n+1}$, in particular $\sigmaF \in \MgguqR{2n+1}$, and $\suo{2n+1}{\sigmaF}=\su{2n+1}-\gammau{\trauo{F}{2n+1}{s}}$. From  $F \in \RqKJ{2n}{-1}$ we see $\gammaF=-\su{-1}$ and that $\sigmaF \in \MggqRag{\seq{\su{j}}{j}{0}{2n}}$. Since  $\trauo{F}{2n+1}{s}$ belongs to $\RqK{-1}$, \rprop{P1401} and \eqref{N17-1} yield $\gammau{\trauo{F}{2n+1}{s}}=\lim_{y \to + \infty} \trauoa{F}{2n+1}{s}{\I y}=\Oqq$. Thus, we get $\suo{2n+1}{\sigmaF}=\su{2n+1}$. Consequently, $\sigmaF$ belongs to $\MggqRag{\seq{\su{j}}{j}{0}{2n+1}}$. Hence, $F \in \RqKJ{2n+1}{-1}$.
\eproof

\bprop
	Let $n\in\NO$ and let  $\seq{\su{j}}{j}{-1}{2n+1}$ be a sequence of Hermitian complex \tqqa{matrices}. Then 
	\begin{multline*}
		\RqoS{2n+1}\\
		=\SetaA{F\in \RqP}{\trauo{F}{2n+1}{s} \in \RqK{-1}\text{ and }\lim_{y\to+\infty}\trauoa{F}{2n+1}{s}{\I y}=\Oqq}.
	\end{multline*}
\eprop
\bproof
	Combine \rpart{R1640.a} of \rprop{R1640}, \eqref{N15-1}, and \rtheo{P0611-2}.
\eproof

\section{On a Schur Type Algorithm for Sequences of Complex \tpqa{Matrices}}\label{S1505}
In this section, we recall some essential facts on a Schur type algorithm for sequences from $\Cpq$, which was introduced and investigated in~\cite{103}. The elementary step of this algorithm is based on the use of the construction of the reciprocal sequence of a finite or infinite sequence from $\Cpq$. For this reason, we first remember the definition of the reciprocal sequence.

Let $\kappa \in \NO \cup \{+\infty\}$ and let $\seq{\su{j}}{j}{0}{\kappa}$ be a sequence of complex \tpqa{matrices}. Then the sequence $\seq{\su{j}^\rez }{j}{0}{\kappa}$\index{^#@$\seq{\su{j}^\rez }{j}{0}{\kappa}$} of complex \tqpa{matrices}, which is given by $\su{0}^\rez \defg \su{0}^\MP$ and, for all $k \in \mn{1}{\kappa}$, recursively by
\[
	\su{k}^\rez 
	\defg - \su{0}^\MP\sum_{j=0}^{k-1}\su{k-j}\su{j}^\rez ,
\]
is called the \emph{reciprocal sequence corresponding to  $\seq{\su{j}}{j}{0}{\kappa}$}. For a detailed treatment of the concept of reciprocal sequences, we refer the reader to~\cite{101}.

Now we explain the elementary step of the Schur type algorithm under consideration. Let $\kappa\in\mn{2}{+\infty}\cup\set{+\infty}$ and let $\seq{\su{j}}{j}{0}{\kappa}$ be a sequence from $\Cpq$ with reciprocal sequence $\seq{\su{j}^\rez}{j}{0}{\kappa}$. Then the sequence  $\seq{\su{j}^\seins }{j}{0}{\kappa-2}$\index{^(1)@$\seq{\su{j}^\seins }{j}{0}{\kappa-2}$} defined for all $j \in \mn{0}{\kappa-2}$ by 
\bgl{seins}
	\su{j}^\seins 
	\defg -\su{0}\su{j+2}^\rez \su{0}
\eg
is said to be the \emph{\tfirstSta{\seq{\su{j}}{j}{0}{\kappa}}}.

\bremal{R1038}
	Let $\kappa\in\mn{2}{+\infty}\cup\set{+\infty}$ and let $\seq{\su{j}}{j}{0}{\kappa}$ be a sequence of complex \tpqa{matrices}. Let $\seq{\su{j}^\seins}{j}{0}{\kappa-2}$ be the \tfirstSta{\seq{\su{j}}{j}{0}{\kappa}}. Let $m\in\mn{2}{\kappa}$. Then from \eqref{seins} it is obvious that $\seq{\su{j}^\seins}{j}{0}{m-2}$ is the \tfirstSta{\seq{\su{j}}{j}{0}{m}}.
\erema

The repeated application of the \tfirstSt{} generates in a natural way a corresponding algorithm for (finite or infinite) sequences of complex \tpqa{matrices}:

Let $\kappa \in \NO \cup \{+\infty\}$ and let $\seq{\su{j}}{j}{0}{\kappa}$ be a sequence of complex \tpqa{matrices}. Then we call the sequence $\seq{\su{j}^\sta{0}}{j}{0}{\kappa}$\index{(0)@$\seq{\su{j}^\sta{0}}{j}{0}{\kappa}$}, given by $\su{j}^\sta{0}\defg \su{j}$ for all $j \in \mn{0}{\kappa}$, the \emph{$0$\nobreakdash-th Schur transform of $\seq{\su{j}}{j}{0}{\kappa}$}. If $\kappa \geq 2$, then we define recursively, the $k$\nobreakdash-th Schur transform: For all $k \in \N$ with $2k \leq \kappa$, the first Schur transform $\seq{\su{j}^\sta{k}}{j}{0}{\kappa-2k}$ of $\seq{\su{j}^\sta{k-1}}{j}{0}{\kappa-2(k-1)}$ is called the \emph{$k$\nobreakdash-th Schur transform of $\seq{\su{j}}{j}{0}{\kappa}$}.

One of the central properties of the just introduced Schur type algorithm is that it preserves the \tHnnd{} extendability of sequences of matrices. This is the content of the following result, which is proved in~\cite{103}*{\cpropss{9.4}{9.5}}.

\bpropl{P1531}
	Let $\kappa\in \NO\cup\set{+\infty}$, let $\seq{\su{j}}{j}{0}{\kappa}\in\Hggequ{\kappa}$, and let $k\in\NO$ with $2k\leq\kappa$. Then the \taSt{k} $\seq{\su{j}^\sta{k}}{j}{0}{\kappa-2k}$ of $\seq{\su{j}}{j}{0}{\kappa}$ belongs to $ \Hggequ{\kappa-2k}$.
\eprop
%

In our considerations below, the special parametrization of block Hankel matrices introduced in~\cites{MR2570113,103} will play an essential role, the so-called \tcHp{}. For the convenience of the reader, we recall this notion. Let $\kappa\in\NO\cup\set{+\infty}$ and let $\seq{\su{j}}{j}{0}{\kappa}$ be a sequence in $\Cpq$. For every choice of non-negative integers $l$ and $m$ with $l\leq m\leq\kappa$, let
\begin{align*}
	\yuu{l}{m}&\defg
	\bMat
		\su{l}\\
		\su{l+1}\\
		\vdots\\
		\su{m}
	\eMat&
	&\text{and}&
	\zuu{l}{m}&\defg\brow\su{l},\su{l+1},\dotsc,\su{m}\erow.
\end{align*}
\index{y_,@$\yuu{l}{m}$}\index{z_,@$\zuu{l}{m}$}For all $n\in\NO$ with $2n\leq\kappa$, let
\[
	\Hu{n}
	\defg\matauuo{\su{j+k}}{j,k}{0}{n}
\]
\index{H_@$\Hu{n}$}and, for all $n\in\NO$ with $2n+1\leq\kappa$, let
\[
	\Ku{n}
	\defg\matauuo{\su{j+k+1}}{j,k}{0}{n}.
\]
\index{K_@$\Ku{n}$}Let
\begin{align*}
	\Mu{0}&\defg\Opq&
	&\text{and}&
	\Mu{n}&\defg\zuu{n}{2n-1}\Hu{n-1}^\MP\yuu{n+1}{2n}
\end{align*}
\index{M_@$\Mu{n}$}for all $n\in\N$ with $2n\leq\kappa$. Furthermore, let
\begin{align*}
	\Nu{0}&\defg\Opq&
	&\text{and}&
	\Nu{n}&\defg\zuu{n+1}{2n}\Hu{n-1}^\MP\yuu{n}{2n-1}
\end{align*}
\index{N_@$\Nu{n}$}for all $n\in\N$ with $2n\leq\kappa$. Let
\begin{align*}
	\Sigmau{0}&\defg\Opq&
	&\text{and}&
	\Sigmau{n}&\defg\zuu{n}{2n-1}\Hu{n-1}^\MP\Ku{n-1}\Hu{n-1}^\MP\yuu{n}{2n-1}
\end{align*}
\index{Sigma_@$\Sigmau{n}$}for all $n\in\N$ with $2n-1\leq\kappa$. For all $n\in\NO$ with $2n\leq\kappa$, let
\[
	\Lambdau{n}
	\defg\Mu{n}+\Nu{n}-\Sigmau{n}.
\]
\index{Lambda_@$\Lambdau{n}$}Let
\begin{align}\label{Lu}
	\Lu{0}&\defg\su{0}&
	&\text{and}&
	\Lu{n}&\defg\su{2n}-\zuu{n}{2n-1}\Hu{n-1}^\MP\yuu{n}{2n-1}
\end{align}
\index{L_@$\Lu{n}$}for all $n\in\N$ with $2n\leq\kappa$. 

%

\bdefil{D0944}
	Let $\kappa\in\N\cup\set{+\infty}$ and let $\seq{\su{j}}{j}{0}{2\kappa}$ be a sequence in $\Cpq$. Then the pair of sequences $[\seq{\kpcu{k}}{k}{1}{\kappa},\seq{\kpdu{k}}{k}{0}{\kappa}]$ given by $\kpcu{k}\defg\su{2k-1}-\Lambdau{k-1}$\index{C_@$\kpcu{k}$} for all $k\in\mn{1}{\kappa}$ and by $\kpdu{k}\defg\Lu{k}$\index{D_@$\kpdu{k}$} for all $k\in\mn{0}{\kappa}$ is called the \emph{\tcHpa{\seq{\su{j}}{j}{0}{2\kappa}}}.
\edefi

\brema
	Let $\kappa\in\N\cup\set{+\infty}$ and let $\seq{\kpcu{k}}{k}{1}{\kappa}$ and $\seq{\kpdu{k}}{k}{0}{\kappa}$ be sequences of complex \tpqa{matrices}. Then one can easily see that there is a unique sequence $\seq{\su{j}}{j}{0}{2\kappa}$ of complex \tpqa{matrices} such that $[\seq{\kpcu{k}}{k}{1}{\kappa},\seq{\kpdu{k}}{k}{0}{\kappa}]$ is the \tcHpa{\seq{\su{j}}{j}{0}{2\kappa}}, namely the sequence given by $\su{0}\defg\kpdu{0}$ and for each $k\in\mn{0}{\kappa}$, by $\su{2k-1}\defg\Lambdau{k-1}+\kpcu{k}$ and $\su{2k}\defg\zuu{k}{2k-1}\Hu{k-1}^\MP\yuu{k}{2k-1}+\kpdu{k}$.
\erema

In~\cites{MR2570113,103} several important classes of sequences of complex \tpqa{matrices} were characterized in terms of their \tcHp{}. From the view of this paper, the class $\Hggequ{\kappa}$ of Hankel non-negative definite extendable sequences is of extreme importance (see \rtheo{1-417}, \rprop{P1407}). In the case of a sequence $\seq{\su{j}}{j}{0}{2\kappa}\in\Hggequ{2\kappa}$, the \tcHp{} can be generated by the above constructed Schur type algorithm. This is the content of the following theorem.

\begin{thm}[\cite{103}*{\ctheo{9.15}}]\label{T1512}
	Let $\kappa\in\N\cup\set{+\infty}$ and let $\seq{\su{j}}{j}{0}{2\kappa}\in\Hggequ{2\kappa}$ with \tcHp{} $[\seq{\kpcu{k}}{k}{1}{\kappa},\seq{\kpdu{k}}{k}{0}{\kappa}]$. Then $\kpcu{k}=\su{1}^\sta{k-1}$ for all $k\in\mn{1}{\kappa}$ and $\kpdu{k}=\su{0}^\sta{k}$ for all $k\in\mn{0}{\kappa}$.
\end{thm}

An essential step in the further considerations of this paper can be described as follows. Let $\kappa\in\mn{2}{+\infty}$ and let $\seq{\su{j}}{j}{0}{\kappa}\in\Hggequ{\kappa}$. Denote by $\seq{\su{j}^\seins}{j}{0}{\kappa-2}$ the \tfirstSta{\seq{\su{j}}{j}{0}{\kappa}}. In view of \rprop{P1531}, we get then $\seq{\su{j}^\seins}{j}{0}{\kappa-2}\in\Hggequ{\kappa-2}$. Thus, \rpart{P1407.b} of \rprop{P1407} yields that both sets $\RFuqa{\kappa}{\seq{\su{j}}{j}{0}{\kappa}}$ and $\RFuqa{\kappa-2}{\seq{\su{j}^\seins}{j}{0}{\kappa-2}}$ are non-empty. Then a central aspect of our strategy is based on the construction of a special bijective mapping $\HTuu{\su{0}}{\su{1}}$ with the property
\[
	\HTuuA{\su{0}}{\su{1}}{\RFuqA{\kappa}{\seq{\su{j}}{j}{0}{\kappa}}}
	=\RFuqA{\kappa-2}{\seq{\su{j}^\seins}{j}{0}{\kappa-2}}.
\]
This mapping $\HTuu{\su{0}}{\su{1}}$ will help us to realize the basic step for our construction of a special Schur type algorithm in the class $\NFq$, which stands in a bijective correspondence to the above described Schur type algorithm for Hankel nonnegative definite extendable sequences. In order to verify several interrelations between the algebraic and function theoretic versions of our Schur type algorithm, we will need various properties of sequences of complex \tqqa{matrices} with prescribed Hankel properties. Now we give a short summary of this material, which is mostly taken from~\cite{103}.

\begin{lem}[\cite{103}*{\clemm{3.1}}]\label{L1037}
	Let $\kappa\in\NO\cup\set{+\infty}$ and let $\seq{\su{j}}{j}{0}{\kappa}\in\Hggequ{\kappa}$. Then:
	\benui
		\il{L1037.a} $\su{2k}\in\Cggq$ for all $k\in\mn{0}{\frac{\kappa}{2}}$.
		\il{L1037.b} $\su{l}^\ad=\su{l}$ for all $l\in\mn{0}{\kappa}$.
		\il{L1037.c} $\bigcup_{j=2k}^\kappa\Bilda{\su{j}}\subseteq\Bilda{\su{2k}}$ and $\Kerna{\su{2k}}\subseteq\bigcap_{j=2k}^\kappa\Kerna{\su{j}}$ for all $k\in\mn{0}{\frac{\kappa}{2}}$.
	\eenui
\end{lem}

\begin{lem}\label{L1031}
	Let $\kappa\in\N\cup\set{+\infty}$ and $\seq{\su{j}}{j}{0}{2\kappa}\in\Hggqu{2\kappa}$. Then:
	\benui
		\il{L1031.a} $\su{2k}\in\Cggq$ for all $k\in\mn{0}{\kappa}$.
		\item $\su{l}^\ad=\su{l}$ for all $l\in\mn{0}{2\kappa}$.
		\il{L1031.c} $\bigcup_{j=0}^{2\kappa-1}\Bilda{\su{j}}\subseteq\Bilda{\su{0}}$  and $\Kerna{\su{0}}\subseteq\bigcap_{j=0}^{2\kappa-1}\Kerna{\su{j}}$.
		\il{L1031.d} For all $k\in\mn{0}{\kappa}$, the matrix $\Lu{k}$ defined in \eqref{Lu} is non-negative Hermitian.
	\eenui
\end{lem}
\bproof
	The assertions of~\eqref{L1031.a}--\eqref{L1031.c} follow from~\cite{103}*{\clemm{3.2}}. From \eqref{Lu} and~\eqref{L1031.a} it follows $\Lu{0}=\su{0}\in\Cggq$. If $\kappa\in\N$, then a standard result on the structure of non-negative Hermitian block matrices (see e.\,g.~\cite{MR1152328}*{\clemm{1.1.9}}) implies $\Lu{k}\in\Cggq$ for $k\in\mn{1}{\kappa}$.
\eproof

Now we recall a class of sequences of complex matrices which, as the consideration in~\cite{103} have shown, turned out to be extremely important in the framework of the above introduced Schur algorithm.

\begin{defn}[\cite{101}*{\cdefi{4.3}}]\label{D1043}
	Let $\kappa\in\NO\cup\set{+\infty}$ and $\seq{\su{j}}{j}{0}{\kappa}$ be a sequence of complex \tpqa{matrices}. We then say that $\seq{\su{j}}{j}{0}{\kappa}$ is \emph{dominated by its first term} (or, simply, that it is \emph{\tftd{}}) if
	\begin{align*}
		\Kerna{\su{0}}&\subseteq\bigcap_{j=0}^{\kappa}\Kerna{\su{j}}&
		&\text{and}&
		\bigcup_{j=0}^{\kappa}\Bilda{\su{j}}&\subseteq\Bilda{\su{0}}.
	\end{align*}
	The set of all \tftd{} sequences $\seq{\su{j}}{j}{0}{\kappa}$ of complex \tpqa{matrices} will be denoted by $\Dpqu{\kappa}$.
\end{defn}

For a comprehensive investigation of \tftd{} sequences, we refer the reader to the paper~\cite{101}.

From the view of the Schur type algorithm for sequences of matrices, the following result (see also~\cite{103}*{\cprop{4.24}}) proved to be of central importance.

\blemml{L1048}
	For all $\kappa\in\NO\cup\set{+\infty}$, the inclusion $\Hggequ{\kappa}\subseteq\Dqqkappa$ holds true.
\elemm
\bproof
	Apply \rpart{L1037.c} of \rlemm{L1037}.
\eproof

Now we turn our attention to further important subclasses of the class of all Hankel non-negative definite sequences. Let $n\in\NO$ and let $\seq{\su{j}}{j}{0}{2n}$ be a sequence from $\Cqq$. Then $\seq{\su{j}}{j}{0}{2n}$ is called \emph{\tHpd{}} if the block Hankel matrix $\Hu{n}\defg\matauuo{\su{j+k}}{j,k}{0}{n}$\index{H_@$\Hu{n}$} is positive Hermitian. A sequence $\seq{\su{j}}{j}{0}{\infty}$ from $\Cqq$ is called \emph{\tHpd{}} if for all $n\in\NO$ the sequence $\seq{\su{j}}{j}{0}{2n}$ is \tHpd{}. For all $\kappa\in\NO\cup\set{+\infty}$, we will write $\Hgqu{2\kappa}$\index{H_,2^>@$\Hgqu{2\kappa}$} for the set of all \tHpd{} sequences $\seq{\su{j}}{j}{0}{2\kappa}$ from $\Cqq$.

\bremal{R1600}
	Let $n\in\NO$. Then~\cite{MR2570113}*{\crema{2.8}} shows that $\Hgqu{2n}\subseteq\Hggequ{2n}$.
\erema

\blemml{L1602}
	Let $\kappa\in\NO\cup\set{+\infty}$ and let $\seq{\su{j}}{j}{0}{2\kappa}\in\Hgqu{2\kappa}$. Then:
	\benui
		\item $\su{2k}\in\Cgq$ for all $k\in\mn{0}{\kappa}$.
		\item $\su{l}^\ad=\su{l}$ for all $l\in\mn{0}{2\kappa}$.
		\il{L1602.c} For all $k\in\mn{0}{\kappa}$, the matrix $\Lu{k}$ defined in \eqref{Lu} is positive Hermitian.
	\eenui
\elemm
\bproof
	All assertions are immediate consequences of results on non-negative Hermitian block matrices (see e.\,g.~\cite{MR1152328}*{\clemm{1.1.7}}).
\eproof

Let $n\in\NO$ and let $\seq{\su{j}}{j}{0}{2n}$ be a sequence from $\Cqq$. Then $\seq{\su{j}}{j}{0}{2n}$ is called \emph{\tHpde{}} if there exist matrices $\su{2n+1}$ and $\su{2n+2}$ from $\Cqq$ such that $\seq{\su{j}}{j}{0}{2(n+1)}\in\Hgqu{2(n+1)}$. The symbol $\Hgequ{2n}$\index{H_,2^>,e@$\Hgequ{2n}$} stands for the set of all \tHpde{} sequences $\seq{\su{j}}{j}{0}{2n}$ from $\Cqq$.

\begin{prop}[\cite{MR2570113}*{\cprop{2.24}}]\label{P1615}
	Let $n\in\NO$ and let $\seq{\su{j}}{j}{0}{2(n+1)}$ be a sequence from $\Cqq$. Then the following statements are equivalent:
	\baeqi{0}
		\item $\seq{\su{j}}{j}{0}{2(n+1)}\in\Hgqu{2(n+1)}$.
		\item $\seq{\su{j}}{j}{0}{2n}\in\Hgqu{2n}$, $\su{2n+1}\in\CHq$ and there exists a matrix $D\in\Cgq$ such that
		\bgl{NR38-1}
			\su{2n+2}
			=\zuu{n+1}{2n+1}\Hu{n}^\MP\yuu{n+1}{2n+1}+D.
		\eg
	\eaeqi
\end{prop}

\bcorol{C1619}
	Let $n\in\NO$. Then $\Hgequ{2n}=\Hgqu{2n}$.
\ecoro
\bproof
	Use \rprop{P1615}.
\eproof

It should be mentioned that \rrema{R1600} is also a consequence of \rcoro{C1619}.

Let $n\in\NO$ and let $\seq{\su{j}}{j}{0}{2n+1}$ be a sequence from $\Cqq$. Then $\seq{\su{j}}{j}{0}{2n+1}$ is called \emph{\tHpde{}} if there exists a matrix $\su{2n+2}\in\Cqq$ such that $\seq{\su{j}}{j}{0}{2(n+1)}\in\Hgqu{2(n+1)}$. The symbol $\Hgequ{2n+1}$\index{H_,2+1^>,e@$\Hgequ{2n+1}$} stands for the set of all \tHpde{} sequences $\seq{\su{j}}{j}{0}{2n+1}$ from $\Cqq$.

\bpropl{P1642}
	Let $n\in\NO$ and let $\seq{\su{j}}{j}{0}{2n+1}$ be a sequence from $\Cqq$.
	\benui
		\il{P1642.a} The following statements are equivalent:
		\baeqii{0}
			\item $\seq{\su{j}}{j}{0}{2n+1}\in\Hgequ{2n+1}$.
			\item $\seq{\su{j}}{j}{0}{2n}\in\Hgqu{2n}$ and $\su{2n+1}\in\CHq$.
		\eaeqii
		\item Let $\seq{\su{j}}{j}{0}{2n+1}\in\Hgequ{2n+1}$ and $\su{2n+2}\in\Cqq$. Then the following statements are equivalent:
		\baeqii{2}
			\item $\seq{\su{j}}{j}{0}{2(n+1)}\in\Hgqu{2(n+1)}$.
			\item There exists a matrix $D\in\Cgq$ such that \eqref{NR38-1}.
		\eaeqii
	\eenui
\eprop
\bproof
	All assertions follow immediately from \rprop{P1615}.
\eproof

\section{On a Coupled Pair of Schur Type Transforms}\label{S1515}

The main goal of this section is to prepare the elementary step of our Schur type algorithm for the class $\NFq$. We will be led to a situation which, roughly speaking, looks as follows: Let $A,B\in\Cpq$, let $\mG$ be a non-empty subset of $\C$, and let $F\colon\mG\to\Cpq$. Then the matrix-valued functions $\HTaooo{F}{1}{A}{B}\colon\mG\to\Cpq$\index{^(+;,)@$\HTaooo{F}{1}{A}{B}$} and $\HTiaooo{F}{1}{A}{B}\colon\mG\to\Cpq$\index{^(-;,)@$\HTiaooo{F}{1}{A}{B}$} which are defined by
\bgl{F+ooa}
	\HTaoooa{F}{1}{A}{B}{z}
	\defg-A\lrk z\Iq+\lek F(z)\rek^\MP A\rrk+B
\eg
and
\bgl{F-ooa}
	\HTiaoooa{F}{1}{A}{B}{z}
	\defg-A\lrk z\Iq+A^\MP\lek F(z)-B\rek\rrk^\MP,
\eg
respectively, will be central objects in our further considerations. The special case $B=\Opq$ will occupy a particular role. For this reason, we set
\begin{align}\label{F+-o}
	\HTaoo{F}{1}{A}&\defg\HTaooo{F}{1}{A}{\Opq},&
	\HTiaoo{F}{1}{A}&\defg\HTiaooo{F}{1}{A}{\Opq}.
\end{align}
\index{^(+;)@$\HTaoo{F}{1}{A}$}\index{^(-;)@$\HTiaoo{F}{1}{A}$}Against to the background of our later considerations, the matrix-valued functions $\HTaooo{F}{1}{A}{B}$ and $\HTiaooo{F}{1}{A}{B}$ are called the \emph{\taaSta{A}{B}{F}} and the \emph{\tiaaSta{A}{B}{F}}.

The generic case studied here concerns the situation where $p=q$, $A$ and $B$ are matrices from $\Cqq$ with later specified properties, $\mG=\ohe$, and $F\in\NFq$.

The use of the transforms introduced in \eqref{F-ooa} and \eqref{F+-o} was inspired by some considerations in the paper Chen/Hu~\cite{MR1624548}. In particular, we mention~\cite{MR1624548}*{\clemm{2.6} and its proof, fromula~(2.7)}. Before treating more general aspects we state some relevant concrete examples for the constructions given by fromulas \eqref{F+ooa} and \eqref{F-ooa}. First we illustrate the transformations given in \eqref{F+ooa} and \eqref{F-ooa} by some examples.

\bexaml{E1436}
	Let $A,B\in\Cqq$. In view of \eqref{F+ooa}, then:
	\benui
		\item Let $\alpha\in\Cqq$ and let $F\colon\ohe\to\Cqq$ be defined by $F(z)\defg\alpha$. Then $\HTaoooa{F}{1}{A}{B}{z}=B-A\alpha^\MP A+z(-A)$ for all $z\in\ohe$.
		\item Let $\beta\in\Cqq$ and let $F\colon\ohe\to\Cqq$ be defined by $F(z)\defg z\beta$. Then $\HTaoooa{F}{1}{A}{B}{z}=B+z(-A)-\frac{1}{z}A\beta^\MP A$ for all $z\in\ohe$.
		\item Let $M\in\Cqq$, let $\tau\in\R$, and let $F\colon\ohe\to\Cqq$ be defined by $F(z)\defg\frac{1+\tau z}{\tau-z}M$. For all $z\in\ohe$, then
	\[
		\HTaoooa{F}{1}{A}{B}{z}
		=
		\begin{cases}
			B+z(AM^\MP A-A)\incase{\tau=0}\\
			B+z(-A)+\frac{1+(-\frac{1}{\tau})z}{(-\frac{1}{\tau})-z}AM^\MP A\incase{\tau\neq0}
		\end{cases}.
	\]
	\eenui
\eexam

For all $\tau\in\R$, let $\Kronu{\tau}$ be the Dirac measure defined on $\BAR$ with unit mass at $\tau$. Furthermore, let $\zmq\colon\BAR\to\Cggq$ be defined by $\zmqa{B}\defg\Oqq$.

\bexam
	Let $A\in\Cqq$ with $-A\in\Cggq$ and let $B\in\CHq$. In view of \rtheo{T1554} and \rexam{E1436}, then:
	\benui
		\item Let $\alpha\in\CHq$ and let $F\colon\ohe\to\Cqq$ be defined by $F(z)\defg\alpha$. Then $F\in\NFq$ and $(\alphaF,\betaF,\nuF)=(\alpha,\Oqq,\zmq)$. Furthermore, $\HTaooo{F}{1}{A}{B}\in\NFq$ and $(\alphau{\HTaooo{F}{1}{A}{B}},\betau{\HTaooo{F}{1}{A}{B}},\nuu{\HTaooo{F}{1}{A}{B}})=(B-A\alpha^\MP A,-A,\zmq)$.
		\item Let $\beta\in\Cggq$ and let $F\colon\ohe\to\Cqq$ be defined by $F(z)\defg z\beta$. Then $F\in\NFq$ and $(\alphaF,\betaF,\nuF)=(\Oqq,\beta,\zmq)$. Furthermore, $\HTaooo{F}{1}{A}{B}\in\NFq$ and $(\alphau{\HTaooo{F}{1}{A}{B}},\betau{\HTaooo{F}{1}{A}{B}},\nuu{\HTaooo{F}{1}{A}{B}})=(B,-A,\Kronu{0}A\beta^\MP A)$.
		\item Let $M\in\Cggq$, let $\tau\in\R\setminus\set{0}$ and let $F\colon\ohe\to\Cqq$ be defined by $F(z)\defg\frac{1+\tau z}{\tau-z}M$. Then $F\in\NFq$ and $(\alphaF,\betaF,\nuF)=(\Oqq,\Oqq,\Kronu{\tau}M)$. Furthermore, $\HTaooo{F}{1}{A}{B}\in\NFq$ and 
		\[
			(\alphau{\HTaooo{F}{1}{A}{B}},\betau{\HTaooo{F}{1}{A}{B}},\nuu{\HTaooo{F}{1}{A}{B}})
			=(B,-A,\Kronu{-\frac{1}{\tau}}AM^\MP A).
		\]
	\eenui
\eexam

\bexam
	Let $A\in\Cqq$, let $B\in\CHq$, let $M\in\Cggq$ with $AM^\MP A\geq A$ and let $F\colon\ohe\to\Cqq$ be defined by $F(z)\defg-\frac{1}{z}M$. In view of \rtheo{T1554} and \rexam{E1436}, then $F\in\NFq$ and $(\alphaF,\betaF,\nuF)=(\Oqq,\Oqq,\Kronu{0}M)$, and, furthermore, $\HTaooo{F}{1}{A}{B}\in\NFq$ and $(\alphau{\HTaooo{F}{1}{A}{B}},\betau{\HTaooo{F}{1}{A}{B}},\nuu{\HTaooo{F}{1}{A}{B}})=(B,AM^\MP A-A,\zmq)$.
\eexam

\bexaml{E1548}
	Let $A,B,\beta\in\Cqq$ and let $F\colon\ohe\to\Cqq$ be defined by $F(z)\defg B+z\beta$. In view of \eqref{F-ooa}, then $\HTiaoooa{F}{1}{A}{B}{z}=-\frac{1}{z}A(\Iq+A^\MP\beta)^\MP$ for all $z\in\ohe$.
\eexam

\bexam
	Let $A\in\Cqq$, let $B\in\CHq$, let $\beta\in\Cggq$ with $A(\Iq+A^\MP\beta)^\MP\in\Cggq$ and let $F\colon\ohe\to\Cqq$ be defined by $F(z)\defg B+z\beta$. In view of \rtheo{T1554} and \rexam{E1548}, then $F\in\NFq$ and $(\alphaF,\betaF,\nuF)=(B,\beta,\zmq)$, and, furthermore, $\HTiaooo{F}{1}{A}{B}\in\NFq$ and
	\[
		(\alphau{\HTiaooo{F}{1}{A}{B}},\betau{\HTiaooo{F}{1}{A}{B}},\nuu{\HTiaooo{F}{1}{A}{B}})
		=\lrk\Oqq,\Oqq,\Kronu{0}A(\Iq+A^\MP\beta)^\MP\rrk.
	\]
\eexam

A central theme of this paper is to choose, for a given function $F\in\NFq$, special matrices $A$ and $B$ from $\Cqq$ such that the function $\HTaooo{F}{1}{A}{B}$ and $\HTiaooo{F}{1}{A}{B}$, respectively, belongs to $\NFq$ or to special subclasses of $\NFq$. The following result provides a first contribution to this topic.

\bprop
	Let $F\in\NFq$. Further, let $A\in\Cqq$ be such that $-A\in\Cggq$ and let $B\in\CHq$. Then $\HTaooo{F}{1}{A}{B}\in\NFq$.
\eprop
\bproof
	Let $G\colon\ohe\to\Cqq$ be defined by
	\bgl{P1144.1}
		G(z)
		\defg B-zA.
	\eg
	In view of $B\in\CHq$, $-A\in\Cggq$, and \eqref{P1144.1}, we see from \rpart{T1554.b} of \rtheo{T1554} then
	\bgl{P1144.2}
		G
		\in\NFq.
	\eg
	Taking $A^\ad=A$ into account, we infer from \rprop{P1511} and \rrema{R0824} then
	\bgl{P1144.3}
		-AF^\MP A
		\in\NFq.
	\eg
	Because of \eqref{F+ooa} and \eqref{P1144.1}, we have
	\bgl{P1144.4}
		\HTaooo{F}{1}{A}{B}
		=G-AF^\MP A.
	\eg
	Using \eqref{P1144.2}, \eqref{P1144.3}, and  \eqref{P1144.4}, we get $\HTaooo{F}{1}{A}{B}\in\NFq$.
\eproof

Furthermore, we will show that under appropriate conditions the equations
\begin{align}\label{F-+=F}
	\HTiaooo{(\HTaooo{F}{1}{A}{B})}{1}{A}{B}
	&=F&
	&\text{and}&
	\HTaooo{(\HTiaooo{F}{1}{A}{B})}{1}{A}{B}
	&=F
\end{align}
hold true. The formulas in \eqref{F-+=F} show that the functions $\HTaooo{F}{1}{A}{B}$ and $\HTiaooo{F}{1}{A}{B}$ form indeed a coupled pair of transformations. Furthermore, it will be clear now our terminologies ``\taaSt{A}{B}'' and ``\tiaaSt{A}{B}''. If all Moore-Penrose inverses in \eqref{F+ooa} and \eqref{F-ooa} would be indeed inverse matrices, then the equations in \eqref{F-+=F} could be confirmed by straightforward direct computations. Unfortunately, this is not the case in more general situations which are of interest for us. So we have to look for a convenient way to prove the equations in \eqref{F-+=F} for situations which will be relevant for us. Now we verify that in important cases the formulas \eqref{F+ooa} and \eqref{F-ooa} can be rewritten as linear fractional transformations with appropriately chosen generating matrix-valued functions. The role of these generating functions will be played by the matix polynomials $\mHTuu{A}{B}$ and $\mHTiuu{A}{B}$ which are studied in \rAppe{A0835}.

\blemml{L1428}
	Let $A,B\in\Cpq$ be such that $\Kerna{A}\subseteq\Kerna{B}$ and let $\mHTuu{A}{B}$ be defined in \eqref{Wuua}. Furthermore, let $\mG$ be a non-empty subset of $\C$ and let $F\colon\mG\to\Cpq$ be a mapping which satisfies for all $z\in\mG$ the conditions $\Kerna{A}\subseteq\Kerna{F(z)}$ and $\Bilda{F(z)}\subseteq\Bilda{A}$. For all $z\in\mG$, then the relations
	\begin{align*}
		F(z)&\in\dblftruu{-A^\MP}{\Iq-A^\MP A}&
		&\text{and}&
		\lftroouA{p}{q}{\mHTuua{A}{B}{z}}{F(z)}&=\HTaoooa{F}{1}{A}{B}{z}
	\end{align*}
	hold true, where $\HTaooo{F}{1}{A}{B}$ is defined in \eqref{F+ooa}.
\elemm
\bproof
	Taking \eqref{F+ooa} and \eqref{Wuua} into account \rlemm{L1213} yields all assertions.
\eproof

\bremal{R0742}
	Let $A\in\Cpq$ and let $\mHTu{A}$ be defined by \eqref{WVu}. Furthermore let $\mG$ be a non-empty subset of $\C$ and let $F\colon\mG\to\Cpq$ be a mapping which satisfies, for all $z\in\mG$, the conditions $\Kerna{A}\subseteq\Kerna{F(z)}$ and $\Bilda{F(z)}\subseteq\Bilda{A}$. Setting $B=\Opq$ in \rlemm{L1428}, then, for all $z\in\mG$, the relations
	\begin{align*}
		F(z)&\in\dblftruu{-A^\MP}{\Iq-A^\MP A}&
		&\text{and}&
		\lftroouA{p}{q}{\mHTua{A}{z}}{F(z)}&=\HTaooa{F}{1}{A}{z}
	\end{align*}
	hold true, where $\HTaoo{F}{1}{A}$ is defined in \eqref{F+-o}.
\erema

The follwing application of \rrema{R0742} and \rlemm{L1428} is important for our further considerations.

\blemml{L1239}
	Let $\kappa\in\NO\cup\set{+\infty}$, let $\seq{\su{j}}{j}{0}{\kappa}\in\Hggequ{\kappa}$ and let $F\in\RFuqa{\kappa}{\seq{\su{j}}{j}{0}{\kappa}}$. Further, let $z\in\ohe$. Then:
	\benui
		\il{L1239.a} $F(z)\in\dblftruu{-\su{0}^\MP}{\Iq-\su{0}^\MP\su{0}}$.
		\il{L1239.b} $\lftrooua{q}{q}{\mHTua{\su{0}}{z}}{F(z)}=\HTaooa{F}{1}{\su{0}}{z}$.
		\il{L1239.c} If $\kappa\in\N\cup\set{+\infty}$, then
		\[
			\lftroouA{q}{q}{\mHTuua{\su{0}}{\su{1}}{z}}{F(z)}
			=\HTaoooa{F}{1}{\su{0}}{\su{1}}{z}.
		\]
	\eenui
\elemm
\bproof
	In view of \rpart{L0908.a} of \rprop{L0908} we have
	\begin{align}\label{L1239.1}
		\KernA{F(z)}&=\Kerna{\su{0}}&
		&\text{and}&
		\BildA{F(z)}&=\Bilda{\su{0}}.
	\end{align}
	Taking \eqref{L1239.1} into account, we infer from \rrema{R0742} the assertions of~\eqref{L1239.a} and~\eqref{L1239.b}.
	
	Suppose now that $\kappa\in\N\cup\set{+\infty}$. Then \rpart{L1037.c} of \rlemm{L1037} gives $\Kerna{\su{0}}\subseteq\Kerna{\su{1}}$. Combining this with \eqref{L1239.1} we see from \rlemm{L1428} that~\eqref{L1239.c} holds.
\eproof

Now we state an important situation where formula \eqref{F-ooa} can be rewritten as linear fractional transformation.

\blemml{L1344}
	Let $A\in\Cggq$, $B\in\Cqq$, and let $\mHTiuu{A}{B}$ be given by \eqref{Vuua}. Let $G\in\NFq$ be such that its \tNp{} $(\alphau{G},\betau{G},\nuu{G})$ satisfies the inclusion $\Kerna{A}\subseteq[\Kerna{\alphau{G}}\cap\Kerna{\betau{G}}\cap\Kerna{\nuua{G}{\R}}]$. Let $F\defg G+B$. For all $z\in\ohe$, then 
	\begin{align}\label{17-2}
			F(z)&\in\dblftruu{A^\MP}{z\Iq-A^\MP B}&
			&\text{and}&
			\HTiaoooa{F}{1}{A}{B}{z}&=\lftroouA{q}{q}{\mHTiuu{A}{B}(z)}{F(z)}.
	\end{align}
\elemm
\bproof
	Let $z\in\ohe$. Since $G\in\NFq$ is supposed, we infer from \rlemm{L1030} that
	\bgl{L1344.1}
		\Kerna{A}
		\subseteq\KernA{G(z)}.
	\eg
	Since $A\in\Cggq$ we have $A\in\CEPq$. Thus, \rpart{L1030.b} of \rlemm{L1030} yields
	\bgl{L1344.4}
		\BildA{G(z)}
		\subseteq\Bilda{A}.
	\eg
	In view of $G\in\NFq$, we have $F(z)-B\in\Iqgg$, whereas the relations \eqref{L1344.1} and \eqref{L1344.4} yield the inclusions $\Kerna{A}\subseteq\Kerna{F(z)-B}$ and $\Bilda{F(z)-B}\subseteq\Bilda{A}$. Thus, the application of \rpart{L1058.d} of \rlemm{L1058} provides us $F(z)\in\dblftruu{A^\MP}{z\Iq-A^\MP B}$. Now the application of \rrema{R1105} yields finally
	\[
		\HTiaoooa{F}{1}{A}{B}{z}
		=\lftroouA{q}{q}{\mHTiuu{A}{B}(z)}{F(z)}.\qedhere
	\]
\eproof

Now we are going to consider the following situation which will turn out to be typical for larger parts of our future considerations. Let $A\in\Cggq$ and $B\in\CHq$ be such that $\Kerna{A}\subseteq\Kerna{B}$. Further, let $F\in\NFq$ be such that
\bgl{86-1}
	\Kerna{A}
	\subseteq\Kerna{\alphaF}\cap\Kerna{\betaF}\cap\KernA{\nuFa{\R}}.
\eg
Then our aim is to investigate the function $\HTiaooo{F}{1}{A}{B}$ given by \eqref{F-ooa}. We start with an auxiliary result.

\blemml{R0930}
	Let  $A\in\Cggq$ and let  $B\in\CHq$ be such that $\Kerna{A}\subseteq\Kerna{B}$. Further, let $F\in \RqP$ be such that \eqref{86-1} holds. Then $G\defg F-B$ belongs to $\RqP$ as well and 
	\begin{equation} \label{18Z}
		(\alphau{G},\betau{G},\nuu{G})
		=(\alphaF-B,\betaF,\nuF).
	\end{equation}
	Furthermore, for all $z\in\ohe$, the relations in \eqref{17-2} hold true. 
\elemm
\bproof
	Since the matrix $B$ is Hermitian and since $F$ belongs to $\RqP$, we see from \rtheo{T1554} that $G\defg F-B$ belongs to $\RqP$ as well and that \eqref{18Z} holds. Because of $\Kerna{A}\subseteq\Kerna{B}$, $\Kerna{\alphaF}\cap \Kerna{B}\subseteq \Kerna{\alphaF-B} = \Kerna{\alphau{G}}$, and \eqref{86-1}, we then have
	\bsp
		\Kerna{A}
		&=\Kerna{A}\cap\Kerna{B}
		\subseteq\Kerna{\alphaF}\cap\Kerna{\betaF}\cap\KernA{\nuFa{\R}}\cap\Kerna{B}\\
		&=\Kerna{\alphau{G}}\cap\Kerna{\betaF}\cap\KernA{\nuFa{\R}}
		=\Kerna{\alphau{G}}\cap\Kerna{\betau{G}}\cap\KernA{\nuua{G}{\R}}.
	\esp
	Applying \rlemm{L1344} completes the proof.
\eproof

The following two specifications of \rlemm{R0930} play an essential role in our subsequent considerations.

\blemml{L1300}
	Let $\seq{\su{j}}{j}{0}{1}\in\Hggequ{1}$ and let $F\in\PFoddqa{\su{0}}$. For each $z\in\ohe$, then
	\ba
		F(z)
		&\in\dblftruu{\su{0}^\MP}{z\Iq-\su{0}^\MP\su{1}}&
		&\text{and}&
		\HTiaoooa{F}{1}{\su{0}}{\su{1}}{z}
		&=\lftroouA{q}{q}{\mHTiuua{\su{0}}{\su{1}}{z}}{F(z)}.
	\ea
\elemm
\bproof
	From \rlemm{L1037} we infer
	\bal{L1300.1}
		\su{0}&\in\Cggq,&
		\su{1}&\in\CHq,&
		&\text{and}&
		\Kerna{\su{0}}&\subseteq\Kerna{\su{1}}.
	\ea
	Let $z\in\ohe$. In view of $F\in\PFoddqa{\su{0}}$, \rpart{L1407.b} of \rlemm{L1407} yields
	\bgl{L1300.2}
		\Kerna{\su{0}}
		\subseteq\KernA{F(z)}.
	\eg
	From \eqref{Poddua} and \rrema{R1433} we see the inclusion $\PFoddqa{\su{0}}\subseteq\NFq$. Thus, from \eqref{L1300.2} and \rpart{L1030.a} of \rlemm{L1030} we obtain
	\bgl{L1300.3}
		\Kerna{\su{0}}
		\subseteq\Kerna{\alphaF}\cap\Kerna{\betaF}\cap\KernA{\nuFa{\R}}.
	\eg
	In view of \eqref{L1300.1} and \eqref{L1300.3}, the application of \rlemm{R0930} yields all assertions.
\eproof

\blemml{L1342}
	Let $\seq{\su{j}}{j}{0}{0}\in\Hggequ{0}$ and let $F\in\PFevenqa{\su{0}}$. Further, let $z\in\ohe$. Then
	\ba
		F(z)
		&\in\dblftruu{\su{0}^\MP}{z\Iq}&
		&\text{and}&
		\HTiaooa{F}{1}{\su{0}}{z}
		&=\lftroouA{q}{q}{\mHTiua{\su{0}}{z}}{F(z)}.
	\ea
\elemm
\bproof
	From \rpart{L1037.a} of \rlemm{L1037} we infer $\su{0}\in\Cggq$. Clearly
	\bal{L1342.2}
		\Oqq&\in\CHq&
		&\text{and}&
		\Kerna{\su{0}}&\subseteq\Kerna{\Oqq}.
	\ea
	In view of $F\in\PFevenqa{\su{0}}$, we see from \eqref{Pevenua} that $F\in\NFuq{-2}$ and $\Kerna{\su{0}}\subseteq\Kerna{\alphaF}\cap\Kerna{\muFa{\R}}$. Thus, \rlemm{L1028} implies
	\bgl{L1342.3}
		\Kerna{\su{0}}
		\subseteq\Kerna{\alphaF}\cap\Kerna{\betaF}\cap\KernA{\nuFa{\R}}.
	\eg
	In view of $\su{0}\in\Cggq$, \eqref{L1342.2}, \eqref{L1342.3}, and \eqref{F+-o}, the application of \rlemm{R0930} completes the proof.
\eproof

Now we formulate the first main result of this section. Assuming the situation of \rlemm{R0930}, we will obtain useful insights into the structure of the \tiaaSta{A}{B}{F}.

\bpropl{L0800}
	Let $A\in\Cggq$ and $B\in\CHq$ be such that $\Kerna{A}\subseteq\Kerna{B}$. Further, let  $F\in \RqP$ be such that \eqref{86-1} holds, and let $H\colon\ohe\to\Cqq$ be defined by $H(z)\defg-B+zA+F(z)$. Then:
	\benui
		\il{L0800.a} $H\in \RqP$ and $(\alphau{H},\betau{H},\nuu{H})=(\alphaF-B,\betaF+A,\nuF)$.	
		\il{L0800.b} For each   $z\in\ohe$,
		\begin{equation} \label{18-2}
			\KernA{H(z)}
			=\Kerna{\betau{H}}\subseteq\Kerna{A}\subseteq\KernA{F(z)-B}
		\end{equation}
		and
		\begin{equation} \label{18-3}
			\det\lrk z\Iq+A^\MP\lek F(z)-B\rek\rrk
			\neq0.
		\end{equation}
		\il{L0800.c} $\HTiaooo{F}{1}{A}{B}=A^\ad(-H^\MP) A$.
		\il{L0800.d} $\HTiaooo{F}{1}{A}{B}\in\Ruua{0}{q}{\seq{t_j}{j}{0}{0}}$ where  $t_0\defg A(A+\betaF)^\MP A$. If $F\in\NFuq{-2}$, then $t_0=A$.
		\il{L0800.e} $\Bilda{\HTiaoooa{F}{1}{A}{B}{z}}=\Bilda{A}$ and $\Kerna{\HTiaoooa{F}{1}{A}{B}{z}}=\Kerna{A}$ for all $z\in\ohe$.
		\il{L0800.f} $\Kerna{\sigmaua{\HTiaooo{F}{1}{A}{B}}{\R}}=\Kerna{A}$ and $\Bilda{\sigmaua{\HTiaooo{F}{1}{A}{B}}{\R}}=\Bilda{A}$.
		\il{L0800.g} If $\det A\neq0$, then $\det\HTiaoooa{F}{1}{A}{B}{z}\neq0$ for all $z\in\ohe$.
	\eenui
\eprop
\bproof
	\eqref{L0800.a} Taking $B\in\CHq$, $A\in\Cggq$, and  $F\in \RqP$ into account, we see from \rtheo{T1554} that $H$ belongs to $\RqP$ as well and that
	\bgl{L0800.1}
		(\alphau{H},\betau{H},\nuu{H})
		=(\alphaF-B,\betaF+A,\nuF).
	\eg
	
	\eqref{L0800.b} From \eqref{L0800.1} we infer
	\begin{equation} \label{19-1}
		\Kerna{\alphaF}\cap\Kerna{B}\subseteq\Kerna{\alphaF-B}
		= \Kerna{\alphau{H}}
	\end{equation}
	and 
	\begin{equation} \label{19-2}
		\Kerna{\betaF}\cap\Kerna{A}\subseteq\Kerna{\betaF+A}
		= \Kerna{\betau{H}}.
	\end{equation}
	Consequently, from $\Kerna{A}\subseteq\Kerna{B}$, \eqref{86-1}, \eqref{19-1}, and \eqref{19-2} then 
	\bsp
			\Kerna{A}
			=\Kerna{A}\cap\Kerna{B} 
			&\subseteq\Kerna{\alphaF}\cap\Kerna{B}\cap\Kerna{\betaF}\cap\Kerna{A}\cap\KernA{\nuFa{\R}}\\
			&\subseteq\Kerna{\alphau{H}}\cap\Kerna{\betau{H}}\cap\KernA{\nuu{H}{(\R)}}
	\esp
	follows. Let $G\defg F-B$. Thus, \rlemm{R0930} shows that $G \in \RqP$ and \eqref{18Z} hold true. \rprop{P1356} yields  
	\begin{equation} \label{19-7-1}
	 \KernA{G(z)}
	 =\Kerna{\alphau{G}}\cap\Kerna{\betau{G}}\cap\KernA{\nuua{G}{\R}}
	\end{equation}
	and 
	\begin{equation} \label{19-4}
	 \KernA{H(z)}
	 =\Kerna{\alphau{H}}\cap\Kerna{\betau{H}}\cap\KernA{\nuu{H}{(\R)}}
	\end{equation}
	for all $z \in \ohe$. Since $A$ and $\betaF$ are non-negative Hermitian matrices, we have $\betau{H}=\betaF +A \geq A \geq \Oqq$ and, consequently,
	\begin{equation} \label{19-8-1}
	  \Kerna{\betau{H}}
	  \subseteq \Kerna{A}.
	\end{equation}
	The inclusion \eqref{19-8-1} and the assumptions $\Kerna{A}\subseteq \Kerna{B}$ and \eqref{86-1} imply together with \eqref{L0800.1} the relations
	\bsp
	 \Kerna{\betau{H}}
	 &\subseteq \Kerna{A}\subseteq \Kerna{B}\cap\Kerna{\alphaF}\cap\KernA{\nuFa{\R}} \\
	 &\subseteq \Kerna{\alphaF-B}\cap\KernA{\nuFa{\R}} 
	 = \Kerna{\alphau{H}}\cap\KernA{\nuua{H}{\R}}.
	\esp
	Thus, we see that
	\begin{equation} \label{19-5}
	 \Kerna{\alphau{H}}\cap\Kerna{\betau{H}}\cap\KernA{\nuu{H}{(\R)}}
	 =\Kerna{\betau{H}}
	\end{equation}
	is true. For all $z \in \ohe$, from \eqref{19-4}, \eqref{19-5}, and \eqref{19-8-1} we know that $\Kerna{H(z)}=\Kerna{\betau{H}} \subseteq \Kerna{A} $ holds for all $z \in \ohe$. Since $\Kerna{A} \subseteq \Kerna{B}$ is assumed, we get from \eqref{86-1}, \eqref{18Z}, and \eqref{19-7-1} that 
	\[
		\begin{split}
			\Kerna{A}
			&\subseteq\Kerna{B} \cap \Kerna{\alphaF}\cap\Kerna{\betaF} \cap\KernA{\nuFa{\R}}\\
			&\subseteq \Kerna{\alphaF-B}\cap\Kerna{\betaF} \cap\KernA{\nuFa{\R}}\\
			&=\Kerna{\alphau{G}}\cap\Kerna{\betau{G}}\cap\KernA{\nuua{G}{\R}}
			=\KernA{G(z)}
			=\KernA{F(z)-B}
		\end{split} 
	\]
	is fulfilled  for all $z \in \ohe$. Thus, \eqref{18-2} is true for all $z \in \ohe$.
	
	In view of $A\in\Cggq$, we have $A\in\CEPq$. Hence, taking \eqref{86-1} into account, we infer from \rlemm{L1030} then
 	\begin{equation}\label{19-7}
	 AA^\MP F(z)
	 =F(z)
	\end{equation}
	for all $z\in\ohe$. Because of $A^\ad=A$, $B^\ad=B$ and $\Kerna{A}\subseteq\Kerna{B}$ we get from \rrema{L0818} then
 	\begin{equation}\label{19-8}
	 AA^\MP B
	 =B.
	\end{equation}
	Let $z \in \ohe$. In order to check that 
 	\begin{equation} \label{75-1}
	 	\KernA{z \Iq +A^\MP\lek F(z)-B\rek}
	 	=\set{0_{\xx{q}{1}}},
	\end{equation}
	we now consider an arbitrary $v \in \Kerna{z \Iq +A^\MP[F(z)-B]}$. From \eqref{19-8} and \eqref{19-7} then
	\[
	 	H(z) v
	 	=\lek-AA^\MP B+zA +AA^\MP F(z)\rek v
	 	=A\lrk z \Iq +A^\MP \lek F(z)-B\rek\rrk v
	 	=0_{\xx{q}{1}}
	\]
	follows. Consequently, \eqref{18-2} implies $[G(z)]v = 0_{\xx{q}{1}}$ and, thus,
	\[
		v
		= \frac{1}{z}(zv +A^\MP 0_{\xx{q}{1}})
		=\frac{1}{z}\lek z \Iq +A^\MP G(z)\rek v
		=\frac{1}{z}\lrk z \Iq +A^\MP \lek F(z) -B\rek\rrk v
		=0_{\xx{q}{1}}.
	\]
	This shows that \eqref{75-1} is true. In other words, \eqref{18-3} holds.

	\eqref{L0800.c} Let us now consider again an arbitrary $z \in \ohe$. From \eqref{18-2} and \rpart{L1622.a} of \rrema{L1622} we conclude $A H^\MP(z)H(z)=A$. Using this, \eqref{18-3}, \eqref{19-7}, \eqref{19-8}, and the assumption $A^\ad=A$, we then obtain
	\begin{equation*}
		\begin{split}
			\HTiaoooa{F}{1}{A}{B}{z}
			&=-A\lek H(z)\rek^\MP H(z)\lrk z \Iq +A^\MP\lek F(z)-B\rek\rrk^\MP \\
			&=-A\lek H(z)\rek^\MP\lek zA +F(z)-B\rek\lrk z \Iq +A^\MP\lek F(z)-B\rek\rrk^\inv \\
			&=-A\lek H(z)\rek^\MP  A\lrk z\Iq +A^\MP\lek F(z)-B\rek\rrk\lrk z \Iq +A^\MP\lek F(z)-B\rek\rrk^\inv \\
			&= A \lek-H(z)\rek^\MP  A
			=A^\ad\lek-H^\MP(z)\rek A.
		\end{split} 
	\end{equation*}
	
	\eqref{L0800.d} Because of $H \in \RqP$ and \rprop{P1511}, we see that $-H^\MP$ belongs to $\RqP$. Thus, \eqref{L0800.c} and~\cite{112}*{\crema{3.4}} show that $\HTiaooo{F}{1}{A}{B}$ belongs to $\RqP$ as well. Since the matrices $A$ and $\betaF$ are non-negative Hermitian, the matrix $t_0$ is non-negative Hermitian as well. From \rprop{P1513} we know that 
  \[ 
	 	\lim_{y \to + \infty}\frac{1}{\I y} H(\I y)
	 	= \betau{H}.
	\]
	Furthermore, for all $y \in [1,+\infty)$, we get from \eqref{18-2} that
  \[ 
	 	\rank \lek\frac{1}{\I y} H(\I y)\rek
	 	=  q -\dim \KernA{ H(\I y)}
	 	=q- \dim \Kerna{\betau{H}}
	 	=\rank \betau{H}
	\]
	holds. Thus, the application of \rlemm{L0921} yields 
  \begin{equation} \label{20-3}
	 	\lim_{y \to + \infty} \lek\frac{1}{\I y} H(\I y)\rek^\MP
	 	= \betau{H}^\MP.
	\end{equation}
	Now, we see from \eqref{20-3}, \rpartss{L0800.a}{L0800.c}, and \rrema{R1013} that
	\[
		\begin{split}
			\Oqq
			&= A^\ad\lrk-\lim_{y \to + \infty} \lek\frac{1}{\I y} H(\I y)\rek^\MP \rrk A+A^\ad \betau{H}^\MP A\\
			&= \lim_{y \to + \infty} \lek\I y A^\ad \lrk-\lek H(\I y)\rek^\MP \rrk A+A( A+ \betaF)^\MP A \rek\\
			&= \lim_{y \to + \infty}  \I y \lek\HTiaooo{F}{1}{A}{B} (\I y)+(\I y)^{-1} t_0 \rek.
		\end{split} 
	\]
	Consequently, in view of $\HTiaooo{F}{1}{A}{B}\in\NFq$, \rpart{T1515.a} of \rtheo{T1515} implies $\HTiaooo{F}{1}{A}{B} \in\Ruua{0}{q}{\seq{t_j}{j}{0}{0}}$. If $F\in\NFuq{-2}$ we see from \rrema{R1538} that $\betaF=\Oqq$. Thus, $t_0=AA^\MP A=A$.
	
	\eqref{L0800.e} In view of \eqref{18-3}, we have
	\[
		\HTiaoooa{F}{1}{A}{B}{z}
		=-A\lrk z\Iq+A^\MP\lek F(z)-B\rek\rrk^\MP
		=-A\lrk z\Iq+A^\MP\lek F(z)-B\rek\rrk^\inv.
	\]
	Thus,
	\ba
		\BildA{\HTiaoooa{F}{1}{A}{B}{z}}
		&=\Bilda{A}&
		&\text{and}&
		\KernA{\lek\HTiaoooa{F}{1}{A}{B}{z}\rek^\ad}
		&=\Kerna{A^\ad}.
	\ea
	In view of $\HTiaooo{F}{1}{A}{B}\in\NFq$, \rlemm{L1444} yields
	\[
		\KernA{\lek\HTiaoooa{F}{1}{A}{B}{z}\rek^\ad}
		=\KernA{\HTiaoooa{F}{1}{A}{B}{z}}.
	\]
	Taking $A^\ad=A$ into account we get
	\[
		\KernA{\HTiaoooa{F}{1}{A}{B}{z}}
		=\Kerna{A}.
	\]
	
	\eqref{L0800.f} In view of~\eqref{L0800.d} and~\eqref{L0800.e}, the application of \rlemm{R1414} yields~\eqref{L0800.f}.
	
	\eqref{L0800.g} This follows immediately from~\eqref{L0800.e}.
\eproof

%

\bcoro
	Let $A\in\Cggq$ and $B\in\CHq$ be such that $\Kerna{A}\subseteq\Kerna{B}$. Let $s\in\Cqq$ with $\Oqq\leq s\leq A$ and $\rank s=\rank A$. Further, let $F\colon\ohe\to\Cqq$ be defined by $F(z)\defg z(A-s)s^\MP A$. Then $F$ belongs to $\NFq$ and fulfills $\Kerna{A}\subseteq\Kerna{\alphaF}\cap\Kerna{\betaF}\cap\Kerna{\nuFa{\R}}$. Moreover, $\HTiaooo{F}{1}{A}{B}\in\RFOq$ and $\suo{0}{\HTiaooo{F}{1}{A}{B}}=s$.
\ecoro
\bproof
	Because of $\Oqq\leq s\leq A$ and $\rank s=\rank A$, we have $\Kerna{s}=\Kerna{A}$ and $\Bilda{s}=\Bilda{A}$, which in view of \rlemm{L1633} implies
	\bal{C1514A.2}
		ss^\MP
		&=AA^\MP&
		&\text{and}&
		s^\MP s
		&=A^\MP A.
	\ea
	Hence, from \eqref{C1514A.2} and \rlemm{L1315} we obtain
	\bsp
		(A-s)s^\MP A
		&=As^\MP A-ss^\MP A
		=A^\ad s^\MP A-ss^\MP A\\
		&\geq A^\ad s^\MP sA^\MP ss^\MP A-ss^\MP A
		=AA^\MP AA^\MP AA^\MP A-AA^\MP A
		=\Oqq.
	\esp
	Thus, \rtheo{T1554} yields $F\in\NFq$ and $(\alphaF,\betaF,\nuF)=(\Oqq,(A-s)s^\MP A,\zmq)$, where $\zmq\colon\BAR\to\Cqq$ is given by $\zmqa{B}\defg\Oqq$. In particular, we have then
	\[
		\Kerna{A}\subseteq\KernA{(A-s)s^\MP A}
		=\Kerna{\alphaF}\cap\Kerna{\betaF}\cap\KernA{\nuFa{\R}}
	\]
	and, taking \eqref{C1514A.2} into account, furthermore
	\bgl{C1514A.3}
		A+\betaF
		=A+(A-s)s^\MP A
		=A+As^\MP A-ss^\MP A
		=As^\MP A.
	\eg
	Thanks to \rpart{L0800.d} of \rprop{L0800} we get then $\HTiaooo{F}{1}{A}{B}\in\Ruua{0}{q}{\seq{t_j}{j}{0}{0}}$, where $t_0\defg A(A+\betaF)^\MP A$. In view of \eqref{Ruua}, this implies $\HTiaooo{F}{1}{A}{B}\in\RFOq$ and, because of \eqref{C1514A.3} and \eqref{C1514A.2}, furthermore
	\bsp
		\suo{0}{\HTiaooo{F}{1}{A}{B}}
		&=t_0
		=A(A+\betaF)^\MP A
		=A(As^\MP A)^\MP A\\
		&=AA^\MP AA^\MP A(As^\MP A)^\MP AA^\MP AA^\MP A
		=ss^\MP ss^\MP A(As^\MP A)^\MP As^\MP ss^\MP s\\
		&=sA^\MP As^\MP A(As^\MP A)^\MP As^\MP AA^\MP s
		=sA^\MP As^\MP AA^\MP s
		=ss^\MP ss^\MP ss^\MP s
		=s.\qedhere
	\esp
\eproof

\bcoro
	Let $A\in\Cggq$ and $B\in\CHq$ such that $\Kerna{A}\subseteq\Kerna{B}$. Let $F\in\NFq$ be such that \eqref{86-1} holds. Then $\HTiaooo{F}{1}{A}{B}$ belongs to $\RFOq$, and $\Oqq\leq\suo{0}{\HTiaooo{F}{1}{A}{B}}\leq A$ and $\rank\suo{0}{\HTiaooo{F}{1}{A}{B}}=\rank A$ hold true.
\ecoro
\bproof
	The application of \rpart{L0800.d} of \rprop{L0800} yields $\HTiaooo{F}{1}{A}{B}\in\Ruua{0}{q}{\seq{t_j}{j}{0}{0}}$, where $t_0\defg A(A+\betaF)^\MP A$. In view of \eqref{Ruua}, this implies $\HTiaooo{F}{1}{A}{B}\in\RFOq$ and
	\[
		\suo{0}{\HTiaooo{F}{1}{A}{B}}
		=t_0
		=A(A+\betaF)^\MP A.
	\]
	Because of $\betaF\in\Cggq$, we have $A+\betaF\geq A\geq\Oqq$, which in view of \rlemm{L1315} implies $A\geq A(A+\betaF)^\MP A\geq\Oqq$ and $\Bilda{A(A+\betaF)^\MP A}=\Bilda{A}$. Thus, we obtain $\Oqq\leq\suo{0}{\HTiaooo{F}{1}{A}{B}}\leq A$ and $\rank\suo{0}{\HTiaooo{F}{1}{A}{B}}=\rank A$.
\eproof

Now we indicate some generic situations in which the formulas in \eqref{F-+=F}, respectively, are satisfied. We start with formula \eqref{F-+=F}.

\bpropl{L1348}
	Let  $A,B\in\Cqq$ and let $F\in \RqP$  be such that \eqref{86-1} and
	\begin{equation} \label{86-11}
		\Bilda{A}
		=\Bilda{\alphaF}+\Bilda{\betaF}+\BildA{\nuFa{\R}}
	\end{equation}
	are fulfilled. Then $\HTiaooo{(\HTaooo{F}{1}{A}{B})}{1}{A}{B}=F$.
\eprop
\bproof
	Let $z \in \ohe$. From \eqref{86-1}, \eqref{86-11}, and \rprop{P1356} we know that 
	\begin{equation} \label{21-1}
		\Kerna{A}
		\subseteq\Kerna{\alphaF}\cap\Kerna{\betaF}\cap\KernA{\nuFa{\R}}
		=\KernA{F(z)}
	\end{equation}
	and 
	\begin{equation}  \label{21-111}
		\Bilda{A}
		=\Bilda{\alphaF}+\Bilda{\betaF}+\BildA{\nuFa{\R}}
		=\BildA{F(z)}.
	\end{equation}
	In view of \eqref{21-1} and \eqref{21-111}, we infer from \rlemm{L1633} then
	\begin{align}\label{L1348.1}
		\lek F(z)\rek\lek F(z)\rek^\MP&=AA^\MP&
		&\text{and}&
		\lek F(z)\rek^\MP\lek F(z)\rek&=A^\MP A.
	\end{align}
	This implies
	\bgl{L1348.3}
		A^\MP A\lek F(z)\rek^\MP
		=\lek F(z)\rek^\MP
	\eg
	and
	\bgl{L1348.4}
		\lek F(z)\rek^\MP AA^\MP\lek F(z)\rek
		=\lek F(z)\rek^\MP\lek F(z)\rek
		=A^\MP A.
	\eg
	Using $(\Iq-A^\MP A)A^\MP=\Oqq$, $A(\Iq-A^\MP A)=\Oqq$, and \eqref{L1348.4}, we obtain
	\bsp
		&\lrk\Iq-A^\MP A+\frac{1}{z}\lek F(z)\rek^\MP A\rrk\lrk\Iq-A^\MP A+zA^\MP\lek F(z)\rek\rrk\\
		&=\Iq-A^\MP A+(\Iq-A^\MP A)\lrk zA^\MP\lek F(z)\rek\rrk\\
		&\quad+\frac{1}{z}\lek F(z)\rek^\MP A(\Iq-A^\MP A)+\lek F(z)\rek^\MP AA^\MP\lek F(z)\rek\\
		&=\Iq-A^\MP A+A^\MP A
		=\Iq.
	\esp
	Thus,
	\bgl{L1348.5}
		\det\lrk\Iq-A^\MP A+\frac{1}{z}\lek F(z)\rek^\MP A\rrk
		\neq0
	\eg
	and
	\[
		\lrk\Iq-A^\MP A+\frac{1}{z}\lek F(z)\rek^\MP A\rrk^\inv
		=\Iq-A^\MP A+zA^\MP\lek F(z)\rek.
	\]
	Using \eqref{F-ooa}, \eqref{F+ooa}, \eqref{L1348.3}, \eqref{L1348.5}, and \eqref{L1348.1}, we get finally
	\bsp
		&\HTiaoooa{(\HTaooo{F}{1}{A}{B})}{1}{A}{B}{z}
		=-A\lrk z\Iq+A^\MP\lek\HTaoooa{F}{1}{A}{B}{z}-B\rek\rrk^\MP\\
		&=-A\lrk z\Iq+A^\MP\lek-A\lrk z\Iq+\lek F(z)\rek^\MP A\rrk\rek\rrk^\MP\\
		&=-A\lrk z(\Iq-A^\MP A)+\lek F(z)\rek^\MP A\rrk^\MP
		=-\frac{1}{z}A\lrk\Iq-A^\MP A+\frac{1}{z}\lek F(z)\rek^\MP A\rrk^\MP\\
		&=-\frac{1}{z}A\lrk\Iq-A^\MP A+zA^\MP\lek F(z)\rek\rrk
		=AA^\MP\lek F(z)\rek
		=F(z).\qedhere
	\esp
\eproof

\bcorol{L0845}
	Let $\kappa\in\NO\cup\set{+\infty}$, $\seq{\su{j}}{j}{0}{\kappa}\in\Hggequ{\kappa}$, and $F\in\RFuqa{\kappa}{\seq{\su{j}}{j}{0}{\kappa}}$.
	\benui
		\il{L0845.a} $\HTiaooo{(\HTaooo{F}{1}{\su{0}}{B})}{1}{\su{0}}{B}=F$ for each $B\in\Cqq$.
		\il{L0845.b} If $\kappa\in\N\cup\set{+\infty}$, then $\HTiaooo{(\HTaooo{F}{1}{\su{0}}{\su{1}})}{1}{\su{0}}{\su{1}}=F$.
		\il{L0845.c} $\HTiaoo{(\HTaoo{F}{1}{\su{0}})}{1}{\su{0}}=F$.
	\eenui
\ecoro
\bproof
	According to \rpart{L0908.c} of \rprop{L0908}, the function $F$ belongs to $\NFq$ and \eqref{N18-1} and \eqref{N18-2} hold. Thus, the application of \rprop{L1348} yields $\HTiaooo{(\HTaooo{F}{1}{\su{0}}{B})}{1}{\su{0}}{B}=F$ for each $B\in\Cqq$. Hence,~\eqref{L0845.a} is proved. Choosing $B=\su{1}$ and $B=\Oqq$ in~\eqref{L0845.a}, we get the assertions of~\eqref{L0845.b} and~\eqref{L0845.c}, respectively.
\eproof


Now we turn our attention to formula \eqref{F-+=F}.

\bpropl{L1438}
	Let $A\in\Cggq$ and $B\in\CHq$ be such that $\Kerna{A}\subseteq\Kerna{B}$. Further, let $F\in\RqP$ be such that \eqref{86-1} holds. Then $\HTaooo{(\HTiaooo{F}{1}{A}{B})}{1}{A}{B}=F$. 
\eprop
\bproof
	Let $z\in\ohe$. In view of $A\in\Cggq$ we have $A\in\CEPq$. Thus taking \eqref{86-1} into account, we infer from \rpart{L1030.b} of \rlemm{L1030} then
	\bgl{Nr.BM1}
		AA^\MP F(z)
		=F(z).
	\eg
	\rPart{L0800.e} of \rprop{L0800} yields $\Kerna{\HTiaoooa{F}{1}{A}{B}{z}}=\Kerna{A}$. Thus, from \rpart{L1622.a} of \rrema{L1622} we infer
	\bgl{21-25}
		A\lek \HTiaoooa{F}{1}{A}{B}{z}\rek^\MP\HTiaoooa{F}{1}{A}{B}{z}
		=A.
	\eg
	Since the matrices $A$ and $B$ are both Hermitian, from $\Kerna{A}\subseteq\Kerna{B}$ the inclusion $\Bilda{B} \subseteq \Bilda{A}$ follows. Consequently, \rpart{L1622.b} of \rrema{L1622} yields \eqref{19-8}. \rPart{L0800.b} of \rprop{L0800} yields \eqref{18-3}. From \eqref{18-3} and \eqref{F-ooa} we see
	\bgl{L1438.2}
		\HTiaoooa{F}{1}{A}{B}{z}
		=-A\lrk z\Iq+A^\MP\lek F(z)-B\rek\rrk^\inv.
	\eg
	Using \eqref{F+ooa}, \eqref{F-ooa}, \eqref{18-3}, \eqref{L1438.2}, \eqref{21-25}, \eqref{Nr.BM1} and \eqref{19-8} we get 
	\begin{equation*}
		\begin{split}
			&\HTaoooa{(\HTiaooo{F}{1}{A}{B})}{1}{A}{B}{z}
			= -A\lrk z \Iq +\lek\HTiaoooa{F}{1}{A}{B}{z}\rek^\MP  A \rrk+B\\
			&= -z A -A  \lek\HTiaoooa{F}{1}{A}{B}{z} \rek^\MP A+B\\
			&= -z A +A  \lek \HTiaoooa{F}{1}{A}{B}{z} \rek^\MP\lek-A\lrk z \Iq + A^\MP\lek F(z)-B\rek\rrk^\inv\rek\\
			&\quad\quad\quad\quad\quad\times\lrk z \Iq + A^\MP\lek F(z)-B\rek\rrk+B\\
			&= -z A +A  \lek\HTiaoooa{F}{1}{A}{B}{z} \rek^\MP 
			\HTiaoooa{F}{1}{A}{B}{z}\lrk z \Iq + A^\MP\lek F(z)-B\rek\rrk+B\\
			&= -z A +A \lrk z \Iq + A^\MP \lek F(z)-B\rek\rrk+B
			= AA^\MP F(z) -A A^\MP B+B
			=F(z).\qedhere
		\end{split} 
	\end{equation*}
\eproof

\bcorol{C1445}
	Let $A\in\Cggq$ and let $F\in\NFq$ be such that \eqref{86-1} holds true. Then $\HTaoo{(\HTiaoo{F}{1}{A})}{1}{A}=F$. 
\ecoro
\bproof
	Choose $B=\Oqq$ in \rprop{L1438}.
\eproof

\bcorol{C1519}
	Let $m\in\NO\cup\set{+\infty}$ and let $\seq{\su{j}}{j}{0}{m}\in\Hggequ{m}$. Then:
	\benui
		\il{C1519.a} $\HTaooo{(\HTiaooo{F}{1}{\su{0}}{\su{1}})}{1}{\su{0}}{\su{1}}=F$ for each $m\in\N\cup\set{+\infty}$. 
		\il{C1519.b} $\HTaoo{(\HTiaoo{F}{1}{\su{0}})}{1}{\su{0}}=F$.
	\eenui
\ecoro
\bproof
	\rPart{L1037.a} of \rlemm{L1037} yields $\su{0}\in\Cggq$. In the case $m\in\N\cup\set{+\infty}$ we see from \rpart{L1037.b} of \rlemm{L1037} that $\su{1}\in\CHq$ and from \rpart{L1037.c} of \rlemm{L1037} that $\Kerna{\su{0}}\subseteq\Kerna{\su{1}}$. Thus, taking \rpart{L0908.c} of \rprop{L0908} into account, the application of \rprop{L1438} and \rcoro{C1445} yields~\eqref{C1519.a} and~\eqref{C1519.b}, respectively.
\eproof

\section{On the $(\su{0},\su{1})$\protect\nobreakdash-Schur Transform for the Classes $\RFuqa{\kappa}{\seq{\su{j}}{j}{0}{\kappa}}$}\label{S1326}
The central topic of this section can be described as follows. Let $m\in\NO$ and let $\seq{\su{j}}{j}{0}{m}\in\Hggequ{m}$. Then \rpart{P1407.b} of \rprop{P1407} tells us that the class $\RFuqa{m}{\seq{\su{j}}{j}{0}{m}}$ is non-empty. If $F\in\RFuqa{m}{\seq{\su{j}}{j}{0}{m}}$, then our interest is concentrated on the \taaSt{\su{0}}{\su{1}} $\HTaooo{F}{1}{\su{0}}{\su{1}}$ of $F$ in the case $m\in\N$ and on the \taaSt{\su{0}}{\Oqq} $\HTaoo{F}{1}{\su{0}}$ of $F$ in the case $m=0$. We will obtain a complete description of these objects. In the case $m=0$, we will show that $\HTaoo{F}{1}{\su{0}}$ belongs to $\PFevenqa{\su{0}}$ (see \rtheo{L1013}). In the case $m=1$, the function $\HTaooo{F}{1}{\su{0}}{\su{1}}$ belongs to $\PFoddqa{\su{0}}$ (see \rtheo{L1015}). The proof of the latter result is mainly based on \rcoro{C1014*} and \rprop{R1640}. Let us now consider the case $m\in\mn{2}{+\infty}$. If $\seq{\su{j}^\seins}{j}{0}{m-2}$ denotes the \tfirstSta{\seq{\su{j}}{j}{0}{m}}, then it will turn out (see \rtheo{L1135}, \rcoro{C0910} and \rtheo{L1314}) that $\HTaooo{F}{1}{\su{0}}{\su{1}}$ belongs to $\RFuqa{m-2}{\seq{\su{j}^\seins}{j}{0}{m-2}}$. Our strategy to prove this is based on the application of Hamburger-Nevanlinna type results for the class $\NFuq{-1}$, which were developed in \rsect{S1510}. Realizing the proofs, we will observe that there is an essential difference between the case of even and odd numbers $m$. In the even case, we will rely on \rtheo{P0611-1}. The main tool in the odd case (, which requires much more work,) is \rtheo{P0611-2}.

Now we start with the detailed treatment of the cases $m=0$ and $m=1$.

\btheol{L1013}
	Let  $\su{0} \in \Cqq$ and let $F\in\Ruqs{0}$. Then $\HTaoo{F}{1}{\su{0}}$ belongs to $\PFevenqa{\su{0}}$.
\etheo
\bproof
	Since $F$ belongs to $\Ruqs{0}$, we have $F\in\ROq$ and $\sigmaF$ belongs to $\MggqRag{\seq{\su{j}}{j}{0}{0}}$. In particular, $\su{0}=\suo{0}{\sigmaF}=\sigmaFa{\R} \in \CHq$ and \rrema{R1433} shows that $F\in\RqP$. From \rprop{P1511} we see that $ -F^\MP$ belongs to  $\RqP$ as well. According to \rrema{R0824}, then $ \su{0}( -F^\MP)\su{0} \in \RqP$. Using \rprop{P1513}, we get
	\begin{equation} \label{21-1a}
		\betau{ \su{0}( -F^\MP)\su{0}}
		= \lim_{y\to+\infty}\lrk\frac{1}{\I y} \su{0}\lek-F(\I y)\rek^\MP\su{0}\rrk.
	\end{equation}
	In view of $F\in\RFuqa{0}{\seq{\su{j}}{j}{0}{0}}$, we conclude from \rpart{T1515.a} of \rtheo{T1515} that
	\begin{equation} \label{21-3}
		\su{0}
		=\lim_{y\to+\infty}\lek-\I yF(\I y)\rek.
	\end{equation}
	From $F\in\ROq$ and \rlemm{R1414} we know that $\Bilda{F(z)}=\Bilda{\sigmaFa{\R}}$ and, hence $\rank F(z) =\rank\sigmaFa{\R} $ for all $z \in \ohe$. Thus, for all $y \in [1,+\infty)$, we have 
	\begin{equation} \label{21-4}
		\rank\lek-\I y F(\I y)\rek
		=\rank\sigmaFa{\R}
		=\rank \su{0}.
	\end{equation}
	Because of \eqref{21-3} and \eqref{21-4}, we see from \rlemm{L0921} that
	\begin{equation} \label{22-1-1}
		\lim_{y\to+\infty}\lek- \I y F(\I y)\rek^\MP
		=\su{0}^\MP.
	\end{equation}
	Combining \eqref{22-1-1} and \eqref{21-1a}, we infer
	\[
		\su{0}
		=\su{0}\lrk\lim_{y\to+\infty}\lek- \I y F(\I y)\rek\rrk^\MP\su{0}
		=\lim_{y\to+\infty}\lrk\frac{1}{\I y} \su{0}\lek F(\I y)\rek^\MP\su{0}\rrk 
		=\betau{ \su{0}( -F^\MP)\su{0}} 
	\]
	and, consequently, $\betau{ \su{0}( -F^\MP)\su{0}} -\su{0}=\Oqq \in \Cggq$. For all $z \in \ohe$, from \eqref{F+ooa} and \eqref{F+-o} we see that 
	\begin{equation} \label{22-3}
		\HTaooa{F}{1}{\su{0}}{z}
		=-\su{0}\lrk z\Iq+\lek F(z)\rek^\MP\su{0}\rrk
			=\lek\su{0}(-F)^\MP\su{0}\rek(z)+z(-\su{0}).
	\end{equation}
	Since $\su{0} (-F)^\MP\su{0}$ belongs to $\RqP$, we get from \eqref{22-3} and \rrema{R1355} that then $\HTaoo{F}{1}{\su{0}}$ belongs to $\RqP$ as well and that $\betau{\HTaoo{F}{1}{\su{0}}}=\Oqq$. Thus,
	\bgl{L1013.2}
		\HTaoo{F}{1}{\su{0}}
		\in\NFuq{-2}
	\eg
	follows from \rrema{R1538}. Taking \eqref{22-3} into account, we conclude
	\[ 
		\HTaooa{F}{1}{\su{0}}{\I}
		=\lek\su{0} (-F)^\MP\su{0}\rek(\I)-\I\su{0}
		=\lrk-\su{0}\lek F(\I)\rek^\MP -\I\Iq \rrk\su{0}
	\]
	and, in particular, $\Kerna{\su{0}}\subseteq \Kerna{\HTaooa{F}{1}{\su{0}}{\I}}$. Thus from \eqref{L1013.2} and \rprop{L1608} we get $\Kerna{\su{0}}\subseteq\Kerna{\alphau{\HTaoo{F}{1}{\su{0}}} }\cap\Kerna{\nuua{\HTaoo{F}{1}{\su{0}}}{\R}}$. Consequently, $\HTaoo{F}{1}{\su{0}} \in \PFevenqa{\su{0}}$.
\eproof

\bcorol{C0835}
	Let $\kappa\in\NO\cup\set{+\infty}$, $\seq{\su{j}}{j}{0}{\kappa}\in\Hggequ{\kappa}$ and $F\in\RFuqa{\kappa}{\seq{\su{j}}{j}{0}{\kappa}}$. Then $\HTaoo{F}{1}{\su{0}}\in\PFevenqa{\su{0}}$.
\ecoro
\bproof
	From \rrema{R0938*} we get $F\in\RFuqa{0}{\seq{\su{j}}{j}{0}{0}}$. Thus, the application of \rtheo{L1013} completes the proof.
\eproof

\btheol{L1015}
	Let $\seq{\su{j}}{j}{0}{1}\in\Hggequ{1}$ and let $F\in\Ruqs{1}$. Then $\HTaooo{F}{1}{\su{0}}{\su{1}}$ belongs to $\PFoddqa{\su{0}}$.
\etheo
\bproof
	For all $z \in \ohe$, we see from \eqref{F+ooa} that
	\begin{equation} \label{23-1}
		\HTaoooa{F}{1}{\su{0}}{\su{1}}{z}
		=\HTaooa{F}{1}{\su{0}}{z} + \su{1}.
	\end{equation}
	Because of $F\in\Ruqs{1}$, we have $F \in \Rkq{1}$ and 
	\begin{equation} \label{23-2}
		\sigmaF
		\in \MggqRag{\seq{\su{j}}{j}{0}{1}}.
	\end{equation}
	Furthermore, we see from \eqref{Ruua} that the function $F$ belongs to $\Ruqs{0}$ as well. Thus, \rtheo{L1013} yields $\HTaoo{F}{1}{\su{0}}\in\PFevenqa{\su{0}}$. In particular, $\HTaoo{F}{1}{\su{0}}\in \RqP$. From \eqref{23-2} we see that the matrices $\su{0}$ and $\su{1}$ are Hermitian. Thus, we conclude from $\HTaoo{F}{1}{\su{0}}\in \RqP$, \eqref{23-1} and \rtheo{T1554} that $\HTaooo{F}{1}{\su{0}}{\su{1}}$ also belongs to $\RqP$. Especially, we get then that $ \Phi \colon [1,+\infty) \to [0, +\infty)$ given by 
	\begin{equation*}
		\Phi(y)
		\defg \frac{\norma{\im \HTaoooa{F}{1}{\su{0}}{\su{1}}{\I y}}}{y}
	\end{equation*}
	is continuous and, in particular, Borel measurable. Since $F$ belongs to $\Rkq{1}$, we have $F \in \RqP$. Consequently, \rprop{P1511} shows that $-F^\MP$ belongs to $\RqP$ as well. In particular, $F^\MP$ is continuous. Furthermore, from $F \in \RqP$ and \eqref{Fuoa} we see that $\trauo{F}{1}{s}$ is continuous. Thus, $\Theta \colon  [1,+\infty) \to [0, +\infty)$ defined by 
	\begin{equation} \label{23-7}
		\Theta (y)
		\defg \normA{\re \lrk-(\I y)^{-1} \su{0} F^\MP(\I y)\lek\trauoa{F}{1}{s}{\I y}-\su{1}\rek \su{0}^\MP\lek\trauoa{F}{1}{s}{\I y} -\su{1}\rek\rrk}
	\end{equation}
	is continuous. Let $\su{-1}\defg\Oqq$. In view of \rrema{R1404}, we then conclude that $F$ belongs to the class $\Rqos{1}$. According to \rcoro{C1014*}, we get
	\begin{equation} \label{23}
		\lim_{y \to +\infty} \trauoa{F}{1}{s}{\I y}
		=\Oqq
	\end{equation}
	and
	\begin{equation} \label{23-3}
			\su{0}
			=\lim_{y \to +\infty}\lek-\I y F(\I y)\rek.
	\end{equation}
	Since $F$ belongs to $\Rkq{1}$ we get from \rlemm{R1414} and \eqref{23-2} that 
	\begin{equation} \label{23-9}
			\KernA{F(z)}
			=\KernA{\sigmaFa{\R}}
			=\Kerna{\su{0}}
	\end{equation}
	and
	\begin{equation} \label{23-8}
			\BildA{F(z)}
			=\BildA{\sigmaFa{\R}}
			=\Bilda{\su{0}}
	\end{equation}
	for all $z \in \ohe$ and, hence, $\rank F(\I y)=\rank \su{0}$ for all $y \in [1,+\infty)$. Using this fact, we see from \eqref{23-3} and \rlemm{L0921} that 
	\begin{equation} \label{24-1}
		\su{0}^\MP
		=\lrk\lim_{y \to +\infty}\lek-\I y F(\I y)\rek\rrk^\MP
		=\lim_{y \to +\infty} \lrk(-\I y)^{-1}\lek F(\I y)\rek^\MP\rrk.
	\end{equation}
	Because of \eqref{23} and \eqref{24-1}, we obtain
	\begin{equation} \label{24-2}
		\begin{split}
			\su{0}\su{0}^\MP\su{1} \su{0}^\MP \su{1}
			&=\su{0}\su{0}^\MP ( \Oqq-\su{1}) \su{0}^\MP (\Oqq-\su{1}) \\
			&=\lim_{y \to +\infty} \su{0} \lrk(-\I y)^{-1}\lek F(\I y)\rek^\MP\rrk\lek\trauoa{F}{1}{s}{\I y}-\su{1}\rek\su{0}^\MP\lek\trauoa{F}{1}{s}{\I y} -\su{1}\rek
		\end{split}
	\end{equation}
	and, in view of \eqref{23-7}, consequently,
	\begin{equation} \label{24-3}
		\normA{- \re(\su{0}\su{0}^\MP\su{1} \su{0}^\MP \su{1})}
		=\lim_{y \to +\infty} \Theta (y).
	\end{equation}
	Since $\Theta$ is continuous, \eqref{24-3} implies the existence of a real number $c$ such that $\Theta (y) \leq c$ for all $y \in [1,+\infty)$. According to \rpart{R1640.a} of \rprop{R1640}, the function $\trauo{F}{1}{s}$ belongs to  $\Rkq{-1}$, i.\,e., $ \trauo{F}{1}{s} \in \RqK{1} $ and $\gammaF = \Oqq$ hold true. In particular,
	\[ 
		\int_{[1,+\infty)}\Abs{\frac{\norma{\im \trauoa{F}{1}{s}{\I y}}}{y}}\tilde\Leb (\dif y)
		< +\infty,
	\]
	where $\tilde\Leb$\index{lambda^~@$\tilde\Leb$} is again the restriction of the Lebesgue measure on $\BA{[1,+\infty)}$. Thus, we get
	\begin{equation}	 \label{CL2}
			\int_{[1,+\infty)}\Abs{\frac{\norma{\im \trauoa{F}{1}{s}{\I y}}}{y} + \frac{c}{y^2}}\tilde\Leb (\dif y)
			< +\infty.
	\end{equation}
	From \eqref{23-9}, \eqref{23-8}, and \rrema{L1622} we conclude that 
	\begin{align} \label{24-2a}
		\su{0}\lek F(z)\rek^\MP F(z)&=\su{0}&
		&\text{and}&
		\su{0} \su{0}^\MP F(z)&=F(z)
	\end{align}
	hold true for all $z \in \ohe$. In view of \eqref{Fuoa}, we have
	\bgl{L1015.1}
		\trauoa{F}{1}{s}{z}
		=z^2F(z)+z\su{0}+\su{1}
	\eg
	for all $z\in\ohe$. Using \eqref{L1015.1}, $\su{-1}=\Oqq$, \eqref{24-2a}, and \eqref{F+ooa}, we get
	\begin{equation} \label{25-1}
		\begin{split}
			&\trauoa{F}{1}{s}{z}+ \frac{1}{z} \su{0} \lek-z F(z)\rek^\MP
			\lek\trauoa{F}{1}{s}{z}-\su{1}\rek\su{0}^\MP\lek\trauoa{F}{1}{s}{z} -\su{1}\rek\\
			&=\trauoa{F}{1}{s}{z} 
			- \frac{1}{z^2} \su{0}\lek F(z)\rek^\MP
			\lek z^2 F(z)+z\su{0}\rek\su{0}^\MP\lek z^2 F(z)+z\su{0}\rek\\
			&=\trauoa{F}{1}{s}{z} 
			- z^2 \su{0} \lek F(z)\rek^\MP F(z) \su{0}^\MP F(z)
			- z \su{0} \lek F(z)\rek^\MP F(z) \su{0}^\MP \su{0} \\
			&\quad- z \su{0} \lek F(z)\rek^\MP \su{0} \su{0}^\MP F(z)
			-	\su{0} \lek F(z)\rek^\MP  \su{0}\su{0}^\MP \su{0} \\
			&=\trauoa{F}{1}{s}{z} 
			- z^2 \su{0}  \su{0}^\MP F(z) - z \su{0}  \su{0}^\MP \su{0}
			- z \su{0} \lek F(z)\rek^\MP F(z) -\su{0} \lek F(z)\rek^\MP   \su{0} \\
			&=z^2 F(z) + z \su{0} + \su{1} -z^2 F(z) -z \su{0}-z \su{0} 
			- \su{0} \lek F(z)\rek^\MP \su{0} \\
			&=-\su{0}\lrk z \Iq + \lek F(z)\rek^\MP\su{0}\rrk+\su{1}
			=\HTaoooa{F}{1}{\su{0}}{\su{1}}{z}.
		\end{split}
	\end{equation}
	Because of \eqref{25-1}, \eqref{23-7}, and $\Theta (y) \leq c$ for all $y \in [1,+\infty)$, we then get
	\bspl{25-2}
		&\Abs{\frac{\norma{\im \HTaoooa{F}{1}{\su{0}}{\su{1}}{\I y}}}{y}}\\
		&\leq \frac{1}{y} \Biggl[\normA {\im \trauoa{F}{1}{s}{\I y}}\\
		&\quad\quad+ \normA{\im \lrk\frac{1}{\I y}  \su{0}\lek-\I y F(z)\rek^\MP\lek\trauoa{F}{1}{s}{\I y}-\su{1}\rek\su{0}^\MP\lek\trauoa{F}{1}{s}{\I y} -\su{1}\rek\rrk} \Biggr]\\
		&= \frac{1}{y}\biggl[\normA{\im \trauoa{F}{1}{s}{\I y}}\\
		&\quad\quad+  \frac{1}{y}\normA{\re\lrk  - \su{0}\lek-\I y F(z)\rek^\MP\lek\trauoa{F}{1}{s}{\I y}-\su{1}\rek\su{0}^\MP\lek\trauoa{F}{1}{s}{\I y} -\su{1}\rek\rrk}\biggr]\\
		&=\frac{1}{y}\normA{\im\trauoa{F}{1}{s}{\I y}}+\frac{1}{y^2}\Theta(y)
		\leq \frac{1}{y} \normA{\im \trauoa{F}{1}{s}{\I y}} +  \frac{c}{y^2}
	\esp
	for all $y \in [1, +\infty)$. Thus, \eqref{25-2} and \eqref{CL2} imply
	\begin{equation}	 \label{25-3}
		\int_{[1,+\infty)}\Abs{\frac{\norma{\im \HTaoooa{F}{1}{\su{0}}{\su{1}}{\I y}}}{y}}\tilde\Leb (\dif y)
		< +\infty.
	\end{equation}
	Since $\HTaooo{F}{1}{\su{0}}{\su{1}}$ belongs to $\RqP$, inequality \eqref{25-3} shows that $\HTaooo{F}{1}{\su{0}}{\su{1}} $ belongs to $\RqK{-1}$. From \rprop{P1401} we then know that
	\begin{equation}	 \label{25}
		\gammau{ \HTaooo{F}{1}{\su{0}}{\su{1}}}
		= \lim_{y \to +\infty}  \HTaoooa{F}{1}{\su{0}}{\su{1}}{\I y}.
	\end{equation}
	Because of \eqref{23}, \eqref{24-2}, \eqref{25-1}, and \eqref{25}, we then have 
	\begin{equation}	 \label{25-26}
		\begin{split}
			\Oqq
			&=\Oqq-0\cdot\su{0}\su{0}^\MP\su{1}\su{0}^\MP\su{1}\\
			&=\lim_{y \to +\infty}\trauoa{F}{1}{s}{\I y}\\
			&\quad+ \lrk\lim_{y \to +\infty} \frac{1}{\I y}\rrk\lim_{y \to +\infty} \lrk\su{0}\lek-\I y F(\I y)\rek^\MP\lek\trauoa{F}{1}{s}{\I y}-\su{1}\rek\su{0}^\MP\lek\trauoa{F}{1}{s}{\I y} -\su{1}\rek\rrk\\
			&= \lim_{y \to +\infty}  \HTaoooa{F}{1}{\su{0}}{\su{1}}{\I y}
			=	\gammau{ \HTaooo{F}{1}{\su{0}}{\su{1}}}.
		\end{split}
	\end{equation}
	Since $\HTaooo{F}{1}{\su{0}}{\su{1}}$ belongs to $\RqK{-1}$, we see from \eqref{25-26} and \eqref{Ruu} that $\HTaooo{F}{1}{\su{0}}{\su{1}}$ belongs to $\Rkq{-1}$. In view of \eqref{23-2}, the function $\eu{1}\colon\R\to\R$ defined by $\eua{1}{t}\defg t$ belongs to $\LaaaR{\R}{\BAR}{\sigmaF}$. Taking~\cite{112}*{\clemm{B.2}} into account, we infer then that $\Kerna{\sigmaFa{\R}} \subseteq \Kerna{\int_{\R} \eu{1} \dif\sigmaF}$, i.\,e., $\Kerna{\su{0}}\subseteq \Kerna{\suo{1}{\sigmaF}}$. Using \eqref{23-2} and \rpart{L1622.a} of \rrema{L1622}, we then get $\su{1}\su{0}^\MP\su{0}=\su{1}$. Consequently,
	\[	 
		\begin{split}	
			\HTaoooa{F}{1}{\su{0}}{\su{1}}{\I}
			&=-\su{0} \lrk\I\Iq +\lek F(\I)\rek^\MP \su{0}\rrk +\su{1}
			=-\I\su{0} - \su{0}\lek F(\I)\rek^\MP \su{0} +\su{1}\su{0}^\MP \su{0} \\
			&= \lrk-\I\Iq - \su{0}\lek F(\I)\rek^\MP+ \su{1}\su{0}^\MP\rrk\su{0}.
		\end{split}
	\]
	In particular, $\Kerna{\su{0}}\subseteq \Kerna{  \HTaoooa{F}{1}{\su{0}}{\su{1}}{\I}}$. Since $ \HTaooo{F}{1}{\su{0}}{\su{1}}$ belongs to $\Rkq{-1}$, from \rlemm{R1410} it follows $\Kerna{\su{0}}\subseteq \Kerna{\muua{\HTaooo{F}{1}{\su{0}}{\su{1}}}{\R}}$. Thus, $\HTaooo{F}{1}{\su{0}}{\su{1}}$ belongs to $\PFoddqa{\su{0}}$.
\eproof

\bcorol{C0856}
	Let $\kappa\in\N\cup\set{+\infty}$, let $\seq{\su{j}}{j}{0}{\kappa}\in\Hggequ{\kappa}$ and let $F\in\RFuqa{\kappa}{\seq{\su{j}}{j}{0}{\kappa}}$. Then $\HTaooo{F}{1}{\su{0}}{\su{1}}\in\PFoddqa{\su{0}}$.
\ecoro
\bproof
	From \rrema{R0938*} we get $F\in\RFuqa{1}{\seq{\su{j}}{j}{0}{1}}$. Thus, the application of \rtheo{L1015} completes the proof.
\eproof

Now we turn our attention to the case $m\in\mn{2}{+\infty}$. As already mentioned we will apply several Hamburger-Nevanlinna type results from \rsect{S1510}. In order to prepare the application of this material, we still need some auxiliary results. First we will compute the functions introduced in \rrema{R1013} (in particular, see \eqref{Fuoa}) for the case that the function $F$ is replaced by $\HTaooo{F}{1}{\su{0}}{\su{1}}$, whereas the role of the sequence is occupied by the \tfirstSt{} $\seq{\su{j}^\seins}{j}{0}{m-2}$ of the original sequence $\seq{\su{j}}{j}{0}{m}$.

\blemml{L1705}
	Let $\mG$ be a non-empty subset of  $\C\setminus\set{0}$, let $F\colon \mG \to \Cpq$, let  $\kappa\in\N\cup\set{+\infty}$, and let $\seq{\su{j}}{j}{0}{\kappa}$ be a sequence of complex \tpqa{matrices}. In the case $\kappa\geq2$, let $\seq{\su{j}^\sta{1}}{j}{0}{\kappa-2}$ be the \tfirstSta{\seq{\su{j}}{j}{0}{\kappa}}. Let $\su{-1}\defg\Opq$, let $\su{-1}^\seins\defg\Opq$, let  $m\in\mn{1}{\kappa}$, and let $\Deltau{m}\colon\mG\to\Cpq$\index{Delta_@$\Deltau{m}$} be defined by
	\begin{equation} \label{Deltaua}
		\Deltaua{m}{z} \defg
		\begin{cases}
			[\trauoa{F}{1}{s}{z}-\su{1}]\su{0}^\MP[\trauoa{F}{1}{s}{z}-\su{1}]\incase{  m=1}\\
			\\
			\begin{aligned}
				&[\trauoa{F}{1}{s}{z}-\su{1}] \su{0}^\MP \trauoa{F}{m}{s}{z}\\
				& -\trauoa{F}{m}{s}{z}\su{0}^\MP(\su{1}+ \sum\nolimits_{j=1}^{m-1}z^{-j}\su{j-1}^\seins)\\
				&+\su{m}\su{0}^\MP\su{1}+\sum\nolimits_{k=1}^{m-1} \su{m-k}\su{0}^\MP\su{k-1}^\seins\\
				&+\sum\nolimits_{j=1}^{m-1}(z^{-j}\sum\nolimits_{k=j}^{m-1} \su{m+j-k}\su{0}^\MP\su{k-1}^\seins)
			\end{aligned}
		 \incase{  m\geq2}
		\end{cases}.
	\end{equation}
	Suppose that $\seq{\su{j}}{j}{0}{m}$ belongs to $\Dpqu{m}$ given in \rdefi{D1043} and that $z\in\mG$ is such that $\Kerna{F(z)}=\Kerna{\su{0}}$ and $\Bilda{F(z)}\subseteq\Bilda{\su{0}}$ are fulfilled.  Then
	\begin{equation} \label{LA}
		\trauoa{(\HTaooo{F}{1}{\su{0}}{\su{1}})}{m-2}{s^\seins}{z}
		=\trauoa{F}{m}{s}{z}+\frac{1}{z}\su{0}\lek-zF(z)\rek^\MP \Deltaua{m}{z}. 
	\end{equation}
\elemm
\bproof 
	Because of $\Kerna{F(z)}=\Kerna{\su{0}}$ and \rpart{L1622.a} of \rrema{L1622}, we have 
	\begin{align}\label{I-5}
		\su{0}\lek F(z)\rek^\MP F(z)&=\su{0}&
		&\text{and}&
		F(z) \su{0}^\MP\su{0}&=F(z),
	\end{align}
	whereas in view of \rpart{L1622.b} of \rrema{L1622} the assumption $\Bilda{F(z)} \subseteq \Bilda{\su{0}} $ yields $\su{0}\su{0}^\MP F(z)=F(z)$. Since $\seq{\su{j}}{j}{0}{m}$ belongs to $\Dpqu{m}$ and since $\su{-1}=\Oqq$, the equations 
	\begin{align} \label{I-1}
		\su{0}\su{0}^\MP \su{j}&=\su{j}&
		&\text{and}&
		\su{j}\su{0}^\MP \su{0}&=\su{j}
	\end{align}
	are valid for all $j\in\mn{-1}{m}$. From~\cite{103}*{\crema{8.5}} we know that
	\begin{equation} \label{I-4}
		\su{0}\su{0}^\MP \su{j}^\seins
		=\su{j}^\seins
	\end{equation}
	is true for all $j\in\mn{-1}{m-2}$. Using \eqref{Fuoa}, \eqref{I-5}, and \eqref{I-1}, we easily check that
	\begin{align} \label{I-11}
		\su{0}\su{0}^\MP \trauoa{F}{m}{s}{z}&=\trauoa{F}{m}{s}{z}&
		&\text{and}&
		\trauoa{F}{m}{s}{z}\su{0}^\MP\su{0}&=\trauoa{F}{m}{s}{z}
	\end{align}
	are fulfilled. Because of \eqref{Fuoa}, \eqref{F+ooa}, \eqref{I-1}, \eqref{I-4}, and \eqref{I-5}, we get 
	\begin{equation} \label{CF1}
		\begin{split}
			&\trauoa{( \HTaooo{F}{1}{\su{0}}{\su{1}})}{m-2}{s^\seins}{z}
			=z^{m-1}\lek-\su{0}\lrk z\Iq +\lek F(z)\rek^\MP\su{0}\rrk+\su{1} +\sum_{k=0}^{m-1} z^{-k}\su{k-1}^\seins\rek\\
			&= z^{m-1}\lrk-z\su{0}\su{0}^\MP \su{0}-\su{0}\lek F(z)\rek^\MP\su{0}+\su{0}\su{0}^\MP\su{1}+\sum_{k=0}^{m-1} z^{-k} \su{0}\su{0}^\MP\su{k-1}^\seins\rrk\\
			&= z^{m-1}\Biggl(-z\su{0}\lek F(z)\rek^\MP F(z)\su{0}^\MP \su{0}-\su{0}\lek F(z)\rek^\MP\su{0}+\su{0}\lek F(z)\rek^\MP F(z)\su{0}^\MP\su{1}\\
			&\quad+\sum_{k=0}^{m-1} z^{-k} \su{0}\lek F(z)\rek^\MP F(z) \su{0}^\MP\su{k-1}^\seins\Biggr)\\
			&= - z^{m-1} \su{0}\lek F(z)\rek^\MP\lek z F(z)\su{0}^\MP \su{0}+\su{0}
			-F(z)\su{0}^\MP\su{1} -\sum_{k=0}^{m-1} z^{-k} F(z) \su{0}^\MP\su{k-1}^\seins\rek\\
			&=  z^{m} \su{0}\lek-zF(z)\rek^\MP\lek \su{0} + F(z)\su{0}^\MP \lrk z\su{0} -\su{1} -\sum_{k=0}^{m-1} z^{-k} \su{k-1}^\seins \rrk\rek.
		\end{split} 
	\end{equation}
	Taking into account $\su{-1}=\Oqq$ and \eqref{Fuoa}, we see that 
	\begin{equation} \label{I-6}
		z^{m+1} F(z)
		= \trauoa{F}{m}{s}{z} -\sum_{j=1}^{m+1} z^{-j} \su{j-1}
	\end{equation}
	holds true. Combing \eqref{CF1}, \eqref{I-6}, \eqref{I-1}, and using \eqref{I-4}, we conclude
	\begin{equation} \label{I-8}
		\begin{split}
			&\trauoa{( \HTaooo{F}{1}{\su{0}}{\su{1}}}{m-2}{s^\seins}{z}\\
			&= \su{0}\lek-zF(z)\rek^\MP\lek z^{m}  \su{0} + \trauoa{F}{m}{s}{z}\su{0}^\MP\lrk\su{0} -z^{-1}\su{1} -\sum _{k=0}^{m-1} z^{-(k+1)} \su{k-1}^\seins\rrk\rek\\
			&= \su{0}\lek-zF(z)\rek^\MP 
			\Biggl[ \trauoa{F}{m}{s}{z}\su{0}^\MP\lrk\su{0} -z^{-1}\su{1} -\sum _{k=0}^{m-1} z^{-(k+1)} \su{k-1}^\seins\rrk+z^{m}  \su{0}  \\
			&\quad\quad\quad\quad\quad\quad\quad-\sum _{j=1}^{m+1} z^{m+1-j} \su{j-1} \su{0}^\MP \su{0} 
			+ \sum _{j=1}^{m+1} z^{m-j} \su{j-1} \su{0}^\MP \su{1}\\
			&\quad\quad\quad\quad\quad\quad\quad+ \sum _{j=1}^{m+1} \sum _{k=0}^{m-1} z^{m-(j+k)} \su{j-1} \su{0}^\MP \su{k-1}^\seins \Biggr] \\
			&= \su{0}\lek-zF(z)\rek^\MP 
			\Biggl[ \trauoa{F}{m}{s}{z}\su{0}^\MP\lrk\su{0} -z^{-1}\su{1} -\sum _{k=0}^{m-1} z^{-(k+1)} \su{k-1}^\seins\rrk\\
			&\quad-\sum _{j=2}^{m+1} z^{m+1-j} \su{j-1} 
			+z^{m-1}  \su{1}  
			+\sum _{j=2}^{m+1} z^{m-j} \su{j-1} \su{0}^\MP \su{1} \\
			&\quad+ \sum _{k=0}^{m-1} z^{m-1-k} \su{k-1}^\seins
			+ \sum _{j=2}^{m+1} \sum _{k=0}^{m-1} z^{m-(j+k)} \su{j-1} \su{0}^\MP \su{k-1}^\seins \Biggr]	.
		\end{split} 
	\end{equation}
	From \eqref{Fuoa}, $\su{-1}=\Oqq$, and the first equation in \eqref{I-11} we get 
	\begin{equation} \label{I-16}
		\begin{split}
			\lek\trauoa{F}{1}{s}{z}-\su{1}\rek\su{0}^\MP \trauoa{F}{m}{s}{z} 
			&= \lek z^2 F(z) +z\su{0} \rek\su{0}^\MP \trauoa{F}{m}{s}{z} 
			= z \trauoa{F}{m}{s}{z}+z^2 F(z) \\
			&=z \trauoa{F}{m}{s}{z}+z^2 F(z)  \su{0}^\MP \trauoa{F}{m}{s}{z}.
		\end{split} 
	\end{equation}
	First we now consider the case $m=1$. Thanks to \eqref{I-8}, $\su{-1}=\Oqq$, and $\su{-1}^\seins=\Oqq$, we  then see that
	\begin{equation} \label{I-10}
		\trauoa{( \HTaooo{F}{1}{\su{0}}{\su{1}})}{m-2}{s^\seins}{z} 
		= \su{0}\lek-zF(z)\rek^\MP \lek\trauoa{F}{m}{s}{z}\su{0}^\MP (\su{0} -z^{-1}) +z^{-1}\su{1} \su{0}^\MP \su{1}\rek. 
	\end{equation}
	Since the first equation in \eqref{I-11} and \eqref{I-5} yield
	\begin{equation} \label{I-181}
		\trauoa{F}{m}{s}{z}
		=\su{0}\su{0}^\MP \trauoa{F}{m}{s}{z}
		= \su{0}\lek F(z)\rek^\MP F(z)\su{0}^\MP \trauoa{F}{m}{s}{z}, 
	\end{equation}
	we see from \eqref{I-10}, the second equation in \eqref{I-11}, \eqref{I-16}, \eqref{Deltaua}, and $m=1$ that 
	\bsp
			&z\lek\trauoa{( \HTaooo{F}{1}{\su{0}}{\su{1}})}{m-2}{s^\seins}{z}-  \trauoa{F}{m}{s}{z}\rek \\
			&= \su{0}\lek-zF(z)\rek^\MP \lek\trauoa{F}{m}{s}{z}\su{0}^\MP (z\su{0} -\su{1}) +\su{1} \su{0}^\MP \su{1}\rek-z\su{0}\lek F(z)\rek^\MP F(z)\su{0}^\MP \trauoa{F}{m}{s}{z} \\
			&= \su{0}\lek-zF(z)\rek^\MP \lek z \trauoa{F}{m}{s}{z} - \trauoa{F}{m}{s}{z}\su{0}^\MP \su{1} +\su{1} \su{0}^\MP \su{1} +z^2 F(z)\su{0}^\MP \trauoa{F}{m}{s}{z} \rek\\ 
			&= \su{0}\lek-zF(z)\rek^\MP \lrk\lek\trauoa{F}{1}{s}{z} - \su{1}\rek\su{0}^\MP \trauoa{F}{m}{s}{z} -\trauoa{F}{m}{s}{z}\su{0}^\MP \su{1} +\su{1} \su{0}^\MP \su{1} \rrk\\ 	 
			&= \su{0}\lek-zF(z)\rek^\MP\lek\trauoa{F}{1}{s}{z} - \su{1}\rek\su{0}^\MP\lek\trauoa{F}{1}{s}{z}- \su{1}\rek
			= \su{0}\lek-zF(z)\rek^\MP \Deltaua{m}{z}.  
	\esp
	Thus, \eqref{LA} is proved if $m=1$. Now we consider the case $m \geq 2$. By virtue of \eqref{I-8} and $\su{-1}^\seins=\Oqq$, we then conclude that 
	\begin{equation} \label{I-19}
		\begin{split}
			&\trauoa{( \HTaooo{F}{1}{\su{0}}{\su{1}})}{m-2}{s^\seins}{z}
			= \su{0}\lek-zF(z)\rek^\MP 
			\Biggl[ \trauoa{F}{m}{s}{z}\su{0}^\MP \lrk\su{0} -z^{-1}\su{1} -\sum _{k=1}^{m-1} z^{-(k+1)} \su{k-1}^\seins \rrk\\
			&\quad -\sum _{j=3}^{m+1} z^{m+1-j} \su{j-1} 
				+\sum _{j=2}^{m+1} z^{m-j} \su{j-1} \su{0}^\MP \su{1} \\
			&\quad + \sum _{k=1}^{m-1} z^{m-1-k} \su{k-1}^\seins
			+ \sum _{j=2}^{m+1} \sum _{k=1}^{m-1} z^{m-(j+k)} \su{j-1} \su{0}^\MP \su{k-1}^\seins \Biggr] \\
			&= \su{0}\lek-zF(z)\rek^\MP 
			\Biggl[\trauoa{F}{m}{s}{z}\su{0}^\MP \lrk\su{0} -z^{-1}\su{1} -\sum _{j=2}^{m} z^{-j} \su{j-2}^\seins \rrk+ z^{-1} \su{m} \su{0}^\MP \su{1}\\
			& \quad +\sum _{j=2}^{m} z^{m-j}(\su{j-2}^\seins + \su{j-1} \su{0}^\MP \su{1} -\su{j})+ \sum _{j=2}^{m+1} \sum _{k=1}^{m-1} z^{m-(j+k)} \su{j-1} \su{0}^\MP \su{k-1}^\seins \Biggr].
		\end{split} 
	\end{equation}
	From~\cite{103}*{\cprop{8.23}} we know that 	
	\begin{equation} \label{I-19A}
		\su{j}^\seins=
		\begin{cases}
			\su{2}-\su{1}\su{0}^\MP\su{1}\incase{ j=0}\\
			\su{j+2} - \su{j+1}\su{0}^\MP\su{1} -\sum_{l=0}^{j-1}\su{j-l}\su{0}^\MP\su{l}^\seins\incase{ j\in\mn{1}{m-2}}
		\end{cases}
	\end{equation}
	holds. If $m=2$, then from  \eqref{I-19} and \eqref{I-19A} we can easily check that 
	\bsp
		\trauoa{( \HTaooo{F}{1}{\su{0}}{\su{1}})}{m-2}{s^\seins}{z}
		&= \su{0}\lek-zF(z)\rek^\MP\Bigl[\trauoa{F}{m}{s}{z}\su{0}^\MP (\su{0} -z^{-1}\su{1} -z^{-2} \su{0}^\seins)\\
		&\quad +z^{-1} (\su{2}\su{0}^\MP \su{1}+ \su{1} \su{0}^\MP \su{0}^\seins+ z^{-1} \su{2} \su{0}^\MP \su{0}^ \seins)\Bigr]
	\esp
	and, consequently, in view of \eqref{I-181}, the second equation in \eqref{I-11}, \eqref{I-16}, and \eqref{Deltaua}, then 
	\[ 
		\begin{split}
			&z\lek\trauoa{( \HTaooo{F}{1}{\su{0}}{\su{1}})}{m-2}{s^\seins}{z} -  \trauoa{F}{m}{s}{z}\rek \\
			&= \su{0}\lek-zF(z)\rek^\MP \Bigl[ z \trauoa{F}{m}{s}{z}\su{0}^\MP \su{0} -\trauoa{F}{m}{s}{z}\su{0}^\MP (\su{1} +z^{-1}\su{0}^ \seins) +  \su{2} \su{0}^\MP \su{1} \\
			& \quad +(\su{1} +z^{-1}\su{2})\su{0}^\MP\su{0}^ \seins
			\Bigr]
			-z\su{0}\lek F(z)\rek^\MP F(z)\su{0}^\MP \trauoa{F}{m}{s}{z} \\	 
			&= \su{0}\lek-zF(z)\rek^\MP \Bigl[ z \trauoa{F}{m}{s}{z} - \trauoa{F}{m}{s}{z}\su{0}^\MP (\su{1} +z^{-1}\su{0}^ \seins) +  \su{2} \su{0}^\MP \su{1} \\
			& \quad +(\su{1} +z^{-1}\su{2})\su{0}^\MP\su{0}^ \seins
			+z^2 F(z)\su{0}^\MP \trauoa{F}{m}{s}{z} \Bigr]\\ 
			&= \su{0}\lek-zF(z)\rek^\MP \biggl(\lek\trauoa{F}{1}{s}{z} - \su{1}\rek\su{0}^\MP \trauoa{F}{m}{s}{z} -\trauoa{F}{m}{s}{z}\su{0}^\MP (\su{1} +z^{-1}\su{0}^ \seins) \\
			& \quad + \su{2} \su{0}^\MP \su{1}
			+(\su{1} +z^{-1}\su{2})\su{0}^\MP\su{0}^ \seins
			\biggr)\\ 	 
			&= \su{0}\lek-zF(z)\rek^\MP \Deltaua{2}{z}
			= \su{0}\lek-zF(z)\rek^\MP \Deltaua{m}{z} 
		\end{split} 
	\]
	holds as well. Thus, \eqref{LA} is also proved in the case $m=2$. Now let $m \geq 3$. Then
	\begin{equation} \label{I-23}
		\begin{split}
			&\sum _{j=2}^{m+1} \sum _{k=1}^{m-1} z^{m-(j+k)} \su{j-1} \su{0}^\MP \su{k-1}^\seins 
			=	\sum _{l=3}^{2m} \sum _{r=\max\set{m+2,l}-(m+1)}^{\min\set{m-1,l-2}} z^{m-l} \su{l-r-1} \su{0}^\MP \su{r-1}^\seins \\
			&=	\sum _{l=3}^{m} \sum _{r=1}^{l-2} z^{m-l} \su{l-r-1} \su{0}^\MP \su{r-1}^\seins
			+\sum _{r=1}^{m-1} z^{-1} \su{m-r} \su{0}^\MP \su{r-1}^\seins \\
			& \quad +	\sum _{l=m+2}^{2m} \sum _{r=l-(m+1)}^{m-1} z^{m-l} \su{m-r-1} \su{0}^\MP \su{r-1}^\seins
			\\
			&=	\sum _{l=3}^{m} \sum _{k=0}^{l-3} z^{m-l} \su{l-k-2} \su{0}^\MP \su{k}^\seins
			+\sum _{r=1}^{m-1} z^{-1} \su{m-r} \su{0}^\MP \su{r-1}^\seins \\
			& \quad+	\sum _{l=m+2}^{2m} \sum _{k=l-m}^{2m} z^{m-l} \su{m-k} \su{0}^\MP \su{k-2}^\seins
			\\ 
			&=	\sum _{j=3}^{m} z^{m-j}\sum _{k=0}^{(j-2)-1} \su{j-2-k} \su{0}^\MP \su{k}^\seins
			+z^{-1}\sum _{r=1}^{m-1}  \su{m-1} \su{0}^\MP \su{r-1}^\seins \\
			& \quad+	\sum _{j=2}^{m} z^{-j} \sum _{k=j}^{m}  \su{m+j-k} \su{0}^\MP \su{k-2}^\seins.
		\end{split} 
	\end{equation}
	Consequently, because of \eqref{I-19}, $\su{-1}=\Oqq$, \eqref{I-23}, and \eqref{I-19A}, we obtain 
	\begin{equation} \label{I-24}
		\begin{split}
			&\trauoa{( \HTaooo{F}{1}{\su{0}}{\su{1}})}{m-2}{s^\seins}{z}\\
			&= \su{0}\lek-zF(z)\rek^\MP 
			\Biggl[ \trauoa{F}{m}{s}{z}\su{0}^\MP \lrk\su{0} -z^{-1}\su{1} 
			-\sum _{j=2}^{m} z^{-j} \su{j-2}^\seins\rrk
			\\
			&\quad + z^{m-2} (\su{1}^\seins  + \su{1}\su{0}^\MP \su{1} - \su{2})
			+\sum _{j=3}^{m} z^{m-j} (\su{j-2}^\seins  + \su{j-1}\su{0}^\MP \su{1} - \su{j}) \\
			& \quad+z^{-1} \su{m}\su{0}^\MP \su{1}
			+\sum _{j=3}^{m} z^{m-j} \sum _{k=0}^{(j-2)-1} \su{(j-2)-k} \su{0}^\MP \su{k}^\seins 
				+z^{-1}\sum _{r=1}^{m-1} \su{m-1} \su{0}^\MP \su{k-2}^\seins 
				\Biggr]\\ 
			&= \su{0}\lek-zF(z)\rek^\MP 
			\Biggl(\trauoa{F}{m}{s}{z}\su{0}^\MP \lrk\su{0} -z^{-1}\su{1} 
			-\sum _{j=2}^{m} z^{-j} \su{j-2}^\seins \rrk
			\\
			& \quad + z^{m-2}\lek\su{1}^\seins  -(\su{2}- \su{1}\su{0}^\MP \su{1})\rek\\
			& \quad +\sum _{j=3}^{m} z^{m-j} \lek\su{j-2}^\seins  -\lrk\su{j}- \su{j-1}\su{0}^\MP \su{1} - \sum _{k=0}^{(j-2)-1} \su{(j-2)-k} \su{0}^\MP \su{k}^\seins\rrk\rek \\
			& \quad +z^{-1}\lrk\su{m}\su{0}^\MP \su{1} +  \sum _{k=1}^{m-1} \su{m-k} \su{0}^\MP \su{k-1}^\seins\rrk
			+\sum _{j=2}^{m} z^{-j} \sum _{k=j}^{m} \su{m+j-k} \su{0}^\MP \su{k-2}^\seins
				\Biggr)
				\\
			&= \su{0}\lek-zF(z)\rek^\MP 
			\Biggl[ \trauoa{F}{m}{s}{z}\su{0}^\MP \lrk\su{0} -z^{-1}\su{1} 
			-\sum _{j=2}^{m} z^{-j} \su{j-2}^\seins\rrk
			\\
			& \quad+z^{-1} \lrk\su{m}\su{0}^\MP \su{1} +  \sum _{k=1}^{m-1} \su{m-k} \su{0}^\MP \su{k-1}^\seins \rrk
			+\sum _{j=2}^{m} z^{-j} \sum _{k=j}^{m} \su{m+j-k} \su{0}^\MP \su{k-2}^\seins
				\Biggr].
		\end{split} 
	\end{equation}
	Using \eqref{I-24}, the first equation in \eqref{I-11}, \eqref{I-5}, and \eqref{Deltaua}, we get
	\[ 
		\begin{split}
			&z\lek\trauoa{(\HTaooo{F}{1}{\su{0}}{\su{1}})}{m-2}{s^\seins}{z} -  \trauoa{F}{m}{s}{z}\rek\\
			&= \su{0}\lek-zF(z)\rek^\MP \Biggl[  \trauoa{F}{m}{s}{z}\su{0}^\MP\lrk z \su{0}-\su{1} 
			- \sum _{j=2}^{m}z^{j+1} \su{j-2}^ \seins\rrk
			+\su{m} \su{0}^\MP \su{1} \\
			&\quad\quad\quad\quad\quad\quad\quad+ \sum _{k=1}^{m-1} \su{m-k}\su{0}^\MP\su{k-1}^ \seins
			+\sum _{j=2}^{m} z^{-j+1} \sum _{k=j}^{m}\su{m+j-k}\su{0}^\MP\su{k-2}^ \seins
			\Biggr]\\
			&\quad-z\su{0}\lek F(z)\rek^\MP F(z)\su{0}^\MP \trauoa{F}{m}{s}{z} 
			\\	 
			&= \su{0}\lek-zF(z)\rek^\MP \Biggl[ z \trauoa{F}{m}{s}{z}\su{0}^\MP\su{0} - \trauoa{F}{m}{s}{z}\su{0}^\MP \lrk\su{1} +\sum _{j=1}^{m-1}z^{-j} \su{j-1}^ \seins \rrk  
				+\su{m} \su{0}^\MP \su{1} \\
			&\quad+ \sum _{k=1}^{m-1} \su{m-k}\su{0}^\MP\su{k-1}^ \seins
			+\sum _{j=1}^{m-1} z^{-j} \sum _{k=j+1}^{m}\su{m+j+1-k}\su{0}^\MP\su{k-2}^ \seins
			-z^2 F(z)\su{0}^\MP \trauoa{F}{m}{s}{z} \Biggr]\\ 
			&= \su{0}\lek-zF(z)\rek^\MP \Biggl(\lek\trauoa{F}{1}{s}{z} - \su{1}\rek\su{0}^\MP \trauoa{F}{m}{s}{z} -\trauoa{F}{m}{s}{z}\su{0}^\MP \lrk\su{1}+\sum _{j=1}^{m-1}z^{-j} \su{j-1}^ \seins \rrk\\
			&\quad+\su{m} \su{0}^\MP \su{1}
			+ \sum _{k=1}^{m-1} \su{m-k}\su{0}^\MP\su{k-1}^ \seins
			+\sum _{j=1}^{m-1} z^{-j} \sum _{k=j}^{m-1}\su{m+j-k}\su{0}^\MP\su{k-1}^ \seins
			\Biggr)\\ 	 
			& = \su{0}\lek-zF(z)\rek^\MP \Deltaua{m}{z}.  
		\end{split} 
	\]
	Thus, \eqref{LA} is also proved in the case $m \geq 3$.
\eproof

Now we study the asymptotic behaviour of the function $\Deltau{m}$, which was introduced in \eqref{Deltaua}.

\blemml{L1627}
	Let $\theta\in[0,2\pi)$ and let $\mG$ be a subset of  $\C\setminus\set{0}$ with $\setaa{\e^{\I  \theta}y}{y \in [1,+\infty)}\subseteq\mathcal{G}$. Let $F\colon \mG \to \Cpq$ be a matrix-valued function, let  $\kappa\in\N\cup\set{+\infty}$, and let $\seq{\su{j}}{j}{0}{\kappa}$ be a sequence of complex \tpqa{matrices}. In the case $\kappa\geq2$ let$\seq{\su{j}^\seins}{j}{0}{\kappa-2}$ be the \tfirstSta{\seq{\su{j}}{j}{0}{\kappa}}. Let $\su{-1}\defg\Opq$,  let  $m\in\mn{1}{\kappa}$, and let $\Deltau{m}\colon\mG\to\Cpq$ be defined by \eqref{Deltaua}. Suppose that 
	\begin{equation} \label{20-1}
		\lim_{r\to+\infty} \trauoa{F}{m}{s}{\e^{\I \theta}r}
		=\Opq 
	\end{equation}
	holds true. Then
	\[ 
			\lim_{r\to+\infty}\Deltaua{m}{\e^{\I  \theta}r}
			=
			\begin{cases}
				\su{1}\su{0}^\MP\su{1}\incase{ m=1}\\
				\su{m}\su{0}^\MP\su{1} +\sum_{k=1}^{m-1}\su{m-k}\su{0}^\MP\su{k-1}^\seins\incase{m\geq2}
			\end{cases}.
	\]
\elemm
\bproof
	In the case $m=1$, in view of \eqref{20-1} and \eqref{Deltaua}, we get immediately 
	\bsp
		\su{m}\su{0}^\MP\su{1}
		&=\lrk\lek\lim_{r\to+\infty}\trauoa{F}{m}{s}{\e^{\I \theta}r}\rek-\su{1}\rrk\su{0}^\MP \lrk\lek\lim_{r\to+\infty} \trauoa{F}{m}{s}{\e^{\I \theta}r}\rek-\su{1}\rrk\\
		&= \lim_{r\to+\infty}\Deltaua{m}{\e^{\I  \theta}r}.
	\esp
	Now assume that $m\geq2$. For all $r \in [1,+\infty)$, from \rrema{R1013} we conclude 
	\begin{equation} \label{20-3a}
		\begin{split} 
			\trauoa{F}{1}{s}{\e^{\I \theta}r}
			&= (\e^{\I \theta}r)^{1-m}\lek\trauoa{F}{m}{s}{\e^{\I \theta}r} - \sum_{j=0}^{m-1-1} (\e^{\I \theta}r)^j\su{m-j}\rek\\
			& =(\e^{\I \theta}r)^{-(m-1)} \trauoa{F}{m}{s}{\e^{\I \theta}r} - \sum_{j=0}^{m-2} (\e^{\I \theta}r)^{j-(m-1)}\su{m-j}.
		\end{split} 
	\end{equation}
	In view of $m\geq2$, the assumption \eqref{20-1}, and \eqref{20-3a}, we have 
	\begin{equation} \label{20-4}
		\begin{split} 
			&\Opq\\
			&=\lek\lim_{r\to+\infty}(\e^{\I \theta}r)^{-(m-1)}\rek\lek\lim_{r\to+\infty} \trauoa{F}{m}{s}{\e^{\I \theta}r}\rek - \sum_{j=0}^{m-2}\lek\lim_{r\to+\infty}(\e^{\I \theta}r)^{j-(m-1)}\rek\su{m-j} \\
			&= \lim_{r\to+\infty}\trauoa{F}{1}{s}{\e^{\I \theta}r}.
		\end{split} 
	\end{equation}
	Consequently, \eqref{20-4}, \eqref{20-1}, and \eqref{Deltaua} yield 
	\[ 
		\begin{split} 
			&\su{m} \su{0}^\MP \su{1} +\sum_{k=1}^{m-1}\su{m-k}\su{0}^\MP\su{k-1}^\seins \\
			& = (\Opq-\su{1}) \su{0}^\MP\cdot\Opq -\Opq\cdot\su{0}^\MP\lrk\su{1}+ \sum_{j=1}^{m-1} 0\cdot \su{j-1}^\seins\rrk+ \su{m}\su{0}^\MP\su{1} \\
			& \quad +\sum_{k=1}^{m-1}\su{m-k}\su{0}^\MP\su{k-1}^\seins + 
			\sum_{j=1}^{m-1}\lrk0\cdot\sum_{l=j}^{m-1}\su{m+j-l}\su{0}^\MP\su{l-1}^\seins\rrk\\
			& =\lrk\lek\lim_{r\to+\infty} \trauoa{F}{1}{s}{\e^{\I \theta}r}\rek-\su{1}\rrk\su{0}^\MP\lek\lim_{r\to+\infty}\trauoa{F}{m}{s}{\e^{\I \theta}r}\rek\\
			&\quad-\lek\lim_{r\to+\infty} \trauoa{F}{m}{s}{\e^{\I \theta}r}\rek\su{0}^\MP\lrk \su{1}+ \sum_{j=1}^{m-1}\lek\lim_{r\to+\infty}{(\e^{\I \theta}r)}^{-j}\rek\su{j-1}^\seins\rrk+ \su{m}\su{0}^\MP\su{1} \\
			&\quad+\sum_{k=1}^{m-1}\su{m-k}\su{0}^\MP\su{k-1}^\seins + 
			\sum_{j=1}^{m-1}\lrk\lek\lim_{r\to+\infty}{(\e^{\I \theta}r)}^{-j}\rek\sum_{l=j}^{m-1}\su{m+j-l}\su{0}^\MP\su{l-1}^\seins\rrk\\	
			&= \lim_{r\to+\infty}\Deltaua{m}{\e^{\I  \theta}r}.\qedhere
		\end{split} 
	\]
\eproof


\btheol{L1135}
	Let $n\in\N$, let  $\seq{\su{j}}{j}{0}{2n} \in \Hggeuu{q}{2n}$  and let $F\in\Ruqs{2n}$. Further, let $\seq{\su{j}^\seins}{j}{0}{2n-2}$ be the \tfirstSta{\seq{\su{j}}{j}{0}{2n}}. Then $\HTaooo{F}{1}{\su{0}}{\su{1}}$ belongs to $\RkqSEjj{2n-2}{0}{2n-2}$.
\etheo
\bproof
	Since $F$ belongs to $\Rkqjj{2n}{0}{2n}$, from \eqref{Ruua} we see that $F$ belongs to $\Rkq{2n}$ and that
 	\begin{equation*}
 		\sigmaF
 		\in \MggqRag{\seq{\su{j}}{j}{0}{2n}}
 	\end{equation*}
	hold true. Furthermore, \rrema{R0938*} shows that $F$ belongs to $\Rkqjj{1}{0}{1}$. Thus, we see from \rtheo{L1015} that $\HTaooo{F}{1}{\su{0}}{\su{1}} \in \PFoddqa{\su{0}}$. In view of \eqref{Ruu}, this in particular means that $\HTaooo{F}{1}{\su{0}}{\su{1}}$ belongs to $\Rkq{-1}$. Consequently, \rprop{P0601} yields $\HTaooo{F}{1}{\su{0}}{\su{1}} \in \RqP$. Letting $\su{-1}\defg\Oqq$, from $\seq{\su{j}}{j}{0}{2n} \in \Hggeuu{q}{2n}$ we get in view of \rpart{L1037.b} of \rlemm{L1037} that $\su{j}^\ad=\su{j}$ for all $j\in \mn{-1}{2n-1}$ and, because of \rprop{P1531}, that the \tfirstSt{} $\seq{\su{j}^\seins}{j}{0}{2n-2}$ of $\seq{\su{j}}{j}{0}{2n}$ belongs to $\Hggeuu{q}{2n-2}$. In particular, setting $\su{-1}^\seins\defg \Oqq$, we see again from \rpart{L1037.b} of \rlemm{L1037} that
 	\begin{equation}	\label{L1014.3}
 		(\su{j}^\seins) ^\ad
 		=\su{j}^\seins 
 	\end{equation}
	for all $j\in \mn{-1}{2n-2}$. Now we verify that the function $F$ satisfies the assumptions of \rlemm{L1705}. From $\seq{\su{j}}{j}{0}{2n} \in \Hggeuu{q}{2n}$ and~\cite{103}*{\clemm{3.1}} we know that $\seq{\su{j}}{j}{0}{2n}\in \Duuu{q}{q}{2n}$. In view of $F\in\RFuqa{2n}{\seq{\su{j}}{j}{0}{2n}}$, we infer from \rpart{L0908.a} of \rprop{L0908} that 
	\bal{L1014.5}
		\KernA{F(z)}
		&=\Kerna{\su{0}}&
		&\text{and}&
		\BildA{F(z)}
		&=\Bilda{\su{0}}
	\ea
	for all $z \in \ohe$. Because of $\seq{\su{j}}{j}{0}{2n}\in \Duuu{q}{q}{2n}$ and \eqref{L1014.5}, \rlemm{L1705}, implies
	\begin{equation} \label{L1014.6}
		\begin{split}
			\trauoa{( \HTaooo{F}{1}{\su{0}}{\su{1}})}{2n-2}{s^\seins}{z}
			=\trauoa{F}{2n}{s}{z}+\frac{1}{z}\su{0}\lek-zF(z)\rek^\MP \Deltaua{2n}{z}
		\end{split} 
	\end{equation}
	for all $z \in \ohe$. In view of $\su{-1}=\Oqq$, the application of \rrema{R1404} yields $F \in \RqKJ{2n}{-1}$. \rcoro{C1014*} then yields 
	\begin{align} \label{L1014.7}
			\lim_{y \to +\infty}  \trauoa{F}{1}{s}{\I y}&=\Oqq,&
			\lim_{y \to +\infty}  \trauoa{F}{2n}{s}{\I y}&=\Oqq
	\end{align}
	and 
	\begin{equation} \label{L1014.8}
			\su{0}
			=\lim_{y \to +\infty}  (-\I y)\lek F(\I y)+\su{-1}\rek
			=\lim_{y \to +\infty}\lek-\I yF(\I y)\rek.
	\end{equation}
	Since \eqref{L1014.5} implies $\rank [-\I yF(\I y)]=\rank \su{0}$ for all $y \in [1,+\infty)$, from \eqref{L1014.8} and \rlemm{L0921} we get that 
	\begin{equation} \label{L1014.9}
		\lim_{y \to +\infty}\lrk\lek-\I yF(\I y)\rek^\MP\rrk
		=\su{0}^\MP.
	\end{equation}
	Because of \eqref{L1014.7}, $\lim_{y \to +\infty}\frac{1}{\I y}=0$, and \eqref{Deltaua}, we have 
	\begin{equation} \label{L1014.10}
		\lim_{y \to +\infty} \Deltaua{2n}{\I y}
		=\su{2n}\su{0}^\MP \su{1} +\sum_{k=1}^{2n-1} \su{2n-k} \su{0}^\MP \su{k-1}^{\seins}.
	\end{equation}
	Consequently, keeping in mind \eqref{L1014.7} again as well as equation \eqref{L1014.9}, \eqref{L1014.10}, and \eqref{L1014.6}, we obtain
	\begin{equation} \label{L1014.11}
		\begin{split}
			&\Oqq \\
			&=\lim_{y \to +\infty}  \trauoa{F}{2n}{s}{\I y} + \lrk\lim_{y \to +\infty}\frac{1}{\I y}\rrk\su{0}  \lek\lim_{y \to +\infty}\lrk\lek-\I yF(\I y)\rek^\MP \rrk\rek\lek\lim_{y \to +\infty} \Deltaua{2n}{\I y}\rek\\
			&= \lim_{y \to +\infty}\lrk\trauoa{F}{2n}{s}{\I y} + \frac{1}{\I y} \su{0}\lek-\I yF(\I y)\rek^\MP \Deltaua{2n}{\I y}\rrk
			= \lim_{y \to +\infty}\trauoa{( \HTaooo{F}{1}{\su{0}}{\su{1}})}{2n-2}{s^\seins}{\I y}.
		\end{split}
	\end{equation}
	From $\HTaooo{F}{1}{\su{0}}{\su{1}} \in \RqP$, \eqref{L1014.3}, \eqref{L1014.11}, and \rtheo{P0611-1} we conclude that  $\HTaooo{F}{1}{\su{0}}{\su{1}}$ belongs to the class $\NFuqa{2n-2}{\seq{\su{j}^\seins}{j}{-1}{2n-2}}$. Because of  $\su{-1}^{\seins} =\Oqq$ we see from \rrema{R1404} that $\HTaooo{F}{1}{\su{0}}{\su{1}}\in\RFuqa{2n-2}{\seq{\su{j}^\seins}{j}{0}{2n-2}}$.
\eproof

\bcorol{C0910}
	Let $\seq{\su{j}}{j}{0}{\infty}\in\Hggequ{\infty}$ and let $F\in\RFuqa{\infty}{\seq{\su{j}}{j}{0}{\infty}}$. Further, let $\seq{\su{j}^\seins}{j}{0}{\infty}$ be the \tfirstSta{\seq{\su{j}}{j}{0}{\infty}}. Then $\HTaooo{F}{1}{\su{0}}{\su{1}}$ belongs to the class $\RFuqa{\infty}{\seq{\su{j}^\seins}{j}{0}{\infty}}$.
\ecoro
\bproof
	Combine \rremass{R0938*}{R1038} and \rtheo{L1135}.
\eproof

\btheol{L1314}
	Let $n\in\N$, let $\seq{\su{j}}{j}{0}{2n+1} \in  \Hggeuu{q}{2n+1}$ and let  $F\in \Rkqjj{2n+1}{0}{2n+1}$. Further, let $\seq{\su{j}^\seins}{j}{0}{2n-1}$ be the \tfirstSta{\seq{\su{j}}{j}{0}{2n+1}}. Then $\HTaooo{F}{1}{\su{0}}{\su{1}}$ belongs to $ \RkqSEjj {2n-1}{0}{2n-1}$.
\etheo
\bproof
	Since $F$ belongs to $\Rkqjj{2n+1}{0}{2n+1}$, we get $F \in \Rkq{2n+1}$ and
	\begin{equation*}
		\sigmaF
		\in \MggqRag{\seq{\su{j}}{j}{0}{2n+1}}.
	\end{equation*}
	Furthermore, \rrema{R0938*} shows that $F$ belongs to $\Rkqjj{1}{0}{1}$. Thus, we see from \rtheo{L1015} that $\HTaooo{F}{1}{\su{0}}{\su{1}} \in \PFoddqa{\su{0}}$. This in particular means that $\HTaooo{F}{1}{\su{0}}{\su{1}}$ belongs to $\Rkq{-1}$. Consequently, $\HTaooo{F}{1}{\su{0}}{\su{1}} \in \RqP$. Letting $\su{-1} \defg \Oqq$, from  $\seq{\su{j}}{j}{0}{2n+1} \Hggeuu{q}{2n+1}$ we get in view of \rpart{L1037.b} of \rlemm{L1037} that $\su{j}^\ad=\su{j}$, holds for all $j \in \mn{-1}{2n}$ and, because of \rprop{P1531}, that the sequence $\seq{\su{j}^\seins}{j}{0}{2n-1}$ belongs to $\Hggequ{2n-1}$. In particular, setting $\su{-1}^\seins \defg \Oqq$, we see again from \rpart{L1037.b} of \rlemm{L1037} that \eqref{L1014.3}, i.\,e.\ $(\su{j}^\seins)^\ad=\su{j}^\seins$, holds for all $j \in \mn{-1}{2n-1}$. We verify now that the function $F$ satisfies the assumptions of \rlemm{L1705}. From $\seq{\su{j}}{j}{0}{2n+1} \in \Hggequ{2n+1}$ and \rlemm{L1048} we know that
	\bgl{L1314.4}
		\seq{\su{j}}{j}{0}{2n+1}
		\in\Duuu{q}{q}{2n+1}.
	\eg
	In view of $F\in\RFuqa{2n+1}{\seq{\su{j}}{j}{0}{2n+1}}$, we infer from \rpart{L0908.a} of \rprop{L0908} that 
	\bal{LK2}
		\KernA{F(z)}
		&=\Kerna{\su{0}}&
		&\text{and}&
		\BildA{F(z)}
		&=\Bilda{\su{0}}
	\ea
	for all $z \in \ohe$. Because of \eqref{L1314.4} and \eqref{LK2}, \rlemm{L1705} implies
	\begin{equation} \label{FS2}
		\trauoa{(\HTaooo{F}{1}{\su{0}}{\su{1}})}{2n-1}{s^\seins}{z} 
		= \trauoa{F}{2n+1}{s}{z}+\frac{1}{z}\su{0}\lek-zF(z)\rek^\MP \Deltaua{2n+1}{z}
	\end{equation}
	for all $z \in \ohe$. In view of $\su{-1}=\Oqq$ and the choice of $F$, \rrema{R1404} yields
	\bgl{L1314.1}
		F
		\in\NFuqA{2n+1}{\seq{\su{j}}{j}{-1}{2n+1}}.
	\eg
	In view of \eqref{L1314.1}, \rcoro{C1014*} then provides us 
	\begin{align} \label{LF1}
		\lim_{y \to +\infty}  \trauoa{F}{1}{s}{\I y}&=\Oqq,&
		\lim_{y \to +\infty}  \trauoa{F}{2n+1}{s}{\I y}&=\Oqq
	\end{align}
	and 
	\begin{equation} \label{LOK}
		\su{0}
		=\lim_{y \to +\infty}  (-\I y)\lek F(\I y)+\su{-1}\rek
		=\lim_{y \to +\infty}\lek-\I yF(\I y)\rek.
	\end{equation}
	Since \eqref{LK2} implies $\rank [-\I yF(\I y)]=\rank \su{0}$ for all $y \in [1,+\infty)$, from \eqref{LOK} and \rlemm{L0921} we get that 
	\begin{equation} \label{LK3}
		\lim_{y \to +\infty}\lek-\I yF(\I y)\rek^\MP
		=\su{0}^\MP.
	\end{equation}
	Because of \eqref{LF1}, the limit relation $\lim_{y\to+\infty}\frac{1}{\I y}=0$ and \eqref{Deltaua}, we have 
	\begin{equation} \label{LIC}
		\lim_{y \to +\infty} \Deltaua{2n+1}{\I y}
		=\su{2n+1}\su{0}^\MP \su{1} +\sum_{k=1}^{2n} \su{2n+1-k} \su{0}^\MP \su{k-1}^{\seins}.
	\end{equation}
	Consequently, keeping in mind \eqref{LF1}, the limit relation $\lim_{y \to +\infty}\frac{1}{\I y}=0$, equation  \eqref{LIC}, \eqref{LK3}, and \eqref{FS2}, we obtain
	\begin{equation} \label{FCS}
		\begin{split}
			&\Oqq \\
			&=\lim_{y \to +\infty}  \trauoa{F}{2n+1}{s}{\I y} + \lrk\lim_{y \to +\infty}\frac{1}{\I y}\rrk\su{0}\lek\lim_{y \to +\infty}\lrk\lek-\I yF(\I y)\rek^\MP\rrk\rek\lek\lim_{y \to +\infty} \Deltaua{2n+1}{\I y}\rek\\
			&= \lim_{y \to +\infty}\lrk\trauoa{F}{2n+1}{s}{\I y} + \frac{1}{\I y} \su{0}\lek-\I yF(\I y)\rek^\MP \Deltaua{2n+1}{\I y}\rrk\\
			&= \lim_{y \to +\infty} \trauoa{(\HTaooo{F}{1}{\su{0}}{\su{1}})}{2n-1}{s^\seins}{\I y}.
		\end{split}
	\end{equation}
	Thus, \rtheo{P0611-2} provides us $\trauo{(\HTaooo{F}{1}{\su{0}}{\su{1}})}{2n-1}{s^\seins} \in \RqP$. In particular, the function $\trauo{(\HTaooo{F}{1}{\su{0}}{\su{1}})}{2n-1}{s^\seins}$ is holomorphic in $\ohe$. This shows us that $\Phi\colon[1,+\infty) \to [0, +\infty)$ defined by
	\begin{equation} \label{BY}
		\Phi(y)
		\defg \frac{1}{y}\normA{\trauoa{(\HTaooo{F}{1}{\su{0}}{\su{1}})}{2n-1}{s^\seins}{\I y}}
	\end{equation}
	is a continuous function. In view of \eqref{L1314.1}, from \rpart{R1640.a} of \rprop{R1640} we get $\trauo{F}{2n+1}{s} \in \Rkq{-1}$. Thus, \eqref{Ruu} shows that $\trauo{F}{2n+1}{s}$ belongs to $\RqK{-1}$ and that $\gammau{\trauo{F}{2n+1}{s}}=\Oqq$. In view of \eqref{Ruo[-1]}, this implies that $\Omega\colon[1,+\infty) \to \R$ given by
	\bgl{L1314.2}
		\Omega(y)
		\defg\frac{1}{y}\normA{\im\trauoa{F}{2n+1}{s}{\I y}}
	\eg
	belongs to $\Loaaaa{1}{[1,+\infty)}{\BA{[1,+\infty)}}{ \tilde{\lambda}}{\R}$, where $\tilde{\lambda}$\index{lambda^~@$\tilde\Leb$} is again the Lebesgue measure defined on $\BA{[1,+\infty)}$. Because of $F \in \Rkqjj{2n+1}{0}{2n+1}$, \rprop{R0947*} and \rtheo{T0933}, we have $F  \in \RqP$. Thus, we see from \eqref{Fuoa} that the functions $\trauo{F}{1}{s}$ and $\trauo{F}{2n+1}{s}$ are both holomorphic in $\ohe$. Therefore, \eqref{Deltaua} shows that $\Deltau{2n+1}$ is holomorphic in $\ohe$. Since $F$ belongs to $\RqP$, we also see from \rprop{P1511} that $F^\MP$ is holomorphic in $\ohe$. Consequently, $\Theta\colon[1, +\infty) \to \R$ defined by
	\begin{equation} \label{Z1}
		\Theta(y)
		\defg\normA{\re\lrk(\I y)^{-1} \su{0}\lek F(\I y)\rek^\MP \Deltaua{2n+1}{\I y}\rrk}
	\end{equation}
	is a continuous function. Using \eqref{LK3}, \eqref{LIC}, and \eqref{Z1}, we get that
	\[ 
		\begin{split}
			&\normA{-\re\lek\su{0} \su{0}^\MP\lrk\su{2n+1}\su{0}^\MP\su{1} + \sum_{k=1}^{2n} \su{2n+1-k} \su{0}^\MP \su{k-1}^\seins \rrk\rek}\\
			&=\normA{-\re\lrk\su{0}\lek\lim_{y \to +\infty}\lrk\lek-\I y F(\I y)\rek^\MP\rrk\rek \lek\lim_{y \to +\infty} \Deltaua{2n+1}{\I y}\rek\rrk}
			= \lim_{y \to +\infty} \Theta (y).
		\end{split}
	\]
	This implies that there exists a non-negative real number $c$ such that  $\abs{\Theta(y)}\leq c$ for all $y \in [1,+\infty)$. Since $\Omega$ belongs to $\Loaaaa{1}{[1,+\infty)}{\BA{[1,+\infty)}}{ \tilde{\lambda}}{\R}$ the function $\Psi\colon[1,+\infty) \to \R$ given by $\Psi(y) \defg \Omega(y)+ \frac{c}{y^2}$ also belongs to the space $\Loaaaa{1}{[1,+\infty)}{\BA{[1,+\infty)}}{ \tilde{\lambda}}{\R}$. Keeping in mind \eqref{BY}, \eqref{L1314.2} and  \eqref{FS2}, we infer
	\[ 
		\begin{split}
			\Abs{\Phi(y)}
			&\leq\frac{1}{y}\normA{-\im\lek\trauoa{F}{2n+1}{s}{\I y}\rek}+\frac{1}{y}\normA{\im\lrk\frac{1}{\I y}\su{0}\lek-\I y F(\I y)\rek^\MP  \Deltaua{2n+1}{\I y}  \rrk}\\
			&= \Omega(y)+ \frac{1}{y^2}\normA{ \re\lrk\su{0}\lek\I y F(\I y)\rek^\MP  \Deltaua{2n+1}{\I y}\rrk}
			=\Omega(y) +\Theta (y)
			\leq \Psi(y)
		\end{split}
	\]
	for all $y \in [1,+\infty)$, and, consequently, $\Phi$ belongs to $\Loaaaa{1}{[1,+\infty)}{\BA{[1,+\infty)}}{ \tilde{\lambda}}{\R}$. Since $ \trauo{(\HTaooo{F}{1}{\su{0}}{\su{1}})}{2n-1}{s^\seins}$ belongs to $\RqP$, we thus see from \eqref{Ruo[-1]} and \eqref{BY} that
	\bgl{L1314.3}
		\trauo{(\HTaooo{F}{1}{\su{0}}{\su{1}})}{2n-1}{s^\seins}
		\in\RqK{-1}.
	\eg
	From $\HTaooo{F}{1}{\su{0}}{\su{1}} \in \RqP$, \eqref{L1014.3}, \eqref{FCS}, \eqref{L1314.3} and \rtheo{P0611-2} we conclude that $\HTaooo{F}{1}{\su{0}}{\su{1}}$ belongs to the class   $\NFuqa{2n-1}{\seq{\su{j}^\seins}{j}{-1}{2n-1}}$. Because of $\su{-1}^{\seins} =\Oqq$, then \rrema{R1404} shows that the function $\HTaooo{F}{1}{\su{0}}{\su{1}}$ belongs to  $\RFuqa{2n-1}{\seq{\su{j}^\seins}{j}{0}{2n-1}}$.
\eproof


\bpropl{P1323*}
	Let $\kappa\in\mn{2}{+\infty}\cup\set{+\infty}$, let $\seq{\su{j}}{j}{0}{\kappa}\in\Hggequ{\kappa}$, and let $F$ belong to $\Ruqs{\kappa}$. Then $\HTaooo{F}{1}{\su{0}}{\su{1}}\in\RFuqa{\kappa-2}{\seq{\su{j}^\seins}{j}{0}{\kappa-2}}$, where $\seq{\su{j}^\seins}{j}{0}{\kappa-2}$ denotes the \tfirstSta{\seq{\su{j}}{j}{0}{\kappa}}.
\eprop
\bproof
	Because of $\seq{\su{j}}{j}{0}{\kappa}\in\Hggequ{\kappa}$, we have $\seq{\su{j}}{j}{0}{m}\in\Hggequ{m}$ for each $m\in\mn{0}{\kappa}$. \rrema{R0938*} yields $F\in\bigcap_{m=0}^\kappa\RFuqa{m}{\seq{\su{j}}{j}{0}{m}}$. Thus, \rtheoss{L1135}{L1314} show then that $\HTaooo{F}{1}{\su{0}}{\su{1}}$ belongs to $\bigcap_{m=0}^{\kappa-2}\RFuqa{m}{\seq{\su{j}^\seins}{j}{0}{m}}$. Consequently, from \rrema{R0938*} then $\HTaooo{F}{1}{\su{0}}{\su{1}}\in\RFuqa{\kappa-2}{\seq{\su{j}^\seins}{j}{0}{\kappa-2}}$ follows.
\eproof


\section{On the Inverse $(\su{0},\su{1})$\protect\nobreakdash-Schur Transform for Special Subclasses of $\NFq$}\label{S1021}
Against to the background of \rpropss{L1348}{L1438} the considerations of \rsect{S1326} lead us to study three inverse problems which will be formulated now. Let $m\in\NO$, $\seq{\su{j}}{j}{0}{m}\in\Hggequ{m}$ and $G\in\RFuqa{m}{\seq{\su{j}}{j}{0}{m}}$. In the case $m=0$, it was shown in \rtheo{L1013} that $\HTaoo{G}{1}{\su{0}}$ belongs to $\PFevenqa{\su{0}}$. Now we start with a function $F\in\PFevenqa{\su{0}}$ and will show in \rprop{L1046} that $\HTiaoo{F}{1}{\su{0}}\in\RFuqa{0}{\seq{\su{j}}{j}{0}{0}}$. In the case $m=1$, \rtheo{L1013} yields that $\HTaooo{G}{1}{\su{0}}{\su{1}}\in\PFoddqa{\su{0}}$. For this reason, we will now consider a function $F\in\PFoddqa{\su{0}}$ and show in \rprop{L1513} that $\HTiaooo{F}{1}{\su{0}}{\su{1}}\in\RFuqa{1}{\seq{\su{j}}{j}{0}{1}}$. Finally, we investigate the case $m\in\mn{2}{+\infty}$. Let $\seq{\su{j}^\seins}{j}{0}{m-2}$ be the \tfirstSta{\seq{\su{j}}{j}{0}{m}}, then we know from \rtheoss{L1135}{L1314} that $\HTaooo{G}{1}{\su{0}}{\su{1}}\in\RFuqa{m-2}{\seq{\su{j}^\seins}{j}{0}{m-2}}$. So we are going to verify now that, for a function $F\in\RFuqa{m-2}{\seq{\su{j}^\seins}{j}{0}{m-2}}$, its \tiaaSt{\su{0}}{\su{1}} $\HTiaooo{F}{1}{\su{0}}{\su{1}}$ belongs to $\RFuqa{m}{\seq{\su{j}}{j}{0}{m}}$ (see \rtheoss{L1241}{L0925}).    

Now we turn our attention to a detailed treatment of the case $m=0$. An application of \rprop{L0800} will provide us quickly the desired result.

\bpropl{L1046}
	Let  $\su{0}\in\Cggq$ and let $F \in \PFevenqa{\su{0}}$. Then   	$\HTiaoo{F}{1}{\su{0}}$ belongs to $\ROqjj{0}{0}$.
\eprop
\bproof
	Since $F$ belongs to $\PFevenqa{\su{0}}$, \rpart{L1407.b} of \rlemm{L1407} shows that the function $F$ belongs to $\RqP$ and that the conditions $\betaF=0$ and
	\[
		\Kerna{\su{0}}
		\subseteq \Kerna{\alphaF}\cap\KernA{\nuFa{\R}}
		=\Kerna{\alphaF}\cap\Kerna{\betaF}\cap\KernA{\nuFa{\R}}
	\]
	are fulfilled. Setting $B\defg \Oqq$ and $t_0\defg \su{0}(\su{0}+\betaF)^\MP \su{0}$, we get $t_0=\su{0}$ and \rpart{L0800.d} of \rprop{L0800} shows that $\HTiaoo{F}{-1}{\su{0}}$ belongs to $\ROqjj{0}{0}$.
\eproof

\bcorol{C1103}
	Let $\seq{\su{j}}{j}{0}{0}\in\Hggequ{0}$. For each $F\in\RFuqa{0}{\seq{\su{j}}{j}{0}{0}}$, let
	\[
		\HTua{\su{0}}{F}
		\defg\HTaoo{F}{1}{\su{0}}.
	\]
	\index{F_(+;)()@$\HTua{\su{0}}{F}$}Then $\HTu{\su{0}}$ generates a bijective correspondence between the classes $\RFuqa{0}{\seq{\su{j}}{j}{0}{0}}$ and $\PFevenqa{\su{0}}$. The inverse mapping $(\HTu{\su{0}})^\inv$ is given for each $G\in\PFevenqa{\su{0}}$ by
	\[
		(\HTu{\su{0}})^\inv(G)
		\defg\HTiaoo{G}{1}{\su{0}}.
	\]
\ecoro
\bproof
	Taking \rpart{L0845.c} of \rcoro{L0845} and \rpart{C1519.b} of \rcoro{C1519} into account, the combination of \rtheo{L1013} and \rprop{L1046} yields all assertions.
\eproof

Now we consider the case $m=1$. First we state a more general result, which is of own interest.

\bpropl{0918}
	Let  $A\in\Cggq$ and  $B\in\CHq$ be such that $\Kerna{A}\subseteq\Kerna{B}$. Let $F\in \RqK{-1}$ be such that 
	\begin{equation} \label{MM1}
		\Kerna{A}
		\subseteq \Kerna{\gammaF}\cap\KernA{\muFa{\R}}.
	\end{equation}
	Then  $\HTiaooo{F}{1}{A}{B} \in \Ruua{1}{q}{\seq{t_j}{j}{0}{1}}$, where $t_0\defg A$ and where $t_1\defg B-\gammaF$. 
\eprop
\bproof
	Since $F$ belongs to $\NFuq{-1}$, we get from \eqref{Ruo[-1]} that
	\bgl{0918.1}
		F
		\in\NFq
	\eg
	and from \eqref{MM1} and \rlemm{L1406} that
	\begin{align}\label{0918.2}
		\Kerna{A}&\subseteq\KernA{F(z)}&\text{for each }z&\in\ohe.
	\end{align}
	Using \eqref{0918.1}, \eqref{0918.2} and \rprop{P1356}, we obtain
	\bgl{0918.3}
		\Kerna{A}
		\subseteq\Kerna{\alphaF}\cap\Kerna{\betaF}\cap\KernA{\nuFa{\R}}.
	\eg
	In view of \eqref{0918.1} and \eqref{0918.3}, we infer from \rpart{L0800.d} of \rprop{L0800} then
	\[
		\HTiaooo{F}{1}{A}{B}
		\in\NFq
	\]
	and from \eqref{18-3} that
	\begin{align}\label{0918.5}
		\det\lrk z\Iq+A^\MP\lek F(z)-B\rek\rrk&\neq0&\text{for each }z&\in\ohe.
	\end{align}
	Using \eqref{0918.1}, \eqref{0918.3}, and \rlemm{L1030}, we get
	\bgl{0918.6}
		A^\MP AF
		=F.
	\eg
	Since $A$ and $B$ are Hermitian matrices with $\Kerna{A}\subseteq\Kerna{B}$, the application of \rrema{L0818} yields
	\bgl{0918.7}
		AA^\MP B
		=B.
	\eg
	Setting $t_{-1}=\Oqq$ and using \eqref{Fuoa}, \eqref{F-ooa}, \eqref{0918.5}, \eqref{0918.6}, and \eqref{0918.7}, we get
	\[ 
		\begin{split}
			&z\lrk\trauoa{(\HTiaooo{F}{1}{A}{B})}{1}{t}{z}-\lek F(z)-\gammaF\rek\rrk\\
			&=z^3\HTiaoooa{F}{1}{A}{B}{z}+z^3t_{-1}+z^{2} t_0+z t_1-zF(z)+z\gammaF \\
			&=-z^3 A\lrk z \Iq + A^\MP\lek F(z)-B\rek\rrk^\MP +z^2 A +z\lek B-F(z)\rek\\
			&=-z^3 A\lrk z \Iq + A^\MP\lek F(z)-B\rek\rrk^\inv+z^2 A +z\lek B-F(z)\rek\\
			&=-z^2 A\lrk\Iq +\frac{1}{z}A^\MP\lek F(z)-B\rek\rrk^\inv+z^2 A +z\lek B-F(z)\rek\\
			&=\lek-z^2A+\lrk z^2 A +z\lek B-F(z)\rek\rrk\lrk\Iq + \frac{1}{z} A^\MP\lek F(z)-B\rek\rrk\rek\\
			&\quad\times\lrk\Iq + \frac{1}{z}A^\MP\lek F(z)-B\rek\rrk^{-1} \\
			&=\lrk z\lek B-F(z)\rek+zAA^\MP\lek F(z)-B\rek+\lek B-F(z)\rek A^\MP\lek F(z)-B\rek\rrk\\
			&\quad\times\lrk\Iq + \frac{1}{z}A^\MP\lek F(z)-B\rek\rrk^{-1} \\
			&=-\lek F(z)-B\rek A^\MP\lek F(z)-B\rek\lrk\Iq + \frac{1}{z}A^\MP\lek F(z)-B\rek\rrk^{-1}.
		\end{split} 
	\]
	This implies
	\begin{multline}\label{M14}
		\trauoa{(\HTiaooo{F}{1}{A}{B})}{1}{t}{z} \\
		=F(z)-\gammaF-\frac{1}{z}\lek F(z)-B\rek A^\MP\lek F(z)-B\rek\lrk\Iq + \frac{1}{z}A^\MP\lek F(z)-B\rek\rrk^{-1}.
	\end{multline}	
	Since $F$ belongs to $\RqK{-1}$, from \rprop{P1401} we see that \eqref{N15-3} is true. From \eqref{N15-3} we get
	\[ 
		\begin{split}
			\Iq
			&=\Iq+0\cdot A^\MP (\gammaF-B)
			=\Iq +\lrk\lim_{y\to+\infty} \frac{1}{\I y}\rrk\lrk A^\MP\lek\lim_{y\to+\infty}F(\I y)\rek-B\rrk\\
			&=\lim_{y\to+\infty} \lrk\Iq +\frac{1}{\I y}A^\MP\lek F(\I y)-B \rek\rrk
		\end{split} 
	\]
	and, consequently, 
	\begin{equation} 	\label{M3}
		\lim_{y\to+\infty}\lek\lrk\Iq +\frac{1}{\I y} A^\MP\lek F(\I y) -B \rek\rrk^{-1}\rek
		=\Iq
	\end{equation}	
	Using \eqref{N15-3}, \eqref{M3}, and \eqref{M14}, it follows
	\begin{equation} 	\label{Q1}
		\begin{split}
			\Oqq
			&=\gammaF -\gammaF - 0\cdot(\gammaF-B)A^\MP(\gammaF-B)\cdot \Iq \\
			&=\lek\lim_{y\to+\infty} F(\I y)\rek-\gammaF -\lrk\lim_{y\to+\infty} \frac{1}{\I y} \rrk\lek\lim_{y\to+\infty} F(\I y) -B \rek A^\MP\\
			&\quad\times\lek\lim_{y\to+\infty} F(\I y)-B\rek\lrk\lim_{y\to+\infty}\lek\lrk\Iq +\frac{1}{\I y} A^\MP\lek F(\I y)-B\rek\rrk^{-1}\rek\rrk\\
			&=\lim_{y\to+\infty} \Biggl[F(\I y)-\gammaF - \frac{1}{\I y}\lek F(\I y) -B\rek A^\MP\lek F(\I y)-B\rek\\
			&\quad\times\lrk\Iq +\frac{1}{\I y} A^\MP \lek F(\I y) -B\rek\rrk^{-1}\Biggr]\\
			&= \lim_{y\to+\infty} \trauoa{(\HTiaooo{F}{1}{A}{B})}{1}{t}{\I y}.
		\end{split}
	\end{equation}
	In view of $F\in\NFuq{-1}$, we get from \rrema{R1425} that
	\bgl{10.10.A}
		(\gammaF)^\ad
		=\gammaF.
	\eg
	Thus, we see from the definition of the sequence $\seq{t_j}{j}{-1}{1}$ that
	\bal{10.10.B}
		t_j^\ad&=t_j&\text{for each }j&\in\set{-1,0,1}.
	\ea
	Hence, taking \eqref{Q1} and \eqref{10.10.B} into account, we infer from \rtheo{P0611-2} that
	\begin{equation} 	\label{Q2}
		\trauo{(\HTiaooo{F}{1}{A}{B})}{1}{t}
		\in \RqP.
	\end{equation}
	Now we are going to show that
	\[
		\trauo{(\HTiaooo{F}{1}{A}{B})}{1}{t}
		=\NFuq{-1}.
	\]
	Since $F$ belongs to $\RqP$, the function $\Theta\colon[1,+\infty) \to \R$ defined by 
	\begin{equation} 	\label{FM1}
		\Theta(y)
		\defg\normA{\re\lek\lek F(\I y) -B\rek A^\MP\lek F(\I y)-B\rek\lrk\Iq +(\I y)^{-1} A^\MP \lek F(\I y)-B\rek\rrk^{-1}\rek}
	\end{equation}
	is continuous. From \eqref{N15-3} and \eqref{M3} we then get 
	\[ 
		\begin{split}
		 &\normA{\re\lek(\gammaF-B)A^\MP(\gammaF-B) \Iq\rek}\\
		 &=
		 \Biggl\| \re\Biggl[\lek\lim_{y\to+\infty}F(\I y) -B\rek A^\MP\lek\lim_{y\to+\infty}F(\I y) -B\rek \\
		 &\quad\times\lrk\lim_{y\to+\infty}\lek\lrk\Iq +(\I y)^{-1} A^\MP\lek F(\I y)-B\rek\rrk^{-1}\rek\rrk\Biggr] \Biggr\| \\
		 & =\lim_{y\to+\infty}\Theta(\I y).
		\end{split}
	\]
	Consequently, there is a non-negative real number $c$ such that
	\begin{align}\label{0918.12}
		\Abs{\Theta(y)}&\leq c&\text{for each }y&\in [1,+\infty).
	\end{align}
	Since $F \in \RqK{-1}$ is supposed, the function $\Psi\colon[1,+\infty) \to \R$ given by
	\begin{equation} \label{FM2}
		\Psi(y)
		\defg \frac{1}{y}\normA{\im F(\I y)} + \frac{c}{y^2}
	\end{equation}
	fulfills
	\begin{equation} \label{NM}
		\Psi
		\in \LoAAAA{1}{[1,+\infty)}{\BA{[1,+\infty)}}{\tilde{\Leb}}{\R},
	\end{equation}
	where $\tilde{\Leb}$ is the Lebesgue measure defined on $\BA{[1,+\infty)}$. For all $y \in [1,+\infty)$ from \eqref{M14}, \eqref{10.10.A}, \eqref{FM1}, \eqref{0918.12}, and \eqref{FM2} we see that $\Omega\colon[1,+\infty) \to \R $ given by 
 	\[ 
		\Omega(y)
		\defg \frac{1}{y}\normA{\im  \trauoa{(\HTiaooo{F}{1}{A}{B})}{1}{t}{\I y}}
	\]
	satisfies, for $y\in[1,+\infty)$, the inequality
	\begin{equation} 	\label{FM3}
		\begin{split}
			&\Abs{\Omega(y)}\\
			&= \frac{1}{y}\biggl\lVert\im F(\I y)\\
			&\quad\quad+  \frac{1}{y}\im\lek\I\lek F(\I y) -B\rek A^\MP\lek F(\I y)-B\rek\lrk\Iq +(\I y)^{-1} A^\MP\lek F(\I y)-B\rek\rrk^{-1}\rek\biggr\rVert\\
			& \leq  \frac{1}{y}\normA{\im F(\I y)}+ \Theta(y)
			\leq \Psi(y).
		\end{split}
	\end{equation}
	Since \eqref{Q2} holds, the function $\Omega$ is continuous and thus Borel-measurable. Hence, \eqref{FM3} and \eqref{NM} yield that $\Omega$ also belongs to $\Loaaaa{1}{[1,+\infty)}{\BA{[1,+\infty)}}{\tilde{\lambda}}{\R}$. Keeping in mind \eqref{Q2}, this implies $ \trauo{(\HTiaooo{F}{1}{A}{B})}{1}{t} \in \RqK{-1} $. Thus, combining this with \eqref{10.10.B}, \eqref{0918.1}, and \eqref{Q1}, we infer from \rtheo{P0611-2} that $\HTiaooo{F}{1}{A}{B}\in\Ruoa{q}{1}{\seq{t_j}{j}{-1}{1}}$. In view of $t_{-1}=\Oqq$ and \rrema{R1404}, this implies that $\HTiaooo{F}{1}{A}{B}\in\Ruua{1}{q}{\seq{t_j}{j}{0}{1}}$.
\eproof

The following result gives a complete answer to the case $m=1$.

\bpropl{L1513}
	Let  $\seq{\su{j}}{j}{0}{1}\in\Hggeuu{q}{1}$ and let $F \in \PFoddqa{\su{0}}$. Then $\HTiaooo{F}{1}{\su{0}}{\su{1}}$ belongs to $\Rkqjj{1}{0}{1}$.
\eprop
\bproof
	Because of $F \in \PFoddqa{\su{0}}$, we get from \eqref{Poddua} that $F \in \Rkq{-1}$ and $\Kerna{\su{0}}\subseteq \Kerna{\muFa{\R}}$. Hence, \eqref{Ruu} yields $F\in\RqK{-1}$ and $\gammaF=\Oqq$. Thus, we get $\Kerna{\su{0}}\subseteq \Kerna{\gammaF}\cap\Kerna{\muFa{\R}}$. Since $\seq{\su{j}}{j}{0}{1}$ belongs to $\Hggeuu{q}{1}$, we see from \rlemm{L1037} that $\su{0}\in \Cggq$, $\su{1}^\ad=\su{1}$, and $\Kerna{\su{0}}\subseteq \Kerna{\su{1}}$. Applying \rprop{0918} and using $\su{1}-\gammaF=\su{1}$, we conclude $\HTiaooo{F}{1}{\su{0}}{\su{1}} \in \Rkqjj{1}{0}{1}$.
\eproof

\bcorol{C1514}
	Let $\seq{\su{j}}{j}{0}{1}\in\Hggequ{1}$. For $F\in\RFuqa{1}{\seq{\su{j}}{j}{0}{1}}$ let
	\[
		\HTuua{\su{0}}{\su{1}}{F}
		\defg\HTaooo{F}{1}{\su{0}}{\su{1}}.
	\]
	\index{F_(+;,)()@$\HTuua{\su{0}}{\su{1}}{F}$}Then $\HTuu{\su{0}}{\su{1}}$ generates a bijective correspondence between the classes $\RFuqa{1}{\seq{\su{j}}{j}{0}{1}}$ and $\PFoddqa{\su{0}}$. The inverse mapping $(\HTuu{\su{0}}{\su{1}})^\inv$ is given for $G\in\PFoddqa{\su{0}}$ by
	\[
		(\HTuu{\su{0}}{\su{1}})^\inv(G)
		\defg\HTiaooo{G}{1}{\su{0}}{\su{1}}.
	\]
\ecoro
\bproof
	Taking \rpart{L0845.c} of \rcoro{L0845} and \rpart{C1519.b} of \rcoro{C1519} into account, the combination of \rtheo{L1015} and \rprop{L1513} completes the proof.
\eproof

Now we turn our attention to the case $m\in\mn{2}{+\infty}$. Similar as in \rsect{S1326}, we will use various Hamburger-Nevanlinna type results from \rsect{S1510}. In order to prepare the application of this material, we still need some auxiliary results, which can be considered as analogues of \rlemmss{L1705}{L1627}. First we will compute the functions introduced in \rrema{R1013} for the case that the function $F$ is replaced by $\HTiaooo{F}{1}{\su{0}}{\su{1}}$.

\blemml{L1651-1}
	Let $\mG$ be a non-empty subset of  $\C\setminus\set{0}$, let $F\colon \mG \to \Cpq$ be a matrix-valued function, let  $\kappa\in\N\cup\set{+\infty}$, and let $\seq{\su{j}}{j}{0}{\kappa}$ be a sequence of complex \tpqa{matrices}. In the case $\kappa\geq2$ let $\seq{\su{j}^\seins}{j}{0}{\kappa-2}$ be the \tfirstSta{\seq{\su{j}}{j}{0}{\kappa}}. Further, let $\su{-1}\defg\Opq$,  let $\su{-1}^\seins\defg\Opq$, let  $m\in\mn{1}{\kappa}$, and let $\nablau{m}\colon\mG\to\Cpq$\index{nabla_@$\nablau{m}$} be defined by
	\begin{equation} \label{nablaua}
		\nablaua{m}{z}
		\defg
		\begin{cases}
			[F(z)-\su{1}]\su{0}^\MP[F(z)-\su{1}]\incase{m=1}\\
			\\
			\begin{aligned}
				&\trauoa{F}{m-2}{s^\seins}{z} \su{0}^\MP [F(z)-\su{1}]-
				\sum\nolimits_{j=0}^{m-1}z^{-j}\su{j+1}\su{0}^\MP \trauoa{F}{m-2}{s^\seins}{z} \\
				&+\su{m}\su{0}^\MP\su{1}+\sum\nolimits_{k=1}^{m-1} \su{m-k}\su{0}^\MP\su{k-1}^\seins\\
				&+\sum\nolimits_{j=1}^{m-1}z^{-j}\sum\nolimits_{k=j}^{m-1} \su{m+j-k}\su{0}^\MP\su{k-1}^\seins
			\end{aligned}
			\incase{m\geq2}
		\end{cases}.
	\end{equation}
	Suppose that $\seq{\su{j}}{j}{0}{m}\in\Dpqu{m}$ and that $z\in\mG$ is such that $\Bilda{F(z)}\subseteq\Bilda{\su{0}}$ and 
	\begin{equation} \label{KS}
		\det\lrk\Iq +z^{-1} \su{0}^\MP\lek F(z)-\su{1}\rek\rrk
		\neq 0
	\end{equation}
	are fulfilled.  Then
	\begin{equation} \label{TM1}
		\trauoa{( \HTiaooo{F}{1}{\su{0}}{\su{1}})}{m}{s}{z}
		=\trauoa{F}{m-2}{s^\seins}{z}-\frac{1}{z} \nablaua{m}{z}\lrk\Iq +z^{-1} \su{0}^\MP \lek F(z)-\su{1}\rek\rrk^{-1}. 
	\end{equation}
\elemm
\bproof 
	In view of \rrema{R1013}, we have 
	\begin{equation} \label{B3}
		z^{-(m-1)}\trauoa{F}{m-2}{s^\seins}{z} -\sum_{k=0}^{m-1}z^{-k}s_{k-1}^{\seins}
		=F(z). 
	\end{equation}
	The assumption $\Bilda{F(z)}\subseteq\Bilda{\su{0}}$  and \rpart{L1622.b} of \rrema{L1622} yield $ \su{0}\su{0}^\MP F(z)=F(z)$. From~\cite{103}*{\crema{8.5}} we know that \eqref{I-4} holds true for all $j\in\mn{-1}{m-2}$. Thus, we get 
	\begin{equation} \label{B6}
		\su{0}\su{0}^\MP \trauoa{F}{m-2}{s^\seins}{z}
		= \trauoa{F}{m-2}{s^\seins}{z}.
	\end{equation}	
	Because of $\seq{\su{j}}{j}{0}{m} \in \Dpqu{m}$, \rdefi{D1043}, and \rpart{L1622.b} of \rrema{L1622}, the first equation in \eqref{I-1} is fulfilled for all $j\in\mn{0}{m}$. Using $\su{-1}=\Opq$, $\su{-1}^\seins=\Opq$, \eqref{Fuoa}, \eqref{F-ooa}, \eqref{KS}, \eqref{B3}, and \eqref{B6}, we conclude
	\[
		\begin{split}
			&\trauoa{( \HTiaooo{F}{1}{\su{0}}{\su{1}})}{m}{s}{z}
			=z^{m+1}\lek\HTiaooo{F}{1}{\su{0}}{\su{1}}(z)+ \sum_{j=0}^{m+1}z^{-j}s_{j-1}\rek\\
			&=z^{m+1} \Biggl[ -\su{0}\lrk z\Iq + \su{0}^\MP\lek F(z)-\su{1}\rek\rrk^{\MP} \\
			&\quad\quad\quad\quad+\sum_{j=0}^{m+1}z^{-j}s_{j-1} \lrk z\Iq + \su{0}^\MP\lek F(z)-\su{1}\rek\rrk\lrk z\Iq +\su{0}^\MP\lek F(z)-\su{1}\rek\rrk^{-1} \Biggr]
			\\
			&=\lek-z^{m+1}\su{0}+\sum_{j=0}^{m+1}z^{m+1-j}s_{j-1}\lrk z\Iq + \su{0}^\MP \lek F(z)-\su{1}\rek\rrk\rek\lrk z\Iq + \su{0}^\MP\lek F(z)-\su{1}\rek\rrk^{-1} \\
			&=  \Biggl[ -z^{m+1}\su{0}\\
			&\quad\quad+\sum_{j=0}^{m+1}z^{m+1-j}s_{j-1}\lrk z\Iq + \su{0}^\MP\lek z^{-(m-1)}\trauoa{F}{m-2}{s^\seins}{z}-\sum_{k=0}^{m-1}z^{-k}s_{k-1}^{\seins}\rek-\su{0}^\MP \su{1}\rrk\Biggr]\\
			&\quad\times\lrk z\Iq +\su{0}^\MP \lek F(z)-\su{1}\rek\rrk^{-1} \\
			&=  \Biggl[ -z^{m}\su{0} 
			+\sum_{j=0}^{m+1}z^{m+1-j}s_{j-1} +
			\sum_{j=0}^{m+1}z^{-j+1} \su{j-1} \su{0}^\MP \trauoa{F}{m-2}{s^\seins}{z} \\
			&\quad\quad- \sum_{j=0}^{m+1}\sum_{k=0}^{m-1}z^{m-(j+k)}s_{j-1} \su{0}^\MP \su{k-1}^{\seins}
			-\sum_{j=0}^{m+1} z^{m-j}\su{j-1}\su{0}^\MP \su{1}  \Biggr]\\
			&\quad\times\lrk\Iq +z^{-1} \su{0}^\MP \lek F(z)-\su{1}\rek\rrk^{-1}
		\end{split}
	\]
	\begin{equation}\label{9-1}
		\begin{split}
			&=  \Biggl[ -z^{m}\su{0}+z^{m+1}\su{-1}+z^{m}\su{0} 
			+\sum_{j=2}^{m+1}z^{m+1-j}s_{j-1} 
			+z^{1}\su{-1} \su{0}^\MP \trauoa{F}{m-2}{s^\seins}{z} \\
			&\quad\quad+z^{0}\su{0} \su{0}^\MP \trauoa{F}{m-2}{s^\seins}{z}
			+	\sum_{j=2}^{m+1}z^{-j+1} \su{j-1} \su{0}^\MP \trauoa{F}{m-2}{s^\seins}{z} 
			- \sum_{k=0}^{m-1}z^{m-k}\su{-1} \su{0}^\MP \su{k-1}^{\seins} \\
			&\quad\quad-\sum_{k=0}^{m-1}z^{m-1-k}s_{0} \su{0}^\MP \su{k-1}^{\seins}
			-\sum_{j=2}^{m+1}\sum_{k=0}^{m-1} z^{m-(j+k)}\su{j-1}\su{0}^\MP \su{k-1} ^{\seins} 
			-z^{m}\su{-1} \su{0}^\MP \su{1}
				\\
			&\quad\quad-z^{m-1}\su{0} \su{0}^\MP \su{1}-\sum_{j=2}^{m+1} z^{m-j}\su{j-1}\su{0}^\MP \su{1}  
			\Biggr]\lrk\Iq +z^{-1} \su{0}^\MP \lek F(z)-\su{1}\rek\rrk^{-1}\\
			&=  \Biggl[ \sum_{j=2}^{m+1}z^{m+1-j}s_{j-1} 
			\trauoa{F}{m-2}{s^\seins}{z} 
			+	\sum_{j=2}^{m+1}z^{-j+1} \su{j-1} \su{0}^\MP \trauoa{F}{m-2}{s^\seins}{z} \\
			&\quad\quad- \sum_{k=0}^{m-1}z^{m-1-k} \su{k-1}^{\seins} 
			-\sum_{j=2}^{m+1}\sum_{k=0}^{m-1} z^{m-(j+k)}\su{j-1}\su{0}^\MP \su{k-1} ^{\seins} 
			-z\su{1}  \\
			&\quad\quad-\sum_{j=2}^{m+1} z^{m-j}\su{j-1}\su{0}^\MP \su{1}  
			\Biggr]\lrk \Iq +z^{-1} \su{0}^\MP \lek F(z)-\su{1}\rek\rrk^{-1}.
		\end{split} 
	\end{equation}
	In the case $m=1$, from \eqref{9-1},  $\su{-1}^{\seins}=\Opq$, and   $\trauoa{F}{-1}{s^\seins}{z}=F(z)$, it follows 	
	\[ 
		\begin{split}
			\trauoa{( \HTiaooo{F}{1}{\su{0}}{\su{1}})}{m}{s}{z} 
			&=\trauoa{(\HTiaooo{F}{1}{\su{0}}{\su{1}})}{1}{s}{z}\\
			&=\lek F(z)+z^{-1} \su{1}\su{0}^\MP F(z) - z^{-1} \su{1}\su{0}^\MP \su{1}\rek\lrk\Iq +z^{-1} \su{0}^\MP\lek F(z)-\su{1}\rek\rrk^{-1}
		\end{split} 
	\]
	and, consequently,
	\bsp
		&z\lek\trauoa{F}{m-2}{s^\seins}{z} -	\trauoa{( \HTiaooo{F}{1}{\su{0}}{\su{1}})}{m}{s}{z}\rek\lrk\Iq +z^{-1} \su{0}^\MP \lek F(z)-\su{1}\rek\rrk\\
		&= z\lrk F(z) -\lek F(z)+z^{-1} \su{1}\su{0}^\MP F(z) - z^{-1} \su{1}\su{0}^\MP \su{1}\rek\lrk\Iq +z^{-1} \su{0}^\MP \lek F(z)-\su{1}\rek\rrk^{-1}\rek\\
		&\quad\times\lrk\Iq +z^{-1} \su{0}^\MP \lek F(z)-\su{1}\rek\rrk\\
		&= F(z) \su{0}^\MP \lek F(z) -\su{1}\rek-\su{1}\su{0}^\MP \lek F(z) -\su{1}\rek 
		= \nablaua{1}{z}
		= \nablaua{m}{z}.
	\esp
	Thus, \eqref{TM1} is proved in the case $m=1$.  If $m \geq 2$, then \eqref{9-1} and $\su{-1}^{\seins}=\Opq$ yield
	\begin{equation} \label{15-1}
		\begin{split}
			&\trauoa{( \HTiaooo{F}{1}{\su{0}}{\su{1}})}{m}{s}{z} \\
			&=  \Biggl[ \sum_{j=3}^{m+1}z^{m+1-j}s_{j-1} 
			\trauoa{F}{m-2}{s^\seins}{z} 
			+	\sum_{j=0}^{m-1}z^{-j+1} \su{j+1} \su{0}^\MP  \trauoa{F}{m-2}{s^\seins}{z} \\
			&\quad- \sum_{k=1}^{m-1}z^{m-1-k} \su{k-1}^{\seins} 
			-\sum_{j=2}^{m+1}\sum_{k=1}^{m-1} z^{m-(j+k)}\su{j-1}\su{0}^\MP \su{k-1} ^{\seins} 
			-z^{m-2}\su{1}\su{0}^\MP \su{1}  \\
			&\quad-\sum_{j=3}^{m+1} z^{m-j}\su{j-1}\su{0}^\MP \su{1}  
			\Biggl]
			\lrk\Iq +z^{-1} \su{0}^\MP \lek F(z)-\su{1}\rek\rrk^{-1}.
		\end{split} 
	\end{equation}
	Furthermore, if $m \geq2$, then~\cite{103}*{\cprop{8.23}} shows that \eqref{I-19A} is valid. Thus, in the case $m=2$, from \eqref{15-1} and \eqref{I-19A} we get 
	\[ 
		\begin{split}
			&\trauoa{( \HTiaooo{F}{1}{\su{0}}{\su{1}})}{m}{s}{z}
			=\trauoa{( \HTiaooo{F}{1}{\su{0}}{\su{1}})}{2}{s}{z} \\
			&=  \Biggl[ \su{2}+  \trauoa{F}{m-2}{s^\seins}{z}  
			+	\sum_{j=0}^{m-1}z^{-j-1} \su{j+1} \su{0}^\MP  \trauoa{F}{m-2}{s^\seins}{z} - \su{0}^{\seins}\\
			&\quad- z^{-1} \su{1}\su{0}^\MP\su{0}^{\seins}
			-  z^{-2} \su{2}\su{0}^\MP\su{0}^{\seins}
			-\su{1}\su{0}^\MP \su{1}
			-z^{-1}\su{2}\su{0}^\MP \su{1} 
			\Biggr]\lrk \Iq +z^{-1} \su{0}^\MP\lek F(z)-\su{1}\rek\rrk^{-1} \\
			&=  \Biggl[ \trauoa{F}{m-2}{s^\seins}{z}  
			+	\sum_{j=0}^{m-1}z^{-j-1} \su{j+1} \su{0}^\MP  \trauoa{F}{m-2}{s^\seins}{z} \\
			&\quad-  z^{-1}( \su{2}\su{0}^\MP\su{1} + \su{1}\su{0}^\MP\su{0}^{\seins})
			-z^{-2}\su{2}\su{0}^\MP \su{0}^{\seins} 
			\Biggr] 
			\lrk\Iq +z^{-1} \su{0}^\MP \lek F(z)-\su{1}\rek\rrk^{-1} \\
			&=  \Biggl[ \trauoa{F}{m-2}{s^\seins}{z}  
			+	\sum_{j=0}^{m-1}z^{-j-1} \su{j+1} \su{0}^\MP  \trauoa{F}{m-2}{s^\seins}{z} 
			-  z^{-1}\lrk\su{2}\su{0}^\MP\su{1}  
			+  \sum_{k=1}^{1} \su{2-k} \su{0}^\MP \su{k-1}^{\seins}\rrk\\
			&\quad- \sum_{j=1}^{1}  z^{-j-1} \sum_{k=j}^{1} \su{2+j-k} \su{0}^\MP \su{k-1}^{\seins} 
			\Biggr]
			\lrk\Iq +z^{-1} \su{0}^\MP \lek F(z)-\su{1}\rek\rrk^{-1}
		\end{split} 
	\]
	and, consequently,
	\begin{equation} \label{16-1}
		\begin{split}
			&\trauoa{( \HTiaooo{F}{1}{\su{0}}{\su{1}})}{m}{s}{z} \\
			&=  \Biggl[ \trauoa{F}{m-2}{s^\seins}{z}  
			+	\sum_{j=0}^{m-1}z^{-j-1} \su{j+1} \su{0}^\MP  \trauoa{F}{m-2}{s^\seins}{z} 
			-  z^{-1}  \lrk\su{m}\su{0}^\MP\su{1}  
			+  \sum_{k=1}^{m-1} \su{m-k} \su{0}^\MP \su{k-1}^{\seins}\rrk\\
			&\quad- \sum_{j=1}^{m-1}z^{-j-1} \sum_{k=j}^{m-1} \su{m+j-k} \su{0}^\MP \su{k-1}^{\seins} 
			\Biggr]\lrk\Iq +z^{-1} \su{0}^\MP \lek F(z)-\su{1}\rek\rrk^{-1}.
		\end{split} 
	\end{equation}
	In case $m \geq 3$, we have
 	\begin{equation} \label{19-1-1}
		\begin{split}
			& \sum_{j=2}^{m+1} \sum_{k=1}^{m-1} z^{m-(j+k)} \su{j-1} \su{0}^\MP \su{k-1}^{\seins} \\
			&= \sum_{j=2}^{m-1} \sum_{k=1}^{m-1} z^{m-(j+k)} \su{j-1} \su{0}^\MP \su{k-1}^{\seins} 
			+  \sum_{k=1}^{m-1} z^{-k} \su{m-1} \su{0}^\MP \su{k-1}^{\seins}
			+  \sum_{k=1}^{m-1} z^{-(k+1)} \su{m} \su{0}^\MP \su{k-1}^{\seins}\\ 
			&= \sum_{l=3}^{2m-2} \sum_{k=\max\left\{1,l-(m-1)\right\}}^{\min\left\{m-1,l-2\right\}} z^{m-l} \su{l-k-1} \su{0}^\MP \su{k-1}^{\seins} 
			+z^{-1}\su{m-1} \su{0}^\MP \su{0}^{\seins} \\
			&\quad+  \sum_{k=2}^{m-1} z^{-k} \su{m-1} \su{0}^\MP \su{k-1}^{\seins}
			+  \sum_{k=1}^{m-2} z^{-(k+1)} \su{m} \su{0}^\MP
				\su{k-1}^{\seins}
				+z^{-m} \su{m} \su{0}^\MP  \su{m-2}^{\seins} \\
			&= \sum_{l=3}^{2m-2} \sum_{k=\max\left\{m,l\right\}-(m-1)}^{\min\left\{m+1,l\right\}-2} z^{m-l} \su{l-k-1} \su{0}^\MP \su{k-1}^{\seins} 
			+z^{-1}\su{m-1} \su{0}^\MP \su{0}^{\seins} \\
			&\quad+  \sum_{k=2}^{m-1} z^{-k} \su{m-1} \su{0}^\MP \su{k-1}^{\seins}
			+  \sum_{k=2}^{m-1} z^{-k} \su{m} \su{0}^\MP
				\su{k-2}^{\seins}
				+z^{-m} \su{m} \su{0}^\MP \su{m-2}^{\seins}  \\
			&= \sum_{l=3}^{m} \sum_{k=1}^{l-2} z^{m-l} \su{l-k-1} \su{0}^\MP \su{k-1}^{\seins} 
			+\sum_{l=m+1}^{2m-2} \sum_{k=l-(m-1)}^{m-1} z^{m-l} \su{l-k-1} \su{0}^\MP \su{k-1}^{\seins}
				\\
			&\quad+z^{-1}\su{m-1} \su{0}^\MP \su{0}^{\seins}
			+  \sum_{k=2}^{m-1} z^{-k}( \su{m-1} \su{0}^\MP \su{k-1}^{\seins} + \su{m} \su{0}^\MP  \su{k-2}^{\seins})
				+z^{-m} \su{m} \su{0}^\MP \su{m-2}^{\seins}  \\
			&= \sum_{l=3}^{m}\lrk z^{m-l} \sum_{k=0}^{(l-2)-1} \su{(l-2)-k} \su{0}^\MP \su{k}^{\seins}\rrk+\sum_{j=0}^{m-3}\lrk z^{-j-1} \sum_{k=j+2}^{m-1}  \su{j+m-k} \su{0}^\MP \su{k-1}^{\seins}\rrk\\
			&\quad+z^{-1}\su{m-1} \su{0}^\MP \su{0}^{\seins} 
			+  \sum_{j=1}^{m-2} z^{-j-1}( \su{m-1} \su{0}^\MP \su{j}^{\seins} + \su{m} \su{0}^\MP  \su{j-1}^{\seins})
				+z^{-m} \su{m} \su{0}^\MP \su{m-2}^{\seins} 
		\end{split} 
	\end{equation}
	and, taking into account \eqref{15-1} and  \eqref{19-1-1}, furthermore,
 	\begin{equation} \label{20.1}
		\begin{split}
			&\trauoa{( \HTiaooo{F}{1}{\su{0}}{\su{1}})}{m}{s}{z} \\
			&=  \Biggl[z^{m-2} \su{2} + \sum_{j=4}^{m+1}z^{m+1-j} \su{j-1}+ \trauoa{F}{m-2}{s^\seins}{z}  
			+	\sum_{j=0}^{m-1}z^{-j-1} \su{j+1} \su{0}^\MP  \trauoa{F}{m-2}{s^\seins}{z}  \\
			&\quad-  z^{m-2}  \su{0}^{\seins} 
			-	\sum_{k=2}^{m-1}z^{m-1-k} \su{k-1}^{\seins} 
			- \sum_{l=3}^{m} \lrk z^{m-l} \sum_{k=0}^{(l-2)-1} \su{(l-2)-k} \su{0}^\MP \su{k}^{\seins}\rrk\\
			&\quad-\sum_{j=0}^{m-3}\lrk z^{-j-1} \sum_{k=j+2}^{m-1}  \su{j+m-k} \su{0}^\MP \su{k-1}^{\seins}\rrk
			-z^{-1}\su{m-1} \su{0}^\MP \su{0}^{\seins} \\
			&\quad-  \sum_{j=1}^{m-2} z^{-j-1}( \su{m-1} \su{0}^\MP \su{j}^{\seins}+ \su{m} \su{0}^\MP  \su{j-1}^{\seins}) 
			- z^{-m} \su{m} \su{0}^\MP \su{m-2}^{\seins} 
			-z^{m-2}\su{1}\su{0}^\MP \su{1}\\
			&\quad-	\sum_{j=3}^{m}z^{m-j} \su{j-1} \su{0}^\MP \su{1} - z^{-1} \su{m}\su{0}^\MP \su{1}
			\Biggr]\lrk\Iq +z^{-1} \su{0}^\MP \lek F(z)-\su{1}\rek\rrk^{-1}\\
			&=  \Biggl[z^{m-2} \su{2} + \sum_{l=3}^{m}z^{m-l} \su{l}+ \trauoa{F}{m-2}{s^\seins}{z}  
			+	\sum_{j=0}^{m-1}z^{-j-1} \su{j+1} \su{0}^\MP  \trauoa{F}{m-2}{s^\seins}{z}  \\
			&\quad-  z^{m-2}  ( \su{2}- \su{1}\su{0}^\MP\su{1})
			-	\sum_{l=3}^{m}z^{m-l} \su{l-2}^{\seins} 
			- \sum_{l=3}^{m} \lrk z^{m-l} \sum_{k=0}^{(l-2)-1} \su{(l-2)-k} \su{0}^\MP \su{k}^{\seins}\rrk\\
			&\quad-\sum_{j=0}^{m-3}\lrk z^{-j-1} \sum_{k=j+2}^{m-1}  \su{j+m-k} \su{0}^\MP \su{k-1}^{\seins}\rrk
			-z^{-1}(\su{m-1} \su{0}^\MP \su{0}^{\seins} +\su{m}\su{0}^\MP\su{1}) \\
			&\quad-  \sum_{j=1}^{m-2} z^{-j-1}( \su{m-1} \su{0}^\MP \su{j}^{\seins}
			+ \su{m} \su{0}^\MP  \su{j-1}^{\seins}) 
			- z^{-m} \su{m} \su{0}^\MP \su{m-2}^{\seins} 
			- z^{m-2} \su{1} \su{0}^\MP \su{1} \\
			&\quad- 	 \sum_{l=3}^{m} z^{m-l} \su{l-1} \su{0}^\MP \su{1}
			\Biggr]
			\lrk\Iq +z^{-1} \su{0}^\MP \lek F(z)-\su{1}\rek\rrk^{-1}	 
			\\
			&=  \Biggl[\trauoa{F}{m-2}{s^\seins}{z} 
			+\sum_{j=0}^{m-1}z^{-j-1} \su{j+1} \su{0}^\MP  \trauoa{F}{m-2}{s^\seins}{z} \\
			&\quad-\sum_{j=0}^{m-3} z^{-j-1} \sum_{k=j+2}^{m-1}  \su{j+m-k} \su{0}^\MP \su{k-1}^{\seins} 
				-z^{-1} (\su{m-1} \su{0}^\MP \su{0}^{\seins} +\su{m}\su{0}^\MP\su{1} ) \\	 	
			&\quad-  \sum_{j=1}^{m-2} z^{-j-1} ( \su{m-1} \su{0}^\MP \su{j}^{\seins} 	 + \su{m} \su{0}^\MP  \su{j-1}^{\seins} ) 
			- z^{-m} \su{m} \su{0}^\MP \su{m-2}^{\seins}  
			\Biggr] \\
			&\quad\times\lrk\Iq +z^{-1} \su{0}^\MP \lek F(z)-\su{1}\rek\rrk^{-1}	  
		\end{split} 
	\end{equation}
	Thus, in the case $m=3$, from \eqref{20.1} we conclude 
 	\begin{equation} \label{22-1}
		\begin{split}
			&\trauoa{( \HTiaooo{F}{1}{\su{0}}{\su{1}})}{m}{s}{z} 
			=
			\trauoa{( \HTiaooo{F}{1}{\su{0}}{\su{1}})}{3}{s}{z} \\
				&=  \Biggl[\trauoa{F}{m-2}{s^\seins}{z} 
				+\sum_{j=0}^{m-1}z^{-j-1} \su{j+1} \su{0}^\MP  \trauoa{F}{m-2}{s^\seins}{z} 
				-z^{-1}\sum_{k=2}^{m-1} \su{m-k} \su{0}^\MP \su{k-1}^{\seins} \\
				&\quad-z^{-1} (\su{m-1} \su{0}^\MP \su{0}^{\seins} +\su{m}\su{0}^\MP\su{1})
				-z^{-2}(\su{2} \su{0}^\MP \su{1}^{\seins} +\su{3}\su{0}^\MP\su{0}^{\seins}) - z^{-3} \su{3} \su{0}^\MP \su{1}^{\seins}  
			\Biggr] \\
			&\quad\times\lrk \Iq +z^{-1} \su{0}^\MP \lek F(z)-\su{1}\rek\rrk^{-1}	 
			\\
			& 	 =  \Biggl[\trauoa{F}{m-2}{s^\seins}{z} 
				+\sum_{j=0}^{m-1}z^{-j-1} \su{j+1} \su{0}^\MP  \trauoa{F}{m-2}{s^\seins}{z} \\
			&\quad-z^{-1} \lrk\su{m} \su{0}^\MP \su{0}^{\seins} +
				\sum_{k=1}^{m-1} \su{m-k} \su{0}^\MP \su{k-1}^{\seins} \rrk
				-\sum_{j=1}^{2}\lrk z^{-j-1} \sum_{k=j}^{2}\su{3-j-k} \su{0}^\MP \su{k-1}^{\seins} \rrk
			\Biggr] 	\\
			&\quad\times\lrk \Iq +z^{-1} \su{0}^\MP \lek F(z)-\su{1}\rek\rrk^{-1}	 
		\end{split} 
	\end{equation}
	and, consequently, that \eqref{16-1} holds true. If $m\geq 4$, then \eqref{20.1} implies
 	\[ 
		\begin{split}
			&\trauoa{( \HTiaooo{F}{1}{\su{0}}{\su{1}})}{m}{s}{z} 
			\\
			&=  \Biggl[\trauoa{F}{m-2}{s^\seins}{z} 
			+\sum_{j=0}^{m-1}z^{-j-1} \su{j+1} \su{0}^\MP  \trauoa{F}{m-2}{s^\seins}{z} 
			-z^{-1}\sum_{k=2}^{m-1} \su{m-k} \su{0}^\MP \su{k-1}^{\seins} \\
			&\quad-\sum_{j=1}^{m-3}\lrk z^{-j-1} \sum_{k=j+2}^{m-1}\su{m+j-k} \su{0}^\MP \su{k-1}^{\seins}\rrk-z^{-1}(\su{m-1} \su{0}^\MP \su{0}^{\seins} +\su{m}\su{0}^\MP\su{1} ) \\
			&\quad-\sum_{j=1}^{m-2}\lrk z^{-j-1} \sum_{k=j}^{j+1}\su{m+j-k} \su{0}^\MP \su{k-1}^{\seins}\rrk
			-\sum_{j=m-1}^{m-1}\lrk z^{-j-1} \sum_{k=m-1}^{m-1}\su{m+j-k} \su{0}^\MP \su{k-1}^{\seins}\rrk
			\Biggr] \\
			&\quad\times\lrk\Iq +z^{-1} \su{0}^\MP\lek F(z)-\su{1}\rek\rrk^{-1}	 
			\\
			&=  \Biggl[\trauoa{F}{m-2}{s^\seins}{z} 
			+\sum_{j=0}^{m-1}z^{-j-1} \su{j+1} \su{0}^\MP  \trauoa{F}{m-2}{s^\seins}{z} 
			-z^{-1} \lrk\su{m} \su{0}^\MP \su{1} +\sum_{k=1}^{m-1} \su{m-k} \su{0}^\MP \su{k-1}^{\seins}\rrk\\
			&\quad-\sum_{j=1}^{m-3}\lrk z^{-j-1} \sum_{k=j}^{m-1}\su{m+j-k} \su{0}^\MP \su{k-1}^{\seins} \rrk-\sum_{j=m-2}^{m-2}\lrk z^{-j-1} \sum_{k=j}^{m-1}\su{m+j-k} \su{0}^\MP \su{k-1}^{\seins}\rrk\\
			&\quad-\sum_{j=m-1}^{m-1}\lrk z^{-j-1} \sum_{k=j}^{m-1}\su{m+j-k} \su{0}^\MP \su{k-1}^{\seins}\rrk
			\Biggr]\lrk\Iq +z^{-1} \su{0}^\MP\lek F(z)-\su{1}\rek\rrk^{-1}	 
		\end{split} 
	\]
	and, consequently, \eqref{16-1}. Hence, \eqref{22-1} holds true if $m\geq2$. Thus, in the case $m\geq2$, we get from \eqref{16-1} and \eqref{nablaua} that
 	\[ 
		\begin{split}
			&z\lek\trauoa{F}{m-2}{s^\seins}{z} -\trauoa{( \HTiaooo{F}{1}{\su{0}}{\su{1}})}{m}{s}{z}\rek\lrk\Iq +z^{-1} \su{0}^\MP\lek F(z)-\su{1}\rek\rrk\\
			&=z \Biggl[ \trauoa{F}{m-2}{s^\seins}{z}  
			-\Biggl(\trauoa{F}{m-2}{s^\seins}{z} 
			+\sum_{j=0}^{m-1}z^{-j-1} \su{j+1} \su{0}^\MP  \trauoa{F}{m-2}{s^\seins}{z} \\
			&\quad-z^{-1}\lrk\su{m} \su{0}^\MP \su{1} +\sum_{k=1}^{m-1} \su{m-k} \su{0}^\MP \su{k-1}^{\seins}\rrk-\sum_{j=1}^{m-1}\lrk z^{-j-1} \sum_{k=j}^{m-1}\su{m+j-k} \su{0}^\MP \su{k-1}^{\seins}\rrk\Biggr) \\
			&\quad\times\lrk\Iq +z^{-1} \su{0}^\MP\lek F(z)-\su{1}\rek\rrk^{-1}\Biggl]\lrk\Iq +z^{-1} \su{0}^\MP\lek F(z)-\su{1}\rek\rrk\\
			&=z  \trauoa{F}{m-2}{s^\seins}{z}+\trauoa{F}{m-2}{s^\seins}{z}  \su{0}^\MP\lek F(z)-\su{1}\rek-z\trauoa{F}{m-2}{s^\seins}{z}\\
			&\quad-\sum_{j=0}^{m-1}z^{-j} \su{j+1} \su{0}^\MP  \trauoa{F}{m-2}{s^\seins}{z} 
			+\su{m} \su{0}^\MP \su{1} +\sum_{k=1}^{m-1} \su{m-k} \su{0}^\MP \su{k-1}^{\seins} \\
			&\quad-\sum_{j=1}^{m-1}\lrk z^{-j} \sum_{k=j}^{m-1}\su{m+j-k} \su{0}^\MP \su{k-1}^{\seins}\rrk 
			= \nablaua{m}{z}
		\end{split} 
	\]
	and, hence, \eqref{TM1} are fulfilled.
\eproof

Now we study the asymptotic behaviour of the function $\nablau{m}$, which was introduced in \eqref{nablaua}.

\blemml{L1651-2}
	Let $\theta\in[0,2\pi)$ and let $\mG$ be a of $\C\setminus\set{0}$ with $\setaa{\e^{\I  \theta}y}{y \in [1,+\infty)}\subseteq\mathcal{G}$. Let $F\colon \mG \to \Cpq$ be a matrix-valued function, let  $\kappa\in\N\cup\set{+\infty}$,  let $\seq{\su{j}}{j}{0}{\kappa}$ be a sequence of complex \tpqa{matrices}. In the case $\kappa\geq2$, let $\seq{\su{j}^\seins}{j}{0}{\kappa-2}$ be the \tfirstSta{\seq{\su{j}}{j}{0}{\kappa}}. Let $\su{-1}\defg\Opq$, let $\su{-1}^\seins \defg\Opq$, let $m\in\mn{1}{\kappa}$ and let $\nablau{m}\colon\mG\to\Cpq$ be defined by \eqref{nablaua}. Suppose  
	\begin{equation} \label{ML1}
		\lim_{r\to+\infty} \trauoa{F}{m-2}{s^\seins}{\e^{\I\theta}r}
		=\Opq.
	\end{equation}
	Then
	\begin{equation}  \label{ML1-1}
			\lim_{r\to+\infty}\nablaua{m}{\e^{\I  \theta}r}
			=
			\begin{cases}
				\su{1}\su{0}^\MP\su{1}\incase{ m=1}\\
				\su{m}\su{0}^\MP\su{1} +\sum_{k=1}^{m-1}\su{m-k}\su{0}^\MP\su{k-1}^\seins\incase{m\geq2}
			\end{cases}.
	\end{equation}
\elemm
\bproof
	Because of $\su{-1}^\seins \defg\Opq$, from \eqref{Fuoa} we have  $\trauoa{F}{-1}{s^\seins}{z}=F(z)$ for all $z \in \mG$. Thus,  in the case $m=1$, from \eqref{ML1} and \eqref{nablaua} we see that
	\[
		\begin{split}
			\su{1}\su{0}^\MP\su{1}
			&=\lrk\lek\lim_{r\to+\infty}\trauoa{F}{m-2}{s^\seins}{\e^{\I\theta}r}\rek-\su{1}\rrk\su{0}^\MP\lrk\lek\lim_{r\to+\infty} \trauoa{F}{m-2}{s^\seins}{\e^{\I\theta}r}\rek-\su{1}\rrk\\
			&=  \lim_{r\to+\infty}\lek F(\e^{\I \theta}r)-\su{1}\rek\su{0}^\MP\lek F(\e^{\I \theta}r)-\su{1}\rek
			=\lim_{r\to+\infty}\nablaua{m}{\e^{\I  \theta}r}
		\end{split}
	\]
	holds true. Consequently, if $m=1$, then \eqref{ML1-1} is checked. Now assume that $m\geq2$. Then we first observe that \rrema{R1013} shows that 
	\begin{equation} \label{ML2}
		F(z)
		= 	\trauoa{F}{-1}{s^\seins}{z}
		=z^{1-m}\trauoa{F}{m-2}{s^\seins}{z} - \sum_{j=0}^{m-2} z^{j-m+1}\su{m-2-j}	
	\end{equation}
	for all $z\in\mG$. Using \eqref{ML1} and \eqref{ML2}, we conclude
	\begin{equation} \label{ML3}
		\begin{split} 
			\Opq
			&=\lek\lim_{r\to+\infty}(\e^{\I \theta}r)^{1-m}\rek\lek\lim_{r\to+\infty} \trauoa{F}{m-2}{s^\seins}{\e^{\I\theta}r}\rek - \sum_{j=0}^{m-2}\lek\lim_{r\to+\infty}(\e^{\I \theta}r)^{j-m+1}\rek\su{m-2-j} \\
			&= \lim_{r\to+\infty}F(\e^{\I \theta}r).
		\end{split} 
	\end{equation}
	Thus, in view of \eqref{ML1}, \eqref{ML3}, and \eqref{nablaua}, we get then 
	\[ 
		\begin{split} 
			&\su{m} \su{0}^\MP\su{1} +\sum_{k=1}^{m-1}\su{m-k}\su{0}^\MP\su{k-1}^\seins\\
			&=\lek\lim_{r\to+\infty} \trauoa{F}{m-2}{s^\seins}{\e^{\I\theta}r}\rek\su{0}^\MP\lek \lim_{r\to+\infty}F{(\e^{\I \theta}r)}-\su{1}\rek- \su{1}\su{0}^\MP\lek\lim_{r\to+\infty} \trauoa{F}{m-2}{s^\seins}{\e^{\I\theta}r}\rek\\	
			& =\lek\lim_{r\to+\infty} \trauoa{F}{m-2}{s^\seins}{\e^{\I\theta}r}\rek\su{0}^\MP\lek\lim_{r\to+\infty}F{(\e^{\I \theta}r)}-\su{1}\rek- \su{1}\su{0}^\MP\lek\lim_{r\to+\infty} \trauoa{F}{m-2}{s^\seins}{\e^{\I\theta}r}\rek\\
			&\quad-	\sum_{j=1}^{m-1}\lek\lim_{r\to+\infty}{(\e^{\I \theta}r)}^{-j}\rek\su{j+1}\su{0}^\MP\lek\lim_{r\to+\infty} \trauoa{F}{m-2}{s^\seins}{\e^{\I\theta}r}\rek\\	
			&\quad+ \su{m}\su{0}^\MP\su{1} 
			+\sum_{k=1}^{m-1}\su{m-k}\su{0}^\MP\su{k-1}^\seins + 
			\sum_{j=1}^{m-1}\lek\lim_{r\to+\infty}{(\e^{\I \theta}r)}^{-j}\rek \sum_{k=j}^{m-1}\su{m+j-k}\su{0}^\MP\su{k-1}^\seins \\	
			& =\lim_{r\to+\infty}\Biggl( \trauoa{F}{m-2}{s^\seins}{\e^{\I\theta}r}\su{0}^\MP\lek F{(\e^{\I \theta}r)}-\su{1}\rek- \su{1}\su{0}^\MP \trauoa{F}{m-2}{s^\seins}{\e^{\I\theta}r}\\
			&\quad-	\sum_{j=1}^{m-1} {(\e^{\I \theta}r)}^{-j} \su{j+1}\su{0}^\MP   \trauoa{F}{m-2}{s^\seins}{\e^{\I\theta}r} 
			+ \su{m}\su{0}^\MP\su{1} \\
			&\quad+\sum_{k=1}^{m-1}\su{m-k}\su{0}^\MP\su{k-1}^\seins + 
			\sum_{j=1}^{m-1}\lek{(\e^{\I \theta}r)}^{-j} \sum_{k=j}^{m-1}\su{m+j-k}\su{0}^\MP\su{k-1}^\seins\rek\Biggl)\\
			& =\lim_{r\to+\infty}\Biggl( \trauoa{F}{m-2}{s^\seins}{\e^{\I\theta}r}\su{0}^\MP\lek F{(\e^{\I \theta}r)}-\su{1}\rek-\sum_{j=1}^{m-1} {(\e^{\I \theta}r)}^{-j} \su{j+1}\su{0}^\MP   \trauoa{F}{m-2}{s^\seins}{\e^{\I\theta}r}\\
			&\quad+ \su{m}\su{0}^\MP\su{1}+\sum_{k=1}^{m-1}\su{m-k}\su{0}^\MP\su{k-1}^\seins + 
			\sum_{j=1}^{m-1}\lek{(\e^{\I \theta}r)}^{-j} \sum_{k=j}^{m-1}\su{m+j-k}\su{0}^\MP\su{k-1}^\seins\rek\Biggl)\\		
			&= \lim_{r\to+\infty}\nablaua{m}{\e^{\I  \theta}r}.\qedhere
		\end{split} 
	\]
\eproof

In preparing the next result, we remember to \rpart{P1407.b} of \rprop{P1407}.

\blemml{L1522}
	Let $m\in\mn{2}{+\infty}$, let $\seq{\su{j}}{j}{0}{m}\in\Hggequ{m}$, let $\seq{\su{j}^\seins}{j}{0}{m-2}$ be the \tfirstSta{\seq{\su{j}}{j}{0}{m}} and let $F\in\RFuqa{m-2}{\seq{\su{j}^\seins}{j}{0}{m-2}}$. Then:
	\benui
		\il{L1522.a} $F\in\PFoddqa{\su{0}}$.
		\il{L1522.b} $\HTiaooo{F}{1}{\su{0}}{\su{1}}\in\RFuq{1}$.
		\il{L1522.c} For each $z\in\ohe$ the inclusion $\Bilda{F(z)}\subseteq\Bilda{\su{0}}$ holds.
	\eenui
\elemm
\bproof
	\eqref{L1522.a} \rlemm{L1553} yields
	\bgl{L1522.1}
		F
		\in\PFoddqa{\su{0}^\seins}.
	\eg
	From \eqref{seins} we get $\Kerna{\su{0}}\subseteq\Kerna{\su{0}^\seins}$. Thus, \rrema{R1543} yields $\PFoddqa{\su{0}^\seins}\subseteq\PFoddqa{\su{0}}$. Combining this with \eqref{L1522.1} we get~\eqref{L1522.a}.
	
	\eqref{L1522.b} From $\seq{\su{j}}{j}{0}{m}\in\Hggequ{m}$ and $m\geq2$ we get $\seq{\su{j}}{j}{0}{1}\in\Hggequ{1}$. Hence, \rprop{L1513} implies $\HTiaooo{F}{1}{\su{0}}{\su{1}}\in\RFuqa{1}{\seq{\su{j}}{j}{0}{1}}$. Thus, from \eqref{Ruua} and \rrema{R1433} we obtain~\eqref{L1522.b}.
	
	\eqref{L1522.c} In view of~\eqref{L1522.a} we infer from \rpart{L1407.b} of \rlemm{L1407} then
	\begin{align}\label{L1522.2}
		\Kerna{\su{0}}&\subseteq\KernA{F(z)},&
		z&\in\ohe.
	\end{align}
	From \rpart{L1037.b} of \rlemm{L1037} we obtain
	\bgl{L1522.3}
		\su{0}^\ad
		=\su{0}.
	\eg
	Because of \eqref{Ruua} and \rrema{R1433} we have $F\in\NFq$. Combining this with \eqref{L1522.2} and \eqref{L1522.3} the application of \rlemm{L1030} yields~\eqref{L1522.c}.
\eproof

\btheol{L1241}
	Let $n\in\N$, let $\seq{\su{j}}{j}{0}{2n}\in\Hggeuu{q}{2n}$, let $\seq{\su{j}^\seins}{j}{0}{2n-2}$ be the \tfirstSta{\seq{\su{j}}{j}{0}{2n}}, and let $F \in \RkqSEjj{2n-2}{0}{2n-2}$. Then $\HTiaooo{F}{1}{\su{0}}{\su{1}}$ belongs to $\Rkqjj{2n}{0}{2n}$.
\etheo
\bproof
	The strategy of our proof is based on an application of \rtheo{P0611-1}. Let $m\defg 2n$ and let
	\bgl{35A}
		\su{-1}
		\defg\Oqq.
	\eg
	From \eqref{35A}, $\seq{\su{j}}{j}{0}{m}\in\Hggequ{m}$ and \rpart{L1037.b} of \rlemm{L1037} we get
	\begin{align}\label{L1241.2}
		\su{j}^\ad&=\su{j}&\text{for each }j&\in\mn{-1}{m}.
	\end{align}
	\rPart{L1522.b} of \rlemm{L1522} yields
	\bgl{L1241.3}
		\HTiaooo{F}{1}{\su{0}}{\su{1}}
		\in\NFq.
	\eg
	Let $\su{-1}^\seins\defg\Oqq$. Then, in view of $F\in\RFuqa{m-2}{\seq{\su{j}^\seins}{j}{0}{m-2}}$, we infer from \rrema{R1404} that
	\bgl{10-42-1}
		F
		\in\NFuqA{m-2}{\seq{\su{j}^\seins}{j}{-1}{m-2}}.
	\eg
	Thus, \rcoro{C1014*} shows that 
	\begin{equation} \label{2157.3}
		\lim_{y\to+\infty} \trauoa{F}{m-2}{s^\seins}{\I y}
		=\Oqq.
	\end{equation} 
	Setting
	\bgl{10-42-2}
		G
		\defg \su{m}\su{0}^\MP \su{1}+ \sum_{k=1}^{m-1} \su{m-k} \su{0}^\MP \su{k-1}^{\seins},
	\eg
	we get from \eqref{2157.3} and \rlemm{L1651-2} then
	\begin{equation} \label{2157.4}
		G
		=\lim_{r \to +\infty} \nablaua{m}{\e^{\I \frac{\pi}{2}}r}
		= \lim_{y \to +\infty} \nablaua{m}{\I y},
	\end{equation} 
	where $\nablau{m}$ is defined in \eqref{nablaua}. From the choice of $F$ and \eqref{Ruua} we conclude that $F\in\RFuq{m-2}$. Thus, \rrema{R1433} yields \eqref{R1437.B1}. Consequently,
	\bgl{85-1}
		\Iq
		=\Iq+0 \cdot\su{0}^\MP (\Oqq-\su{1})
		= \lim_{y\to+\infty}\lrk\Iq+(\I y)^{-1} \su{0}^\MP\lek F(\I y)-\su{1}\rek\rrk,
	\eg
	which implies 
	\begin{equation}	\label{2157.6}
		\lim_{y\to+\infty}\lek\lrk\Iq+(\I y)^{-1} \su{0}^\MP\lek F(\I y)-\su{1}\rek\rrk^{-1}\rek
		=\Iq.
	\end{equation}
	Now we are going to apply \rlemm{L1651-1}. In view of $\seq{\su{j}}{j}{0}{m}\in\Hggequ{m}$ and \rlemm{L1048}, we get $\seq{\su{j}}{j}{0}{m}\in\Dqqu{m}$. \rPart{L1522.c} of \rlemm{L1522} yields
	\bal{L1241.5}
		\BildA{F(z)}&\subseteq\Bilda{\su{0}}&\text{for each }z&\in\ohe.
	\ea
	In view of $\seq{\su{j}}{j}{0}{m}\in\Hggequ{m}$ we infer from \rpartsss{L1037.a}{L1037.b}{L1037.c} of \rlemm{L1037} then $\su{0}\in\Cggq$, $\su{1}^\ad=\su{1}$ and $\Kerna{\su{0}}\subseteq\Kerna{\su{1}}$. Thus, \rpart{L0800.b} of \rprop{L0800} yields
	\bgl{L1241.6}
		\det\lrk\Iq+z^{-1} \su{0}^\MP\lek F(\I y)-\su{1}\rek\rrk
		=z^{-q}\det\lrk z\Iq+ \su{0}^\MP\lek F(\I y)-\su{1}\rek\rrk
		\neq 0
	\eg
	for all $z \in \ohe$. Since $\seq{\su{j}}{j}{0}{m}$ belongs to $\Dqqu{m}$, from \eqref{L1241.5}, \eqref{L1241.6} and \rlemm{L1651-1} we then conclude that
	\begin{equation}	\label{2157.8}
		\trauoa{(\HTiaooo{F}{1}{\su{0}}{\su{1}})}{m}{s}{z}
		=\trauoa{F}{m-2}{s^\seins}{z}-\frac{1}{z} \nablaua{m}{z}\lrk\Iq+z^{-1} \su{0}^\MP\lek F(z)-\su{1}\rek\rrk^{-1}
	\end{equation}
	for all $z \in \ohe$. By virtue of \eqref{2157.3}, \eqref{2157.4}, and \eqref{2157.6}, from \eqref{2157.8} we see that 
	\begin{equation}	\label{2157.9}
		\begin{split}
			\Oqq
			&=\Oqq -0\cdot G \cdot \Iq \\
			&=\lek\lim_{y \to +\infty} \trauoa{F}{m-2}{s^\seins}{\I y}\rek-\lrk\lim_{y \to +\infty} \frac{1}{\I y}\rrk\lek\lim_{y \to +\infty}\nablaua{m}{\I y}\rek\\
			&\quad\times\lrk\lim_{y \to +\infty}\lek\lrk\Iq+(\I y)^{-1} \su{0}^\MP\lek F(\I y)-\su{1}\rek\rrk^{-1}\rek\rrk\\
			&= \lim_{y \to +\infty} \trauoa{(\HTiaooo{F}{1}{\su{0}}{\su{1}})}{m}{s}{\I y}
		\end{split}
	\end{equation}
	holds true. Taking \eqref{L1241.2}, \eqref{L1241.3}, \eqref{2157.9}, and $m=2n$ into account, \rtheo{P0611-1} then yields that  $\HTiaooo{F}{1}{\su{0}}{\su{1}}$  belongs to  $\RqKJ{2n}{-1}$. In view of \eqref{35A}, then \rrema{R1404} shows that the function $\HTiaooo{F}{1}{\su{0}}{\su{1}}$ also  belongs to $\Rkqjj{2n}{0}{2n}$.
\eproof

\bcorol{C1028}
	Let $\seq{\su{j}}{j}{0}{\infty}\in\Hggequ{\infty}$, let $\seq{\su{j}^\seins}{j}{0}{\infty}$ be the \tfirstSta{\seq{\su{j}}{j}{0}{\infty}} and let $F\in\RFuqa{\infty}{\seq{\su{j}^\seins}{j}{0}{\infty}}$. Then $\HTiaooo{F}{1}{\su{0}}{\su{1}}$ belongs to $\RFuqa{\infty}{\seq{\su{j}}{j}{0}{\infty}}$.
\ecoro
\bproof
	Combine \rremass{R0938*}{R1038} and \rtheo{L1241}.
\eproof

\btheol{L0925}
	Let $n\in\N$, let $\seq{\su{j}}{j}{0}{2n+1} \in\Hggeuu{q}{2n+1}$, let $\seq{\su{j}^\seins}{j}{0}{2n-1}$ be the \tfirstSta{\seq{\su{j}}{j}{0}{2n+1}} and let $F \in \RkqSEjj{2n-1}{0}{2n-1}$. Then $\HTiaooo{F}{1}{\su{0}}{\su{1}}$ belongs to $\Rkqjj{2n+1}{0}{2n+1}$.
\etheo
\bproof
	The strategy of our proof is based on an application of \rtheo{P0611-2}. Let $m\defg2n+1$ and $\su{-1}\defg\Oqq$. Then, from  $\seq{\su{j}}{j}{0}{m}\in\Hggequ{m}$ and \rpart{L1037.b} of \rlemm{L1037} we get \eqref{L1241.2}.
	\rPart{L1522.b} of \rlemm{L1522} yields \eqref{L1241.3}.
	Let $\su{-1}^\seins\defg\Oqq$. In view of $F\in\RFuqa{m-2}{\seq{\su{j}^\seins}{j}{0}{m-2}}$, we infer then from \rrema{R1404} that \eqref{10-42-1} holds.
	In view of \eqref{10-42-1}, the applictaion of
	\rcoro{C1014*} shows that \eqref{2157.3} is true.
	Let $G$ be defined by \eqref{10-42-2}. Then
	we get from \eqref{2157.3} and \rlemm{L1651-2} then \eqref{2157.4},
	where $\nablau{m}$ is defined in \eqref{nablaua}. From the choice of $F$ and \eqref{Ruua} we know that $F\in\RFuq{m-2}$. Thus, \rrema{R1437} implies that \eqref{R1437.B1}
	holds true. Consequently, \eqref{85-1} is true,
	which implies \eqref{2157.6}. 
	Now we are going to apply \rlemm{L1651-1}. In view of  $\seq{\su{j}}{j}{0}{m}\in\Hggequ{m}$ and \rlemm{L1048} we get $\seq{\su{j}}{j}{0}{m}\in\Dqqu{m}$. \rPart{L1522.c} of \rlemm{L1522} yields \eqref{L1241.5}.
	In view of $\seq{\su{j}}{j}{0}{m}\in\Hggequ{m}$, we infer from \rlemm{L1037} then $\su{0}\in\Cggq$, $\su{1}^\ad=\su{1}$ and $\Kerna{\su{0}}\subseteq\Kerna{\su{1}}$. Thus, \rpart{L0800.b} of \rprop{L0800} yields \eqref{L1241.6}
	for all $z \in \ohe$. Because of $\seq{\su{j}}{j}{0}{m}\in\Dqqu{m}$, \eqref{L1241.5} and \eqref{L1241.6}, we then conclude from \rlemm{L1651-1} that
	\begin{equation}	\label{IW}
		\trauoa{(\HTiaooo{F}{1}{\su{0}}{\su{1}})}{2n+1}{s}{z}
		=\trauoa{F}{2n-1}{s^\seins}{z}-\frac{1}{z} \nablaua{2n+1}{z}\lrk\Iq+z^{-1} \su{0}^\MP \lek F(z)-\su{1}\rek\rrk^{-1}
	\end{equation}
	for all $z \in \ohe$. By virtue of \eqref{2157.3}, \eqref{2157.4} and \eqref{2157.6}, we see that \eqref{2157.9}
	holds true. Taking \eqref{L1241.2}, \eqref{L1241.3}, \eqref{2157.9}, and $m=2n+1$ into account, we get from \rtheo{P0611-2} then $\trauo{(\HTiaooo{F}{1}{\su{0}}{\su{1}})}{2n+1}{s} \in \RqP$ follows. In particular, $\trauo{(\HTiaooo{F}{1}{\su{0}}{\su{1}})}{2n+1}{s}$ is holomorphic in $\ohe$. This shows us the function $\Omega\colon[1,+\infty) \to \R $ defined by 
	\[ 
		\Omega(y)
		\defg \frac{1}{y}\normA{\im  \trauoa{(\HTiaooo{F}{1}{\su{0}}{\su{1}})}{2n+1}{s}{\I y}}
	\]
	is continuous. Because of the continuity of the holomorphic function $F$, \eqref{Fuoa}, \eqref{nablaua}, \eqref{L1241.6} and $m=2n+1$, the function $\Theta\colon[1,+\infty) \to \R$ given by 
	\begin{equation} 	\label{6-1}
		\Theta(y)
		\defg\normA{\re\lek\nablaua{2n+1}{\I y}\lrk\Iq+(\I y)^{-1} \su{0}^\MP\lek F(\I y)-\su{1}\rek\rrk^{-1}\rek}
	\end{equation}
	is continuous. In view of \eqref{2157.4}, \eqref{2157.6}, $m=2n+1$, and \eqref{6-1}, we conclude 
	\bsp
		&\norma{\re G}
		=\normA{\re ( G\cdot\Iq)}\\
		&=\normA{\re\lek\lek\lim_{y \to +\infty}\nablaua{2n+1}{\I y}\rek\lrk\lim_{y \to +\infty}\lek\lrk\Iq+(\I y)^{-1} \su{0}^\MP\lek F(\I y)-\su{1}\rek\rrk^{-1}\rek\rrk\rek}\\
		& = \lim_{y\to+\infty}\Theta(y).
	\esp
	Consequently, there is a real number $c$ such that $\abs{\Theta(y)}\leq c$ for all $y \in [1,+\infty)$. In view of \eqref{10-42-1} and $m=2n+1$,  we see from \rpart{R1640.a} of \rprop{R1640} that $ \trauo{F}{2n-1}{s^\seins}$ belongs to $\Rkq{-1}$ and (see \eqref{Ruu}), in particular, to $\RqK{-1}$. This means that $\trauo{F}{2n-1}{s^\seins}\in \RqP$ and that $\Phi\colon[1,+\infty) \to \R$ defined by $\Phi(y) \defg \frac{1}{y}\norma{\im \trauoa{F}{2n-1}{s^\seins}{\I y}}$  belongs to $ \Loaaaa{1}{[1,+\infty)}{\BA{[1,+\infty)}}{ \tilde{\lambda}}{\R}$, where $\tilde\Leb$\index{lambda^~@$\tilde\Leb$} is the Lebesgue measure defined on $\BA{[1,+\infty)}$. Hence, $\Psi\colon[1,+\infty) \to \R$ given by   $\Psi(y) \defg \Phi(y) + \frac{c}{y^2}$ belongs to $\Loaaaa{1}{[1,+\infty)}{\BA{[1,+\infty)}}{ \tilde{\lambda}}{\R}$. Since $\trauo{(\HTiaooo{F}{1}{\su{0}}{\su{1}})}{2n+1}{s}$ belongs to the class  $\RqP$, the function $\Lambda\colon[1,+\infty) \to \R$ given by
	\bgl{L0925.7}
		\Lambda(y)
		\defg \frac{1}{y}\normA{\im \trauoa{(\HTiaooo{F}{1}{\su{0}}{\su{1}})}{2n+1}{s}{\I y}}
	\eg
	is continuous and, because of \eqref{IW}, for all $y \in [1,+\infty)$ it satisfies
	\[
		\begin{split}
			\Abs{\Lambda(y)}
			&\leq \frac{1}{y}\lrk\normA{\im\trauoa{F}{2n-1}{s^\seins}{\I y}}+\normA{\im\lek\frac{1}{\I y} \nablaua{2n+1}{\I y}\lrk\Iq+\frac{1}{\I y} \su{0}^\MP\lek F(\I y)-\su{1}\rek\rrk^{-1}\rek}\rrk\\
			&= \frac{1}{y}\normA{\im\trauoa{F}{2n-1}{s^\seins}{\I y}}+ \frac{1}{y^2}\normA{\re\lek\nablaua{2n+1}{\I y}\lrk\Iq+\frac{1}{\I y} \su{0}^\MP\lek F(\I y)-\su{1}\rek\rrk^{-1}\rek}\\
			&= \Phi(y)+\frac{1}{y^2} \Theta(y)
			\leq\Psi(y).
		\end{split}
	\]
	Thus, the function $\Lambda$ also belongs to   $\Loaaaa{1}{[1,+\infty)}{\BA{[1,+\infty)}}{ \tilde{\lambda}}{\R}$. Since  $\trauo{(\HTiaooo{F}{1}{\su{0}}{\su{1}})}{2n+1}{s}$ is a member of the class   $\RqP$, we then see from \eqref{L0925.7} that
	\[
		\trauo{(\HTiaooo{F}{1}{\su{0}}{\su{1}})}{2n+1}{s}
		\in\NFuq{-1}.
	\]
	Taking \eqref{L1241.2}, \eqref{L1241.3}, \eqref{2157.9}, and $m=2n+1$ into account, \rtheo{P0611-2} implies
	\bgl{L0925.9}
		\HTiaooo{F}{1}{\su{0}}{\su{1}}
		\in\NFuqa{2n+1}{\seq{\su{j}}{j}{-1}{2n+1}}.
	\eg
	In view of $\su{-1}\defg\Oqq$ and \eqref{L0925.9}, \rrema{R1404} yields finally
	\[
		\HTiaooo{F}{1}{\su{0}}{\su{1}}
		\in\RFuqa{2n+1}{\seq{\su{j}}{j}{0}{2n+1}}.\qedhere
	\]
\eproof

\bcorol{C1332}
	Let $m\in\mn{2}{+\infty}\cup\set{+\infty}$, let $\seq{\su{j}}{j}{0}{m}\in\Hggequ{m}$ and let $\seq{\su{j}^\seins}{j}{0}{m-2}$ be the \tfirstSta{\seq{\su{j}}{j}{0}{m}}. Denote $\HTuu{\su{0}}{\su{1}}$ the bijective mapping defined in \rcoro{C1514}. Then $\HTuu{\su{0}}{\su{1}}$ generates a bijective correspondence between the sets $\RFuqa{m}{\seq{\su{j}}{j}{0}{m}}$ and $\RFuqa{m-2}{\seq{\su{j}^\seins}{j}{0}{m-2}}$.  The inverse mapping $(\HTuu{\su{0}}{\su{1}})^\inv$ is given for $G\in\RFuqa{m-2}{\seq{\su{j}^\seins}{j}{0}{m-2}}$ by $(\HTuu{\su{0}}{\su{1}})^\inv(G)=\HTiaooo{G}{1}{\su{0}}{\su{1}}$.
\ecoro
\bproof
	We will consider the cases of even $m$, odd $m$ and $m=\infty$. In any of these cases, we can apply \rcoro{L0845} and \rpart{C1519.b} of \rcoro{C1519} to verify the shape of the inverse mapping. If $m$ is even, then the assertion follows from \rtheoss{L1135}{L1241}. If $m$ is odd, then the application of \rtheoss{L1314}{L0925} yields the desired result. Finally, if $m=\infty$, then one has to apply \rcoross{C0910}{C1028}.
\eproof

We mention that the investigations in \rsectss{S1326}{S1021} were influenced to some extent by considerations in Chen/Hu~\cite{MR1624548}. In particular,~\cite{MR1624548}*{\clemm{2.6}} played an essential role for the choice of our strategy. This concerns first the development of a Schur type algorithm for sequences of complex matrices and then the construction of an interrelated Schur type algorithm for functions belonging to special subclasses of $\NFq$. The contents of~\cite{MR1624548}*{\clemm{2.6}} are covered by \rtheoss{L1135}{L1241}. The method of Chen/Hu to prove~\cite{MR1624548}*{\clemm{2.6}} strongly differs from our approach. It uses results on generalized Bezoutians to Anderson/Jury~\cite{MR0444175} and Gekhtman/Shmoish~\cite{MR1340693}.

\section{A Schur-Nevanlinna Type Algorithm for the Class $\RFuq{\kappa}$}\label{S1610}

The results of \rsect{S1326} suggest the construction of a Schur-Nevanlinna type algorithm for the class $\RFuq{\kappa}$ with $\kappa\in\NO\cup\set{+\infty}$. The main theme of this section is to work out the details of this algorithm. In \rsect{S1326}, we fixed an $m\in\NO$ and a sequence $\seq{\su{j}}{j}{0}{m}\in\Hggequ{m}$. \rprop{P1407} tells us then that $\RFuqa{m}{\seq{\su{j}}{j}{0}{m}}\neq\emptyset$. Let $F\in\RFuqa{m}{\seq{\su{j}}{j}{0}{m}}$. Then the central theme of \rsect{S1326} was to study the \taaSt{\su{0}}{\Oqq} $\HTaoo{F}{1}{\su{0}}$ of $F$ and furthermore, in case $m\geq1$, the \taaSt{\su{0}}{\su{1}} $\HTaooo{F}{1}{\su{0}}{\su{1}}$ of $F$. The following observation shows that in the case $\kappa\in\NO\cup\set{+\infty}$ the results of \rsect{S1326} can be applied to arbitrary functions belonging to $\RFuq{\kappa}$.

\blemml{L1355}
	Let $\kappa\in\NO\cup\set{+\infty}$ and let  $F\in\RFuq{\kappa}$. Then $\sigmaF\in\MgguqR{\kappa}$, $\seq{\suo{j}{\sigmaF}}{j}{0}{\kappa}\in\Hggequ{\kappa}$ and
	\[
		F
		\in\RFuqA{\kappa}{\seq{\suo{j}{\sigmaF}}{j}{0}{\kappa}}.
	\]
\elemm
\bproof
	From \eqref{Ruu} and \eqref{Ruo} we see that $\sigmaF\in\MgguqR{\kappa}$. Thus, from \eqref{Ruua} we infer now $F\in\RFuqa{\kappa}{\seq{\suo{j}{\sigmaF}}{j}{0}{\kappa}}$. Hence, \rpart{P1407.b} of \rprop{P1407} yields $\seq{\suo{j}{\sigmaF}}{j}{0}{\kappa}\in\Hggequ{\kappa}$.
\eproof

In view of \rlemm{L1355}, we introduce the following construction of \tSNt{s} (shortly \tsSNt{s}) for functions belonging to the class $\RFuq{\kappa}$ with $\kappa\in\NO\cup\set{+\infty}$.

\bdefil{D1416}
	Let $\kappa\in\NO\cup\set{+\infty}$ and let  $F\in\RFuq{\kappa}$.
	\benui
		\item The function $F^\Fa{0}\defg F$\index{^(0)@$F^\Fa{0}$} is called the \emph{\tFaa{0}{F}}.
		\il{D1416.b} The function $F^\Fa{1}\defg \HTaoo{F}{1}{\suo{0}{\sigmaF}}$\index{^(1)@$F^\Fa{1}$} is called the \emph{\tFeinsa{F}}.
		\il{D1416.c} In the case $\kappa\in\N\cup\set{+\infty}$, the function $F^\Fa{2}\defg \HTaooo{F}{1}{\suo{0}{\sigmaF}}{\suo{1}{\sigmaF}}$\index{^(2)@$F^\Fa{2}$} is called the \emph{\tFzweia{F}}.
	\eenui
\edefi

We relate now the just introduced notions to the \tiSt{s} studied in \rsect{S1515}.

\brema
	Let $\kappa\in\NO\cup\set{+\infty}$ and let  $F\in\RFuq{\kappa}$. From \rlemm{L1355} and \rpartss{L0845.b}{L0845.c} of \rcoro{L0845}, one can see then that $\HTiaoo{(F^\Fa{1})}{1}{\suo{0}{\sigmaF}}=F$ and, in the case $\kappa\geq1$, furthermore
	\[
		\HTiaooo{(F^\Fa{2})}{1}{\suo{0}{\sigmaF}}{\suo{1}{\sigmaF}}
		=F.
	\]
\erema

\bpropl{P1433}
	Let $\kappa\in\NO\cup\set{+\infty}$ and let  $F\in\RFuq{\kappa}$. Then $F^\Fa{1}\in\PFevenqa{\suo{0}{\sigmaF}}$ and, in the case $\kappa\in\N\cup\set{+\infty}$, $F^\Fa{2}\in\PFoddqa{\suo{0}{\sigmaF}}$.
\eprop
\bproof
	In view of \rlemm{L1355}, the application of \rcoross{C0835}{C0856} yields the first and second assertion, respectively.
\eproof

\bpropl{P1438}
	Let $\kappa\in\mn{2}{+\infty}\cup\set{+\infty}$ and let  $F\in\RFuq{\kappa}$. Let $\seq{\su{j}^\seins}{j}{0}{\kappa-2}$ be the \tfirstSta{\seq{\suo{j}{\sigmaF}}{j}{0}{\kappa}}. Then $F^\Fa{2}\in\RFuqa{\kappa-2}{\seq{\su{j}^\seins}{j}{0}{\kappa-2}}$, $\sigmau{F^\Fa{2}}\in\MggqRag{\seq{\su{j}^\seins}{j}{0}{\kappa-2}}$, and $\seq{\su{j}^\seins}{j}{0}{\kappa-2}=\seq{\suo{j}{\sigmau{F^\Fa{2}}}}{j}{0}{\kappa-2}$.
\eprop
\bproof
	In view of \rlemm{L1355}, we have $\sigmaF\in\MgguqR{\kappa}$, $\seq{\suo{j}{\sigmaF}}{j}{0}{\kappa}\in\Hggequ{\kappa}$, and $F\in\RFuqa{\kappa}{\seq{\suo{j}{\sigmaF}}{j}{0}{\kappa}}$. Thus, we get from \rrema{P1323*} in the case $\kappa\in\mn{2}{+\infty}$ and \rcoro{C0910} in the case $\kappa=+\infty$ that $F^\Fa{2}\in\RFuqa{\kappa-2}{\seq{\su{j}^\seins}{j}{0}{\kappa-2}}$. Hence, \eqref{Ruua} implies $F^\Fa{2}\in\RFuq{\kappa-2}$ and $\sigmau{F^\Fa{2}}\in\MggqRag{\seq{\su{j}^\seins}{j}{0}{\kappa-2}}$. Consequently, $\seq{\suo{j}{\sigmau{F^\Fa{2}}}}{j}{0}{\kappa-2}=\seq{\su{j}^\seins}{j}{0}{\kappa-2}$.
\eproof

\blemml{R1300}
	Let $\kappa\in\N\cup\set{+\infty}$ and let $F\in\RFuq{\kappa}$. Then $F^\Fa{2}\in\RFuq{\kappa-2}$.
\elemm
\bproof
	In the case $\kappa=1$, the combination of \rprop{P1433} and \eqref{Poddua} yields
	\[
		F^\Fa{2}
		\in\PFoddqa{\suo{0}{\sigmaF}}
		\subseteq\RFuq{-1}
		=\RFuq{1-2}.
	\]
	In the case $\kappa\in\mn{2}{+\infty}\cup\set{+\infty}$, the application of \rprop{P1438} and \eqref{Ruua} brings $F^\Fa{2}\in\RFuq{\kappa-2}$.
\eproof

In view of \rlemm{R1300}, we introduce in recursive way the following notions.

\bdefil{D1609}
	Let $\kappa\in\N\cup\set{+\infty}$ and let $F\in\RFuq{\kappa}$. For all $k\in\NO$ with $2k+1\leq\kappa$, we will call the \tFzwei{} $F^\Fa{2(k+1)}$\index{^(2k)@$F^\Fa{2k}$} of the \tFa{2k} $F^\Fa{2k}$ the \emph{\tFaa{2(k+1)}{F}}.
\edefi

\bdefil{D1028}
	Let $\kappa\in\NO\cup\set{+\infty}$ and let $F\in\RFuq{\kappa}$. For all $k\in\NO$ with $2k\leq\kappa$, we will call the \tFeins{} $F^\Fa{2k+1}$\index{^(2k+1)@$F^\Fa{2k+1}$} of the \tFa{2k} $F^\Fa{2k}$ the \emph{\tFaa{(2k+1)}{F}}.
\edefi

\brema
	Let $\kappa\in\mn{-1}{+\infty}\cup\set{+\infty}$, let $F\in\RFuq{\kappa}$, and let $k\in\NO$ with $2k-1\leq\kappa$. From \rdefisss{D1416}{D1609}{D1028} and \rlemm{R1300}, we see that:
	\benui
		\item $F^\Fa{2k}\in\RFuq{\kappa-2k}$.
		\item $(F^\Fa{2k})^\Fa{l}=F^\Fa{2k+l}$ for each $l\in\NO$ with $l-1\leq\kappa-2k$.
	\eenui
\erema
%

The content of our next considerations can be described as follows. Let $\kappa\in\mn{2}{+\infty}\cup\set{+\infty}$ and let $\seq{\su{j}}{j}{0}{\kappa}\in\Hggequ{\kappa}$. Then we are going to study the above introduced Schur-Nevanlinna type algorithm particularly for functions which belong to the class $\RFuqa{\kappa}{\seq{\su{j}}{j}{0}{\kappa}}$. The next results contain essential information about the algorithm applied to functions belonging to the class $\RFuqa{\kappa}{\seq{\su{j}}{j}{0}{\kappa}}$. For every choice of $F\in\NFq$ and $\su{0},\su{1}\in\Cqq$, we set
\ba
	\HTuua{\su{0}}{\su{1}}{F}&\defg\HTaooo{F}{1}{\su{0}}{\su{1}},&
	\HTiuua{\su{0}}{\su{1}}{F}&\defg\HTiaooo{F}{1}{\su{0}}{\su{1}}\\
	\intertext{and}
	\HTua{\su{0}}{F}&\defg\HTaoo{F}{1}{\su{0}},&
	\HTiua{\su{0}}{F}&\defg\HTiaoo{F}{1}{\su{0}}.
\ea\index{F_(+;,)()@$\HTuua{\su{0}}{\su{1}}{F}$}\index{F_(-;,)()@$\HTiuua{\su{0}}{\su{1}}{F}$}\index{F_(+;)()@$\HTua{\su{0}}{F}$}\index{F_(-;)()@$\HTiua{\su{0}}{F}$}

\bpropl{P0852}
	Let $\kappa\in\N\cup\set{+\infty}$, let $\seq{\su{j}}{j}{0}{\kappa}\in\Hggequ{\kappa}$, and let $F$ belong to $\RFuqa{\kappa}{\seq{\su{j}}{j}{0}{\kappa}}$. Further, let $n\in\NO$ be such that $2n+1\leq\kappa$ and let $F^\Fa{2(n+1)}$ be the \tFaa{2(n+1)}{F}. For each integer $m$ with $0\leq m\leq\min\set{n+1,\frac{\kappa}{2}}$, let $\seq{\su{j}^\sta{m}}{j}{0}{\kappa-2m}$ be the \taSta{m}{\seq{\su{j}}{j}{0}{\kappa}}. Then:
	\benui
		\il{P0852.a}
		\bgl{P0852.1}
			F^\Fa{2(n+1)}
			=[\HTuu{\su{0}^\sta{n}}{\su{1}^\sta{n}}\circ\HTuu{\su{0}^\sta{n-1}}{\su{1}^\sta{n-1}}\circ\dotsb\circ\HTuu{\su{0}^\seins}{\su{1}^\seins}\circ\HTuu{\su{0}^\sta{0}}{\su{1}^\sta{0}}](F)
		\eg
		and
		\bgl{P0852.2}
			F^\Fa{2(n+1)}
			\in
			\begin{cases}
				\Ruua{\kappa-2(n+1)}{q}{\seq{\su{j}^\sta{n+1}}{j}{0}{\kappa-2(n+1)}}\incase{2n+1<\kappa}\\
				\PFoddqa{\su{0}^\sta{n}}\incase{2n+1=\kappa}
			\end{cases}.
		\eg
		\il{P0852.b} $F=[\HTiuu{\su{0}^\sta{0}}{\su{1}^\sta{0}}\circ\HTiuu{\su{0}^\sta{1}}{\su{1}^\sta{1}}\circ\dotsb\circ\HTiuu{\su{0}^\sta{n}}{\su{1}^\sta{n}}](F^\Fa{2(n+1)})$.
	\eenui
\eprop
\bproof
	\eqref{P0852.a} Formula~\eqref{P0852.1} follows immediately from the definition of $F^\Fa{2(n+1)}$ (see \rdefi{D1609} and \rpart{D1416.c} of \rdefi{D1416}). Using \rcoro{C0856}, \rtheo{L1135}, \rcoro{C0910}, and \rtheo{L1314} formula~\eqref{P0852.2} follows by induction.
	
	\eqref{P0852.b} In view of \rcoross{C1514}{C1332}, \rpart{P0852.b} is an immediate consequence of~\eqref{P0852.a}.
\eproof

\bpropl{P0921}
	Let $\kappa\in\NO\cup\set{+\infty}$, let $\seq{\su{j}}{j}{0}{\kappa}\in\Hggequ{\kappa}$ and let $F$ belong to $\RFuqa{\kappa}{\seq{\su{j}}{j}{0}{\kappa}}$. Further, let $n\in\NO$ be such that $2n\leq\kappa$ and denote by $F^\Fa{2n+1}$ the \tFaa{(2n+1)}{F}. For each $m\in\mn{0}{n}$ let $\seq{\su{j}^\sta{m}}{j}{0}{\kappa-2m}$ be the \taSta{m}{\seq{\su{j}}{j}{0}{\kappa}}. Then the following statements hold true:
	\benui
		\il{P0921.a}
		\bgl{P0921.1}
			F^\Fa{2n+1}
			=
			\begin{cases}
				[\HTu{\su{0}^\sta{n}}\circ\HTuu{\su{0}^\sta{n-1}}{\su{1}^\sta{n-1}}\circ\dotsb\circ\HTuu{\su{0}^\sta{0}}{\su{1}^\sta{0}}](F)\incase{n\in\N}\\
				\HTua{\su{0}^\sta{n}}{F}\incase{n=0}
			\end{cases}
		\eg
		and
		\bgl{P0921.2}
			F^\Fa{2n+1}
			\in\PFevenqa{\su{0}^\sta{n}}.
		\eg
		\il{P0921.b}
		\[
			F=
			\begin{cases}
				\HTiua{\su{0}^\sta{n}}{F^\Fa{2n+1}}\incase{n=0}\\
				[\HTiuu{\su{0}^\sta{0}}{\su{1}^\sta{0}}\circ\dotsb\circ\HTiuu{\su{0}^\sta{n-1}}{\su{1}^\sta{n-1}}\circ\HTiu{\su{0}^\sta{n}}](F^\Fa{2n+1})\incase{n\in\N}
			\end{cases}.
		\]
	\eenui
\eprop
\bproof
	\eqref{P0921.a} Formula~\eqref{P0921.1} follows immediately from the definition of $F^\Fa{2n+1}$ (see \rdefi{D1028} and \rpartss{D1416.b}{D1416.c} of \rdefi{D1416}). Using \rcoro{C0835}, \rtheo{L1135}, \rcoro{C0910}, and \rtheo{L1314}, formula~\eqref{P0921.2} follows by induction.
	
	\eqref{P0921.b} In view of \rcoross{C1103}{C1332}, \rpart{P0921.b} is an immediate consequence of~\eqref{P0921.a}.
\eproof

Our next considerations are aimed to study the inversion of the Schur-Nevanlinna algorithm. The following result can be considered as an inverse statement with respect to \rprop{P0852}.

\bpropl{P1035}
	Let $n\in\NO$ and $\seq{\su{j}}{j}{0}{2n+1}\in\Hggequ{2n+1}$. For each $m\in\mn{0}{n}$, let $\seq{\su{j}^\sta{m}}{j}{0}{2(n-m)+1}$ be the \taSta{m}{\seq{\su{j}}{j}{0}{2n+1}}. Further, let $G\in\PFoddqa{\su{0}^\sta{n}}$ and let
	\[
		F
		\defg[\HTiuu{\su{0}^\sta{0}}{\su{1}^\sta{0}}\circ\HTiuu{\su{0}^\seins}{\su{1}^\seins}\circ\dotsb\circ\HTiuu{\su{0}^\sta{n}}{\su{1}^\sta{n}}](G).
	\]
	Then the following statements hold true:
	\benui
		\il{P1035.a} $F\in\RFuqa{2n+1}{\seq{\su{j}}{j}{0}{2n+1}}$.
		\il{P1035.b} The \tFa{2(n+1)} $F^\Fa{2(n+1)}$ of $F$ satisfies $F^\Fa{2(n+1)}=G$.
	\eenui
\eprop
\bproof
	\eqref{P1035.a} This follows by combining \rprop{L1513} and \rtheo{L0925}.
	
	\eqref{P1035.b} Taking the definition of $F$ and~\eqref{P1035.a} into account, the application of \rcoross{C1514}{C1332} yields
	\[
		G=[\HTuu{\su{0}^\sta{n}}{\su{1}^\sta{n}}\circ\HTuu{\su{0}^\sta{n-1}}{\su{1}^\sta{n-1}}\circ\dotsb\circ\HTuu{\su{0}^\seins}{\su{1}^\seins}\circ\HTuu{\su{0}^\sta{0}}{\su{1}^\sta{0}}](F).
	\]
	Combining this with~\eqref{P1035.a}, we infer from \rpart{P0852.a} of \rprop{P0852} that
	\[
		F^\Fa{2(n+1)}
		=G.\qedhere
	\]
\eproof

The combination of \rpropss{P0852}{P1035} gives us now a complete description of the SN~algorithm in the class $\RFuqa{2n+1}{\seq{\su{j}}{j}{0}{2n+1}}$.

\btheol{T1107}
	Let $n\in\NO$ and $\seq{\su{j}}{j}{0}{2n+1}\in\Hggequ{2n+1}$. For each $m\in\mn{0}{n}$, let $\seq{\su{j}^\sta{m}}{j}{0}{2(n-m)+1}$ be the \taSta{m}{\seq{\su{j}}{j}{0}{2n+1}}. Let
	\bgl{cF-2n+1}
		\HTiU{\seq{\su{j}}{j}{0}{2n+1}}
		\defg\HTiuu{\su{0}^\sta{0}}{\su{1}^\sta{0}}\circ\HTiuu{\su{0}^\seins}{\su{1}^\seins}\circ\dotsb\circ\HTiuu{\su{0}^\sta{n}}{\su{1}^\sta{n}}
	\eg
	\index{F_(-;2n+1)@$\HTiU{\seq{\su{j}}{j}{0}{2n+1}}$}and let
	\bgl{cF+2n+1}
		\HTU{\seq{\su{j}}{j}{0}{2n+1}}
		\defg\HTuu{\su{0}^\sta{n}}{\su{1}^\sta{n}}\circ\HTuu{\su{0}^\sta{n-1}}{\su{1}^\sta{n-1}}\circ\dotsb\circ\HTuu{\su{0}^\seins}{\su{1}^\seins}\circ\HTuu{\su{0}^\sta{0}}{\su{1}^\sta{0}}.
	\eg
	\index{F_(+;2n+1)@$\HTU{\seq{\su{j}}{j}{0}{2n+1}}$}Then the following statements hold true:
	\benui
		\item The mapping $\HTiU{\seq{\su{j}}{j}{0}{2n+1}}$ generates a bijective correspondence between $\PFoddqa{\su{0}^\sta{n}}$ and $\RFuqa{2n+1}{\seq{\su{j}}{j}{0}{2n+1}}$.
		\item Let $(\HTiU{\seq{\su{j}}{j}{0}{2n+1}})^\inv$ be the inverse mapping of $\HTiU{\seq{\su{j}}{j}{0}{2n+1}}$. For each matrix-valued function $F\in\RFuqa{2n+1}{\seq{\su{j}}{j}{0}{2n+1}}$, then
		\[
			(\HTiU{\seq{\su{j}}{j}{0}{2n+1}})^\inv(F)
			=\HTUa{\seq{\su{j}}{j}{0}{2n+1}}{F}
			=F^\Fa{2(n+1)},
		\]
		where $F^\Fa{2(n+1)}$ is the \tFaa{2(n+1)}{F}.
	\eenui
\etheo
\bproof
	Combine \rpropss{P0852}{P1035}.
\eproof

The following result can be considered as an inverse statement with respect to \rprop{P0921}.

\bpropl{P1309}
	Let $n\in\NO$ and let $\seq{\su{j}}{j}{0}{2n}\in\Hggequ{2n}$. For each $m\in\mn{0}{n}$, let $\seq{\su{j}^\sta{m}}{j}{0}{2(n-m)}$ be the \taSta{m}{\seq{\su{j}}{j}{0}{2n}}. Further, let $G\in\PFevenqa{\su{0}^\sta{n}}$ and let
	\[
		F
		\defg
		\begin{cases}
			\HTiua{\su{0}^\sta{n}}{G}\incase{n=0}\\
			[\HTiuu{\su{0}^\sta{0}}{\su{1}^\sta{0}}\circ\dotsb\circ\HTiuu{\su{0}^\sta{n-1}}{\su{1}^\sta{n-1}}\circ\HTiu{\su{0}^\sta{n}}](G)\incase{n\in\N}
		\end{cases}.
	\]
	Then the following statements hold true:
	\benui
		\il{P1309.a} $F\in\RFuqa{2n}{\seq{\su{j}}{j}{0}{2n}}$.
		\il{P1309.b} The \tFa{(2n+1)} $F^\Fa{2n+1}$ of $F$ satisfies $F^\Fa{2n+1}=G$.
	\eenui
\eprop
\bproof
	\eqref{P1309.a} This follows by combining \rprop{L1046} and \rtheo{L1241}.
	
	\eqref{P1309.b} Taking the definition of $F$ and~\eqref{P1309.a} into account, the application of \rcoross{C1103}{C1332} yields
	\[
		G
		=
		\begin{cases}
			\HTua{\su{0}^\sta{n}}{F}\incase{n=0}\\
			[\HTu{\su{0}^\sta{n}}\circ\HTuu{\su{0}^\sta{n-1}}{\su{1}^\sta{n-1}}\circ\dotsb\circ\HTuu{\su{0}^\sta{0}}{\su{1}^\sta{0}}](F)\incase{n\in\N}
		\end{cases}.
	\]
	Combining this with~\eqref{P1309.a}, we infer from \rpart{P0921.a} of \rprop{P0921} that
	\[
		F^\Fa{2n+1}
		=G.\qedhere
	\]
\eproof

The combination of \rpropss{P0921}{P1309} gives us now a complete description of the SN~algorithm in the class $\RFuqa{2n}{\seq{\su{j}}{j}{0}{2n}}$.

\btheol{T1421}
	Let $n\in\NO$ and $\seq{\su{j}}{j}{0}{2n}\in\Hggequ{2n}$. For each $m\in\mn{0}{n}$ let $\seq{\su{j}^\sta{m}}{j}{0}{2(n-m)}$ be the \taSta{m}{\seq{\su{j}}{j}{0}{2n}}. Let
	\bgl{cF-2n}
		\HTiU{\seq{\su{j}}{j}{0}{2n}}
		\defg
		\begin{cases}
			\HTiu{\su{0}^\sta{n}}\incase{n=0}\\
			\HTiuu{\su{0}^\sta{0}}{\su{1}^\sta{0}}\circ\dotsb\circ\HTiuu{\su{0}^\sta{n-1}}{\su{1}^\sta{n-1}}\circ\HTiu{\su{0}^\sta{n}}\incase{n\in\N}
		\end{cases}
	\eg
	\index{F_(-;2n)@$\HTiU{\seq{\su{j}}{j}{0}{2n}}$}and
	\bgl{cF+2n}
		\HTU{\seq{\su{j}}{j}{0}{2n}}
		\defg
		\begin{cases}
			\HTu{\su{0}^\sta{n}}\incase{n=0}\\
			\HTu{\su{0}^\sta{n}}\circ\HTuu{\su{0}^\sta{n-1}}{\su{1}^\sta{n-1}}\circ\dotsb\circ\HTuu{\su{0}^\sta{0}}{\su{1}^\sta{0}}\incase{n\in\N}
		\end{cases}.
	\eg
	\index{F_(+;2n)@$\HTU{\seq{\su{j}}{j}{0}{2n}}$}Then the following statements hold true:
	\benui
		\item The mapping $\HTiU{\seq{\su{j}}{j}{0}{2n}}$ generates a bijective correspondence between $\PFevenqa{\su{0}^\sta{n}}$ and $\RFuqa{2n}{\seq{\su{j}}{j}{0}{2n}}$.
		\item Let $(\HTiU{\seq{\su{j}}{j}{0}{2n}})^\inv$ be the inverse mapping of $\HTiU{\seq{\su{j}}{j}{0}{2n}}$. For each matrix-valued function $F\in\RFuqa{2n}{\seq{\su{j}}{j}{0}{2n}}$, then
		\[
			(\HTiU{\seq{\su{j}}{j}{0}{2n}})^\inv(F)
			=\HTUa{\seq{\su{j}}{j}{0}{2n}}{F}
			=F^\Fa{2n+1}
		\]
		where $F^\Fa{2n+1}$ is the \tFaa{(2n+1)}{F}.
	\eenui
\etheo
\bproof
	Combine \rpropss{P0921}{P1309}.
\eproof

Our next considerations are aimed to rewrite the mappings introduced in \rtheoss{T1107}{T1421} as linear fractional transformations of matrices. The essential tool in realizing this goal will be the matrix polynomials introduced in \rAppe{A0835}. More precisely, we will use finite products of such matrix polynomials.

Let $\kappa\in\NO\cup\set{+\infty}$ and let $\seq{\su{j}}{j}{0}{\kappa}$ be a sequence of complex \tpqa{matrices}. Let $n\in\NO$ be such that $2n\leq\kappa$. For all $m\in\mn{0}{n}$, let $\seq{\su{j}^\sta{m}}{j}{0}{\kappa-2m}$ be the \taSta{m}{\seq{\su{j}}{j}{0}{\kappa}}. For all $n\in\NO$ with $2n\leq\kappa$, let
\bgl{VVo2n}
	\rmiupou{s}{2n}
	\defg
	\begin{cases}
		\mHTiu{\su{0}^\sta{0}}\incase{n=0}\\
		\mHTiuu{\su{0}^\sta{0}}{\su{1}^\sta{0}}\mHTiuu{\su{0}^\sta{1}}{\su{1}^\sta{1}}\cdot\dotso\cdot\mHTiuu{\su{0}^\sta{n-1}}{\su{1}^\sta{n-1}}\mHTiu{\su{0}^\sta{n}}\incase{n\geq1}
	\end{cases}
\eg
\index{V^(^2n)@$\rmiupou{s}{2n}$}and, for all $n\in\NO$ with $2n+1\leq\kappa$, let
\bgl{VVo2n+1}
	\rmiupou{s}{2n+1}
	\defg\mHTiuu{\su{0}^\sta{0}}{\su{1}^\sta{0}}\mHTiuu{\su{0}^\sta{1}}{\su{1}^\sta{1}}\cdot\dotso\cdot\mHTiuu{\su{0}^\sta{n}}{\su{1}^\sta{n}}.
\eg
\index{V^(^2n+1)@$\rmiupou{s}{2n+1}$}Furthermore, for all $m\in\mn{0}{\kappa}$, let
\[
	\rmiupou{s}{m}
	=
	\bMat
		\rmiupnwou{s}{m}&\rmiupneou{s}{m}\\
		\rmiupswou{s}{m}&\rmiupseou{s}{m}
	\eMat
\]
\index{v_11^(^m)@$\rmiupnwou{s}{m}$}\index{v_12^(^m)@$\rmiupneou{s}{m}$}\index{v_21^(^m)@$\rmiupswou{s}{m}$}\index{v_22^(^m)@$\rmiupseou{s}{m}$}be the block representation of $\rmiupou{s}{m}$ with \tppa{block} $\rmiupnwou{s}{m}$.

\bremal{R1043}
	Let $\kappa\in\mn{2}{+\infty}\cup\set{+\infty}$ and let $\seq{\su{j}}{j}{0}{\kappa}$ be a sequence of complex \tpqa{matrices}. For all $n\in\N$ with $2n\leq\kappa$ and all $k\in\mn{0}{n-1}$, one can see then from \eqref{VVo2n}, \eqref{VVo2n+1}, and~\cite{103}*{\crema{9.2}} that
	\[
		\rmiupou{s^\sta{k}}{2(n-k)}
		=\rmiupou{s^\sta{k}}{2(n-k)-1}\mHTiu{\su{0}^\sta{n}}
	\]
	and
	\[
		\rmiupou{s^\sta{k}}{2(n-k)}
		=\mHTiuu{\su{0}^\sta{k}}{\su{1}^\sta{k}}\rmiupou{t}{2(n-k-1)},
	\]
	where $t_j\defg\su{j}^\sta{k+1}$ for all $j\in\mn{0}{2(n-k-1)}$.
\erema

\brema
	Let $\kappa\in\mn{3}{+\infty}\cup\set{+\infty}$ and let $\seq{\su{j}}{j}{0}{\kappa}$ be a sequence of complex \tpqa{matrices}. Then, for all $n\in\N$ with $2n+1\leq\kappa$ and all $k\in\mn{0}{n-1}$, one can see from \eqref{VVo2n+1} and~\cite{103}*{\crema{9.2}} that
	\[
		\rmiupou{s^\sta{k}}{2(n-k)+1}
		=\rmiupou{s^\sta{k}}{2(n-k)-1}\mHTiuu{\su{0}^\sta{n}}{\su{1}^\sta{n}}
	\]
	and
	\[
		\rmiupou{s^\sta{k}}{2(n-k)+1}
		=\mHTiuu{\su{0}^\sta{k}}{\su{1}^\sta{k}}\rmiupou{t}{2(n-k)-1},
	\]
	where $t_j\defg\su{j}^\sta{k+1}$ for all $j\in\mn{0}{2(n-k)-1}$.
\erema

Let $\kappa\in\NO\cup\set{+\infty}$ and let $\seq{\su{j}}{j}{0}{\kappa}$ be a sequence of complex \tpqa{matrices}. Let $n\in\NO$ be such that $2n\leq\kappa$. For all $m\in\mn{0}{n}$, let $\seq{\su{j}^\sta{m}}{j}{0}{\kappa-2m}$ be the \taSta{m}{\seq{\su{j}}{j}{0}{\kappa}}. Let
\bgl{WWo2n}
	\rmupou{s}{2n}
	\defg
	\begin{cases}
		\mHTu{\su{0}^\sta{0}}\incase{n=0}\\
		\mHTu{\su{0}^\sta{n}}\mHTuu{\su{0}^\sta{n-1}}{\su{1}^\sta{n-1}}\cdot\dotso\cdot\mHTuu{\su{0}^\sta{1}}{\su{1}^\sta{1}}\mHTuu{\su{0}^\sta{0}}{\su{1}^\sta{0}}\incase{n\geq1}
	\end{cases}.
\eg
\index{W^(^2n)@$\rmupou{s}{2n}$}For all $n\in\NO$ with $2n+1\leq\kappa$, let
\bgl{WWo2n+1}
	\rmupou{s}{2n+1}
	\defg\mHTuu{\su{0}^\sta{n}}{\su{1}^\sta{n}}\mHTuu{\su{0}^\sta{n-1}}{\su{1}^\sta{n-1}}\cdot\dotso\cdot\mHTuu{\su{0}^\sta{0}}{\su{1}^\sta{0}}
\eg
\index{W^(^2n+1)@$\rmupou{s}{2n+1}$}Furthermore, for all $m\in\mn{0}{\kappa}$, let
\[
	\rmupou{s}{m}
	=
	\bMat
		\rmupnwou{s}{m}&\rmupneou{s}{m}\\
		\rmupswou{s}{m}&\rmupseou{s}{m}
	\eMat
\]
\index{w_11^(^m)@$\rmupnwou{s}{m}$}\index{w_12^(^m)@$\rmupneou{s}{m}$}\index{w_21^(^m)@$\rmupswou{s}{m}$}\index{w_22^(^m)@$\rmupseou{s}{m}$}be the block representation of $\rmupou{s}{m}$ with \tppa{block} $\rmupnwou{s}{m}$.

\bremal{R1024}
	Let $\kappa\in\mn{2}{+\infty}\cup\set{+\infty}$ and let $\seq{\su{j}}{j}{0}{\kappa}$ be a sequence of complex \tpqa{matrices}. For all $n\in\N$ with $2n\leq\kappa$ and all $k\in\mn{0}{n-1}$, one can then see from \eqref{WWo2n}, \eqref{WWo2n+1}, and~\cite{103}*{\crema{9.2}} that
	\[
		\rmupou{s^\sta{k}}{2(n-k)}
		=\mHTu{\su{0}^\sta{n}}\rmupou{s^\sta{k}}{2(n-k)-1}
	\]
	and
	\[
		\rmupou{s^\sta{k}}{2(n-k)}=\rmupou{t}{2(n-k-1)}\mHTuu{\su{0}^\sta{k}}{\su{1}^\sta{k}},
	\]
	where $t_j\defg\su{j}^\sta{k+1}$ for all $j\in\mn{0}{2(n-k-1)}$.
\erema

\brema
	Let $\kappa\in\mn{3}{+\infty}\cup\set{+\infty}$ and let $\seq{\su{j}}{j}{0}{\kappa}$ be a sequence of complex \tpqa{matrices}. Then, for all $n\in\N$ with $2n+1\leq\kappa$ and all $k\in\mn{0}{n-1}$, one can see from \eqref{WWo2n+1} and~\cite{103}*{\crema{9.2}} that
	\[
		\rmupou{s^\sta{k}}{2(n-k)+1}
		=\mHTuu{\su{0}^\sta{n}}{\su{1}^\sta{n}}\rmupou{s^\sta{k}}{2(n-k)-1}
	\]
	and
	\[
		\rmupou{s^\sta{k}}{2(n-k)+1}
		=\rmupou{t}{2(n-k)-1}\mHTuu{\su{0}^\sta{k}}{\su{1}^\sta{k}},
	\]
	where $t_j\defg\su{j}^\sta{k+1}$ for all $j\in\mn{0}{2(n-k)-1}$.
\erema

\bpropl{P1746}
	Let $n\in\NO$, let $\seq{\su{j}}{j}{0}{2n+1}\in\Hggequ{2n+1}$ and let $\seq{\su{j}^\sta{n}}{j}{0}{1}$ be the \taSta{n}{\seq{\su{j}}{j}{0}{2n+1}}. Further, let $G\in\PFoddqa{\su{0}^\sta{n}}$. Then:
	\benui
		\il{P1746.a} For all $z\in\ohe$,
		\[
			\det\lek\rmiupswoua{s}{2n+1}{z}G(z)+\rmiupseoua{s}{2n+1}{z}\rek
			\neq0.
		\]
		\il{P1746.b} Let $\HTiU{\seq{\su{j}}{j}{0}{2n+1}}$ be given via \eqref{cF-2n+1}. Then
		\[
			\HTiUa{\seq{\su{j}}{j}{0}{2n+1}}{G}
			=\lftrfua{\rmiupou{s}{2n+1}}{G}.
		\]
	\eenui
\eprop
\bproof
	For each $m\in\mn{0}{n}$ let $\seq{\su{j}^\sta{m}}{j}{0}{2(n-m)+1}$ be the \taSta{m}{\seq{\su{j}}{j}{0}{2n+1}}. For each $m\in\mn{0}{n}$, in view of \rprop{P1531}, then
	\[
		\seq{\su{j}^\sta{m}}{j}{0}{2(n-m)+1}
		\in\Hggequ{2(n-m)+1}.
	\]
	Thus, for each $m\in\mn{0}{n}$, we have
	\bgl{P1746.1}
		\seq{\su{j}^\sta{m}}{j}{0}{1}
		\in\Hggequ{1}.
	\eg
	Because of \eqref{cF-2n+1}, we have
	\bgl{P1746.2}
		\HTiU{\seq{\su{j}}{j}{0}{2n+1}}
		=\HTiuu{\su{0}^\sta{0}}{\su{1}^\sta{0}}\circ\HTiuu{\su{0}^\seins}{\su{1}^\seins}\circ\dotsb\circ\HTiuu{\su{0}^\sta{n}}{\su{1}^\sta{n}}.
	\eg
	For each $k\in\mn{0}{n}$ we infer from \rprop{L1513} and \rtheo{L0925} inductively
	\begin{multline*}
		[\HTiuu{\su{0}^\sta{k}}{\su{1}^\sta{k}}\circ\HTiuu{\su{0}^\sta{k+1}}{\su{1}^\sta{k+1}}\circ\dotsb\circ\HTiuu{\su{0}^\sta{n}}{\su{1}^\sta{n}}](G)\\
		\in\RFuqA{2(n-k)+1}{\seq{\su{j}^\sta{k}}{j}{0}{2(n-k)+1}}.
	\end{multline*}
	For each $k\in\mn{0}{n}$, then \rlemm{L1553} shows that
	\bgl{P1746.3}
		[\HTiuu{\su{0}^\sta{k}}{\su{1}^\sta{k}}\circ\HTiuu{\su{0}^\sta{k+1}}{\su{1}^\sta{k+1}}\circ\dotsb\circ\HTiuu{\su{0}^\sta{n}}{\su{1}^\sta{n}}](G)
		\in\PFoddqa{\su{0}^\sta{k}}.
	\eg
	Taking \eqref{P1746.1}, \eqref{P1746.2}, and \eqref{P1746.3} into account, we get from \rlemm{L1300} that
	\bgl{P1746.4}
		\HTiUa{\seq{\su{j}}{j}{0}{2n+1}}{G}
		=(\lftrfu{\mHTiuu{\su{0}^\sta{0}}{\su{1}^\sta{0}}}\circ\lftrfu{\mHTiuu{\su{0}^\seins}{\su{1}^\seins}}\circ\dotsb\circ\lftrfu{\mHTiuu{\su{0}^\sta{n}}{\su{1}^\sta{n}}})(G).
	\eg
	In view of \eqref{P1746.4} and \eqref{VVo2n+1}, now \rprop{P1509} yields~\eqref{P1746.a} and~\eqref{P1746.b}.
\eproof

\bpropl{P1814}
	Let $n\in\NO$ and let $\seq{\su{j}}{j}{0}{2n+1}\in\Hggequ{2n+1}$. Further, let $F$ belong to $\RFuqa{2n+1}{\seq{\su{j}}{j}{0}{2n+1}}$. Then the following statements hold true:
	\benui
		\il{P1814.a} For all $z\in\ohe$,
		\[
			\det\lek\rmupswoua{s}{2n+1}{z}F(z)+\rmupseoua{s}{2n+1}{z}\rek
			\neq0.
		\]
		\il{P1814.b} Let $\HTU{\seq{\su{j}}{j}{0}{2n+1}}$ be given by \eqref{cF+2n+1}. Then
		\[
			\HTUa{\seq{\su{j}}{j}{0}{2n+1}}{F}
			=\lftrfua{\rmupou{s}{2n+1}}{F}.
		\]
	\eenui
\eprop
\bproof
	For each $m\in\mn{0}{n}$ let $\seq{\su{j}^\sta{m}}{j}{0}{2(n-m)+1}$ be the \taSta{m}{\seq{\su{j}}{j}{0}{2n+1}}. For each $m\in\mn{0}{n}$ in view of \rprop{P1531} then
	\[
		\seq{\su{j}^\sta{m}}{j}{0}{2(n-m)+1}
		\in\Hggequ{2(n-m)+1}.
	\]
	Thus, for every choice of $m\in\mn{0}{n}$, we have
	\bgl{P1814.1}
		\seq{\su{j}^\sta{m}}{j}{0}{1}
		\in\Hggequ{1}.
	\eg
	Because of \eqref{cF+2n+1}, we have
	\bgl{P1814.2}
		\HTU{\seq{\su{j}}{j}{0}{2n+1}}
		=\HTuu{\su{0}^\sta{n}}{\su{1}^\sta{n}}\circ\HTuu{\su{0}^\sta{n-1}}{\su{1}^\sta{n-1}}\circ\dotsb\circ\HTuu{\su{0}^\sta{0}}{\su{1}^\sta{0}}.
	\eg
	For $n=0$, the assertions of~\eqref{P1814.a} and~\eqref{P1814.b} are an immediate consequence of \eqref{P1814.2}, \eqref{VVo2n+1}, and \rpart{L1239.c} of \rlemm{L1239}. Now let $n\in\N$. For each $k\in\mn{0}{n-1}$, we infer from \rtheo{L1314} inductively that
	\bgl{P1814.3}
		[\HTuu{\su{0}^\sta{k}}{\su{1}^\sta{k}}\circ\dotsb\circ\HTuu{\su{0}^\sta{0}}{\su{1}^\sta{0}}](F)
		\in\RFuqA{2(n-k)+1}{\seq{\su{j}^\sta{k}}{j}{0}{2(n-k)+1}}.
	\eg
	Taking \eqref{P1814.1}, \eqref{P1814.2}, and \eqref{P1814.3} into account, \rpart{L1239.c} of \rlemm{L1239} yields
	\bgl{P1814.4}
		\HTUa{\seq{\su{j}}{j}{0}{2n+1}}{F}
		=(\lftrfu{\mHTuu{\su{0}^\sta{n}}{\su{1}^\sta{n}}}\circ\lftrfu{\mHTuu{\su{0}^\sta{n-1}}{\su{1}^\sta{n-1}}}\circ\dotsb\circ\lftrfu{\mHTuu{\su{0}^\sta{0}}{\su{1}^\sta{0}}})(F).
	\eg
	In view of \eqref{P1814.4} and \eqref{WWo2n+1}, now \rprop{P1509} yields~\eqref{P1814.a} and~\eqref{P1814.b}.
\eproof

\bpropl{P1245}
	Let $n\in\NO$, let $\seq{\su{j}}{j}{0}{2n}\in\Hggequ{2n}$, and let $\seq{\su{j}^\sta{n}}{j}{0}{0}$ be the \taSta{n}{\seq{\su{j}}{j}{0}{2n}}. Further, let $G\in\PFevenqa{\su{0}^\sta{n}}$. Then:
	\benui
		\il{P1245.a} For all $z\in\ohe$,
		\[
			\det\lek\rmiupswoua{s}{2n}{z}G(z)+\rmiupseoua{s}{2n}{z}\rek
			\neq0.
		\]
		\il{P1245.b} Let $\HTiU{\seq{\su{j}}{j}{0}{2n}}$ be given by \eqref{cF-2n}. Then
		\[
			\HTiUa{\seq{\su{j}}{j}{0}{2n}}{G}
			=\lftrfua{\rmiupou{s}{2n}}{G}.
		\]
	\eenui
\eprop
\bproof
	From \eqref{cF-2n} and \eqref{cF-2n+1} we infer that
	\bgl{P1245.1}
		\HTiU{\seq{\su{j}}{j}{0}{2n}}
		=
		\begin{cases}
			\HTiu{\su{0}^\sta{n}}\incase{n=0}\\
			\HTiU{\seq{\su{j}}{j}{0}{2n-1}}\circ\HTiu{\su{0}^\sta{n}}\incase{n\in\N}
		\end{cases}.
	\eg
	First we consider the case $n=0$. Taking \eqref{VVo2n} and \eqref{P1245.1} into account, \rlemm{L1342} yields the assertions of~\eqref{P1245.a} and~\eqref{P1245.b}. Now let $n\in\N$. Then clearly
	\bgl{P1245.2}
		\seq{\su{j}}{j}{0}{2n-1}
		\in\Hggequ{2n-1}.
	\eg
	Furthermore, \rprop{P1531} implies $\seq{\su{j}^\sta{n}}{j}{0}{0}\in\Hggequ{0}$. Then \rlemm{L1342} provides us
	\bgl{P1245.4}
		\HTiua{\su{0}^\sta{n}}{G}
		=\lftrfua{\mHTiu{\su{0}^\sta{n}}}{G}.
	\eg
	From \rcoro{C1103} we see that
	\[
		\HTiua{\su{0}^\sta{n}}{G}
		\in\RFuqA{0}{\seq{\su{j}^\sta{n}}{j}{0}{0}}.
	\]
	Thus, \rlemm{L1553} yields
	\bgl{P1245.5}
		\HTiua{\su{0}^\sta{n}}{G}
		\in\PFoddqa{\su{0}^\sta{n}}.
	\eg
	Let $\seq{\su{j}^\sta{n-1}}{j}{0}{2}$ be the \taSta{(n-1)}{\seq{\su{j}}{j}{0}{2n}}. In view of \eqref{seins}, then $\su{0}^\sta{n}=-\su{0}^\sta{n-1}(\su{2}^\sta{n-1})^\MP\su{0}^\sta{n-1}$. In particular, $\Kerna{\su{0}^\sta{n-1}}\subseteq\Kerna{\su{0}^\sta{n}}$ and \rrema{R1543} implies
	\bgl{P1245.6}
		\PFoddqa{\su{0}^\sta{n}}
		\subseteq\PFoddqa{\su{0}^\sta{n-1}}.
	\eg
	From \eqref{P1245.5} and \eqref{P1245.6} we get $\HTiua{\su{0}^\sta{n}}{G}\in\PFoddqa{\su{0}^\sta{n-1}}$, and taking \eqref{P1245.2} into account, \rprop{P1746} yields
	\bgl{P1245.8}
		\HTiUA{\seq{\su{j}}{j}{0}{2n-1}}{\HTiua{\su{0}^\sta{n}}{G}}
		=\lftrfuA{\rmiupou{s}{2n-1}}{\HTiua{\su{0}^\sta{n}}{G}}.
	\eg
	Using \eqref{P1245.1}, \eqref{P1245.8}, and \eqref{P1245.4}, we get
	\bgl{P1245.9}
		\HTiUa{\seq{\su{j}}{j}{0}{2n}}{G}
		=(\lftrfu{\rmiupou{s}{2n-1}}\circ\lftrfu{\mHTiu{\su{0}^\sta{n}}})(G).
	\eg
	In view of \rrema{R1043}, we have
	\bgl{P1245.10}
		\rmiupou{s}{2n}
		=\rmiupou{s}{2n-1}\mHTiu{\su{0}^\sta{n}}.
	\eg
	In view of \eqref{P1245.9} and \eqref{P1245.10}, \rprop{P1509} yields then~\eqref{P1245.a} and~\eqref{P1245.b}.
\eproof

\bprop
	Let $n\in\NO$, let $\seq{\su{j}}{j}{0}{2n}\in\Hggequ{2n}$ and let $F\in\RFuqa{2n}{\seq{\su{j}}{j}{0}{2n}}$. Then the following statements hold true:
	\benui
		\il{P1320.a} For all $z\in\ohe$,
		\[
			\det\lek\rmupswoua{s}{2n}{z}F(z)+\rmupseoua{s}{2n}{z}\rek
			\neq0.
		\]
		\il{P1320.b} Let $\HTU{\seq{\su{j}}{j}{0}{2n}}$ be given by \eqref{cF+2n}. Then
		\[
			\HTUa{\seq{\su{j}}{j}{0}{2n}}{F}
			=\lftrfua{\rmupou{s}{2n}}{F}.
		\]
	\eenui
\eprop
\bproof
	Let $\seq{\su{j}^\sta{n}}{j}{0}{0}$ be the \taSta{n}{\seq{\su{j}}{j}{0}{2n}}. From \eqref{cF+2n} and \eqref{cF+2n+1} we infer that
	\bgl{P1320.1}
		\HTU{\seq{\su{j}}{j}{0}{2n}}
		=
		\begin{cases}
			\HTu{\su{0}^\sta{n}}\incase{n=0}\\
			\HTu{\su{0}^\sta{n}}\circ\HTU{\seq{\su{j}}{j}{0}{2n-1}}\incase{n\in\N}
		\end{cases}.
	\eg
	First we consider the case $n=0$. Taking \eqref{WWo2n} and \eqref{P1320.1} into account the application of \rpartss{L1239.a}{L1239.b} of \rlemm{L1239} yields the assertions of~\eqref{P1320.a} and~\eqref{P1320.b}. Now let $n\in\N$. Then clearly
	\bgl{P1320.2}
		\seq{\su{j}}{j}{0}{2n-1}
		\in\Hggequ{2n-1}.
	\eg
	From \eqref{Ruua} we see that
	\bgl{P1320.3}
		F
		\in\RFuqA{2n-1}{\seq{\su{j}}{j}{0}{2n-1}}.
	\eg
	In view of \eqref{P1320.2} and \eqref{P1320.3}, we conclude from \rprop{P1814} that
	\bgl{P1320.4}
		\HTUa{\seq{\su{j}}{j}{0}{2n-1}}{F}
		=\lftrfua{\rmupou{s}{2n-1}}{F}
	\eg
	and from \rpart{P0852.a} of \rprop{P0852} that
	\bgl{P1320.5}
		\HTUa{\seq{\su{j}}{j}{0}{2n-1}}{F}
		\in\RFuqA{0}{\seq{\su{j}^\sta{n}}{j}{0}{0}}.
	\eg
	In view of \rprop{P1531}, we have $\seq{\su{j}^\sta{n}}{j}{0}{0}\in\Hggequ{0}$. Then from \eqref{P1320.5} and \rpartss{L1239.a}{L1239.b} of \rlemm{L1239} we get
	\bgl{P1320.7}
		\HTua{\su{0}^\sta{n}}{\HTUA{\seq{\su{j}}{j}{0}{2n-1}}{F}}
		=\lftrfuA{\mHTu{\su{0}^\sta{n}}}{\HTUa{\seq{\su{j}}{j}{0}{2n-1}}{F}}.
	\eg
	Using \eqref{P1320.1}, \eqref{P1320.7}, and \eqref{P1320.4}, we obtain
	\bgl{P1320.8}
		\HTUa{\seq{\su{j}}{j}{0}{2n}}{F}
		=(\lftrfu{\mHTu{\su{0}^\sta{n}}}\circ\lftrfu{\rmupou{s}{2n-1}})(F).
	\eg
	In view of \rrema{R1024}, we have
	\bgl{P1320.9}
		\rmupou{s}{2n}
		=\mHTu{\su{0}^\sta{n}}\rmupou{s}{2n-1}.
	\eg
	Taking \eqref{P1320.8} and \eqref{P1320.9} into account, the application of \rprop{P1509} yields the assertions of~\eqref{P1320.a} and~\eqref{P1320.b}.
\eproof

\section{Descriptions of the Sets $\RFuqa{m}{\seq{\su{j}}{j}{0}{m}}$}\label{S1413}
Let $m\in\NO$ and let $\seq{\su{j}}{j}{0}{m}\in\Hggequ{m}$. Then we know from \rtheo{1-417} that the solution set $\MggqRag{\seq{\su{j}}{j}{0}{m}}$ of \rprob{\mproblem{\R}{m}{=}} is non-empty. Furthermore, \rcoro{R1103} tells us that $\RFuqa{m}{\seq{\su{j}}{j}{0}{m}}$ is exactly the set of \tStit{s} of all the measures belonging to $\MggqRag{\seq{\su{j}}{j}{0}{m}}$. On the basis of our Schur-Nevanlinna-type algorithm, which was introduced in \rsect{S1610}, we have obtained important insights into the structure of the set $\RFuqa{m}{\seq{\su{j}}{j}{0}{m}}$. Using \rtheo{T1512}, we will now rewrite the descriptions of $\RFuqa{m}{\seq{\su{j}}{j}{0}{m}}$, which were obtained in \rsect{S1610}, in a form which is better adapted to the original data sequence $\seq{\su{j}}{j}{0}{m}$. This is our first main theorem.

\btheol{T1425}
	Let $n\in\NO$, let $\seq{\su{j}}{j}{0}{2n}\in\Hggequ{2n}$, let $\rmiupou{s}{2n}$ be defined by \eqref{VVo2n} and let $\Lu{n}$ be given by \eqref{Lu}. Then the following statements hold true:
	\benui
		\il{T1425.a} $\RFuqa{2n}{\seq{\su{j}}{j}{0}{2n}}=\lftrfuA{\rmiupou{s}{2n}}{\PFevenqa{\Lu{n}}}$.
		\il{T1425.b} For each $F\in\RFuqa{2n}{\seq{\su{j}}{j}{0}{2n}}$, there is a unique $G\in\PFevenqa{\Lu{n}}$ which satisfies
		\[
			\lftrfuA{\rmiupou{s}{2n}}{G}
			=F,
		\]
		namely $G=F^\Fa{2n+1}$, where $F^\Fa{2n+1}$ stands for the \tFaa{(2n+1)}{F}.
	\eenui
\etheo
\bproof
	Let $\seq{\su{j}^\sta{n}}{j}{0}{0}$ be the \taSta{n}{\seq{\su{j}}{j}{0}{2n}}. In view of \rtheo{T1512} and \rdefi{D0944}, then $\su{0}^\sta{n}=\Lu{n}$. Hence, the application of \rtheo{T1421} and \rprop{P1245} completes the proof.
\eproof

It should be mentioned that a result similar to \rpart{T1425.a} of \rtheo{T1425} is contained in~\cite{MR1624548}*{\ctheo{4.1}}.

\bcorol{C1439}
	Let $n\in\NO$, let $\seq{\su{j}}{j}{0}{2n}\in\Hgqu{2n}$ and let $\rmiupou{s}{2n}$ be defined by \eqref{VVo2n}. Then
	\[
		\RFuqA{2n}{\seq{\su{j}}{j}{0}{2n}}
		=\lftrfuA{\rmiupou{s}{2n}}{\NFuq{-2}}.
	\]
\ecoro
\bproof
	From \rrema{R1600} we get
	\bgl{C1439.1}
		\seq{\su{j}}{j}{0}{2n}
		\in\Hggequ{2n}.
	\eg
	Because of the assumption $\seq{\su{j}}{j}{0}{2n}\in\Hgqu{2n}$, we infer from \rpart{L1602.c} of \rlemm{L1602} that the matrix $\Lu{n}$ defined in \eqref{Lu} is positive Hermitian. In particular, $\det\Lu{n}\neq0$. Consequently, \rrema{R0827} implies $\PFevenqa{\Lu{n}}=\NFuq{-2}$. Thus, taking \eqref{C1439.1} into account, the application of \rpart{T1425.a} of \rtheo{T1425} completes the proof.
\eproof

Observe, that under the assumptions of \rcoro{C1439}, an alternative description of the class $\RFuqa{2n}{\seq{\su{j}}{j}{0}{2n}}$ was obtained in~\cite{MR1018213}*{\ctheo{8.2}} with the aid of the RKHS~method.


Let $n\in\NO$ and let $\seq{\su{j}}{j}{0}{2n}\in\Hggqu{2n}$. Then $\seq{\su{j}}{j}{0}{2n}$ is called \emph{\tcd{}}, if $\Lu{n}=\Oqq$, where $\Lu{n}$ is defined in \eqref{Lu}. Observe that the set $\Hggdqu{2n}$\index{H_,2^>=,cd@$\Hggdqu{2n}$} of all \tcd{} sequences belonging to $\Hggqu{2n}$ is a subclass of $\Hggequ{2n}$ (see~\cite{MR2570113}*{\ccoro{2.14}}). This class is connected to the case of a unique solution, which was already discussed in~\cite{MR1624548}*{\ccoro{3.5}}.


\btheol{T1250*}
	Let $n\in\NO$ and let $\seq{\su{j}}{j}{0}{2n}\in\Hggequ{2n}$. Then:
	\benui
		\item The following statements are equivalent:
		\baeqii{0}
			\il{T1250*.i} The set $\RFuqa{2n}{\seq{\su{j}}{j}{0}{2n}}$ consists of exactly one element.
			\il{T1250*.ii} $\seq{\su{j}}{j}{0}{2n}\in\Hggdqu{2n}$.
		\eaeqii
		\il{T1250*.b} If~\ref{T1250*.i} holds true, then $\det\rmiupseoua{s}{2n}{z}\neq0$ for all $z\in\ohe$ and
		\[
			\RuuA{2n}{q}{\seq{\su{j}}{j}{0}{2n}}\\
			=\Set{\rmiupneou{s}{2n}(\rmiupseou{s}{2n})^\inv}.
		\]
	\eenui
\etheo
\bproof
	\bimp{T1250*.i}{T1250*.ii}
		Use \rtheo{T0933}, \rcoro{C1403}, \eqref{Ruua}, \rcoro{R1103}, and~\cite{MR2570113}*{\ctheo{8.7}}.
	\eimp
	
	\bimp{T1250*.ii}{T1250*.i}
		Let $\Lu{n}$ be defined in \eqref{Lu}. Then~\ref{T1250*.ii} means $\Lu{n}=\Oqq$. Thus, \rpart{L1249*.a} of \rprop{L1249*}, implies that the set $\PFevenqa{\Lu{n}}$ consists of exactly one element, namely the constant function defined on $\ohe$ with value $\Oqq$. Using \rpart{T1425.a} of \rtheo{T1425}, then~\ref{T1250*.i} follows, and we see moreover that~\eqref{T1250*.b} holds true.
	\eimp
\eproof

Now we formulate our second main theorem.

\btheol{T1452}
	Let $n\in\NO$, let $\seq{\su{j}}{j}{0}{2n+1}\in\Hggequ{2n+1}$, let $\rmiupou{s}{2n+1}$ be defined by \eqref{VVo2n+1}, and let $\Lu{n}$ be given by \eqref{Lu}. Then the following statements hold:
	\benui
		\il{T1452.a} $\RFuqa{2n+1}{\seq{\su{j}}{j}{0}{2n+1}}=\lftrfuA{\rmiupou{s}{2n+1}}{\PFoddqa{\Lu{n}}}$.
		\il{T1452.b} For each $F\in\RFuqa{2n+1}{\seq{\su{j}}{j}{0}{2n+1}}$, there is a unique $G\in\PFoddqa{\Lu{n}}$ which satisfies
		\[
			\lftrfuA{\rmiupou{s}{2n+1}}{G}
			=F,
		\]
		namely $G=F^\Fa{2(n+1)}$, where $F^\Fa{2(n+1)}$ stands for the \tFaa{2(n+1)}{F}.
	\eenui
\etheo
\bproof
	Let $\seq{\su{j}^\sta{n}}{j}{0}{1}$ be the \taSta{n}{\seq{\su{j}}{j}{0}{2n+1}}. From the choice of $\seq{\su{j}}{j}{0}{2n+1}$ we get $\seq{\su{j}}{j}{0}{2n}\in\Hggequ{2n}$. Thus, we infer from \rtheo{T1512} and \rdefi{D0944} then $\su{0}^\sta{n}=\Lu{n}$. Hence, the application of \rtheo{T1107} and \rprop{P1746} completes the proof.
\eproof

In the scalar case $q=1$, the following result goes back to~\cite{2011arXiv1101.0162D}*{\ccoro{5.2}}.

\bcoro
	Let $n\in\NO$, let $\seq{\su{j}}{j}{0}{2n+1}\in\Hgequ{2n+1}$ and let $\rmiupou{s}{2n+1}$ be defined by \eqref{VVo2n+1}. Then
	\[
		\RFuqA{2n+1}{\seq{\su{j}}{j}{0}{2n+1}}
		=\lftrfuA{\rmiupou{s}{2n+1}}{\RFuq{-1}}.
	\]
\ecoro
\bproof
	From $\seq{\su{j}}{j}{0}{2n+1}\in\Hgequ{2n+1}$ we get
	\bgl{C1459.1}
		\seq{\su{j}}{j}{0}{2n+1}
		\in\Hggequ{2n+1}.
	\eg
	Because of \rpart{P1642.a} of \rprop{P1642}, we also get $\seq{\su{j}}{j}{0}{2n}\in\Hgqu{2n}$. Then \rpart{L1602.c} of \rlemm{L1602} yields $\Lu{n}\in\Cgq$. In particular, $\det\Lu{n}\neq0$. Consequently, \rrema{R0827} implies $\PFoddqa{\Lu{n}}=\RFuq{-1}$. In view of \eqref{C1459.1} and \rpart{T1452.a} of \rtheo{T1452}, this completes the proof.
\eproof

Now we derive an analogue of \rtheo{T1250*}.

\btheo
	Let $n\in\NO$ and let $\seq{\su{j}}{j}{0}{2n+1}\in\Hggequ{2n+1}$. Then:
	\benui
		\item The following statements are equivalent:
		\baeqii{0}
			\il{T1054*.i} The set $\RFuqa{2n+1}{\seq{\su{j}}{j}{0}{2n+1}}$ consists of exactly one element.
			\il{T1054*.ii} $\seq{\su{j}}{j}{0}{2n}\in\Hggdqu{2n}$.
		\eaeqii
		\il{T1054*.b} If~\ref{T1054*.i} holds true, then $\det\rmiupseoua{s}{2n+1}{z}\neq0$ for all $z\in\ohe$ and
		\[
			\RuuA{2n+1}{q}{\seq{\su{j}}{j}{0}{2n+1}}\\
			=\Set{\rmiupneou{s}{2n+1}(\rmiupseou{s}{2n+1})^\inv}.
		\]
	\eenui
\etheo
\bproof
	\bimp{T1054*.i}{T1054*.ii}
		Use \rtheo{T0933}, \rcoro{C1403}, \eqref{Ruua}, \rcoro{R1103}, and~\cite{MR2570113}*{\ctheo{8.9}}.
	\eimp
	
	\bimp{T1054*.ii}{T1054*.i}
		Let $\Lu{n}$ be defined in \eqref{Lu}. Then~\ref{T1054*.ii} means $\Lu{n}=\Oqq$. Thus, \rpart{L1249*.a} of \rprop{L1249*} implies that the set $\PFoddqa{\Lu{n}}$ consists of exactly one element, namely the constant function defined on $\ohe$ with value $\Oqq$. Using \rpart{T1452.a} of \rtheo{T1452}, then~\ref{T1054*.i} and~\eqref{T1054*.b} follow.
	\eimp
\eproof

By using \rprop{L1249*}, we are able to derive alternate descriptions of the set $\RFuqa{m}{\seq{\su{j}}{j}{0}{m}}$ if $m\in\NO$ and $\seq{\su{j}}{j}{0}{m}\in\Hggequ{m}$ are arbitrarily given. This gives reformulations of our two main results stated in \rtheoss{T1425}{T1452}.

\btheol{T1620}
	Let $n\in\NO$ and let $\seq{\su{j}}{j}{0}{2n}\in\Hggequ{2n}$. Let $\Lu{n}$ be defined in \eqref{Lu}. Suppose that $r\defg\rank\Lu{n}$ fulfills $r\geq1$. Let $u_1,u_2,\dotsc,u_r$ be an orthonormal basis of $\Bilda{\Lu{n}}$ and let $U\defg\brow u_1,u_2,\dotsc,u_r\erow$. Then:
	\benui
		\il{T1620.a} $\RFuqa{2n}{\seq{\su{j}}{j}{0}{2n}}=\lftrfua{\rmiupou{s}{2n}}{\setaa{UfU^\ad}{f\in\NFuu{-2}{r}}}$.
		\il{T1620.b} Let $F\in\RFuqa{2n}{\seq{\su{j}}{j}{0}{2n}}$. Then there is a unique $f\in\NFuu{-2}{r}$ such that
		\[
			\lftrfua{\rmiupou{s}{2n}}{UfU^\ad}
			=F,
		\]
		namely $f\defg U^\ad F^\Fa{2n+1}U$, where $F^\Fa{2n+1}$ is the \tFaa{(2n+1)}{F}.
	\eenui
\etheo
\bproof
	Since $\seq{\su{j}}{j}{0}{2n}\in\Hggequ{2n}$ implies $\seq{\su{j}}{j}{0}{2n}\in\Hggqu{2n}$ we get from \rpart{L1031.d} of \rlemm{L1031} that $\Lu{n}\in\Cggq$. In particular, $\Lu{n}^\ad=\Lu{n}$. Thus, \rpart{L1249*.b} of \rprop{L1249*} yields
	\bgl{T1620.2}
		\PFevenqa{\Lu{n}}
		=\SetaA{UfU^\ad}{f\in\NFuu{-2}{r}},
	\eg
	whereas \rrema{R1526} implies $U^\ad U=\Iu{r}$.
	
	\eqref{T1620.a} Because of \eqref{T1620.2}, we infer~\eqref{T1620.a} from \rpart{T1425.a} of \rtheo{T1425}.
	
	\eqref{T1620.b} Using \eqref{T1620.2}, $U^\ad U=\Iu{r}$, and~\eqref{T1620.a}, we conclude~\eqref{T1620.b} from \rpart{T1425.b} of \rtheo{T1425}.
\eproof

\btheol{T1640}
	Let $n\in\NO$ and let $\seq{\su{j}}{j}{0}{2n+1}\in\Hggequ{2n+1}$. Let $\Lu{n}$ be defined in \eqref{Lu}. Suppose that $r\defg\rank\Lu{n}$ fulfills $r\geq1$. Let $u_1,u_2,\dotsc,u_r$ be an orthonormal basis of $[\Kerna{\su{0}^\sta{n}}]^\bot$ and let $U\defg\brow u_1,u_2,\dotsc,u_r\erow$. Then:
	\benui
		\il{T1640.a} $\RFuqa{2n+1}{\seq{\su{j}}{j}{0}{2n+1}}=\lftrfua{\rmiupou{s}{2n+1}}{\setaa{UfU^\ad}{f\in\RFuu{-1}{r}}}$.
		\il{T1640.b} Let $F\in\RFuqa{2n+1}{\seq{\su{j}}{j}{0}{2n+1}}$. Then there is a unique $f\in\RFuu{-1}{r}$ such that
		\[
			\lftrfua{\rmiupou{s}{2n+1}}{UfU^\ad}
			=F,
		\]
		namely $f\defg U^\ad F^\Fa{2(n+1)}U$, where $F^\Fa{2(n+1)}$ is the \tFaa{2(n+1)}{F}.
	\eenui
\etheo
\bproof
	Since $\seq{\su{j}}{j}{0}{2n+1}\in\Hggequ{2n+1}$ implies $\seq{\su{j}}{j}{0}{2n}\in\Hggqu{2n}$, we get again from \rpart{L1031.d} of \rlemm{L1031} then $\Lu{n}^\ad=\Lu{n}$. Thus, \rpart{L1249*.b} of \rprop{L1249*} yields
	\bgl{T1640.2}
		\PFoddqa{\Lu{n}}
		=\SetaA{UfU^\ad}{f\in\RFuu{-1}{r}},
	\eg
	whereas \rrema{R1526} implies $U^\ad U=\Iu{r}$.
	
	\eqref{T1640.a} In view of \eqref{T1640.2}, we infer~\eqref{T1640.a} from \rpart{T1452.a} of \rtheo{T1452}.
	
	\eqref{T1640.b} Using \eqref{T1640.2}, $U^\ad U=\Iu{r}$, and~\eqref{T1640.a}, we get~\eqref{T1640.b} from \rpart{T1452.b} of \rtheo{T1452}.
\eproof

\appendix\section{Some Results On Moore-Penrose Inverses of Matrices}\label{A1549}
For the convenience of the reader, we state some well-known and some special results on Moore-Penrose inverses of matrices (see e.\,g., Rao/Mitra~\cite{MR0338013} or~\cite{MR1152328}*{\cSect{1}}).
If $A\in\Cpq$, then (by definition) the Moore-Penrose inverse $A^\MP$ of $A$ is the unique matrix $A^\MP\in\Cqp$\index{^+@$A^\MP$} which satisfies the four equations
\begin{align*}
	AA^\MP A&=A,&A^\MP AA^\MP&=A^\MP,&(AA^\MP)^\ad&=AA^\MP,&
	&\text{and}&
	(A^\MP A)^\ad&=A^\MP A.
\end{align*}

\bremal{R1553}
	Let $A\in\Cpq$. Then one can easily check that:
	\benui
		\il{R1553.a} $(A^\MP)^\MP=A$, $(A^\MP)^\ad=(A^\ad)^\MP$, and $\Ip-AA^\MP\in\Cggq$.
		\il{R1553.c} $\Bilda{A^\MP}=\Bilda{A^\ad}$, $\rank(A^\MP)=\rank A$, and $\Kerna{A^\MP}=\Kerna{A^\ad}$.
	\eenui
\erema


\begin{prop}[see, e.\,g.~\cite{MR1152328}*{\ctheo{1.1.1}}]\label{P1619}
	If $A\in\Cpq$ then a matrix $G\in\Cqp$ is the Moore-Penrose inverse of $A$ if and only if $AG=P_A$ and $GA=P_G$, where $P_A$ and $P_G$ are respectively, the matrices associated with the orthogonal projection in $\Cp$ onto $\Bilda{A}$ and the orthogonal projection in $\Cq$ onto $\Bilda{G}$.
\end{prop}

\bremal{L1622}
	Let $A\in\Cpq$. Then it is readily checked that:
	\benui
		\il{L1622.a} Let $r\in\N$ and let $B\in\Coo{r}{q}$. Then $\Kerna{A}\subseteq\Kerna{B}$ if and only if $BA^\MP A=B$. Furthermore, $\Kerna{B}=\Kerna{A}$ if and only if $B^\MP B=A^\MP A$.
		\il{L1622.b} Let $s\in\N$ and let $C\in\Coo{p}{s}$. Then $\Bilda{C}\subseteq\Bilda{A}$ if and only if $AA^\MP C=C$. Furthermore, $\Bilda{C}=\Bilda{A}$ if and only if $CC^\MP=AA^\MP$.
	\eenui
\erema

\bremal{L0818}
	Let $A,B\in\CHq$. From \rremass{L1622}{R1553} one can easily see that the following statements are equivalent:
	\baeqi{0}
		\item $\Kerna{A}\subseteq\Kerna{B}$.
		\item $BA^\MP A=B$.
		\item $AA^\MP B=B$.
		\item $\Bilda{A}\subseteq\Bilda{B}$.
	\eaeqi
\erema

\blemml{L1633}
	Let $A,X\in\Cpq$. Then the following statements are equivalent:
	\baeqi{0}
		\il{L1633.i} $\Kerna{A}\subseteq\Kerna{X}$ and $\Bilda{A}\subseteq\Bilda{X}$.
		\il{L1633.ii} $\Kerna{A}=\Kerna{X}$ and $\Bilda{A}=\Bilda{X}$.
		\il{L1633.iii} $A^\MP A=X^\MP X$ and $AA^\MP=XX^\MP$.
		\il{L1633.iv} $\Kerna{A^\MP}=\Kerna{X^\MP}$ and $\Bilda{A^\MP}=\Bilda{X^\MP}$.
	\eaeqi
\elemm
\bproof
	In view of $\dim[\Kerna{A}]+\dim[\Bilda{A}]=q$ and $\dim[\Kerna{X}]+\dim[\Bilda{X}]=q$, the implication~\impl{L1633.i}{L1633.ii} is true. Otherwise, the implication~\impl{L1633.ii}{L1633.i} is trivial. \rPart{L1622.b} of \rrema{L1622} yields the equivalence of~\ref{L1633.ii} and~\ref{L1633.iii}. Taking \rpart{R1553.a} of \rrema{R1553} into account, we see that the equivalence of~\ref{L1633.iii} and~\ref{L1633.iv} is an immediate consequence of the already verified equivalence of~\ref{L1633.ii} and~\ref{L1633.iii}.
\eproof

\bremal{R1526}
	Let $A\in\Cpq\setminus\set{\Opq}$. Let $r\defg\rank A$, let $u_1,u_2,\dotsc,u_r$ be an orthonormal basis of $\Bilda{A^\ad}$, and let $U\defg\brow u_1,u_2,\dotsc,u_r\erow$. Then $U^\ad U=\Iu{r}$, and, in view of \rprop{P1619} and \rrema{R1553}, furthermore $UU^\ad=A^\MP A$.
\erema

\blemml{L1315}
	Let $A,B\in\CHq$. Then the following statements are equivalent:
	\baeqi{0}
		\il{L1315.i} $A\geq B\geq\Oqq$.
		\il{L1315.ii} $B^\MP\geq B^\MP BA^\MP BB^\MP\geq\Oqq$ and $\Kerna{A}\subseteq\Kerna{B}$.
		\il{L1315.iii} $B\geq BA^\MP B\geq\Oqq$ and $\Kerna{A}\subseteq\Kerna{B}$.
		\il{L1315.iv} $\bmat A&B\\ B&B\emat\geq\Ouu{2q}{2q}$.
 \eaeqi
 If~\ref{L1315.i} is fulfilled, then $\Kerna{BA^\MP B}=\Kerna{B}$ and $\Bilda{BA^\MP B}=\Bilda{B}$.
\elemm
\bproof
	\bimp{L1315.i}{L1315.ii}
		From \rrema{R1553} and~\ref{L1315.i} we get $\Kerna{A}\subseteq\Kerna{B}$ and
		\[
			\Iq
			\geq AA^\MP
			=\sqrt{A^\MP}A\sqrt{A^\MP}
			\geq\sqrt{A^\MP}B\sqrt{A^\MP}
			=(\sqrt{B}\sqrt{A^\MP})^\ad(\sqrt{B}\sqrt{A^\MP}),
		\]
		which implies $\Iq\geq(\sqrt{B}\sqrt{A^\MP})(\sqrt{B}\sqrt{A^\MP})^\ad$. Hence,
		\[
			B^\MP
			=\sqrt{B}^\MP\cdot\Iq\cdot\sqrt{B}^\MP
			\geq\sqrt{B}^\MP(\sqrt{B}\sqrt{A^\MP})(\sqrt{B}\sqrt{A^\MP})^\ad\sqrt{B}^\MP
			\geq\Oqq.
		\]
		In view of $\sqrt{B}^\MP(\sqrt{B}\sqrt{A^\MP})(\sqrt{B}\sqrt{A^\MP})^\ad\sqrt{B}^\MP=B^\MP BA^\MP BB^\MP$, then~\ref{L1315.ii} is fulfilled.
	\eimp
	
	\bimp{L1315.ii}{L1315.iii}
		From~\ref{L1315.ii} we get
		\[
			B
			=BB^\MP B
			\geq BB^\MP BA^\MP BB^\MP B
			\geq\Oqq.
		\]
		In view of $BB^\MP BA^\MP BB^\MP B=BA^\MP B$ and~\ref{L1315.ii}, then~\ref{L1315.iii} is fulfilled.
	\eimp
	
	\bimpp{L1315.iii}{L1315.iv}{L1315.iv}{L1315.i}
		These implications are immediate consequences of a well-known characterization of non-negative Hermitian block matrices (see, e.\,g.,~\cite{MR1152328}*{\clemm{1.1.9}, p.~18}) and the equation $BB^\MP B=B$.
	\eimpp
	
%

	Finally, suppose that~\ref{L1315.i} is fulfilled. Hence, we have $\Kerna{A}\subseteq\Kerna{B}$, which, in view of \rrema{L0818}, implies $AA^\MP B=B$. Thus, we obtain
	\bsp
		\Kerna{B}
		&\subseteq\Kerna{BA^\MP B}
		=\KernA{(\sqrt{A^\MP}B)^\ad(\sqrt{A^\MP}B)}
		=\Kerna{\sqrt{A^\MP}B}\\
		&\subseteq\Kerna{A\sqrt{A^\MP}\sqrt{A^\MP}B}
		=\Kerna{AA^\MP B}
		=\Kerna{B},
	\esp
	i.\,e., $\Kerna{BA^\MP B}=\Kerna{B}$. Because of $A,B\in\CHq$, this yields $\Bilda{BA^\MP B}=\Bilda{B}$.
\eproof

Now we state some basic facts on the class
\[
	\CEPq
	\defg\SetaA{A\in\Cqq}{\Bilda{A^\ad}=\Bilda{A}}.
\]
\index{C_EP^x@$\CEPq$}

\bpropl{P1645}
	Let $A\in\Cqq$. Then the following statements are equivalent:
	\baeqi{0}
		\il{P1645.i} $A\in\CEPq$.
		\il{P1645.ii} $\Kerna{A^\ad}=\Kerna{A}$.
		\il{P1645.iii} $AA^\MP=A^\MP A$.
		\il{P1645.iv} $A^\MP\in\CEPq$.
		\il{P1645.v} $\re(A^\MP)=A^\MP\re( A)(A^\MP)^\ad$ and $\im(A^\MP)=-A^\MP\im( A)(A^\MP)^\ad$.
		\il{P1645.vi} $\re(A^\MP)=(A^\MP)^\ad\re( A)A^\MP$ and $\im(A^\MP)=-(A^\MP)^\ad\im( A)A^\MP$.
		\il{P1645.vii} $\re A=A\re( A^\MP)A^\ad$ and $\im A=-A\im( A^\MP)A^\ad$.
		\il{P1645.viii} $\re A=A^\ad\re( A^\MP)A$ and $\im A=-A^\ad\im( A^\MP)A$.
	\eaeqi
\eprop
\bproof
	The equivalence of~\ref{P1645.i},~\ref{P1645.ii},~\ref{P1645.iii}, and~\ref{P1645.iv} can be found in Cheng/Tian~\cite{MR2013464} or Tian/Wang~\cite{MR2763588}. What concerns the equivalence of~\ref{P1645.i},~\ref{P1645.v},~\ref{P1645.vi},~\ref{P1645.vii}, and~\ref{P1645.viii} we refer to~\cite{FKLR}*{\cprop{A.6}}.
\eproof

%

\begin{lem}[\cite{FKLR}*{\clemm{A.10}}]\label{L1412}
	Let $A\in\Iqgg$ where $\Iqgg$ is given via \eqref{RIugg}. Then $\Kerna{A}\subseteq\Kerna{\im A}$, $\Bilda{\im A}\subseteq\Bilda{A}$, and $A\in\CEPq$.
\end{lem}

At the end of this section, we give a slight generalization of a result due to S.~L.~Campbell and C.~D.~Meyer~Jr. This result can be proved by an obvious modification of the proof given in~\cite{MR1105324}*{\ctheo{10.4.1}}.

\begin{lem}[\cite{MR1105324}*{\ctheo{10.4.1}}]\label{L0921}
	Let $\seq{A_n}{n}{1}{\infty}$ be a sequence of complex \tpqa{matrices} which converges to a complex \tpqa{matrix} $A$. Then $\seq{A_n^\MP}{n}{1}{\infty}$ is convergent if and only if there is a positive integer $m$ such that $\rank A_n=\rank A$ for each integer $n$ with $n\geq m$. In this case, $\seq{A_n^\MP}{n}{1}{\infty}$ converges to $A^\MP$.
\end{lem}

\section{On Linear Fractional Transformations of Matrices}\label{A0853}
In this appendix, we summarize some basic facts on linear fractional transformations of matrices which are needed in the paper. This material is mostly taken from~\cite{MR566141} and~\cite{MR1152328}*{\cSect{1.6}}.

Let $a\in\Coo{p}{p}$, $b\in\Coo{p}{q}$, $c\in \Cqp$, $d\in\Cqq$, and let 
\[ 
	E
	\defg
	\bMat
		a&b\\
		c&d
	\eMat.
\]
If the set $\dblftruu{c}{d}\defg\setaa{x\in\Cpq}{\det(cx+d)\neq0}$\index{Q_[,]@$\dblftruu{c}{d}$}is non-empty then the linear fractional transformation $\lftroou{p}{q}{E}\colon\dblftruu{c}{d}\to\Coo{p}{q}$\index{S_^(,)@$\lftroou{p}{q}{E}$} is defined by
\[
	\lftrooua{p}{q}{E}{x}
	\defg(ax+b)(cx+d)^\inv.
\]


The following well-known result shows that the composition of two linear fractional transformations is again a mapping of this type.


\begin{prop}[see, e.\,g.~\cite{MR1152328}*{\cprop{1.6.3}}]\label{P1509}
	Let $a_1,a_2\in\Cpp$, let $b_1,b_2\in\Cpq$, let $c_1,c_2\in\Cqp$, and let $d_1,d_2\in\Cqq$ be such that $\rank\brow c_1,d_1\erow=\rank\brow c_2,d_2\erow=q$. Furthermore, let $E_1\defg\bmat a_1&b_1\\ c_1&d_1\emat$, $E_2\defg\bmat a_2&b_2\\ c_2&d_2\emat$, $E\defg E_2E_1$, and $E=\bmat a&b\\c&d\emat$ be the block representation of $E$ with \tppa{block} $a$. Then $\mathcal{Q}\defg\setaa{x\in\dblftruu{c_1}{d_1}}{\lftrooua{p}{q}{E_1}{x}\in\dblftruu{c_2}{d_2}}$ is a nonempty subset of the set $\dblftruu{c}{d}$ and $\lftrooua{p}{q}{E_2}{\lftrooua{p}{q}{E_1}{x}}=\lftrooua{p}{q}{E}{x}$ holds true for all $x\in\mathcal{Q}$.
\end{prop}

We make the following convention: If a non-empty subset $\mG$ of $\C$ and a matrix-valued function $V\colon\mG\to\Coo{2q}{2q}$ with \tqqa{block} partition $V=\bmat v_{11}&v_{12}\\v_{21}&v_{22}\emat$ and a matrix-valued function $F\colon\mG\to\Cqq$ with $\det[v_{21}(z)F(z)+v_{22}(z)]\neq0$ for all $z\in\mG$ are given, then we will use the notation $\lftrfua{V}{F}$\index{S_()@$\lftrfua{V}{F}$} for the function $\lftrfua{V}{F}\colon\mG\to\Cqq$ defined by $\lftrfuaa{V}{F}{z}\defg\lftrooua{q}{q}{V(z)}{F(z)}$ for all $z\in\mG$.

\section{The Matrix Polynomials $\mHTiuu{A}{B}$ and $\mHTuu{A}{B}$}\label{A0835}
In this appendix, we study special linear \taaa{(p+q)}{(p+q)}{matrix} polynomials which are intensively used in \rsect{S1515}. Let $A,B\in\Cpq$. Then we define $\mHTuu{A}{B}\colon\C\to\Coo{(p+q)}{(p+q)}$\index{W_,@$\mHTuu{A}{B}$} and $\mHTiuu{A}{B}\colon\C\to\Coo{(p+q)}{(p+q)}$\index{V_,@$\mHTiuu{A}{B}$} by
\bgl{Wuua}
	\mHTuua{A}{B}{z}
	\defg
	\bMat
		z\Ip-BA^\MP&A\\
		-A^\MP&\Iq-A^\MP A
	\eMat
\eg
and
\bgl{Vuua}
	\mHTiuua{A}{B}{z}
	\defg
	\bMat
		\Opp&-A\\
		A^\MP&z\Iq-A^\MP B
	\eMat.
\eg
If $B=\Opq$, then we set
\begin{align}\label{WVu}
	\mHTu{A}&\defg\mHTuu{A}{\Opq}&
	&\text{and}&
	\mHTiu{A}&\defg\mHTiuu{A}{\Opq}.
\end{align}\index{W_@$\mHTu{A}$}\index{V_@$\mHTiu{A}$}

The use of the matrix polynomial $\mHTiuu{A}{B}$ was inspired by some constructions in the paper~\cite{MR1624548}. In particular, we mention~\cite{MR1624548}*{p.~225, formula~(4.12)}. In their constructions Chen and Hu used Drazin inverses instead of Moore-Penrose inverses of matrices. Since both types of generalized inverses coincide for Hermitian matrices (see~\cite{103}*{\cprop{A.2}} we can conclude that in generic case the matrix polynomials $\mHTiuu{A}{B}$ coincide with the objects used in~\cite{MR1624548}.

\brema
	Let $A,B\in\Cpq$ and let $z\in\C$. Then one can easily see that
	\[
		\mHTiuua{A}{B}{z}\mHTuua{A}{B}{z}
		=\diag\bRow AA^\MP, A^\MP[A-B(\Iq-A^\MP A)]+z(\Iq-A^\MP A)\eRow
	\]
	and
	\[
		\mHTuua{A}{B}{z}\mHTiuua{A}{B}{z}
		=\lek
		\begin{array}{c|c}
			AA^\MP&BA^\MP A-AA^\MP B\\
			\hline
			\Oqp&A^\MP A+z(\Iq-A^\MP A)
		\end{array}
		\rek.
	\]
\erema
%
%
%


Now we are going to study the linear fractional transformation generated by the matrix $\mHTuua{A}{B}{z}$ for arbitrarily given $z\in\C$.

\blemml{L0935}
	Let $A\in\Cpq$. Then:
	\benui
		\il{L0935.a} The matrix $-A$ belongs to $\dblftruu{-A^\MP}{\Iq-A^\MP A}$. In particular, $\dblftruu{-A^\MP}{\Iq-A^\MP A}\neq\emptyset$.
		\il{L0935.b} Let $X\in\Cpq$ be such that $\Kerna{A}\subseteq\Kerna{X}$ and $\Bilda{A}\subseteq\Bilda{X}$. Then $X\in\dblftruu{-A^\MP}{\Iq-A^\MP A}$ and $(-A^\MP X +\Iq-A^\MP A)^\inv=-X^\MP A+\Iq-A^\MP A$.
	\eenui
\elemm
\bproof
	\eqref{L0935.a} This follows from $-A^\MP(-A)+\Iq-A^\MP A=\Iq$.
	
	\eqref{L0935.b} In view of \rlemm{L1633}, we have $AA^\MP=XX^\MP$ and $A^\MP A=X^\MP X$. Hence,
	\bsp
		&(-A^\MP X +\Iq-A^\MP A)(-X^\MP A+\Iq-A^\MP A)\\
		&=A^\MP XX^\MP A-A^\MP X+A^\MP XA^\MP A-X^\MP A+\Iq-A^\MP A\\
		&\quad+A^\MP AX^\MP A-A^\MP A+A^\MP AA^\MP A\\
		&=A^\MP AA^\MP A-A^\MP X+A^\MP XX^\MP X-X^\MP A+\Iq-A^\MP A+X^\MP XX^\MP A-A^\MP A+A^\MP A\\
		&=A^\MP A-A^\MP X+A^\MP X-X^\MP A+\Iq-A^\MP A+X^\MP A
		=\Iq.
	\esp
	This completes the proof of~\eqref{L0935.b}. 
\eproof

\blemml{L1213}
	Let $A,B\in\Cpq$ and let $\mHTuu{A}{B}$ be defined via \eqref{Wuua}. Let $X\in\Cpq$ be such that the inclusions $\Kerna{A}\subseteq\Kerna{X}$ and $\Bilda{A}\subseteq\Bilda{X}$ are satisfied. Furthermore, let $z\in\C$. Then:
	\benui
		\il{L1213.a} The matrix $X$ belongs to $\dblftruu{-A^\MP}{\Iq-A^\MP A}$. Furthermore,
		\begin{gather*}
			\lftrooua{p}{q}{\mHTuua{A}{B}{z}}{X}
			=-A(z\Iq+X^\MP A)+BA^\MP A,\\
			\Kerna{A}
			\subseteq\KernA{\lftrooua{p}{q}{\mHTuua{A}{B}{z}}{X}},
		\end{gather*}
		and
		\[
			\BildA{\lftrooua{p}{q}{\mHTuua{A}{B}{z}}{X}-BA^\MP A}
			\subseteq\Bilda{A}.
		\]
		\il{L1213.b} If $\Kerna{A}\subseteq\Kerna{B}$, then
		\ba
			\lftrooua{p}{q}{\mHTuua{A}{B}{z}}{X}
			&=-A(z\Iq+X^\MP A)+B&
			&\text{and}&
			\BildA{\lftrooua{p}{q}{\mHTuua{A}{B}{z}}{X}-B}
			&\subseteq\Bilda{A}.
		\ea
	\eenui
\elemm
\bproof
	In view of \rpart{L0935.b} of \rlemm{L0935}, we have $X\in\dblftruu{-A^\MP}{\Iq-A^\MP A}$ and
	\bgl{L1213.1}
		(-A^\MP X +\Iq-A^\MP A)^\inv
		=-X^\MP A+\Iq-A^\MP A.
	\eg
	Because of the choice of $X$, \rpartss{L1622.a}{L1622.b} of \rrema{L1622} yield $XA^\MP A=X$ and $XX^\MP A=A$, and, in view of \eqref{L1213.1}, then
	\bspl{L1213.4}
		&\lftrooua{p}{q}{\mHTuua{A}{B}{z}}{X}
		=\lek(z\Ip-BA^\MP)X+A\rek(-A^\MP X+\Iq-A^\MP A)^\inv\\
		&=\lek(z\Ip-BA^\MP)X+A\rek(-X^\MP A+\Iq-A^\MP A)\\
		&=(z\Ip-BA^\MP)X(-X^\MP A+\Iq-A^\MP A)+A(-X^\MP A+\Iq-A^\MP A)\\
		&=-(z\Ip-BA^\MP)XX^\MP A+(z\Ip-BA^\MP)X(\Iq-A^\MP A)-AX^\MP A+A-AA^\MP A\\
		&=-(z\Ip-BA^\MP)A+(z\Ip-BA^\MP)X(\Iq-A^\MP A)-AX^\MP A\\
		&=(z\Ip-BA^\MP)\lek X(\Iq-A^\MP A)-A\rek-AX^\MP A
		=-(z\Ip-BA^\MP)A-AX^\MP A\\
		&=-A(z\Iq+X^\MP A)+BA^\MP A.
	\esp
	The remaining assertions of~\eqref{L1213.a} are immediate consequences of \eqref{L1213.4}.
	
	\eqref{L1213.b} follows from~\eqref{L1213.a} and \rpart{L1622.a} of \rrema{L1622}.
\eproof

\bremal{R1105}
	Let $A,B\in\Cpq$ and let $z\in\C$. Further, let $X\in\dblftruu{A^\MP}{z\Iq-A^\MP B}$. Then from \eqref{Vuua} we see that
	\[
		\lftrooua{p}{q}{\mHTiuua{A}{B}{z}}{X}
		=-A\lek z\Iq+A^\MP(X-B)\rek^\MP
	\]
	and, in view of $\det[z\Iq+A^\MP(X-B)]\neq0$, thus $\Bilda{\lftrooua{p}{q}{\mHTiuua{A}{B}{z}}{X}}=\Bilda{A}$.
\erema

\blemml{L1058}
	Let $A\in\Cggq$ and $B\in\Cqq$. Further, let $X\in\Cqq$ be such that
	\begin{align*}
		X-B&\in\Iqgg,&
		\Kerna{A}&\subseteq\Kerna{X-B},&
		&\text{and}&
		\Bilda{X-B}&\subseteq\Bilda{A}
	\end{align*}
	are satisfied. For each $z\in\ohe$, then the following statements hold true:
	\benui
		\il{L1058.a} $\im(X+zA-B)\geq(\im z)A\in\Cggq$.
		\il{L1058.b} $\Kerna{X+zA-B}\subseteq\Kerna{A}$.
		\il{L1058.c} $X+zA-B=A(A^\MP X+z\Iq-A^\MP B)$.
		\il{L1058.d} The matrix $X$ belongs to $\dblftruu{A^\MP}{z\Iq-A^\MP B}$. Furthermore,
		\ba
			\lftrooua{q}{q}{\mHTiuua{A}{B}{z}}{X}
			&=-A(X+zA-B)^\MP A&
			&\text{and}&
			\Kerna{A}
			&\subseteq\KernA{\lftrooua{q}{q}{\mHTiuua{A}{B}{z}}{X}}.
		\ea
	\eenui
\elemm
\bproof
	Let $z\in\C$.
	
	\eqref{L1058.a} Because of $A\in\Cggq$, $\im(X-B)\in\Cggq$, and $\im z\in(0,+\infty)$, we have
	\[
		\im(X+zA-B)
		=\im(X-B)+(\im z)A
		\geq(\im z)A
		\in\Cggq.
	\]
	
	\eqref{L1058.b} From~\eqref{L1058.a} and $\im z\neq0$ we get
	\bgl{L1058.1}
		\KernA{\im(X+zA-B)}
		\subseteq\KernA{(\im z)A}
		=\Kerna{A}.
	\eg
	In view of~\eqref{L1058.a}, we have $X+zA-B\in\Iqgg$. Thus, \rlemm{L1412} gives
	\bgl{L1058.2}
		\Kerna{X+zA-B}
		\subseteq\KernA{\im(X+zA-B)}.
	\eg
	Now the combination of \eqref{L1058.2} and \eqref{L1058.1} yields~\eqref{L1058.b}.
	
	\eqref{L1058.c} In view of $\Bilda{X-B}\subseteq\Bilda{A}$, \rrema{L1622} gives $AA^\MP(X-B)=X-B$. Thus,
	\[
		X+zA-B
		=zA+AA^\MP(X-B)
		=A(A^\MP X+z\Iq-A^\MP B).
	\]
	
	\eqref{L1058.d} We first prove that
	\bgl{L1058.3}
		\Kerna{A^\MP X+z\Iq-A^\MP B}
		=\set{\Ouu{q}{1}}.
	\eg
	In view of~\eqref{L1058.c},~\eqref{L1058.b} and $\Kerna{A}\subseteq\Kerna{X-B}$, we get
	\bgl{L1058.4}
		\Kerna{A^\MP X+z\Iq-A^\MP B}
		\subseteq\Kerna{X+zA-B}
		\subseteq\Kerna{A}
		\subseteq\Kerna{X-B}.
	\eg
	Let
	\bgl{L1058.5}
		v
		\in\Kerna{A^\MP X+z\Iq-A^\MP B}.
	\eg
	From \eqref{L1058.5} and \eqref{L1058.4} we get then $(X-B)v=\Ouu{q}{1}$. Thus, taking again \eqref{L1058.5} into account we see
	\[
		v
		=\frac{1}{z}\lek zv+A^\MP(X-B)v\rek
		=\frac{1}{z}(A^\MP X+z\Iq-A^\MP B)v
		=\frac{1}{z}\cdot\Ouu{q}{1}
		=\Ouu{q}{1}.
	\]
	Thus \eqref{L1058.3} is proved. Hence
	\bgl{L1058.6}
		X
		\in\dblftruu{A^\MP}{z\Iq-A^\MP B}.
	\eg
	In view of~\eqref{L1058.b} we infer from \rpart{L1622.a} of \rrema{L1622} then
	\bgl{L1058.7}
		A(X+zA-B)^\MP(X+zA-B)
		=A.
	\eg
	Using \eqref{Vuua}, \eqref{L1058.6}, \eqref{L1058.7} and~\eqref{L1058.c} we get
	\bsp
		\lftrooua{q}{q}{\mHTiuua{A}{B}{z}}{X}
		&=-A(A^\MP X+z\Iq-A^\MP B)^\inv\\
		&=-A(X+zA-B)^\MP(X+zA-B)(A^\MP X+z\Iq-A^\MP B)^\inv\\
		&=-A(X+zA-B)^\MP A(A^\MP X+z\Iq-A^\MP B)(A^\MP X+z\Iq-A^\MP B)^\inv\\
		&=-A(X+zA-B)^\MP A.
	\esp
	The last equation implies
	\[
		\Kerna{A}
		\subseteq\KernA{\lftrooua{q}{q}{\mHTiuua{A}{B}{z}}{X}}.
	\]
	
	The proof is complete.
\eproof

\end{document}